\documentclass[10pt,a4paper]{article}
\usepackage{graphics}
\usepackage{amsmath,amssymb,amsthm,graphics,graphicx,epic,eepic,multicol,ascmac}
\usepackage[mathscr]{euscript}
\usepackage[all]{xy}
\usepackage{fancybox}
\usepackage{url}
\usepackage{ytableau}
\usepackage{here}
\usepackage[toc,page]{appendix}

\ytableausetup{centertableaux,mathmode,boxsize=1em}
\def\young(#1){\ytableaushort{#1}}
\def\yng(#1){\ydiagram{#1}}

\numberwithin{equation}{section}
\theoremstyle{theorem}
\newtheorem{thm}{Theorem}[section]
\newtheorem{prop}[thm]{Proposition}

\newtheorem{rem}[thm]{Remark}

\theoremstyle{definition}
\newtheorem{defn}[thm]{Definition}
\newtheorem{ex}[thm]{Example}

\def\al{\alpha}

\def\wht(#1){\widehat{\ #1\ }}

\newcommand{\ch}{\mathrm{ch}}

\newcommand{\lbr}{\begin{bmatrix}}
\newcommand{\rbr}{\end{bmatrix}}

\newcommand{\cd}{commutative diagram }

\def\al{\alpha}

\def\beneme{\begin{enumerate}}
\def\beq{\begin{equation}}
\def\beqn{\begin{eqnarray}}
\def\beqnn{\begin{eqnarray*}}

\def\bfii0{{\bf i_0}}

\def\bbra#1,#2,#3{\left\{\begin{array}{c}\hspace{-5pt}
#1;#2\\ \hspace{-5pt}#3\end{array}\hspace{-5pt}\right\}}
\def\cd{\cdots}
\def\ci(#1,#2){c_{#1}^{(#2)}}
\def\Ci(#1,#2){C_{#1}^{(#2)}}
\def\mpp(#1,#2,#3){#1^{(#2)}_{#3}}
\def\bCi(#1,#2){\ovl C_{#1}^{(#2)}}
\def\ch(#1,#2){c_{#2,#1}^{-h_{#1}}}
\def\cc(#1,#2){c_{#2,#1}}

\def\di(#1,#2){D_{#1}^{(#2)}}
\def\dbi(#1,#2){\ovl D_{#1}^{(#2)}}

\def\eneme{\end{enumerate}}

\def\eeq{\end{equation}}
\def\eeqn{\end{eqnarray}}
\def\eeqnn{\end{eqnarray*}}

\def\gau#1,#2{\left[\begin{array}{c}\hspace{-5pt}#1\\
\hspace{-5pt}#2\end{array}\hspace{-5pt}\right]}

\def\ify{\infty}
\def\io{\iota}
\def\ji(#1,#2){j_{#1}^{(#2)}}

\def\lan{\langle}

\def\nd{\noindent}

\def\ovl{\overline}

\def\qed{\hfill\framebox[2mm]{}}
\def\QQ{\mathbb Q}

\def\ran{\rangle}

\def\TY(#1,#2,#3){#1^{(#2)}_{#3}}

\def\vp{\varphi}

\def\xxi(#1,#2,#3){\displaystyle {}^{#1}\Xi^{(#2)}_{#3}}
\def\xsi(#1,#2,#3){\displaystyle {}^{#1}\Sigma^{(#2)}_{#3}}
\def\xE(#1,#2,#3){\displaystyle {}^{#1}E_{#2}[#3]}
\def\xF(#1,#2){\displaystyle {}^{#1}F_{#2}}
\def\xx(#1,#2){\displaystyle {}^{#1}\Xi_{#2}}
\def\W1{W(\varpi_1)}

\def\what{\widehat}

\def\ZZ{\mathbb Z}

\def\m@th{\mathsurround=0pt}
\def\fsquare(#1,#2){
\hbox{\vrule$\hskip-0.4pt\vcenter to #1{\normalbaselines\m@th
\hrule\vfil\hbox to #1{\hfill$\scriptstyle #2$\hfill}\vfil\hrule}$\hskip-0.4pt
\vrule}}

\newcommand{\ba}{\begin{array}}
\newcommand{\ea}{\end{array}}

\newcommand{\eq}{\begin{eqnarray}}
\newcommand{\eneq}{\end{eqnarray}}

\title{\textbf{\large{Polyhedral realizations for crystal bases
and Young walls of classical affine types
}}}
\author{\normalsize{YUKI KANAKUBO\thanks{Max Planck Institute for Mathematics,
Vivatsgasse 7, 53111 Bonn, Germany : {j\_chi\_sen\_you\_ky@eagle.sophia.ac.jp}}}
}
\date{}

\begin{document}

\maketitle
\vspace{-10pt}

\begin{abstract}
For affine Lie algebra $\mathfrak{g}$ of type $A^{(1)}_{n-1}$, $B^{(1)}_{n-1}$, $C^{(1)}_{n-1}$, $D^{(1)}_{n-1}$, $A^{(2)}_{2n-2}$, $A^{(2)}_{2n-3}$ or $D^{(2)}_{n}$,
let $B(\lambda)$ and $B(\infty)$ be the crystal bases
of integrable highest weight representation $V(\lambda)$ and negative part $U_q^-(\mathfrak{g})$ of quantum group $U_q(\mathfrak{g})$.
We consider the polyhedral realizations of crystal bases, which realize 
$B(\lambda)$ and $B(\infty)$ as sets of integer points of some polytopes and cones in $\mathbb{R}^{\infty}$.
It is a natural problem to find explicit forms of the polytopes and cones.
In this paper, 
we introduce pairs of truncated walls, which are defined
as modifications of level $2$-Young walls and
describe inequalities defining the polytopes and cones in terms of level $1$-proper Young walls and pairs of truncated walls.
As an application, we also give
combinatorial
descriptions of $\varepsilon_k^*$-functions on $B(\infty)$ in terms of Young walls and truncated walls.
\end{abstract}

\section{Introduction}

The
crystal base is one of the most important objects in the representation theory of symmetrizable Kac-Moody Lie algebras $\mathfrak{g}$
or quantum groups $U_q(\mathfrak{g})$ \cite{K0}. Realizing crystal bases as objects like as Young tableaux, path models, Laurent monomials,
one can combinatorially study the structure of representations.

In \cite{Kang, Kang-Lee}, when $\mathfrak{g}$ is a classical affine Lie algebra,
crystal bases $B(\lambda)$ for integrable highest weight representations are realized in terms of
Young walls. When $\lambda$ is a fundamental weight with level $1$, the crystal bases are realized via sets of reduced proper Young walls (except for type $C^{(1)}_{n-1}$-case) \cite{Kang}. For a general level $\ell$, $B(\lambda)$ are realized by level $\ell$ Young walls \cite{Kang-Lee}. In particular, when $\lambda=\Lambda_k$ is
a fundamental weight of level $2$, the crystal base is described as a set of pairs of walls, which are obtained from a ground state wall $Y_{\Lambda_k}$ (Definition \ref{YWk-def})
by stacking several blocks.

Polyhedral realizations for $B(\lambda)$ and $B(\infty)$ are introduced for symmetrizable Kac-Moody Lie algebras $\mathfrak{g}$,
which realize crystal bases as sets of integer points of polytopes or cones in $\mathbb{R}^{\infty}$ \cite{N99,NZ}.
Here,
$B(\infty)$ is the crystal base for the negative part $U_q^-(\mathfrak{g})$ of quantum group.
There exist embeddings of crystals $\Psi_{\iota}^{(\lambda)}:B(\lambda)\hookrightarrow \mathbb{Z}^{\infty}_{\iota}\otimes R_{\lambda}$
and $\Psi_{\iota}:B(\infty)\hookrightarrow \mathbb{Z}^{\infty}_{\iota}$
associated with an infinite sequence $\iota$ from the index set $I=\{1,2,\cdots,n\}$ of $\mathfrak{g}$,
where $R_{\lambda}$ is a crystal with a single element (Example \ref{r-ex}) and
\[
\mathbb{Z}^{\infty}_{\iota}
=\{(\cd,a_r,\cd,a_2,a_1)| a_r\in\ZZ
\,\,{\rm and}\,\,a_r=0\,\,{\rm for}\,\,r\gg 0\}
\]
has a structure of crystal associated with $\iota$.
When $\mathfrak{g}$ is a finite dimensional simple Lie algebra, the set ${\rm Im}(\Psi_{\iota})$ can be identified with the set of integer points in string cones in \cite{BZ, Lit}.
It is a natural problem to find explicit forms of ${\rm Im}(\Psi_{\iota})(\cong B(\infty))$ and ${\rm Im}(\Psi_{\iota}^{(\lambda)})(\cong B(\lambda))$.
A bunch of researchers are working on this problem. 
In case that $\mathfrak{g}$ is a finite dimensional simple Lie algebra and $\iota$ is a specific one then
explicit forms of ${\rm Im}(\Psi_{\iota})$ and ${\rm Im}(\Psi_{\iota}^{(\lambda)})$ are
given in \cite{H1, KS, Lit, N99,NZ}. For the finite dimensional simple Lie algebra $\mathfrak{g}$ of type $A$,
combinatorial descriptions of string cones are provided in \cite{GKS16, GP}.
When
$\mathfrak{g}$ is a classical affine Lie algebra and $\iota$ is a specific one,
explicit forms of them are also
given in \cite{H2, HN, N99, NZ}.
In our previous papers \cite{KaN, KaN2},
for finite dimensional classical Lie algebras $\mathfrak{g}$,
assuming that $\iota$ is adapted (see Definition \ref{adapt}), we explicitly describe inequalities defining ${\rm Im}(\Psi_{\iota})$ and ${\rm Im}(\Psi_{\iota}^{(\lambda)})$
in terms of rectangular Young tableaux. 
In \cite{Ka, Ka2}, taking 
$\mathfrak{g}$ as a classical affine Lie algebra of type
$X=A^{(1)}_{n-1}$, $C^{(1)}_{n-1}$, $A^{(2)}_{2n-2}$ or $D^{(2)}_{n}$
and $\iota$ as an adapted one,
we give explicit forms of inequalities for ${\rm Im}(\Psi_{\iota})$ and ${\rm Im}(\Psi_{\iota}^{(\lambda)})$
in terms of extended Young diagrams and Young walls of type $X^L$. Here, $X^L$ means the Langlands dual type of $X$ (see (\ref{Langtype}))
and
extended Young diagrams are introduced
to realize irreducible integrable highest weight modules $V(\Lambda)$ of the quantum affine algebras $U_q(\mathfrak{g})$
as Fock space representations in the case $\mathfrak{g}$ is of type $A^{(1)}_{n-1}$, $C^{(1)}_{n-1}$, $A^{(2)}_{2n-2}$ or $D^{(2)}_{n}$
for almost all fundamental weights $\Lambda$ \cite{Ha, JMMO, KMM}.
We can expect that inequalities defining ${\rm Im}(\Psi_{\iota})$ and ${\rm Im}(\Psi^{(\lambda)}_{\iota})$
can be described in terms of some combinatorial objects related to fundamental representations of $U_q(\mathfrak{g}^L)$ for more general $\mathfrak{g}$.

Based on this philosophy,
in this paper, we assume $\mathfrak{g}$ is of classical affine type $X$, that is, 
$X=A^{(1)}_{n-1}$, $B^{(1)}_{n-1}$, $C^{(1)}_{n-1}$, $D^{(1)}_{n-1}$, $A^{(2)}_{2n-2}$, $A^{(2)}_{2n-3}$ or $D^{(2)}_{n}$ and $\iota$ is adapted
and
introduce the notion of the pairs of proper truncated walls, which are modified notion of proper level-$2$ Young walls.
The main results are as follows: 

(I) We describe inequalities defining ${\rm Im}(\Psi_{\iota})$ in terms of Young walls and truncated walls.

(II) We also describe inequalities defining ${\rm Im}(\Psi^{(\lambda)}_{\iota})$ in terms of Young walls and truncated walls and easy boxes.

(III) As an application, we also give combinatorial
descriptions of $\varepsilon_k^*$-functions on $B(\infty)$ (Theorem \ref{mainthm3}) in terms of Young walls and truncated walls.

Let us explain the detail.
(I) We will define a set of linear maps $\Xi'_{\iota}$ in (\ref{xiprime}), which satisfies
${\rm Im}(\Psi_{\iota})=\{\textbf{a}\in\mathbb{Z}^{\infty} | \varphi(\textbf{a})\geq0\ \text{for any }\varphi\in\Xi'_{\iota}\}$.
The set $\Xi'_{\iota}$ is naturally decomposed 
into certain subsets $\Xi'_{s,k,\iota}$ as $\Xi'_{\iota}=\bigcup_{k\in I, s\in\mathbb{Z}_{\geq1}} \Xi'_{s,k,\iota}$.
As for level $1$ indices $k\in I$ (except for type $C^{(1)}_{n-1}$-case), the set $\Xi'_{s,k,\iota}$ is described by using proper Young walls of type $X^L$
as 
\[
\Xi'_{s,k,\iota}=\{L^{X^L}_{s,k,\iota}(Y) | s\in\mathbb{Z}_{\geq1},\ Y\in {\rm YW}^{X^L}_k\},
\]
where
${\rm YW}^{X^L}_k$ is the set of proper Young walls of ground state $\Lambda_k$ for type $X^L$,
$L^{X^L}_{s,k,\iota}(Y)\in(\mathbb{Q}^{\infty})^*$ is an assigned linear map (see subsection \ref{binf}).
As for other indices $k\in I$, 
the set $\Xi'_{s,k,\iota}$ is similarly described in terms of truncated walls of ground state $\Lambda_k$ (Theorem \ref{mainthm1}).

(II) We also define the set of functions $\Xi'_{\iota}[\lambda]$ in (\ref{xilamdef}) such that
${\rm Im}(\Psi_{\iota}^{(\lambda)})=\{\textbf{a}\in\mathbb{Z}^{\infty} | \varphi(\textbf{a})\geq0\ \text{for any }\varphi\in\Xi'_{\iota}[\lambda]\}$.
The set $\Xi'_{\iota}[\lambda]$ can be written as $\Xi'_{\iota}[\lambda]=\Xi'_{\iota}\cup \bigcup_{k\in I} \Xi'_{k,\iota}[\lambda]$ (see subsection 6.2).
Depending on $\iota$, the set $\Xi'_{k,\iota}[\lambda]$ is described as one of the following three : (1) a set of single or double functions, (2)
a set of functions parametrized by boxes, (3) a set of functions parametrized by proper Young walls or truncated walls (Theorem \ref{mainthm2}).

(III) Using functions in $\Xi'_{k,\iota}[0]$ of (II), one can describe $\varepsilon_k^*$-function on $B(\infty)$.
Thus, as an application of (II),
one obtains a combinatorial description of $\varepsilon_k^*$-function (Theorem \ref{mainthm3}).

The organization
of this paper is as follows:
In Sect.2, we review crystals and polyhedral realizations. In particular, we review
procedures to compute the explicit forms of ${\rm Im}(\Psi_{\iota})$ and ${\rm Im}(\Psi_{\iota}^{(\lambda)})$.
In Sect.3, we recall Young walls and introduce truncated walls.
Sect.4 is devoted to present our main theorems.
In Sect.5, we show key propositions, which imply adding/removing blocks of walls are compatible with procedures of Sect.2 in a certain sense.
We complete all proofs in Sect.6.

\vspace{2mm}

\nd
{\bf Acknowledgements}
The author thanks Max Planck Institute for Mathematics in Bonn for its hospitality and financial support.

\section{Polyhedral realizations and procedures}

\subsection{Notation}

Let $\mathfrak{g}$ be a symmetrizable Kac-Moody algebra over $\mathbb{Q}$,
$A=(a_{i,j})_{i,j\in I}$ be a generalized Cartan matrix
with the index set $I=\{1,2,\cdots,n\}$,
$\mathfrak{h}$ a Cartan subalgebra, $P\subset \mathfrak{h}^*$ a weight lattice,
$\{\alpha_i\}_{i\in I}$ a set of simple roots and $\{h_i\}_{i\in I}$ a set of simple coroots. 
Let
$\langle \cdot,\cdot \rangle : \mathfrak{h} \times \mathfrak{h}^*\rightarrow \mathbb{Q}$
denote the canonical pairing, $P^*:=\{h\in\mathfrak{h} | \langle h,P \rangle\subset\mathbb{Z}\}$ and
$P^+:=\{\lambda\in P | \langle h_i,\lambda \rangle \in\mathbb{Z}_{\geq0} \text{ for all }i\in I \}$.
In particular, we have $\langle h_{i},\alpha_j \rangle=a_{i,j}$. For each $i\in I$, the fundamental weight $\Lambda_i\in P^+$ is defined by
$\langle h_j,\Lambda_i \rangle=\delta_{i,j}$.
Let $U_q(\mathfrak{g})$ be the quantized universal enveloping algebra, which is an associative $\mathbb{Q}(q)$-algebra
with generators $e_i$, $f_i$ ($i\in I$) and $q^h$ ($h\in P^*$) satisfying the usual relations. 
Let $U_q^-(\mathfrak{g})$ denote
the subalgebra of $U_q(\mathfrak{g})$ generated by $f_i$ ($i\in I$).
For $\lambda\in P^+$, 
it is known that the irreducible integrable highest weight module $V(\lambda)$ of $U_q(\mathfrak{g})$ with highest weight $\lambda$
has a crystal base $(L(\lambda),B(\lambda))$. The algebra $U_q^-(\mathfrak{g})$ also has
a crystal base $(L(\infty),B(\infty))$ (see \cite{K0,K1}).

We will use the following
numbering of vertices in affine Dynkin diagrams:
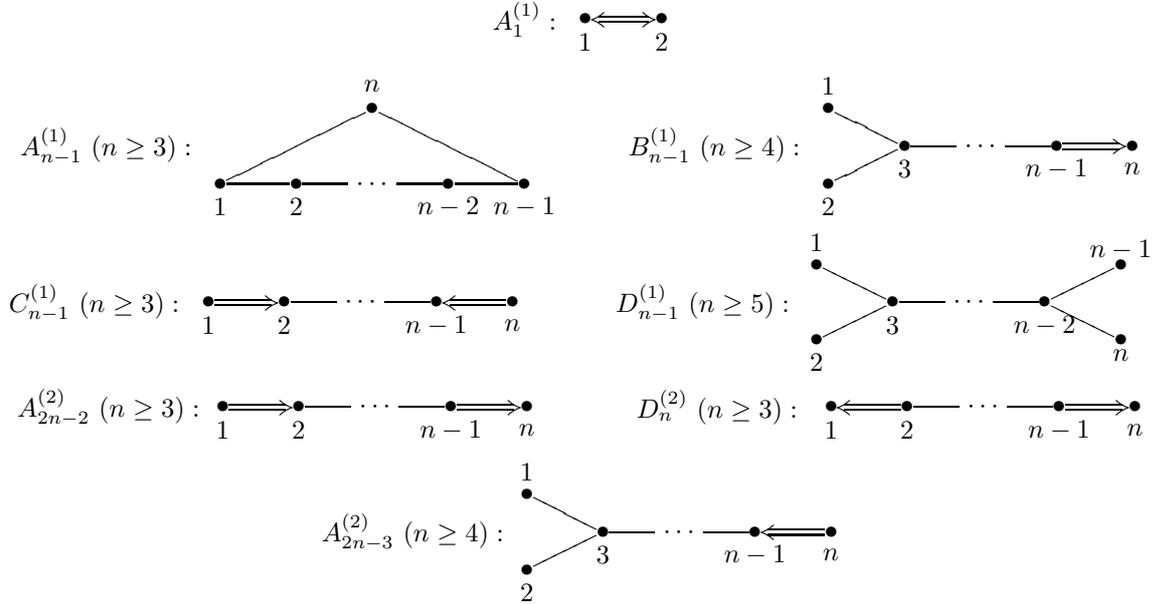
\[
\begin{xy}
(-8,0) *{A_{1}^{(1)} : }="A1",
(0,0) *{\bullet}="1",
(0,-3) *{1}="1a",
(10,0)*{\bullet}="2",
(10,-3) *{2}="2a",
\ar@{<=>} "1";"2"^{}
\end{xy}
\]
\[
\begin{xy}
(-15,5) *{A_{n-1}^{(1)} \ (n\geq 3) : }="A1",
(20,10) *{\bullet}="n",
(20,13) *{n}="na",
(0,0) *{\bullet}="1",
(0,-3) *{1}="1a",
(10,0)*{\bullet}="2",
(10,-3)*{2}="2a",
(20,0)*{\ \cdots\ }="3",
(30,0)*{\bullet}="4",
(30,-3)*{n-2}="4a",
(40,0)*{\bullet}="5",
(40,-3)*{n-1}="5a",
(65,5) *{B_{n-1}^{(1)} \ (n\geq 4) : }="B1",
(80,10) *{\bullet}="11",
(80,0) *{\bullet}="111",
(80,13) *{1}="111a",
(80,-3) *{2}="11a",
(90,5)*{\bullet}="22",
(90,2)*{3}="22a",
(100,5)*{\ \cdots\ }="33",
(110,5)*{\bullet}="44",
(110,2)*{n-1}="44a",
(120,5)*{\bullet}="55",
(120,2)*{n}="55a",
\ar@{-} "1";"2"^{}
\ar@{-} "2";"3"^{}
\ar@{-} "3";"4"^{}
\ar@{-} "4";"5"^{}
\ar@{-} "1";"n"^{}
\ar@{-} "n";"5"^{}
\ar@{-} "11";"22"^{}
\ar@{-} "111";"22"^{}
\ar@{-} "22";"33"^{}
\ar@{-} "33";"44"^{}
\ar@{=>} "44";"55"^{}
\end{xy}
\]
\[
\begin{xy}
(-15,5) *{C_{n-1}^{(1)} \ (n\geq 3) : }="C1",
(0,5) *{\bullet}="11",
(0,2) *{1}="11a",
(10,5)*{\bullet}="22",
(10,2)*{2}="22a",
(20,5)*{\ \cdots\ }="33",
(30,5)*{\bullet}="44",
(30,2)*{n-1}="44a",
(40,5)*{\bullet}="55",
(40,2)*{n}="55a",
(65,5) *{D_{n-1}^{(1)} \ (n\geq 5) : }="D1",
(80,10) *{\bullet}="111",
(80,0) *{\bullet}="1111",
(80,13) *{1}="111a1",
(80,-3) *{2}="11a1",
(90,5)*{\bullet}="221",
(90,2)*{3}="22a1",
(100,5)*{\ \cdots\ }="331",
(110,5)*{\bullet}="441",
(110,2)*{n-2}="44a1",
(120,10)*{\bullet}="551",
(120,0)*{\bullet}="5511",
(120,12)*{n-1}="55a1",
(120,-2)*{n}="55a11",
\ar@{=>} "11";"22"^{}
\ar@{-} "22";"33"^{}
\ar@{-} "33";"44"^{}
\ar@{<=} "44";"55"^{}
\ar@{-} "111";"221"^{}
\ar@{-} "1111";"221"^{}
\ar@{-} "221";"331"^{}
\ar@{-} "331";"441"^{}
\ar@{-} "441";"551"^{}
\ar@{-} "441";"5511"^{}
\end{xy}
\]
\[
\begin{xy}
(-15,5) *{A_{2n-2}^{(2)}\ (n\geq 3) : }="A2",
(0,5) *{\bullet}="1",
(0,2) *{1}="1a",
(10,5)*{\bullet}="2",
(10,2)*{2}="2a",
(20,5)*{\ \cdots\ }="3",
(30,5)*{\bullet}="4",
(30,2)*{n-1}="4a",
(40,5)*{\bullet}="5",
(40,2)*{n}="5a",
(65,5) *{D_{n}^{(2)} \ (n\geq 3) : }="D2",
(80,5) *{\bullet}="11",
(80,2) *{1}="11a",
(90,5)*{\bullet}="22",
(90,2)*{2}="22a",
(100,5)*{\ \cdots\ }="33",
(110,5)*{\bullet}="44",
(110,2)*{n-1}="44a",
(120,5)*{\bullet}="55",
(120,2)*{n}="55a",
\ar@{=>} "1";"2"^{}
\ar@{-} "2";"3"^{}
\ar@{-} "3";"4"^{}
\ar@{=>} "4";"5"^{}
\ar@{<=} "11";"22"^{}
\ar@{-} "22";"33"^{}
\ar@{-} "33";"44"^{}
\ar@{=>} "44";"55"^{}
\end{xy}
\]
\[
\begin{xy}
(65,5) *{A_{2n-3}^{(2)} \ (n\geq 4) : }="A2",
(80,10) *{\bullet}="11",
(80,0) *{\bullet}="111",
(80,13) *{1}="111a",
(80,-3) *{2}="11a",
(90,5)*{\bullet}="22",
(90,2)*{3}="22a",
(100,5)*{\ \cdots\ }="33",
(110,5)*{\bullet}="44",
(110,2)*{n-1}="44a",
(120,5)*{\bullet}="55",
(120,2)*{n}="55a",
\ar@{-} "11";"22"^{}
\ar@{-} "111";"22"^{}
\ar@{-} "22";"33"^{}
\ar@{-} "33";"44"^{}
\ar@{<=} "44";"55"^{}
\end{xy}
\]
Replacing our numbering $1,2,\cdots,n-1,n,n+1$ of $A_{2n}^{(2)}$ with $n,n-1,\cdots,1,0$, one obtains the numbering in \cite{Kang}.
We define the type $A_{2n-2}^{(2)\dagger}$ by the following diagram whose numbering is same as \cite{Kang}:
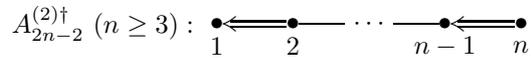
\[
\begin{xy}
(-15,5) *{A_{2n-2}^{(2)\dagger}\ (n\geq 3) : }="A2",
(0,5) *{\bullet}="1",
(0,2) *{1}="1a",
(10,5)*{\bullet}="2",
(10,2)*{2}="2a",
(20,5)*{\ \cdots\ }="3",
(30,5)*{\bullet}="4",
(30,2)*{n-1}="4a",
(40,5)*{\bullet}="5",
(40,2)*{n}="5a",
\ar@{<=} "1";"2"^{}
\ar@{-} "2";"3"^{}
\ar@{-} "3";"4"^{}
\ar@{<=} "4";"5"^{}
\end{xy}
\]
For other types,
replacing our numbering $1,2,\cdots,n-1,n,n+1$
with $0,1,2,\cdots,n$, we obtain the numbering in \cite{Kang, Kang-Lee}.

\subsection{Crystals}

Let us recall the notion of {\it crystals} following \cite{K3}:

\begin{defn}
Let $\mathcal{B}$ be a set and
${\rm wt}:\mathcal{B}\rightarrow P$,
$\varepsilon_k,\varphi_k:\mathcal{B}\rightarrow \mathbb{Z}\sqcup \{-\infty\}$
and $\tilde{e}_k$,$\tilde{f}_k:\mathcal{B}\rightarrow \mathcal{B}\sqcup\{0\}$
($k\in I$) be maps. Then the set $\mathcal{B}$ together with these maps
is said to be a {\it crystal} if 
the following holds: For $b,b'\in\mathcal{B}$ and $k\in I$,
\begin{itemize}
\item $\varphi_k(b)=\varepsilon_k(b)+\langle h_k,{\rm wt}(b)\rangle$,
\item ${\rm wt}(\tilde{e}_kb)={\rm wt}(b)+\alpha_k$ if $\tilde{e}_k(b)\in\mathcal{B}$,
\quad ${\rm wt}(\tilde{f}_kb)={\rm wt}(b)-\alpha_k$ if $\tilde{f}_k(b)\in\mathcal{B}$,
\item $\varepsilon_k(\tilde{e}_k(b))=\varepsilon_k(b)-1,\ \ 
\varphi_k(\tilde{e}_k(b))=\varphi_k(b)+1$\ if $\tilde{e}_k(b)\in\mathcal{B}$, 
\item $\varepsilon_k(\tilde{f}_k(b))=\varepsilon_k(b)+1,\ \ 
\varphi_k(\tilde{f}_k(b))=\varphi_k(b)-1$\ if $\tilde{f}_k(b)\in\mathcal{B}$, 
\item $\tilde{f}_k(b)=b'$ if and only if $b=\tilde{e}_k(b')$,
\item if $\varphi_k(b)=-\infty$ then $\tilde{e}_k(b)=\tilde{f}_k(b)=0$.
\end{itemize}
Here, symbols $0$ and $-\infty$ are additional elements which do not belong to $\mathcal{B}$ and $\mathbb{Z}$, respectively.
\end{defn}
In the above definition, the maps $\tilde{e}_k$,$\tilde{f}_k$ are called {\it Kashiwara operators}. It is well-known that
the sets $B(\lambda)$ $(\lambda\in P^+)$ and $B(\infty)$
have crystal structures.

\begin{ex}\label{r-ex}
For $\lambda\in P$,
let $R_{\lambda}:=\{r_{\lambda}\}$ be the set consisting of a single element $r_{\lambda}$.
One can define a crystal structure on $R_{\lambda}$ as follows:
\[
{\rm wt}(r_{\lambda})=\lambda,\ \ \varepsilon_k(r_{\lambda})= -\lan h_k,\lambda\ran,\ \ \varphi_k(r_{\lambda})=0,\ \ 
\tilde{e}_k(r_{\lambda})=\tilde{f}_k(r_{\lambda})=0.
\]
\end{ex}

\begin{ex}\label{star-ex}
Let $\{\tilde{e}_i\}_{i\in I},\{\tilde{f}_i\}_{i\in I},\{\varepsilon_i\}_{i\in I},\{\varphi_i\}_{i\in I}, {\rm wt}$ denote the crystal
structure on $B(\infty)$ and
$*:U_q(\frak{g})\rightarrow U_q(\frak{g})$ be the antiautomorphism
such that $e_i^*=e_i$, $f_i^*=f_i$ and $(q^h)^*=q^{-h}$ in \cite{K1}.
The antiautomorphism $*$ induces a bijection $*:B(\infty)\rightarrow B(\infty)$ satisfying $*\circ *=id$.
Using this antiautomorphism, we can define another crystal structure on $B(\infty)^{*}$ as follows : The underlying set $B(\infty)^{*}$ equals $B(\infty)$ and
the maps on $B(\infty)^{*}$ are defined by $\tilde{e}_i^*:=*\circ \tilde{e}_i\circ *$, 
$\tilde{f}_i^*:=*\circ \tilde{f}_i\circ *$,
$\varepsilon_i^*:=\varepsilon_i\circ *$,
$\varphi_i^*:=\varphi_i\circ *$ and ${\rm wt}^*:={\rm wt}$.
\end{ex}

\begin{defn}
\begin{enumerate}
\item[$(1)$] For two crystals $\mathcal{B}_1$, $\mathcal{B}_2$,
we say a map $\psi : \mathcal{B}_1\sqcup\{0\}\rightarrow \mathcal{B}_2\sqcup\{0\}$ is a {\it strict morphism} 
and denoted by $\psi : \mathcal{B}_1\rightarrow \mathcal{B}_2$
if $\psi(0)=0$ and the following holds: For $k\in I$ and $b\in\mathcal{B}_1$,
\begin{itemize}
\item if $\psi(b)\in \mathcal{B}_2$ then
\[
{\rm wt}(\psi(b))={\rm wt}(b),\quad
\varepsilon_k(\psi(b))=\varepsilon_k(b),\quad
\varphi_k(\psi(b))=\varphi_k(b),
\]
\item it holds $\tilde{e}_k(\psi(b))=\psi(\tilde{e}_k(b))$ and $\tilde{f}_k(\psi(b))=\psi(\tilde{f}_k(b))$, where
we define $\tilde{e}_k(0)=\tilde{f}_k(0)=0$.
\end{itemize}
\item[$(2)$] If a strict morphism 
$\psi : \mathcal{B}_1\sqcup\{0\}\rightarrow \mathcal{B}_2\sqcup\{0\}$ is injective (resp. bijective) then
$\psi$ is said to be a {\it strict embedding} (resp. {\it isomorphism}) and denoted by $\psi:\mathcal{B}_1 \hookrightarrow \mathcal{B}_2$
(resp. $\psi:\mathcal{B}_1 \overset{\sim}{\rightarrow} \mathcal{B}_2$).
\end{enumerate}
\end{defn}

For two crystals $\mathcal{B}_1$, $\mathcal{B}_2$,
one can define a new crystal called {\it tensor product} $\mathcal{B}_1\otimes\mathcal{B}_2$.
The underlying set of $\mathcal{B}_1\otimes\mathcal{B}_2$ is the direct product
$\mathcal{B}_1\times\mathcal{B}_2$ and a crystal structure is defined in \cite{K3}.

\subsection{Polyhedral realizations of $B(\infty)$}

We fix an infinite sequence
$\io=(\cd,i_r,\cd,i_2,i_1)$ of indices from $I$ such that
\begin{equation}
{\hbox{
$i_r\ne i_{r+1}$ for $r\in\mathbb{Z}_{\geq1}$ and $\sharp\{r\in\mathbb{Z}_{\geq1}| i_r=k\}=\ify$ for all $k\in I$.}}
\label{seq-con}
\end{equation}
Putting
\[
\ZZ^{\ify}
:=\{(\cd,a_r,\cd,a_2,a_1)| a_r\in\ZZ\text{ for }r\in\mathbb{Z}_{\geq1}
\text{ and it holds }a_r=0\,\,{\rm for}\ r\gg 0\},
\]
one can define a crystal structure on the set $\ZZ^{\ify}$ associated with $\iota$ (see subsection 2.4 of \cite{NZ})
and denote it by $\ZZ^{\ify}_{\iota}$.

\begin{thm}\cite{K3,NZ}
\label{emb}
There is the unique strict embedding of crystals
\begin{equation}
\Psi_{\io}:B(\ify)\hookrightarrow \ZZ^{\ify}_{\io}
\label{psi}
\end{equation}
satisfying
$\Psi_{\io} (u_{\ify}) = \textbf{0}$, where
$u_{\ify}\in B(\ify)$ is
the highest weight vector and $\textbf{0}:=(\cd,0,\cd,0,0)\in \ZZ^{\ify}_{\io}$.
\end{thm}

Based on the original paper \cite{NZ}, we gave a procedure
to compute ${\rm Im}(\Psi_{\iota})\subset\mathbb{Z}^{\infty}_{\iota}$ in \cite{Ka}. 
We set
\[
\QQ^{\ify}:=\{\textbf{a}=
(\cd,a_r,\cd,a_2,a_1)| a_r \in \QQ\text{ for }r\in\mathbb{Z}_{\geq1}
\text{ and it holds }
a_r = 0\,\,{\rm for\ } r \gg 0\}.
\]
For $r\in \mathbb{Z}_{\geq1}$, we define $x_r\in(\QQ^{\ify})^*$ as
$x_r(\cd,a_r,\cd,a_2,a_1)=a_r$. One understands $x_r:=0$ when $r\in \mathbb{Z}_{<1}$.
For $r\in \mathbb{Z}_{\geq1}$, let
\[
r^{(+)}:={\rm min}\{\ell\in\mathbb{Z}_{\geq1}\ |\ \ell>r\,\,{\rm and }\,\,i_r=i_{\ell}\},\ \ 
r^{(-)}:={\rm max} \{\ell\in\mathbb{Z}_{\geq1}\ |\ \ell<r\,\,{\rm and }\,\,i_r=i_{\ell}\}\cup\{0\},
\]
and
\begin{equation}
\beta_r:=
x_r+\sum_{r<j<r^{(+)}}\lan h_{i_r},\al_{i_j}\ran x_j+x_{r^{(+)}}\in (\QQ^{\ify})^*.
\label{betak}
\end{equation}
We understand $\beta_0:=0$. 
For $r\in\mathbb{Z}_{\geq1}$,
the map
$S_r'=S_{r,\io}':(\QQ^{\ify})^*\rightarrow (\QQ^{\ify})^*$ is defined as follows:
For $\vp=\sum_{\ell\in\mathbb{Z}_{\geq1}}c_{\ell}x_{\ell}\in(\QQ^{\ify})^*$,
\begin{equation}
S_r'(\vp):=
\begin{cases}
\vp-\beta_r & {\mbox{ if }}\ c_r>0,\\
 \vp+\beta_{r^{(-)}} & {\mbox{ if }}\ c_r< 0,\\
 \vp &  {\mbox{ if }}\ c_r= 0.
\end{cases}
\label{Sk}
\end{equation}
The map $S_r'(\vp)$ is often simply written as $S_r'\vp$. 
Using the map, one defines
\begin{eqnarray}
\Xi_{\io}' &:=  &\{S_{j_{\ell}}'\cd S_{j_2}'S_{j_1}'x_{j_0}\,|\,
\ell\in\mathbb{Z}_{\geq0},j_0,j_1,\cd,j_{\ell}\in\mathbb{Z}_{\geq1}\}, \label{xiprime}\\
\Sigma_{\io}' & := &
\{\textbf{a}\in \ZZ^{\ify}\subset \QQ^{\ify}\,|\,\vp(\textbf{a})\geq0\,\,{\rm for}\,\,
{\rm any}\,\,\vp\in \Xi_{\io}'\}.\nonumber
\end{eqnarray}
If the following condition holds then we say the sequence $\io$ satisfies the $\Xi'$-{\it positivity condition}:
\begin{equation}
{\hbox{for any 
$\vp=\sum_{\ell\in\mathbb{Z}_{\geq1}} c_{\ell}x_{\ell}\in \Xi_{\io}'$, if $\ell^{(-)}=0$ then $c_{\ell}\geq0$}}.
\label{posi}
\end{equation}
We obtain the following procedure:
\begin{thm}\label{inf-thm}\cite{Ka}
We assume that the sequence $\io$ satisfies $(\ref{seq-con})$ 
and $(\ref{posi})$. Then it holds
${\rm Im}(\Psi_{\io})=\Sigma_{\io}'$.
\end{thm}

\subsection{Polyhedral realizations of $B(\lambda)$}

We assume a sequence $\iota$ satisfies (\ref{seq-con}).
For $\lambda\in P^+$,
we set $\mathbb{Z}^{\infty}_{\iota}[\lambda]:=\mathbb{Z}^{\infty}_{\iota}\otimes R_{\lambda}$, where
$R_{\lambda}$ is a crystal defined in Example \ref{r-ex}. Since $R_{\lambda}$ has a single element,
we can naturally identify $\mathbb{Z}^{\infty}_{\iota}[\lambda]$ with $\mathbb{Z}_{\iota}^{\infty}$ as a set.
Just as in Theorem \ref{emb}, the following holds:
\begin{thm}\cite{N99}\label{emb2}
There is the unique strict embedding of crystals
\[
\Psi^{(\lambda)}_{\io}:B(\lambda)\hookrightarrow \ZZ^{\ify}_{\iota}[\lambda]
\]
such that $\Psi^{(\lambda)}_{\io}(u_{\lambda})=(\cdots,0,0,0)$, where
$u_{\lambda}$ is the highest weight vector in $B(\lambda)$ and we identify $\mathbb{Z}^{\infty}_{\iota}[\lambda]$ with $\mathbb{Z}_{\iota}^{\infty}$ as a set.
\end{thm}
In \cite{Ka2}, we gave a procedure to compute ${\rm Im}(\Psi^{(\lambda)}_{\io})\subset \ZZ^{\ify}_{\iota}[\lambda]$ by slightly modifying the method of \cite{N99}. 
Let $\beta^{(\pm)}_r$ be $\beta^{(+)}_r=\beta_{r}$ and
\begin{equation*}
\beta^{(-)}_r = 
\begin{cases}
x_{r^{(-)}}+\sum_{r^{(-)}<j<r} \lan h_{i_r},\alpha_{i_j}\ran x_j +x_r=\sigma_{r^{(-)}}-\sigma_r & {\rm if}\ r^{(-)}>0,\\
- \lan h_{i_r},\lambda\ran + \sum_{1\leq j<r} \lan h_{i_r},\alpha_{i_j}\ran x_j + x_r
=\sigma_0^{(i_r)}-\sigma_r
 & {\rm if}\ r^{(-)}=0.
\end{cases}
\end{equation*}
Note that if $r^{(-)}>0$ then $\beta^{(-)}_r$ coincides with $\beta_{r^{(-)}}$ defined in (\ref{betak}).
For $r\in\mathbb{Z}_{\geq1}$
and $\varphi=c+\sum_{r\geq1} \varphi_rx_r$ with $c$, $\varphi_r\in\mathbb{Q}$,
a new function $\what{S}_r'(\varphi)$ is defined as follows:
\begin{equation}\label{Shat}
\what{S}_r'(\varphi):=
\begin{cases}
\varphi - \beta^{(+)}_r & {\rm if}\ \varphi_r>0, \\
\varphi + \beta^{(-)}_r & {\rm if}\ \varphi_r<0, \\
\varphi & {\rm if}\ \varphi_r=0.
\end{cases}
\end{equation}
It is easy to see that except for the case $\varphi_r<0$ and $r^{(-)}=0$,
it follows 
\begin{equation}\label{twoS}
\what{S}_r'(\varphi)=S_r'(\varphi-c)+c,
\end{equation}
where $S_r'$ is defined in (\ref{Sk}).
For $k\in I$ and $\lambda\in P^+$, we define
\begin{equation}\label{lmi-def}
\iota^{(k)}:={\rm min}\{r\in\mathbb{Z}_{\geq1} | i_r=k\},\qquad
\lambda^{(k)} :=\lan h_k, \lambda \ran - \sum_{1\leq j<\iota^{(k)}} \lan h_k, \alpha_{i_j} \ran x_j - x_{\iota^{(k)}}.
\end{equation}
Let $\Xi'_{\iota}[\lambda]$ be the set of all
functions generated from the functions $x_j$ ($j \geq 1$) and $\lambda^{(k)}$
($k\in I$) by applying $\what{S}'_r$, that is,
\begin{equation}\label{xilamdef}
\begin{array}{l}
\Xi'_{\iota}[\lambda]: =\{\what{S}'_{j_{\ell}}\cdots \what{S}'_{j_1}x_{j_0} | \ell\in\mathbb{Z}_{\geq0},\ j_0,\cdots,j_{\ell}\in \mathbb{Z}_{\geq1} \} \\
\qquad\qquad \cup \{\what{S}'_{j_{\ell}}\cdots \what{S}'_{j_1} \lambda^{(k)} | \ell\in\mathbb{Z}_{\geq0},\ k\in I,\ j_1,\cdots,j_{\ell}\in \mathbb{Z}_{\geq1} \}.
\end{array}
\end{equation}

We also set
\begin{equation}\label{siglamdef}
\Sigma_{\iota}'[\lambda]:=
\{
\mathbf{a}\in\mathbb{Z}^{\infty}_{\iota}[\lambda] | \varphi(\mathbf{a}) \geq0\ \ {\rm for\ any\ }\varphi\in \Xi'_{\iota}[\lambda]
\}
\end{equation}
and consider the following condition:
\begin{defn}\cite{Ka2}\label{ample1}
We say the pair $(\iota,\lambda)$ is $\Xi'$-{\it ample}
if $\textbf{0}:=(\cdots,0,0,0)\in \Sigma'_{\iota}[\lambda]$.
\end{defn}

\begin{thm}\cite{Ka2}\label{Nthm1}
For $\lambda\in P^+$,
we assume the pair $(\iota,\lambda)$ is $\Xi'$-ample. 
Then
the image ${\rm Im}(\Psi^{(\lambda)}_{\iota}) (\cong B(\lambda))$ coincides with
$\Sigma'_{\iota}[\lambda]$.
\end{thm}

\subsection{An application for $\varepsilon_k^*$}

Using the operators $S'_{j}$, one can compute the function $\varepsilon_k^*$
in Example \ref{star-ex}.
For $k\in I$, we define
\[
\xi^{(k)}:=- \sum_{1\leq j<\iota^{(k)}} \lan h_k, \alpha_{i_j} \ran x_j - x_{\iota^{(k)}}\in (\mathbb{Q}^{\infty})^*.
\]
Note that for any $\lambda\in P^+$, it is related with $\lambda^{(k)}$ in (\ref{lmi-def}) by
$\xi^{(k)}=\lambda^{(k)}-\lan h_k, \lambda \ran$. Using the maps in (\ref{Sk}), we define
\[
\Xi'^{(k)}_{\iota}:=\{S'_{j_{\ell}}\cdots S'_{j_1}\xi^{(k)} | \ell\in\mathbb{Z}_{\geq0},\ j_1,\cdots,j_{\ell}\in\mathbb{Z}_{\geq1}\}.
\]
As defined in \cite{Ka2},
we say the sequence $\iota$ satisfies the $\Xi'$-{\it strict positivity condition} if it holds the following:
\begin{equation}\label{strict-cond}
\text{for any }
\vp=\sum_{\ell\geq1} \varphi_{\ell}x_{\ell}\in \left(\bigcup_{k\in I}\Xi'^{(k)}_{\iota}\setminus\{\xi^{(k)}\}\right)\cup\Xi'_{\iota},
\text{ if }\ell^{(-)}=0\text{ then }\varphi_{\ell}\geq0,
\end{equation}
where $\Xi'_{\iota}$ is defined in (\ref{xiprime}).
\begin{thm}\cite{Ka2}\label{thm2}
If $\iota$ satisfies the $\Xi'$-strict positivity condition then for $k\in I$ and $x\in {\rm Im}(\Psi_{\iota})$, it holds
\[
\varepsilon^*_k(x)={\rm max}\{-\varphi(x)| \varphi\in \Xi'^{(k)}_{\iota}\}.
\]
\end{thm}

\section{Young walls}

In this section, let $X=A^{(1)}_{n-1}$, $B^{(1)}_{n-1}$, $C^{(1)}_{n-1}$, $D^{(1)}_{n-1}$, $A^{(2)\dagger}_{2n-2}$, $A^{(2)}_{2n-3}$ or
$D^{(2)}_{n}$ and $k\in I$.

\begin{rem}\label{rem1}
For the sake of simplicity, we often omit the rank of type. For example,
we often simply write $A^{(1)}_{n-1}$ as $A^{(1)}$.
\end{rem}

\subsection{Periodic maps}\label{periodic map}
\begin{defn}
\begin{enumerate}
\item[$(1)$] We define a map $\{1,2,\cdots,n\}\rightarrow \{1,2,\cdots,n\}$ as $\ell \mapsto \ell$ for all $\ell\in \{1,2,\cdots,n\}$
and extend it to the map $\pi'_{A^{(1)}}:\mathbb{Z}_{\geq1}\rightarrow \{1,2,\cdots,n\}=I$ by periodicity $n$. Furthermore,
we extend $\pi_{A^{(1)}}'$ to $\pi_{A^{(1)}}:\mathbb{Z}\rightarrow I$ by periodicity $n$.
\item[$(2)$]
We define a map
$\{1,2,\cdots,2n-2\}\rightarrow \{1,2,\cdots,n\}$
as
\[
\ell\mapsto \ell,\ 2n-\ell\mapsto \ell \qquad (2\leq \ell\leq n-1),
\]
\[
1\mapsto1,\ n\mapsto n,
\]
and extend it to the map 
\begin{equation}\label{prime1}
\pi_{C^{(1)}}':\mathbb{Z}_{\geq1}\rightarrow\{1,2,\cdots,n\}=I
\end{equation}
by periodicity $2n-2$. We also define maps $\pi_{D^{(2)}}'$ and $\pi_{A^{(2)\dagger}_{2n-2}}'$ as
\[
\pi_{D^{(2)}}':=\pi_{C^{(1)}}',\quad \pi_{A^{(2)\dagger}_{2n-2}}':=\pi_{C^{(1)}}'.
\]
\item[$(3)$]
We define a map $\{1,\frac{3}{2},2,3,\cdots,2n-4\}\rightarrow I$ as
\[
\ell\mapsto\ell+1\quad (2\leq \ell\leq n-1),\quad 2n-\ell -1 \mapsto \ell\quad (3\leq \ell\leq n-1),
\]
\[
1\mapsto 1,\ \frac{3}{2}\mapsto2
\]
and extend it to the map 
\begin{equation}\label{prime2}
\pi_{B^{(1)}}':\mathbb{Z}_{\geq1}\cup\left\{\frac{3}{2}+(2n-4)z|z\in\mathbb{Z}_{\geq0}\right\}
\rightarrow I
\end{equation}
by periodicity $2n-4$. We also define $\pi_{A^{(2)}_{2n-3}}'$ as $\pi_{A^{(2)}_{2n-3}}':=\pi_{B^{(1)}}'$.
\item[$(4)$] We define a map $\{1,\frac{3}{2},2,3,\cdots,n-2,n-\frac{3}{2},n-1,n,\cdots, 2n-6\}\rightarrow I$ as
\[
\ell\mapsto\ell+1\quad (2\leq \ell\leq n-2),\quad 2n-\ell -3 \mapsto \ell\quad (3\leq \ell\leq n-2),
\]
\[
1\mapsto 1,\ \frac{3}{2}\mapsto2,\ n-\frac{3}{2}\mapsto n
\]
and extend it to the map 
\begin{equation}\label{prime3}
\pi_{D^{(1)}}':\mathbb{Z}_{\geq1}\cup\left\{\frac{3}{2}+(2n-6)z|z\in\mathbb{Z}_{\geq0}\right\}
\cup\left\{n-\frac{3}{2}+(2n-6)z|z\in\mathbb{Z}_{\geq0}\right\}
\rightarrow I
\end{equation}
by periodicity $2n-6$. 
\end{enumerate}
\end{defn}

\begin{defn}\label{classdef}
For $k\in I$ and $X$ other than $C^{(1)}_{n-1}$, we say $k$ is in class $1$ if the fundamental weight $\Lambda_k$ is level $1$.
We say $k\in I$ is in class $2$ if $k$ is not in class $1$. When $X=C^{(1)}_{n-1}$, we understand all $k\in I$ are in class $2$. 
The list of $k\in I$ in class $1$ is as follows:
\begin{table}[H]
  \begin{tabular}{|c|c|c|c|c|c|c|} \hline
  Type of $\mathfrak{g}$ & $A^{(1)}_{n-1}$ & $B^{(1)}_{n-1}$ & $D^{(1)}_{n-1}$ & $A^{(2)\dagger}_{2n-2}$ & $A^{(2)}_{2n-3}$ & $D^{(2)}_{n}$ \\ \hline
  $k$ & $1,2,\cdots,n$ & $1,2,n$ & $1,2,n-1,n$ & $1$ & $1,2$ & $1,n$  \\ \hline
  \end{tabular}
\end{table}

\end{defn}

Let $D_X$ denote the domain of definition of $\pi'_X$. For example,
\[
D_{B^{(1)}}=\mathbb{Z}_{\geq1}\cup\left\{\frac{3}{2}+(2n-4)z|z\in\mathbb{Z}_{\geq0}\right\}.
\]
For $k\in I$ in class $1$,
let $\ell\in D_X$ be the smallest number such that $\pi'_X(\ell)= k$ and we define $T^X_k:=\lfloor \ell \rfloor$,
where $\lfloor\ \rfloor$ is the floor function.
For $k\in I$,
let $\overline{T}^X_k$ (resp. $\overline{\overline{T}}^X_k$) be the smallest number (resp. second smallest number) $\ell\in D_X$ such that
$k=\pi'_X(\ell)$. Here, we often omit the type $X$ as $T_k$, $\overline{T}_k$, 
$\overline{\overline{T}}_k$.

\subsection{Proper Young walls}\label{pyws}

Let us review the notion of Young walls introduced in \cite{Kang}.
Each Young wall consists of $I$-colored blocks of three different
shapes:
\begin{enumerate}
\item[(1)] block with unit width, unit height and unit thickness:
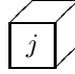
\[
\begin{xy}
(3,3) *{j}="0",
(0,0) *{}="1",
(6,0)*{}="2",
(6,6)*{}="3",
(0,6)*{}="4",
(3,9)*{}="5",
(9,9)*{}="6",
(9,3)*{}="7",
\ar@{-} "1";"2"^{}
\ar@{-} "1";"4"^{}
\ar@{-} "2";"3"^{}
\ar@{-} "3";"4"^{}
\ar@{-} "5";"4"^{}
\ar@{-} "5";"6"^{}
\ar@{-} "3";"6"^{}
\ar@{-} "2";"7"^{}
\ar@{-} "6";"7"^{}
\end{xy}
\]
\item[(2)] block with unit width, unit height and half-unit thickness:
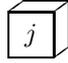
\[
\begin{xy}
(3,3) *{j}="0",
(0,0) *{}="1",
(6,0)*{}="2",
(6,6)*{}="3",
(0,6)*{}="4",
(2,7.5)*{}="5",
(8,7.5)*{}="6",
(8,1.5)*{}="7",
\ar@{-} "1";"2"^{}
\ar@{-} "1";"4"^{}
\ar@{-} "2";"3"^{}
\ar@{-} "3";"4"^{}
\ar@{-} "5";"4"^{}
\ar@{-} "5";"6"^{}
\ar@{-} "3";"6"^{}
\ar@{-} "2";"7"^{}
\ar@{-} "6";"7"^{}
\end{xy}
\]
\item[(3)] block with unit width, half-unit height and unit thickness:
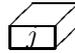
\[
\begin{xy}
(3,1.5) *{j}="0",
(0,0) *{}="1",
(6,0)*{}="2",
(6,3)*{}="3",
(0,3)*{}="4",
(3,6)*{}="5",
(9,6)*{}="6",
(9,3)*{}="7",
\ar@{-} "1";"2"^{}
\ar@{-} "1";"4"^{}
\ar@{-} "2";"3"^{}
\ar@{-} "3";"4"^{}
\ar@{-} "5";"4"^{}
\ar@{-} "5";"6"^{}
\ar@{-} "3";"6"^{}
\ar@{-} "2";"7"^{}
\ar@{-} "6";"7"^{}
\end{xy}
\]
\end{enumerate}
We simply call (1), (2), and (3) a {\it unit cube}, {\it a block with half-unit thickness} and {\it a block with half-unit height}, respectively.
Using these blocks,
we will build a wall whose thickness is less than or
equal to one unit.
Just as in \cite{Kang},
we simply express blocks (1) with color $j\in I$
as
\begin{equation}\label{smpl1}
\begin{xy}
(3,3) *{j}="0",
(0,0) *{}="1",
(6,0)*{}="2",
(6,6)*{}="3",
(0,6)*{}="4",
\ar@{-} "1";"2"^{}
\ar@{-} "1";"4"^{}
\ar@{-} "2";"3"^{}
\ar@{-} "3";"4"^{}
\end{xy}
\end{equation}
blocks (2) with color $j\in I$ as
\begin{equation}\label{smpl12}
\begin{xy}
(1.5,4.5) *{j}="0",
(0,0) *{}="1",
(6,0)*{}="2",
(6,6)*{}="3",
(0,6)*{}="4",
\ar@{-} "1";"4"^{}
\ar@{-} "1";"3"^{}
\ar@{-} "3";"4"^{}
\end{xy}
\end{equation}
or
\begin{equation}\label{smpl11}
\begin{xy}
(4.5,1.5) *{j}="0",
(0,0) *{}="1",
(6,0)*{}="2",
(6,6)*{}="3",
(0,6)*{}="4",
\ar@{-} "1";"2"^{}
\ar@{-} "1";"3"^{}
\ar@{-} "2";"3"^{}
\end{xy}
\end{equation}
and (3) with color $j\in I$ as
\begin{equation}\label{smpl2}
\begin{xy}
(1.5,1.5) *{\ \ j}="0",
(0,0) *{}="1",
(6,0)*{}="2",
(6,3.5)*{}="3",
(0,3.5)*{}="4",
\ar@{-} "1";"2"^{}
\ar@{-} "1";"4"^{}
\ar@{-} "2";"3"^{}
\ar@{-} "3";"4"^{}
\end{xy}
\end{equation}
The blocks (\ref{smpl1}), (\ref{smpl12}), (\ref{smpl11}) and (\ref{smpl2}) are called $j$-blocks.
If colored blocks are stacked as follows 
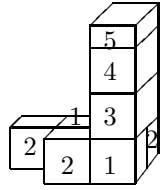
\[
\begin{xy}
(-9,2.5) *{\ \ 2}="02half2",
(-3.5,0) *{\ 2}="00",
(1.5,0) *{\ \ 1}="0",
(7,3.5) *{\ \ 2}="02half",
(-3,6.5) *{\ \ 1}="01half",
(1.5,6.5) *{\ \ 3}="02",
(1.5,12.5) *{\ \ 4}="03",
(1.5,16.5) *{\ \ 5}="04",
(7.5,-0.5)*{}="2-aa",
(7.5,5)*{}="2-aaa",
(9,1.5)*{}="2-a",
(9,6.5)*{}="3-a",
(9,12.5)*{}="3-1-a",
(9,18.5)*{}="3-2-a",
(9,21.5)*{}="3-3-a",
(3,21.5)*{}="4-3-a",
(-3,12.5)*{}="5-1-a",
(0,12.5)*{}="5-1-ab",
(0,5)*{}="4-ab",
(-10.5,5)*{}="4-abc",
(-9,6.5)*{}="7-a",
(0,6.5)*{}="7-ab",
(0,-2.5) *{}="1",
(6,-2.5)*{}="2",
(6,3.5)*{}="3",
(0,3.5)*{}="4",
(-6,3.5)*{}="5",
(-3,6.5)*{}="5-a",
(-6,-2.5)*{}="6",
(-10.5,5)*{}="7",
(-10.5,-0.5)*{}="7-v",
(-6,-0.5)*{}="6-v",
(-12,-2.5)*{}="8",
(6,9.5)*{}="3-1",
(6,15.5)*{}="3-2",
(6,18.5)*{}="3-3",
(0,9.5)*{}="4-1",
(0,15.5)*{}="4-2",
(0,18.5)*{}="4-3",
(-6,9.5)*{}="5-1",
\ar@{-} "6-v";"7-v"^{}
\ar@{-} "7";"7-v"^{}
\ar@{-} "4-ab";"4-abc"^{}
\ar@{-} "5";"5-a"^{}
\ar@{-} "1";"2"^{}
\ar@{-} "1";"4"^{}
\ar@{-} "2";"3"^{}
\ar@{-} "3";"4"^{}
\ar@{-} "5";"6"^{}
\ar@{-} "5";"4"^{}
\ar@{-} "1";"6"^{}
\ar@{-} "3";"3-1"^{}
\ar@{-} "4-1";"3-1"^{}
\ar@{-} "4-1";"4"^{}
\ar@{-} "4-1";"4-2"^{}
\ar@{-} "3-2";"3-1"^{}
\ar@{-} "3-2";"4-2"^{}
\ar@{-} "3-2";"3-3"^{}
\ar@{-} "4-3";"4-2"^{}
\ar@{-} "4-3";"3-3"^{}
\ar@{-} "2";"2-a"^{}
\ar@{-} "3";"3-a"^{}
\ar@{-} "3-1";"3-1-a"^{}
\ar@{-} "3-2";"3-2-a"^{}
\ar@{-} "3-3";"3-3-a"^{}
\ar@{-} "2-a";"3-3-a"^{}
\ar@{-} "2-aa";"2-aaa"^{}
\ar@{-} "4-3";"4-3-a"^{}
\ar@{-} "3-3-a";"4-3-a"^{}
\ar@{-} "7";"7-a"^{}
\ar@{-} "7-ab";"7-a"^{}
\end{xy}
\]
it is simply described as
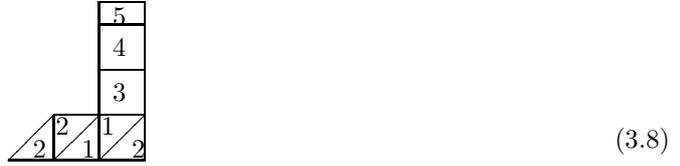
\begin{equation}\label{smpl3}
\begin{xy}
(0,2) *{\ \ 1}="1half",
(4,-1) *{\ \ 2}="2half",
(-2,-1) *{\ 1}="1half2",
(-6,2) *{\ \ 2}="2hal2",
(-9,-1) *{\ \ 2}="2hal3",
(1.5,6.5) *{\ \ 3}="02",
(1.5,12.5) *{\ \ 4}="03",
(1.5,16.5) *{\ \ 5}="04",
(0,-2.5) *{}="1",
(6,-2.5)*{}="2",
(6,3.5)*{}="3",
(0,3.5)*{}="4",
(-6,3.5)*{}="5",
(-6,-2.5)*{}="6",
(-12,-2.5)*{}="8",
(6,9.5)*{}="3-1",
(6,15.5)*{}="3-2",
(6,18.5)*{}="3-3",
(0,9.5)*{}="4-1",
(0,15.5)*{}="4-2",
(0,18.5)*{}="4-3",
(-6,9.5)*{}="5-1",
\ar@{-} "5";"8"^{}
\ar@{-} "1";"3"^{}
\ar@{-} "4";"6"^{}
\ar@{-} "1";"2"^{}
\ar@{-} "1";"4"^{}
\ar@{-} "2";"3"^{}
\ar@{-} "3";"4"^{}
\ar@{-} "5";"6"^{}
\ar@{-} "5";"4"^{}
\ar@{-} "1";"6"^{}
\ar@{-} "6";"8"^{}
\ar@{-} "3";"3-1"^{}
\ar@{-} "4-1";"3-1"^{}
\ar@{-} "4-1";"4"^{}
\ar@{-} "4-1";"4-2"^{}
\ar@{-} "3-2";"3-1"^{}
\ar@{-} "3-2";"4-2"^{}
\ar@{-} "3-2";"3-3"^{}
\ar@{-} "4-3";"4-2"^{}
\ar@{-} "4-3";"3-3"^{}
\end{xy}
\end{equation}
In this way, a block (2) in front is expressed by (\ref{smpl12}) and a block (2) at the back
is expressed by $(\ref{smpl11})$.

Let us recall {\it ground state walls} $Y_{\Lambda_k}$ of type $X$ for $k\in I$ in class $1$ \cite{Kang}:
\begin{enumerate}
\item[(1)] For $X=A^{(1)}_{n-1}$ and $\lambda=\Lambda_k$ ($k\in I$), let $Y_{\lambda}$ be the wall which has no blocks.
\item[(2)]
For $X=A^{(2)\dagger}_{2n-2}$ and $\lambda=\Lambda_1$, one defines
\[
Y_{\Lambda_1}:=
\begin{xy}
(-15.5,1.5) *{\ \cdots}="0000",
(-9.5,1.5) *{\ 1}="000",
(-3.5,1.5) *{\ 1}="00",
(1.5,1.5) *{\ \ 1}="0",
(0,0) *{}="1",
(6,0)*{}="2",
(6,3.5)*{}="3",
(0,3.5)*{}="4",
(-6,3.5)*{}="5",
(-6,0)*{}="6",
(-12,3.5)*{}="7",
(-12,0)*{}="8",
\ar@{-} "1";"2"^{}
\ar@{-} "1";"4"^{}
\ar@{-} "2";"3"^{}
\ar@{-} "3";"4"^{}
\ar@{-} "5";"6"^{}
\ar@{-} "5";"4"^{}
\ar@{-} "1";"6"^{}
\ar@{-} "7";"8"^{}
\ar@{-} "7";"5"^{}
\ar@{-} "6";"8"^{}
\end{xy}
\]
Here, the wall $Y_{\Lambda_1}$ has infinitely many $1$-blocks with half-unit height
and
extends infinitely to the left.
\item[(3)]
For $X=D^{(2)}_{n}$ and $\lambda=\Lambda_1$, $\Lambda_n$, one defines
\[
Y_{\Lambda_1}:=
\begin{xy}
(-15.5,1.5) *{\ \cdots}="0000",
(-9.5,1.5) *{\ 1}="000",
(-3.5,1.5) *{\ 1}="00",
(1.5,1.5) *{\ \ 1}="0",
(0,0) *{}="1",
(6,0)*{}="2",
(6,3.5)*{}="3",
(0,3.5)*{}="4",
(-6,3.5)*{}="5",
(-6,0)*{}="6",
(-12,3.5)*{}="7",
(-12,0)*{}="8",
\ar@{-} "1";"2"^{}
\ar@{-} "1";"4"^{}
\ar@{-} "2";"3"^{}
\ar@{-} "3";"4"^{}
\ar@{-} "5";"6"^{}
\ar@{-} "5";"4"^{}
\ar@{-} "1";"6"^{}
\ar@{-} "7";"8"^{}
\ar@{-} "7";"5"^{}
\ar@{-} "6";"8"^{}
\end{xy}
\]
and
\[
Y_{\Lambda_n}:=
\begin{xy}
(-15.5,1.5) *{\ \cdots}="0000",
(-9.5,1.5) *{\ n}="000",
(-3.5,1.5) *{\ n}="00",
(1.5,1.5) *{\ \ n}="0",
(0,0) *{}="1",
(6,0)*{}="2",
(6,3.5)*{}="3",
(0,3.5)*{}="4",
(-6,3.5)*{}="5",
(-6,0)*{}="6",
(-12,3.5)*{}="7",
(-12,0)*{}="8",
\ar@{-} "1";"2"^{}
\ar@{-} "1";"4"^{}
\ar@{-} "2";"3"^{}
\ar@{-} "3";"4"^{}
\ar@{-} "5";"6"^{}
\ar@{-} "5";"4"^{}
\ar@{-} "1";"6"^{}
\ar@{-} "7";"8"^{}
\ar@{-} "7";"5"^{}
\ar@{-} "6";"8"^{}
\end{xy}
\]
\item[(4)] For $X=A^{(2)}_{2n-3}$ or $B^{(1)}_{n-1}$ and $\lambda=\Lambda_1$, $\Lambda_2$,
\[
Y_{\Lambda_1}:=
\begin{xy}
(3,-1) *{\ \ 2}="2half",
(-2,-1) *{\ 1}="1half2",
(-9,-1) *{\ \ 2}="2hal3",
(-15,-1) *{\ \ 1}="1hal3",
(-22,-1) *{\ \ \cdots}="dot",
(0,-2.5) *{}="1",
(6,-2.5)*{}="2",
(6,3.5)*{}="3",
(0,3.5)*{}="4",
(-6,3.5)*{}="5",
(-6,-2.5)*{}="6",
(-12,-2.5)*{}="8",
(-12,3.5)*{}="8-1",
(-18,-2.5)*{}="9",
\ar@{-} "9";"8-1"^{}
\ar@{-} "8-1";"8"^{}
\ar@{-} "9";"8"^{}
\ar@{-} "5";"8"^{}
\ar@{-} "1";"3"^{}
\ar@{-} "4";"6"^{}
\ar@{-} "1";"2"^{}
\ar@{-} "1";"4"^{}
\ar@{-} "2";"3"^{}
\ar@{-} "5";"6"^{}
\ar@{-} "1";"6"^{}
\ar@{-} "6";"8"^{}
\end{xy}\qquad
Y_{\Lambda_2}:=
\begin{xy}
(3,-1) *{\ \ 1}="2half",
(-2,-1) *{\ 2}="1half2",
(-9,-1) *{\ \ 1}="2hal3",
(-15,-1) *{\ \ 2}="1hal3",
(-22,-1) *{\ \ \cdots}="dot",
(0,-2.5) *{}="1",
(6,-2.5)*{}="2",
(6,3.5)*{}="3",
(0,3.5)*{}="4",
(-6,3.5)*{}="5",
(-6,-2.5)*{}="6",
(-12,-2.5)*{}="8",
(-12,3.5)*{}="8-1",
(-18,-2.5)*{}="9",
\ar@{-} "9";"8-1"^{}
\ar@{-} "8-1";"8"^{}
\ar@{-} "9";"8"^{}
\ar@{-} "5";"8"^{}
\ar@{-} "1";"3"^{}
\ar@{-} "4";"6"^{}
\ar@{-} "1";"2"^{}
\ar@{-} "1";"4"^{}
\ar@{-} "2";"3"^{}
\ar@{-} "5";"6"^{}
\ar@{-} "1";"6"^{}
\ar@{-} "6";"8"^{}
\end{xy}
\]
Here, the walls $Y_{\Lambda_1}$, $Y_{\Lambda_2}$ extend infinitely to the left. For
$X=B^{(1)}_{n-1}$ and $\lambda=\Lambda_n$,
we set
\[
Y_{\Lambda_n}:=
\begin{xy}
(-15.5,1.5) *{\ \cdots}="0000",
(-9.5,1.5) *{\ n}="000",
(-3.5,1.5) *{\ n}="00",
(1.5,1.5) *{\ \ n}="0",
(0,0) *{}="1",
(6,0)*{}="2",
(6,3.5)*{}="3",
(0,3.5)*{}="4",
(-6,3.5)*{}="5",
(-6,0)*{}="6",
(-12,3.5)*{}="7",
(-12,0)*{}="8",
\ar@{-} "1";"2"^{}
\ar@{-} "1";"4"^{}
\ar@{-} "2";"3"^{}
\ar@{-} "3";"4"^{}
\ar@{-} "5";"6"^{}
\ar@{-} "5";"4"^{}
\ar@{-} "1";"6"^{}
\ar@{-} "7";"8"^{}
\ar@{-} "7";"5"^{}
\ar@{-} "6";"8"^{}
\end{xy}
\]

\item[(5)] For $X=D^{(1)}_{n-1}$ and $\lambda=\Lambda_1$, $\Lambda_2$, $\Lambda_{n-1}$ or $\Lambda_n$,
\[
Y_{\Lambda_1}:=
\begin{xy}
(3,-1) *{\ \ 2}="2half",
(-2,-1) *{\ 1}="1half2",
(-9,-1) *{\ \ 2}="2hal3",
(-15,-1) *{\ \ 1}="1hal3",
(-22,-1) *{\ \ \cdots}="dot",
(0,-2.5) *{}="1",
(6,-2.5)*{}="2",
(6,3.5)*{}="3",
(0,3.5)*{}="4",
(-6,3.5)*{}="5",
(-6,-2.5)*{}="6",
(-12,-2.5)*{}="8",
(-12,3.5)*{}="8-1",
(-18,-2.5)*{}="9",
\ar@{-} "9";"8-1"^{}
\ar@{-} "8-1";"8"^{}
\ar@{-} "9";"8"^{}
\ar@{-} "5";"8"^{}
\ar@{-} "1";"3"^{}
\ar@{-} "4";"6"^{}
\ar@{-} "1";"2"^{}
\ar@{-} "1";"4"^{}
\ar@{-} "2";"3"^{}
\ar@{-} "5";"6"^{}
\ar@{-} "1";"6"^{}
\ar@{-} "6";"8"^{}
\end{xy}\qquad
Y_{\Lambda_2}:=
\begin{xy}
(3,-1) *{\ \ 1}="2half",
(-2,-1) *{\ 2}="1half2",
(-9,-1) *{\ \ 1}="2hal3",
(-15,-1) *{\ \ 2}="1hal3",
(-22,-1) *{\ \ \cdots}="dot",
(0,-2.5) *{}="1",
(6,-2.5)*{}="2",
(6,3.5)*{}="3",
(0,3.5)*{}="4",
(-6,3.5)*{}="5",
(-6,-2.5)*{}="6",
(-12,-2.5)*{}="8",
(-12,3.5)*{}="8-1",
(-18,-2.5)*{}="9",
\ar@{-} "9";"8-1"^{}
\ar@{-} "8-1";"8"^{}
\ar@{-} "9";"8"^{}
\ar@{-} "5";"8"^{}
\ar@{-} "1";"3"^{}
\ar@{-} "4";"6"^{}
\ar@{-} "1";"2"^{}
\ar@{-} "1";"4"^{}
\ar@{-} "2";"3"^{}
\ar@{-} "5";"6"^{}
\ar@{-} "1";"6"^{}
\ar@{-} "6";"8"^{}
\end{xy}
\]
\[
Y_{\Lambda_{n-1}}:=
\begin{xy}
(3,-1) *{\ \ n}="2half",
(-2,-1) *{_{n-1}}="1half2",
(-9,-1) *{\ \ n}="2hal3",
(-15,-1) *{\ _{n-1}}="1hal3",
(-22,-1) *{\ \ \cdots}="dot",
(0,-2.5) *{}="1",
(6,-2.5)*{}="2",
(6,3.5)*{}="3",
(0,3.5)*{}="4",
(-6,3.5)*{}="5",
(-6,-2.5)*{}="6",
(-12,-2.5)*{}="8",
(-12,3.5)*{}="8-1",
(-18,-2.5)*{}="9",
\ar@{-} "9";"8-1"^{}
\ar@{-} "8-1";"8"^{}
\ar@{-} "9";"8"^{}
\ar@{-} "5";"8"^{}
\ar@{-} "1";"3"^{}
\ar@{-} "4";"6"^{}
\ar@{-} "1";"2"^{}
\ar@{-} "1";"4"^{}
\ar@{-} "2";"3"^{}
\ar@{-} "5";"6"^{}
\ar@{-} "1";"6"^{}
\ar@{-} "6";"8"^{}
\end{xy}\qquad
Y_{\Lambda_n}:=
\begin{xy}
(3,-1) *{\ \ _{n-1}}="2half",
(-2,-1) *{\ n}="1half2",
(-9,-1) *{\ \ _{n-1}}="2hal3",
(-15,-1) *{\ \ n}="1hal3",
(-22,-1) *{\ \ \cdots}="dot",
(0,-2.5) *{}="1",
(6,-2.5)*{}="2",
(6,3.5)*{}="3",
(0,3.5)*{}="4",
(-6,3.5)*{}="5",
(-6,-2.5)*{}="6",
(-12,-2.5)*{}="8",
(-12,3.5)*{}="8-1",
(-18,-2.5)*{}="9",
\ar@{-} "9";"8-1"^{}
\ar@{-} "8-1";"8"^{}
\ar@{-} "9";"8"^{}
\ar@{-} "5";"8"^{}
\ar@{-} "1";"3"^{}
\ar@{-} "4";"6"^{}
\ar@{-} "1";"2"^{}
\ar@{-} "1";"4"^{}
\ar@{-} "2";"3"^{}
\ar@{-} "5";"6"^{}
\ar@{-} "1";"6"^{}
\ar@{-} "6";"8"^{}
\end{xy}
\]
\end{enumerate}

In \cite{Kang-Lee}, the following stacking patterns are introduced:
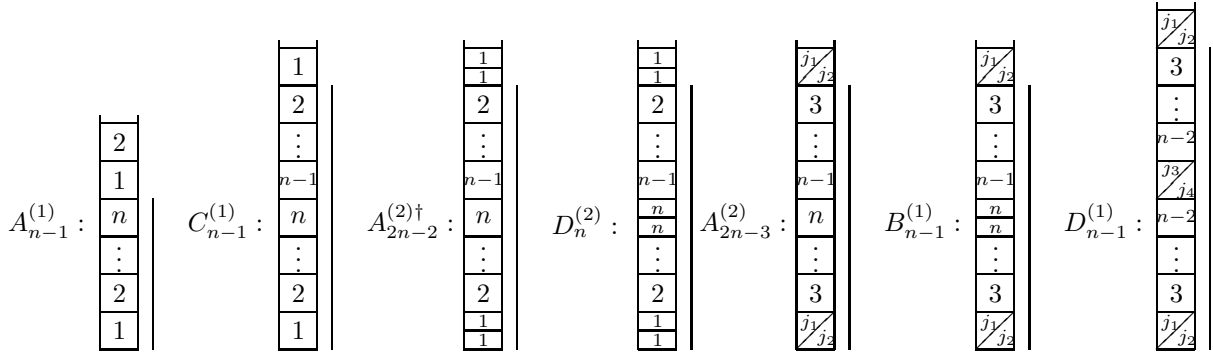
\[
\begin{xy}
(-12,17) *{A^{(1)}_{n-1}:}="0000",
(-3,27.5) *{\ 2}="c6",
(-3,22.5) *{\ 1}="c5",
(-3,17.5) *{\ n}="c4",
(-3,13) *{\ \vdots}="c3",
(-3,7.5) *{\ 2}="c2",
(-3,2.5) *{\ 1}="c1",
(0,0) *{}="1",
(2,0) *{}="1b",
(0,5)*{}="2",
(0,10)*{}="3",
(0,15)*{}="4",
(0,20) *{}="5",
(0,25)*{}="6",
(0,30)*{}="7",
(-5,0) *{}="1a",
(-5,5)*{}="2a",
(-5,10)*{}="3a",
(-5,15)*{}="4a",
(-5,20) *{}="5a",
(-5,25)*{}="6a",
(-5,30)*{}="7a",
(0,31)*{}="top",
(-5,31)*{}="topa",
(2,20)*{}="topb",
\ar@{-} "1";"1a"^{}
\ar@{-} "2";"2a"^{}
\ar@{-} "3";"3a"^{}
\ar@{-} "4";"4a"^{}
\ar@{-} "5";"5a"^{}
\ar@{-} "6";"6a"^{}
\ar@{-} "7";"7a"^{}
\ar@{-} "1";"top"^{}
\ar@{-} "1a";"topa"^{}
\ar@{-} "1b";"topb"^{}
\end{xy}\quad
\begin{xy}
(-12,17) *{C^{(1)}_{n-1}:}="0000",
(-3,37.5) *{\ 1}="c8",
(-3,32.5) *{\ 2}="c7",
(-3,28) *{\ \vdots}="c6",
(-3,22.5) *{\ \scriptstyle{n-1}}="c5",
(-3,17.5) *{\ n}="c4",
(-3,13) *{\ \vdots}="c3",
(-3,7.5) *{\ 2}="c2",
(-3,2.5) *{\ 1}="c1",
(-7,15) *{}="1c",
(2,0) *{}="1b",
(0,0) *{}="1",
(0,5)*{}="2",
(0,10)*{}="3",
(0,15)*{}="4",
(0,20) *{}="5",
(0,25)*{}="6",
(0,30)*{}="7",
(0,35)*{}="8",
(0,40)*{}="9",
(-5,0) *{}="1a",
(-5,5)*{}="2a",
(-5,10)*{}="3a",
(-5,15)*{}="4a",
(-5,17.5)*{}="4.5a",
(-5,20) *{}="5a",
(-5,25)*{}="6a",
(-5,30)*{}="7a",
(-5,35)*{}="8a",
(-5,37.3)*{}="8.5a",
(-5,40)*{}="9a",
(0,41)*{}="top",
(-5,41)*{}="topa",
(2,35)*{}="topb",
(-7,35)*{}="topc",
\ar@{-} "1";"1a"^{}
\ar@{-} "2";"2a"^{}
\ar@{-} "3";"3a"^{}
\ar@{-} "4";"4a"^{}
\ar@{-} "5";"5a"^{}
\ar@{-} "6";"6a"^{}
\ar@{-} "7";"7a"^{}
\ar@{-} "8";"8a"^{}
\ar@{-} "9";"9a"^{}
\ar@{-} "1";"top"^{}
\ar@{-} "1a";"topa"^{}
\ar@{-} "1b";"topb"^{}
\end{xy}\quad
\begin{xy}
(-12,17) *{A^{(2)\dagger}_{2n-2}:}="0000",
(-3,36.2) *{\ \scriptstyle{1}}="c8-1",
(-3,38.7) *{\ \scriptstyle{1}}="c8-2",
(-3,32.5) *{\ 2}="c7",
(-3,28) *{\ \vdots}="c6",
(-3,22.5) *{\ \scriptstyle{n-1}}="c5",
(-3,17.5) *{\ n}="c4",
(-3,13) *{\ \vdots}="c3",
(-3,7.5) *{\ 2}="c2",
(-3,1.3) *{\ \scriptstyle{1}}="c1-1",
(-3,3.8) *{\ \scriptstyle{1}}="c1-2",
(-7,15) *{}="1c",
(2,0) *{}="1b",
(0,0) *{}="1",
(0,2.5) *{}="1.5",
(0,5)*{}="2",
(0,10)*{}="3",
(0,15)*{}="4",
(0,17.5)*{}="4.5",
(0,20) *{}="5",
(0,25)*{}="6",
(0,30)*{}="7",
(0,35)*{}="8",
(0,37.3)*{}="8.5",
(0,40)*{}="9",
(-5,0) *{}="1a",
(-5,2.5) *{}="1.5a",
(-5,5)*{}="2a",
(-5,10)*{}="3a",
(-5,15)*{}="4a",
(-5,20) *{}="5a",
(-5,25)*{}="6a",
(-5,30)*{}="7a",
(-5,35)*{}="8a",
(-5,37.3)*{}="8.5a",
(-5,40)*{}="9a",
(0,41)*{}="top",
(-5,41)*{}="topa",
(2,35)*{}="topb",
(-7,37.5)*{}="topc",
\ar@{-} "1";"1a"^{}
\ar@{-} "1.5";"1.5a"^{}
\ar@{-} "2";"2a"^{}
\ar@{-} "3";"3a"^{}
\ar@{-} "4";"4a"^{}
\ar@{-} "5";"5a"^{}
\ar@{-} "6";"6a"^{}
\ar@{-} "7";"7a"^{}
\ar@{-} "8";"8a"^{}
\ar@{-} "8.5";"8.5a"^{}
\ar@{-} "9";"9a"^{}
\ar@{-} "1";"top"^{}
\ar@{-} "1a";"topa"^{}
\ar@{-} "1b";"topb"^{}
\end{xy}
\quad
\begin{xy}
(-12,17) *{D^{(2)}_{n}:}="0000",
(-3,36.2) *{\ \scriptstyle{1}}="c8-1",
(-3,38.7) *{\ \scriptstyle{1}}="c8-2",
(-3,32.5) *{\ 2}="c7",
(-3,28) *{\ \vdots}="c6",
(-3,22.5) *{\ \scriptstyle{n-1}}="c5",
(-3,16) *{\ \scriptstyle{n}}="c4-1",
(-3,18.5) *{\ \scriptstyle{n}}="c4-2",
(-3,13) *{\ \vdots}="c3",
(-3,7.5) *{\ 2}="c2",
(-3,1.3) *{\ \scriptstyle{1}}="c1-1",
(-3,3.8) *{\ \scriptstyle{1}}="c1-2",
(-7,17) *{}="1c",
(2,0) *{}="1b",
(0,0) *{}="1",
(0,2.5) *{}="1.5",
(0,5)*{}="2",
(0,10)*{}="3",
(0,15)*{}="4",
(0,17.5)*{}="4.5",
(0,20) *{}="5",
(0,25)*{}="6",
(0,30)*{}="7",
(0,35)*{}="8",
(0,37.3)*{}="8.5",
(0,40)*{}="9",
(-5,0) *{}="1a",
(-5,2.5) *{}="1.5a",
(-5,5)*{}="2a",
(-5,10)*{}="3a",
(-5,15)*{}="4a",
(-5,17.5)*{}="4.5a",
(-5,20) *{}="5a",
(-5,25)*{}="6a",
(-5,30)*{}="7a",
(-5,35)*{}="8a",
(-5,37.3)*{}="8.5a",
(-5,40)*{}="9a",
(0,41)*{}="top",
(-5,41)*{}="topa",
(2,35)*{}="topb",
(-7,37)*{}="topc",
\ar@{-} "1";"1a"^{}
\ar@{-} "1.5";"1.5a"^{}
\ar@{-} "2";"2a"^{}
\ar@{-} "3";"3a"^{}
\ar@{-} "4";"4a"^{}
\ar@{-} "4.5";"4.5a"^{}
\ar@{-} "5";"5a"^{}
\ar@{-} "6";"6a"^{}
\ar@{-} "7";"7a"^{}
\ar@{-} "8";"8a"^{}
\ar@{-} "8.5";"8.5a"^{}
\ar@{-} "9";"9a"^{}
\ar@{-} "1";"top"^{}
\ar@{-} "1a";"topa"^{}
\ar@{-} "1b";"topb"^{}
\end{xy}
\begin{xy}
(-12,17) *{A^{(2)}_{2n-3}:}="0000",
(-3.7,38.5) *{\ \scriptstyle{j_1}}="c8-1",
(-1.5,36.5) *{\ \scriptstyle{j_2}}="c8-2",
(-3,32.5) *{\ 3}="c7",
(-3,28) *{\ \vdots}="c6",
(-3,22.5) *{\ \scriptstyle{n-1}}="c5",
(-3,17.5) *{\ n}="c4",
(-3,13) *{\ \vdots}="c3",
(-3,7.5) *{\ 3}="c2",
(-3.7,3.5) *{\ \scriptstyle{j_1}}="c1-1",
(-1.5,1.5) *{\ \scriptstyle{j_2}}="c1-2",
(-7,15) *{}="1c",
(2,0) *{}="1b",
(0,0) *{}="1",
(0,5)*{}="2",
(0,10)*{}="3",
(0,15)*{}="4",
(0,20) *{}="5",
(0,25)*{}="6",
(0,30)*{}="7",
(0,35)*{}="8",
(0,40)*{}="9",
(-5,0) *{}="1a",
(-5,5)*{}="2a",
(-5,10)*{}="3a",
(-5,15)*{}="4a",
(-5,20) *{}="5a",
(-5,25)*{}="6a",
(-5,30)*{}="7a",
(-5,35)*{}="8a",
(-5,40)*{}="9a",
(0,41)*{}="top",
(-5,41)*{}="topa",
(2,35)*{}="topb",
(-7,35)*{}="topc",
\ar@{-} "1";"1a"^{}
\ar@{-} "2";"1a"^{}
\ar@{-} "2";"2a"^{}
\ar@{-} "3";"3a"^{}
\ar@{-} "4";"4a"^{}
\ar@{-} "5";"5a"^{}
\ar@{-} "6";"6a"^{}
\ar@{-} "7";"7a"^{}
\ar@{-} "8";"8a"^{}
\ar@{-} "9";"9a"^{}
\ar@{-} "9";"8a"^{}
\ar@{-} "1";"top"^{}
\ar@{-} "1a";"topa"^{}
\ar@{-} "1b";"topb"^{}
\end{xy}\quad
\begin{xy}
(-12,17) *{B^{(1)}_{n-1}:}="0000",
(-3.7,38.5) *{\ \scriptstyle{j_1}}="c8-1",
(-1.5,36.5) *{\ \scriptstyle{j_2}}="c8-2",
(-3,32.5) *{\ 3}="c7",
(-3,28) *{\ \vdots}="c6",
(-3,22.5) *{\ \scriptstyle{n-1}}="c5",
(-3,16) *{\ \scriptstyle{n}}="c4-1",
(-3,18.5) *{\ \scriptstyle{n}}="c4-2",
(-3,13) *{\ \vdots}="c3",
(-3,7.5) *{\ 3}="c2",
(-3.7,3.5) *{\ \scriptstyle{j_1}}="c1-1",
(-1.5,1.5) *{\ \scriptstyle{j_2}}="c1-2",
(-7,17) *{}="1c",
(2,0) *{}="1b",
(0,0) *{}="1",
(0,5)*{}="2",
(0,10)*{}="3",
(0,15)*{}="4",
(0,17.5)*{}="4.5",
(0,20) *{}="5",
(0,25)*{}="6",
(0,30)*{}="7",
(0,35)*{}="8",
(0,40)*{}="9",
(-5,0) *{}="1a",
(-5,5)*{}="2a",
(-5,10)*{}="3a",
(-5,15)*{}="4a",
(-5,17.5)*{}="4.5a",
(-5,20) *{}="5a",
(-5,25)*{}="6a",
(-5,30)*{}="7a",
(-5,35)*{}="8a",
(-5,40)*{}="9a",
(0,41)*{}="top",
(-5,41)*{}="topa",
(2,35)*{}="topb",
(-7,35)*{}="topc",
\ar@{-} "1";"1a"^{}
\ar@{-} "2";"2a"^{}
\ar@{-} "2";"1a"^{}
\ar@{-} "3";"3a"^{}
\ar@{-} "4";"4a"^{}
\ar@{-} "4.5";"4.5a"^{}
\ar@{-} "5";"5a"^{}
\ar@{-} "6";"6a"^{}
\ar@{-} "7";"7a"^{}
\ar@{-} "8";"8a"^{}
\ar@{-} "9";"9a"^{}
\ar@{-} "9";"8a"^{}
\ar@{-} "1";"top"^{}
\ar@{-} "1a";"topa"^{}
\ar@{-} "1b";"topb"^{}
\end{xy}\quad
\begin{xy}
(-12,17) *{D^{(1)}_{n-1}:}="0000",
(-3.7,43.5) *{\ \scriptstyle{j_1}}="c9-1",
(-1.5,41.5) *{\ \scriptstyle{j_2}}="c9-2",
(-3,37.5) *{\ 3}="c8",
(-3,32.9) *{\ \vdots}="c7",
(-3,28) *{\ \scriptstyle{n-2}}="c6",
(-3.6,24) *{\ \scriptstyle{j_3}}="c5-1",
(-1.6,21.1) *{\ \scriptstyle{j_4}}="c5-2",
(-3,17.5) *{\ \scriptstyle{n-2}}="c4",
(-3,13) *{\ \vdots}="c3",
(-3,7.5) *{\ 3}="c2",
(-3.7,3.5) *{\ \scriptstyle{j_1}}="c1-1",
(-1.5,1.5) *{\ \scriptstyle{j_2}}="c1-2",
(-7,20) *{}="1c",
(2,0) *{}="1b",
(0,0) *{}="1",
(0,5)*{}="2",
(0,10)*{}="3",
(0,15)*{}="4",
(0,20) *{}="5",
(0,25)*{}="6",
(0,30)*{}="7",
(0,35)*{}="8",
(0,40)*{}="9",
(0,45)*{}="10",
(-5,0) *{}="1a",
(-5,5)*{}="2a",
(-5,10)*{}="3a",
(-5,15)*{}="4a",
(-5,20) *{}="5a",
(-5,25)*{}="6a",
(-5,30)*{}="7a",
(-5,35)*{}="8a",
(-5,40)*{}="9a",
(-5,45)*{}="10a",
(0,46)*{}="top",
(-5,46)*{}="topa",
(2,40)*{}="topb",
(-7,40)*{}="topc",
\ar@{-} "1";"1a"^{}
\ar@{-} "2";"1a"^{}
\ar@{-} "2";"2a"^{}
\ar@{-} "3";"3a"^{}
\ar@{-} "4";"4a"^{}
\ar@{-} "5";"5a"^{}
\ar@{-} "6";"5a"^{}
\ar@{-} "6";"6a"^{}
\ar@{-} "7";"7a"^{}
\ar@{-} "8";"8a"^{}
\ar@{-} "9";"9a"^{}
\ar@{-} "10";"10a"^{}
\ar@{-} "10";"9a"^{}
\ar@{-} "1";"top"^{}
\ar@{-} "1a";"topa"^{}
\ar@{-} "1b";"topb"^{}
\end{xy}
\]
The stacking pattern with right-hand line drawn repeated infinitely many times.
Here, $\{j_1,j_2\}=\{1,2\}$ and $\{j_3,j_4\}=\{n-1,n\}$.
For $k\in I$ in class $1$,
let $Pat(X,\geq k)$ denote the pattern obtained from the above pattern for type $X$ by truncating all blocks below the first $k$-block from the bottom.
For $k\in I$ in class $2$,
let $\overline{Pat}(X,\geq k)$
(resp. $\overline{\overline{Pat}}(X,\geq k)$)
 denote the pattern obtained from the above pattern for type $X$ by truncating all blocks below the first (resp. second) $k$-block from the bottom
and changing the lowest $k$-block to two $k$-blocks with half-unit height:
\[
\begin{xy}
(6,9) *{k}="00",
(6,3) *{k}="0",
(0,0) *{}="1",
(12,0)*{}="2",
(12,12)*{}="3",
(0,12)*{}="4",
(0,6)*{}="5",
(12,6)*{}="6",
\ar@{-} "1";"2"^{}
\ar@{-} "1";"4"^{}
\ar@{-} "2";"3"^{}
\ar@{-} "3";"4"^{}
\ar@{-} "5";"6"^{}
\end{xy}
\]
For example,
\[
\begin{xy}
(-22,17) *{Pat(A^{(1)}_{n-1},\geq k)=}="0000",
(-3,27.5) *{\ 2}="c6",
(-3,22.5) *{\ 1}="c5",
(-3,17.5) *{\ n}="c4",
(-3,13) *{\ \vdots}="c3",
(-3,7.5) *{\ \scriptstyle{k+1}}="c2",
(-3,2.5) *{\ k}="c1",
(0,0) *{}="1",
(0,5)*{}="2",
(0,10)*{}="3",
(0,15)*{}="4",
(0,20) *{}="5",
(0,25)*{}="6",
(0,30)*{}="7",
(-5,0) *{}="1a",
(-5,5)*{}="2a",
(-5,10)*{}="3a",
(-5,15)*{}="4a",
(-5,20) *{}="5a",
(-5,25)*{}="6a",
(-5,30)*{}="7a",
(0,31)*{}="top",
(-5,31)*{}="topa",
\ar@{-} "1";"1a"^{}
\ar@{-} "2";"2a"^{}
\ar@{-} "3";"3a"^{}
\ar@{-} "4";"4a"^{}
\ar@{-} "5";"5a"^{}
\ar@{-} "6";"6a"^{}
\ar@{-} "7";"7a"^{}
\ar@{-} "1";"top"^{}
\ar@{-} "1a";"topa"^{}
\end{xy}\quad
\begin{xy}
(-20,17) *{Pat(A^{(2)}_{2n-3},\geq j_1)=}="0000",
(-3.7,38.5) *{\ \scriptstyle{j_1}}="c8-1",
(-1.5,36.5) *{\ \scriptstyle{j_2}}="c8-2",
(-3,32.5) *{\ 3}="c7",
(-3,28) *{\ \vdots}="c6",
(-3,22.5) *{\ \scriptstyle{n-1}}="c5",
(-3,17.5) *{\ n}="c4",
(-3,13) *{\ \vdots}="c3",
(-3,7.5) *{\ 3}="c2",
(-3.7,3.5) *{\ \scriptstyle{j_1}}="c1-1",
(-1.5,1.5) *{\ \scriptstyle{j_2}}="c1-2",
(-7,15) *{}="1c",
(2,0) *{}="1b",
(0,0) *{}="1",
(0,5)*{}="2",
(0,10)*{}="3",
(0,15)*{}="4",
(0,20) *{}="5",
(0,25)*{}="6",
(0,30)*{}="7",
(0,35)*{}="8",
(0,40)*{}="9",
(-5,0) *{}="1a",
(-5,5)*{}="2a",
(-5,10)*{}="3a",
(-5,15)*{}="4a",
(-5,20) *{}="5a",
(-5,25)*{}="6a",
(-5,30)*{}="7a",
(-5,35)*{}="8a",
(-5,40)*{}="9a",
(0,41)*{}="top",
(-5,41)*{}="topa",
(2,20)*{}="topb",
(-7,35)*{}="topc",
\ar@{-} "1";"1a"^{}
\ar@{-} "2";"1a"^{}
\ar@{-} "2";"2a"^{}
\ar@{-} "3";"3a"^{}
\ar@{-} "4";"4a"^{}
\ar@{-} "5";"5a"^{}
\ar@{-} "6";"6a"^{}
\ar@{-} "7";"7a"^{}
\ar@{-} "8";"8a"^{}
\ar@{-} "9";"9a"^{}
\ar@{-} "9";"8a"^{}
\ar@{-} "1";"top"^{}
\ar@{-} "1a";"topa"^{}
\end{xy}\quad
\begin{xy}
(-22,17) *{\overline{Pat}(A^{(2)\dagger}_{2n-2},\geq 2)=}="0000",
(-3,36.2) *{\ \scriptstyle{1}}="c8-1",
(-3,38.7) *{\ \scriptstyle{1}}="c8-2",
(-3,32.5) *{\ 2}="c7",
(-3,28) *{\ \vdots}="c6",
(-3,22.5) *{\ \scriptstyle{n-1}}="c5",
(-3,17.5) *{\ n}="c4",
(-3,13) *{\ \vdots}="c3",
(-3,7.5) *{\ 3}="c2",
(-3,1.3) *{\ \scriptstyle{2}}="c1-1",
(-3,3.8) *{\ \scriptstyle{2}}="c1-2",
(-7,15) *{}="1c",
(2,2.5) *{}="1b",
(0,0) *{}="1",
(0,2.5) *{}="1.5",
(0,5)*{}="2",
(0,10)*{}="3",
(0,15)*{}="4",
(0,17.5)*{}="4.5",
(0,20) *{}="5",
(0,25)*{}="6",
(0,30)*{}="7",
(0,35)*{}="8",
(0,37.3)*{}="8.5",
(0,40)*{}="9",
(-5,0) *{}="1a",
(-5,2.5) *{}="1.5a",
(-5,5)*{}="2a",
(-5,10)*{}="3a",
(-5,15)*{}="4a",
(-5,20) *{}="5a",
(-5,25)*{}="6a",
(-5,30)*{}="7a",
(-5,35)*{}="8a",
(-5,37.3)*{}="8.5a",
(-5,40)*{}="9a",
(0,41)*{}="top",
(-5,41)*{}="topa",
(2,20)*{}="topb",
(-7,37.5)*{}="topc",
\ar@{-} "1";"1a"^{}
\ar@{-} "1.5";"1.5a"^{}
\ar@{-} "2";"2a"^{}
\ar@{-} "3";"3a"^{}
\ar@{-} "4";"4a"^{}
\ar@{-} "5";"5a"^{}
\ar@{-} "6";"6a"^{}
\ar@{-} "7";"7a"^{}
\ar@{-} "8";"8a"^{}
\ar@{-} "8.5";"8.5a"^{}
\ar@{-} "9";"9a"^{}
\ar@{-} "1";"top"^{}
\ar@{-} "1a";"topa"^{}
\end{xy}
\quad
\begin{xy}
(-23,17) *{\overline{\overline{Pat}}(B^{(1)}_{n-1},\geq n-1)=}="0000",
(-3,48) *{\ \scriptstyle{n-1}}="c9",
(-3,41.2) *{\ \scriptstyle{n}}="c8-1",
(-3,43.7) *{\ \scriptstyle{n}}="c8-2",
(-3,38) *{\ \scriptstyle{n-1}}="c77",
(-3,32.9) *{\ \vdots}="c7",
(-3,28) *{\ 3}="c6",
(-4.3,24) *{\ \scriptstyle{j_1}}="c5-1",
(-1.7,21) *{\ \scriptstyle{j_2}}="c5-2",
(-3,17.5) *{\ 3}="c4",
(-3,13) *{\ \vdots}="c3",
(-3,7.5) *{\ \scriptstyle{n-2}}="c2",
(-3,1.3) *{\ \scriptstyle{n-1}}="c1-1",
(-3,3.8) *{\ \scriptstyle{n-1}}="c1-2",
(-7,15) *{}="1c",
(2,2.5) *{}="1b",
(0,0) *{}="1",
(0,2.5) *{}="1.5",
(0,5)*{}="2",
(0,10)*{}="3",
(0,15)*{}="4",
(0,17.5)*{}="4.5",
(0,20) *{}="5",
(0,25)*{}="6",
(0,30)*{}="7",
(0,35)*{}="77",
(0,40)*{}="8",
(0,42.3)*{}="8.5",
(0,45)*{}="9",
(0,50)*{}="10",
(-5,0) *{}="1a",
(-5,2.5) *{}="1.5a",
(-5,5)*{}="2a",
(-5,10)*{}="3a",
(-5,15)*{}="4a",
(-5,20) *{}="5a",
(-5,25)*{}="6a",
(-5,30)*{}="7a",
(-5,35)*{}="77a",
(-5,40)*{}="8a",
(-5,42.3)*{}="8.5a",
(-5,45)*{}="9a",
(-5,50)*{}="10a",
(0,51)*{}="top",
(-5,51)*{}="topa",
(2,20)*{}="topb",
(-7,37.5)*{}="topc",
\ar@{-} "1";"1a"^{}
\ar@{-} "1.5";"1.5a"^{}
\ar@{-} "2";"2a"^{}
\ar@{-} "3";"3a"^{}
\ar@{-} "4";"4a"^{}
\ar@{-} "5";"5a"^{}
\ar@{-} "6";"5a"^{}
\ar@{-} "6";"6a"^{}
\ar@{-} "77";"77a"^{}
\ar@{-} "7";"7a"^{}
\ar@{-} "8";"8a"^{}
\ar@{-} "8.5";"8.5a"^{}
\ar@{-} "9";"9a"^{}
\ar@{-} "10";"10a"^{}
\ar@{-} "1";"top"^{}
\ar@{-} "1a";"topa"^{}
\end{xy}
\]

\begin{defn}\cite{Kang}\label{def-YW1}
Let $k\in I$ in class $1$.
A wall $Y$ is said to be {\it Young wall} of ground state $\Lambda_k$ of type $X$ if the following holds:
\begin{enumerate}
\item The wall $Y$ is obtained from $Y_{\Lambda_k}$
by stacking finitely many colored blocks on $Y_{\Lambda_k}$.
\item 
The colored blocks in $i$th column from the right are stacked following the patterns as below:
\begin{enumerate}
\item For $X=A^{(1)}_{n-1}$, the pattern is $Pat(A^{(1)}_{n-1},\geq \pi_{A^{(1)}}(k-i+1))$.
\item For $X=A^{(2)\dagger}_{2n-2}$ or $D^{(2)}_{n}$, the pattern is $Pat(X,\geq k)$ for all $i\in\mathbb{Z}_{\geq1}$.
\item For $X=B^{(1)}_{n-1}$ or $A^{(2)}_{2n-3}$ and $k\neq n$, the pattern is $Pat(X,\geq k)$ with $j_1=k$ when $i$ is odd
and $Pat(X,\geq k)$ with $j_2=k$ when $i$ is even. If $k=n$ so that $X=B^{(1)}_{n-1}$
 then
the pattern is $Pat(B^{(1)}_{n-1},\geq n)$ with $j_1=2$ when $i$ is odd
and $Pat(B^{(1)}_{n-1},\geq n)$ with $j_1=1$ when $i$ is even. 
\item For $X=D^{(1)}_{n-1}$ and $k=1$ or $2$, the pattern is $Pat(D^{(1)}_{n-1},\geq k)$ with $j_1=k$, $j_3=n$ when $i$ is odd
and $Pat(D^{(1)}_{n-1},\geq k)$ with $j_2=k$, $j_3=n-1$ when $i$ is even.
For $k=n-1$ or $n$, the pattern is $Pat(D^{(1)}_{n-1},\geq k)$ with $j_3=k$, $j_1=2$ when $i$ is odd and $Pat(D^{(1)}_{n-1},\geq k)$ with $j_4=k$, $j_1=1$ when $i$ is even.
\end{enumerate} 
\item No block can be put on the top of a single block with half-unit thickness.
\item Except for the right-most column,
there are no free spaces to the
right of any blocks.
\end{enumerate}
\end{defn}

The patterns mentioned in (ii) are explicitly described in Sect.5 of \cite{Kang}.

\begin{defn}\cite{Kang}\label{def-YW2}
Let $Y$ be a Young wall of ground state $\lambda$ of type $X$.
\begin{enumerate}
\item A column of $Y$ is called a {\it full column} if its height is a multiple of the unit length
and its top has unit thickness.
\item For $X=A^{(2)}_{2n-3}$, $A^{(2)\dagger}_{2n-2}$, $B^{(1)}_{n-1}$, $D^{(1)}_{n-1}$ and $D^{(2)}_n$,
the Young wall $Y$ is said to be {\it proper} if none of two full columns of $Y$ have the same height.
\item For $X=A^{(1)}_{n-1}$, we understand all Young walls are proper.
\end{enumerate}
\end{defn}

\begin{defn}\cite{Kang}\label{def-YW2a}
Let $Y$ be a proper Young wall and $i\in I$.
\begin{enumerate}
\item An $i$-block in $Y$ is called a {\it removable} $i$-{\it block} 
if the wall obtained from $Y$ by removing this $i$-block 
remains a proper Young wall. 
\item If we obtain a proper Young wall by adding
an $i$-block to a place of $Y$ then the place is called an
$i$-{\it admissible} {\it slot}.
\end{enumerate}
\end{defn}

\begin{defn}\cite{Ka}\label{def-YW2b}
Let $Y$ be a proper Young wall and $t=1$ or $t=n$. 
\begin{enumerate}
\item[(i)]
Let $Y'$ be a wall
obtained by adding two $t$-blocks of shape (\ref{smpl2}) to the top of a column of $Y$:
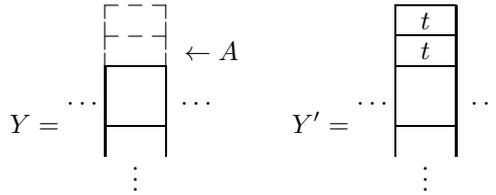
\[
Y=
\begin{xy}
(14,10) *{\leftarrow A}="A",
(-3,3) *{\cdots}="dot1",
(12,3) *{\cdots}="dot2",
(0,-4) *{}="-1",
(8,-4)*{}="-2",
(4,-6)*{\vdots}="dot3",
(0,0) *{}="1",
(8,0)*{}="2",
(8,8)*{}="3",
(0,8)*{}="4",
(8,12)*{}="5",
(0,12)*{}="6",
(8,16)*{}="7",
(0,16)*{}="8",
\ar@{-} "1";"-1"^{}
\ar@{-} "-2";"2"^{}
\ar@{-} "1";"2"^{}
\ar@{-} "1";"4"^{}
\ar@{-} "2";"3"^{}
\ar@{-} "3";"4"^{}
\ar@{--} "3";"5"^{}
\ar@{--} "4";"6"^{}
\ar@{--} "6";"5"^{}
\ar@{--} "7";"5"^{}
\ar@{--} "6";"8"^{}
\ar@{--} "7";"8"^{}
\end{xy}\qquad
Y'=
\begin{xy}
(14,14) *{}="B",
(4,10)*{t}="t1",
(4,14)*{t}="t2",
(-3,3) *{\cdots}="dot1",
(12,3) *{\cdots}="dot2",
(0,-4) *{}="-1",
(8,-4)*{}="-2",
(4,-6)*{\vdots}="dot3",
(0,0) *{}="1",
(8,0)*{}="2",
(8,8)*{}="3",
(0,8)*{}="4",
(8,12)*{}="5",
(0,12)*{}="6",
(8,16)*{}="7",
(0,16)*{}="8",
\ar@{-} "1";"-1"^{}
\ar@{-} "-2";"2"^{}
\ar@{-} "1";"2"^{}
\ar@{-} "1";"4"^{}
\ar@{-} "2";"3"^{}
\ar@{-} "3";"4"^{}
\ar@{-} "3";"5"^{}
\ar@{-} "4";"6"^{}
\ar@{-} "6";"5"^{}
\ar@{-} "7";"5"^{}
\ar@{-} "6";"8"^{}
\ar@{-} "7";"8"^{}
\end{xy}
\]
In the Young wall $Y$, a slot is named $A$ as above. 
The slot $A$ in $Y$ is said to be {\it double} $t$-{\it admissible}
if $Y'$ is a proper Young wall.
\item[(ii)]
Let $Y''$ be a wall
obtained by removing two $t$-blocks of shape (\ref{smpl2}) from the top of a column of $Y$:
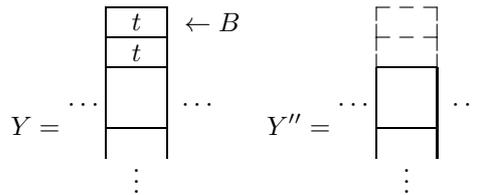
\[
Y=
\begin{xy}
(14,14) *{\leftarrow B}="B",
(4,10)*{t}="t1",
(4,14)*{t}="t2",
(-3,3) *{\cdots}="dot1",
(12,3) *{\cdots}="dot2",
(0,-4) *{}="-1",
(8,-4)*{}="-2",
(4,-6)*{\vdots}="dot3",
(0,0) *{}="1",
(8,0)*{}="2",
(8,8)*{}="3",
(0,8)*{}="4",
(8,12)*{}="5",
(0,12)*{}="6",
(8,16)*{}="7",
(0,16)*{}="8",
\ar@{-} "1";"-1"^{}
\ar@{-} "-2";"2"^{}
\ar@{-} "1";"2"^{}
\ar@{-} "1";"4"^{}
\ar@{-} "2";"3"^{}
\ar@{-} "3";"4"^{}
\ar@{-} "3";"5"^{}
\ar@{-} "4";"6"^{}
\ar@{-} "6";"5"^{}
\ar@{-} "7";"5"^{}
\ar@{-} "6";"8"^{}
\ar@{-} "7";"8"^{}
\end{xy}\quad
Y''=
\begin{xy}
(14,10) *{}="A",
(-3,3) *{\cdots}="dot1",
(12,3) *{\cdots}="dot2",
(0,-4) *{}="-1",
(8,-4)*{}="-2",
(4,-6)*{\vdots}="dot3",
(0,0) *{}="1",
(8,0)*{}="2",
(8,8)*{}="3",
(0,8)*{}="4",
(8,12)*{}="5",
(0,12)*{}="6",
(8,16)*{}="7",
(0,16)*{}="8",
\ar@{-} "1";"-1"^{}
\ar@{-} "-2";"2"^{}
\ar@{-} "1";"2"^{}
\ar@{-} "1";"4"^{}
\ar@{-} "2";"3"^{}
\ar@{-} "3";"4"^{}
\ar@{--} "3";"5"^{}
\ar@{--} "4";"6"^{}
\ar@{--} "6";"5"^{}
\ar@{--} "7";"5"^{}
\ar@{--} "6";"8"^{}
\ar@{--} "7";"8"^{}
\end{xy}
\]
A block in the Young wall $Y$ is named $B$ as above.
The block $B$ in $Y$ is said to be {\it double} $t$-{\it removable}
if $Y''$ is a proper Young wall.
\item[(iii)] Other $i$-admissible slots (resp. removable $i$-blocks)
are said to be single admissible (resp. single removable) for $i\in I$.
\end{enumerate}
\end{defn}

\subsection{Truncated walls}

We take $k\in I$ in class $2$ and 
consider walls obtained from proper Young walls by truncating blocks below a horizontal line.

\begin{defn}\label{def-YW3}
We say a wall $Y$ is a {\it truncated wall} of supporting
(resp. covering)
 ground state $\Lambda_k$ of type $X$ if the following holds:
\begin{enumerate}
\item The wall $Y$ is obtained from
\[
\begin{xy}
(-15.5,1.5) *{\ \cdots}="0000",
(-9.5,1.5) *{\ k}="000",
(-3.5,1.5) *{\ k}="00",
(1.5,1.5) *{\ \ k}="0",
(0,0) *{}="1",
(6,0)*{}="2",
(6,3.5)*{}="3",
(0,3.5)*{}="4",
(-6,3.5)*{}="5",
(-6,0)*{}="6",
(-12,3.5)*{}="7",
(-12,0)*{}="8",
\ar@{-} "1";"2"^{}
\ar@{-} "1";"4"^{}
\ar@{-} "2";"3"^{}
\ar@{-} "3";"4"^{}
\ar@{-} "5";"6"^{}
\ar@{-} "5";"4"^{}
\ar@{-} "1";"6"^{}
\ar@{-} "7";"8"^{}
\ar@{-} "7";"5"^{}
\ar@{-} "6";"8"^{}
\end{xy}
\]
by stacking finitely many colored blocks on it.
\item The colored blocks in $i$th column from the right are stacked following the patterns as below:
\begin{enumerate}
\item For $X=C^{(1)}_{n-1}$, $A^{(2)\dagger}_{2n-2}$ or $D^{(2)}_{n}$, the pattern is $\overline{Pat}(X,\geq k)$ (resp. $\overline{\overline{Pat}}(X,\geq k)$) for all $i\in\mathbb{Z}_{\geq1}$.
\item For $X=B^{(1)}_{n-1}$ or $A^{(2)}_{2n-3}$, the pattern is $\overline{Pat}(X,\geq k)$ (resp. $\overline{\overline{Pat}}(X,\geq k)$) with $j_1=2$ when $i$ is odd
and $\overline{Pat}(X,\geq k)$ (resp. $\overline{\overline{Pat}}(X,\geq k)$) with $j_1=1$ when $i$ is even.
\item For $X=D^{(1)}_{n-1}$, the pattern is $\overline{Pat}(D^{(1)}_{n-1},\geq k)$ (resp. $\overline{\overline{Pat}}(D^{(1)}_{n-1},\geq k)$) with $j_1=2$, $j_3=n$ when $i$ is odd
and $\overline{Pat}(D^{(1)}_{n-1},\geq k)$ (resp. $\overline{\overline{Pat}}(D^{(1)}_{n-1},\geq k)$) with $j_1=1$, $j_3=n-1$ when $i$ is even.
\end{enumerate} 
\item $Y$ satisfies the conditions (iii), (iv) of Definition \ref{def-YW1}.
\end{enumerate}
\end{defn}

We define the notion of {\it full columns} for truncated walls by the same way as in Definition \ref{def-YW2} (i).
Just as in Definition \ref{def-YW2} (ii),
for $X=C^{(1)}_{n-1}$, $A^{(2)}_{2n-3}$, $A^{(2)\dagger}_{2n-2}$, $B^{(1)}_{n-1}$, $D^{(1)}_{n-1}$ and $D^{(2)}_n$,
a truncated wall $Y$ is said to be {\it proper} if none of the two full columns in $Y$ have the same height.
For example, the following $\ovl{Y}_1$ (resp. $\ovl{\ovl{Y}}_1$) is a
proper truncated wall of supporting
(resp. covering) ground state $\Lambda_3$ of type $B^{(1)}_{3}$:
\begin{equation}\label{y1yy1}
\ovl{Y}_1:=
\begin{xy}
(-15.5,-2) *{\ 3}="000-1",
(-9.5,-2) *{\ 3}="00-1",
(-3.5,-2) *{\ 3}="0-1",
(1.5,-2) *{\ \ 3}="0-11",
(1.5,2) *{\ \ 3}="0-2",
(-3.5,2) *{\ 3}="0-22",
(-9.5,2) *{\ 3}="0-222",
(1.5,5.5) *{\ \ 4}="0-3",
(1.5,9) *{\ \ 4}="0-4",
(-3.5,5.5) *{\ 4}="0-33",
(-9.5,5.5) *{\ 4}="0-44",
(1.5,13.5) *{\ \ 3}="0-5",
(0.5,20.6) *{\ \ 2}="0-56",
(-21.5,-2) *{\dots}="00000",
(0,0) *{}="1",
(0,-3.5) *{}="1-u",
(0,3.5) *{}="1-a",
(0,7) *{}="1-aa",
(6,0)*{}="2",
(6,-3.5)*{}="2-u",
(6,3.5)*{}="2-a",
(6,7)*{}="2-aa",
(6,-3.5)*{}="-1-u",
(6,16.5)*{}="3",
(0,16.5)*{}="4",
(6,10.5)*{}="33",
(0,10.5)*{}="44",
(0,22.5)*{}="444",
(6,22.5)*{}="4444",
(-6,3.5)*{}="5",
(-6,0)*{}="6",
(-6,-3.5)*{}="6-u",
(-6,3.5)*{}="6-a",
(-6,7)*{}="6-aa",
(-12,3.5)*{}="7",
(-12,0)*{}="8",
(-12,-3.6)*{}="8-u",
(-12,7)*{}="8-uu",
(-12,3.5)*{}="8-uuu",
(-18,3.5)*{}="9",
(-18,0)*{}="10-a",
(-18,0)*{}="10",
(-18,-3.5)*{}="10-u",
(-18,-3.5)*{}="11-u",
(0,-3.5)*{}="0-u",
\ar@{-} "-1-u";"2-u"^{}
\ar@{-} "2-u";"0-u"^{}
\ar@{-} "10-u";"11-u"^{}
\ar@{-} "8";"10-a"^{}
\ar@{-} "10-u";"10-a"^{}
\ar@{-} "10-u";"2-u"^{}
\ar@{-} "8";"8-u"^{}
\ar@{-} "8";"8-uu"^{}
\ar@{-} "6";"6-u"^{}
\ar@{-} "6";"6-aa"^{}
\ar@{-} "2";"2-u"^{}
\ar@{-} "1";"1-u"^{}
\ar@{-} "1";"2"^{}
\ar@{-} "1";"1-a"^{}
\ar@{-} "1-aa";"1-a"^{}
\ar@{-} "1-aa";"2-aa"^{}
\ar@{-} "8-uu";"2-aa"^{}
\ar@{-} "8-uuu";"1-a"^{}
\ar@{-} "1-aa";"4"^{}
\ar@{-} "2-a";"1-a"^{}
\ar@{-} "2";"3"^{}
\ar@{-} "3";"4"^{}
\ar@{-} "33";"44"^{}
\ar@{-} "444";"44"^{}
\ar@{-} "444";"4444"^{}
\ar@{-} "4";"4444"^{}
\ar@{-} "1";"6"^{}
\ar@{-} "6";"8"^{}
\end{xy}\quad,\quad
\ovl{\ovl{Y}}_1=
\begin{xy}
(-15.5,-2) *{\ 3}="000-1",
(-9.5,-2) *{\ 3}="00-1",
(-3.5,-2) *{\ 3}="0-1",
(1.5,-2) *{\ \ 3}="0-1",
(1.5,2) *{\ \ 3}="0-2",
(-3.5,2) *{\ 3}="0-22",
(-9.5,2) *{\ 3}="0-222",
(3,4.8) *{\ \ 1}="0-3",
(0.5,8) *{\ \ 2}="0-3l",
(-3.5,4.8) *{\ \ 2}="0-4l",
(1.5,13.5) *{\ \ 3}="0-5",
(0,16.5)*{}="4",
(6,9.5)*{}="33",
(0,9.5)*{}="44",
(-21.5,-2) *{\dots}="00000",
(0,0) *{}="1",
(0,3.5) *{}="1-a",
(0,-3.5) *{}="1-u",
(6,0)*{}="2",
(6,-3.5)*{}="2-u",
(6,3.5)*{}="2-a",
(6,9.5)*{}="2-aa",
(6,-3.5)*{}="-1-u",
(6,16.5)*{}="3",
(-6,3.5)*{}="5",
(-6,0)*{}="6",
(-6,-3.5)*{}="6-u",
(-12,3.5)*{}="7",
(-12,0)*{}="8",
(-12,-3.6)*{}="8-u",
(-18,3.5)*{}="9",
(-18,0)*{}="10-a",
(-18,0)*{}="10",
(-18,-3.5)*{}="10-u",
(-18,-3.5)*{}="11-u",
(6,-3.5)*{}="0-u",
\ar@{-} "-1-u";"2-u"^{}
\ar@{-} "2-u";"0-u"^{}
\ar@{-} "10-u";"11-u"^{}
\ar@{-} "8";"10-a"^{}
\ar@{-} "10-u";"10-a"^{}
\ar@{-} "10-u";"2-u"^{}
\ar@{-} "8";"8-u"^{}
\ar@{-} "6";"6-u"^{}
\ar@{-} "2";"2-u"^{}
\ar@{-} "1";"1-u"^{}
\ar@{-} "1";"1-a"^{}
\ar@{-} "2-a";"1-a"^{}
\ar@{-} "2-aa";"1-a"^{}
\ar@{-} "1";"2"^{}
\ar@{-} "1";"4"^{}
\ar@{-} "2";"3"^{}
\ar@{-} "5";"6"^{}
\ar@{-} "1";"6"^{}
\ar@{-} "7";"8"^{}
\ar@{-} "7";"2-a"^{}
\ar@{-} "6";"8"^{}
\ar@{-} "3";"4"^{}
\ar@{-} "33";"44"^{}
\ar@{-} "5";"44"^{}
\end{xy}
\end{equation}
Similarly to Definition \ref{def-YW2a}, one also defines the notion of removable blocks and admissible slots:
\begin{defn}
Let $Y$ be a proper truncated wall of supporting or covering
ground state $\Lambda_k$ and $i\in I$. 
\begin{enumerate}
\item An $i$-block in $Y$
except for
\[
\begin{xy}
(1.5,1.5) *{\ \ k}="0",
(0,0) *{}="1",
(6,0)*{}="2",
(6,3.5)*{}="3",
(0,3.5)*{}="4",
\ar@{-} "1";"2"^{}
\ar@{-} "1";"4"^{}
\ar@{-} "2";"3"^{}
\ar@{-} "3";"4"^{}
\end{xy}
\]
 is called a {\it removable} $i$-{\it block} 
if the wall obtained from $Y$ by removing this $i$-block 
remains a proper truncated wall. 
\item If we obtain a proper truncated wall by adding
an $i$-block
except for
\[
\begin{xy}
(1.5,1.5) *{\ \ k}="0",
(0,0) *{}="1",
(6,0)*{}="2",
(6,3.5)*{}="3",
(0,3.5)*{}="4",
\ar@{-} "1";"2"^{}
\ar@{-} "1";"4"^{}
\ar@{-} "2";"3"^{}
\ar@{-} "3";"4"^{}
\end{xy}
\]
 to a place of $Y$ then the place is called an
$i$-{\it admissible} {\it slot}.
\end{enumerate}
\end{defn}

By the same way as in Definition \ref{def-YW2b}, one defines the notion of 
double removable blocks and double admissible slots for proper truncated walls.

\section{Expression of inequalities by Young walls}

Just as in the previous section, let $X=A^{(1)}_{n-1}$, $B^{(1)}_{n-1}$, $C^{(1)}_{n-1}$, $D^{(1)}_{n-1}$, $A^{(2)\dagger}_{2n-2}$, $A^{(2)}_{2n-3}$ or
$D^{(2)}_{n}$ and $k\in I$.
We define the Langlands dual type $X^L$ as
\begin{equation}\label{Langtype}
X^L=X\ (X=A^{(1)}_{n-1},D^{(1)}_{n-1}),\ (C^{(1)}_{n-1})^L=D^{(2)}_{n},\ (D^{(2)}_{n})^L=C^{(1)}_{n-1},\ 
(A^{(2)}_{2n-3})^L=B^{(1)}_{n-1},\ (B^{(1)}_{n-1})^L=A^{(2)}_{2n-3}.
\end{equation}
We also define $(A^{(2)}_{2n-3})^L=A^{(2)\dagger}_{2n-3}$ and $(A^{(2)\dagger}_{2n-3})^L=A^{(2)}_{2n-3}$. Note that $(X^L)^L=X$.

\subsection{Adaptedness and non-negative integers $P^X_{\ell}$}

\begin{defn}\label{adapt}\cite{KaN}
Let $A=(a_{i,j})_{i,j\in I}$ be the symmetrizable generalized Cartan matrix of $\mathfrak{g}$ and $\io=(\cdots,i_3,i_2,i_1)$
be a sequence
satisfying $(\ref{seq-con})$.
If the following condition holds then we say $\iota$ is {\it adapted} to $A$:
For $i,j\in I$ such that $a_{i,j}<0$, the subsequence of $\iota$ consisting of all $i$, $j$ is
\[
(\cdots,i,j,i,j,i,j,i,j)\quad {\rm or}\quad (\cdots,j,i,j,i,j,i,j,i).
\]
When the Cartan matrix $A$ is fixed, the sequence $\iota$ is shortly said to be {\it adapted}.
\end{defn}

\begin{ex}
Let
$\mathfrak{g}$ be of type $A^{(1)}_2$ and $\iota=(\cdots,3,1,2,3,1,2,3,1,2)$.
It follows $a_{1,2}=a_{2,3}=a_{1,3}=-1$ and 
\begin{itemize}
\item
the subsequence consisting of $1$ and $2$ is $(\cdots,1,2,1,2,1,2)$,
\item
the subsequence consisting of $2$ and $3$ is $(\cdots,3,2,3,2,3,2)$,
\item
the subsequence consisting of $1$ and $3$ is $(\cdots,3,1,3,1,3,1)$.
\end{itemize}
Hence, $\iota$ is an adapted sequence.
\end{ex}

In the rest of paper, we fix a sequence $\iota=(\cdots,i_3,i_2,i_1)$ adapted to $A$.
Let us define a set of integers $(p_{i,j})_{i,j\in I;a_{i,j}<0}$ by
\begin{equation}\label{pij}
p_{i,j}=\begin{cases}
1 & {\rm if}\ {\rm the\ subsequence\ of\ }\iota{\rm\ consisting\ of}\ i,j\ {\rm is}\ (\cdots,j,i,j,i,j,i),\\
0 & {\rm if}\ {\rm the\ subsequence\ of\ }\iota{\rm\ consisting\ of}\ i,j\ {\rm is}\ (\cdots,i,j,i,j,i,j).
\end{cases}
\end{equation}
It is easy to see if $a_{i,j}<0$ then
\begin{equation}\label{pij2}
p_{i,j}+p_{j,i}=1.
\end{equation}
For each $i\in I$, we set $p_{i,i}:=0$.

For $\ell\in\mathbb{Z}$,
one inductively defines integers $P^{A^{(1)}}_{\ell}(t)$ $(t\in\mathbb{Z})$ as $P^{A^{(1)}}_{\ell}(\ell)=0$ and
\begin{equation}\label{pk1}
P^{A^{(1)}}_{\ell}(t):=P^{A^{(1)}}_{\ell}(t-1)+p_{\pi_{A^{(1)}}(t),\pi_{A^{(1)}}(t-1)}\ (\text{for } t>\ell),
\end{equation}
\begin{equation}\label{pk2}
P^{A^{(1)}}_{\ell}(t):=P^{A^{(1)}}_{\ell}(t+1)+p_{\pi_{A^{(1)}}(t),\pi_{A^{(1)}}(t+1)}\ (\text{for } t<\ell).
\end{equation}
For $X=C^{(1)}_{n-1}$, $D^{(2)}_{n}$ or $A^{(2)\dagger}_{2n-2}$ and $\ell\in D_X=\mathbb{Z}_{\geq1}$, we define integers
$P^X_{\ell}(t)$ ($t\in\mathbb{Z}_{\geq \ell}$) as
\[
P^X_{\ell}(\ell):=0,\qquad P^X_{\ell}(t)=P^X_k(t-1)+p_{\pi'_X(t),\pi'_X(t-1)}\ (t>\ell).
\]
For $X=B^{(1)}_{n-1}$ or $A^{(2)}_{2n-3}$ and $\ell\in D_X$, 
we define integers $P^X_{\ell}(t)$ ($t\in D_X$) as follows.
One defines
$P^X_{\ell}(t):=0$ when $t=\ell$ or $t\leq \ell-1$ and
$P^X_{j_1}(j_2):=p_{3,\pi'_X(j_1)}-p_{3,\pi'_X(j_2)}$ 
when $\{j_1,j_2\}=\{1,\frac{3}{2}\}$
and
 for $t>\ell+\frac{1}{2}$,
\[
P^X_{\ell}(t)=
\begin{cases}
P^X_{\ell}(t-\frac{3}{2})+p_{2,3} & \text{if }\pi'_X(t)=2,\\
P^X_{\ell}(t-1)+p_{\pi'_X(t),\pi'_X(t-1)} & \text{otherwise}.\\
\end{cases}
\]
Finally, for $X=D^{(1)}_{n-1}$ and $\ell\in D_X$,
we define integers $P^X_{\ell}(t)$ ($t\in D_X$) as follows.

One defines
$P^X_{\ell}(t):=0$ when $t=\ell$ or $t\leq \ell-1$ and
$P^X_{j_1}(j_2):=p_{3,\pi'_X(j_1)}-p_{3,\pi'_X(j_2)}$ (resp. $P^X_{j_1}(j_2):=p_{n-2,\pi'_X(j_1)}-p_{n-2,\pi'_X(j_2)}$)
when $\{j_1,j_2\}=\{1,\frac{3}{2}\}$ (resp. $\{j_1,j_2\}=\{n-2,n-\frac{3}{2}\}$)
and for $t>\ell+\frac{1}{2}$,
\[
P^X_{\ell}(t)=
\begin{cases}
P^X_{\ell}(t-\frac{3}{2})+p_{2,3} & \text{if }\pi'_X(t)=2,\\
P^X_{\ell}(t-\frac{3}{2})+p_{n,n-2} & \text{if }\pi'_X(t)=n,\\
P^X_{\ell}(t-1)+p_{\pi'_X(t),\pi'_X(t-1)} & \text{otherwise}.\\
\end{cases}
\]
Note that for any type $X$, it holds $P^X_{\ell}(t)\geq0$ except for the case $\{\ell,t\}=\{1,\frac{3}{2}\}$ or $\{\ell,t\}=\{n-2,n-\frac{3}{2}\}$.

\subsection{An identification of single indices and double indices}\label{s-d}

For a fixed sequence $\iota$,
we identify the set of single indices $r\in\mathbb{Z}_{\geq1}$ with the set of double indices $\mathbb{Z}_{\geq1}\times I$ as follows:
Each single index $r\in\mathbb{Z}_{\geq1}$ is identified with a double index $(s,k)\in \mathbb{Z}_{\geq1}\times I$
when $i_r=k$ and $k$ appears $s$ times in $i_r$, $i_{r-1}$, $\cdots,i_1$.
For example, for $\iota=(\cdots,2,1,3,2,1,3,2,1,3)$, single indices
$\cdots,6,5,4,3,2,1$ are identified with double indices
\[
\cdots,(2,2),(2,1),(2,3),(1,2),(1,1),(1,3).
\]
The notation $x_{r}$, $\beta_{r}$, $\beta^{\pm}_{r}$, $S'_{r}$ and $\what{S}'_{r}$ in Sect.2 are also rewritten as
\[
x_r=x_{s,k},\quad \beta_r=\beta_{s,k},\quad \beta^{\pm}_r=\beta^{\pm}_{s,k},\quad S'_r=S'_{s,k},\quad \what{S}'_r=\what{S}'_{s,k}.
\]
If $(s,k)\notin \mathbb{Z}_{\geq1}\times I$ then we define $x_{s,k}:=0$.
By this identification and the ordinary order on $\mathbb{Z}_{\geq1}$, that is, $1<2<3<4<5<6<\cdots$,
one can define an order on $\mathbb{Z}_{\geq1}\times I$. For $\iota=(\cdots,2,1,3,2,1,3,2,1,3)$, the order is
$\cdots>(2,2)>(2,1)>(2,3)>(1,2)>(1,1)>(1,3)$.
Using the notation in (\ref{pij}), if $\iota$ is adapted then $\beta_{s,k}$ can be written as
\[
\beta_{s,k}=x_{s,k}+x_{s+1,k}+\sum_{j\in I; a_{k,j}<0} a_{k,j}x_{s+p_{j,k},j}.
\]

\subsection{Descriptions of inequalities for $B(\infty)$}\label{binf}

In what follows, we draw proper Young walls of ground state $\Lambda_k$ of type $X$ in $\mathbb{R}_{\leq0}\times \mathbb{R}_{\geq T^X_k}$ for $k\in I$ in class $1$.
Similarly, we draw truncated walls of supporting (resp. covering) ground state $\Lambda_k$ of type $X$ in $\mathbb{R}_{\leq0}\times \mathbb{R}_{\geq \overline{T}^X_k}$
(resp. in $\mathbb{R}_{\leq0}\times \mathbb{R}_{\geq \overline{\overline{T}}^X_k}$)  for $k\in I$ in class $2$. 
Here, $T^X_k$, $\overline{T}^X_k$ and $\overline{\overline{T}}^X_k$ are defined in the end of subsection \ref{periodic map}.
For example, when $X=B^{(1)}_2$, the ground state wall $Y_{\Lambda_1}$ is drawn in $\mathbb{R}_{\leq0}\times \mathbb{R}_{\geq 1}$:
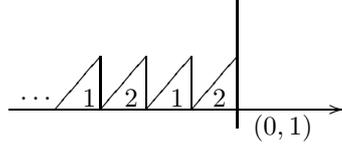
\[
\begin{xy}
(-14,-2) *{\ 1}="000-1",
(-8.5,-2) *{\ 2}="00-1",
(-2.5,-2) *{\ 1}="0-1",
(2.5,-2) *{\ \ 2}="0-11",
(-20.5,-2) *{\dots}="00000",
(12,-6) *{(0,1)}="origin",
(0,0) *{}="1",
(0,-3.5) *{}="1-u",
(6,0)*{}="2",
(6,-3.5)*{}="2-u",
(6,-6)*{}="2-uu",
(6,-7.5)*{}="-1-u",
(6,12)*{}="3",
(-6,3.5)*{}="5",
(-6,0)*{}="6",
(-6,-3.5)*{}="6-u",
(-12,3.5)*{}="7",
(-12,-3.5)*{}="7-u",
(-12,0)*{}="8",
(-12,-3.6)*{}="8-u",
(-18,3.5)*{}="9",
(-18,0)*{}="10-a",
(-18,0)*{}="10",
(-18,-3.5)*{}="10-u",
(-24,-3.5)*{}="11-u",
(20,-3.5)*{}="0-u",
(0,3.5) *{}="1-a",
(6,3.5) *{}="2-a",
\ar@{-} "6-u";"1-a"^{}
\ar@{-} "1-u";"2-a"^{}
\ar@{-} "1-u";"1-a"^{}
\ar@{-} "6-u";"5"^{}
\ar@{-} "7-u";"5"^{}
\ar@{-} "7-u";"7"^{}
\ar@{-} "10-u";"7"^{}
\ar@{->} "2-u";"0-u"^{}
\ar@{-} "10-u";"11-u"^{}
\ar@{-} "10-u";"2-u"^{}
\ar@{-} "2";"2-u"^{}
\ar@{-} "2-uu";"2-u"^{}
\ar@{-} "2";"3"^{}
\end{xy}
\]
When $X=C^{(1)}_2$, a truncated wall of covering ground state $\Lambda_1$ are drawn as
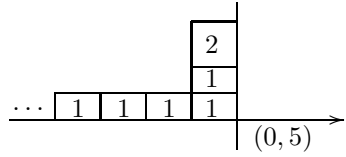
\[
\begin{xy}
(-15.5,-2) *{\ 1}="000-1",
(-9.5,-2) *{\ 1}="00-1",
(-3.5,-2) *{\ 1}="0-1",
(1.5,-2) *{\ \ 1}="0-1",
(1.5,2) *{\ \ 1}="0-2",
(1.5,6.5) *{\ \ 2}="0-3",
(-21.5,-2) *{\dots}="00000",
(12,-6) *{(0,5)}="origin",
(0,0) *{}="1",
(0,3.5) *{}="1-a",
(0,9.5) *{}="1-aa",
(0,-3.5) *{}="1-u",
(6,0)*{}="2",
(6,-3.5)*{}="2-u",
(6,3.5)*{}="2-a",
(6,9.5)*{}="2-aa",
(6,-7.5)*{}="-1-u",
(6,12)*{}="3",
(-6,3.5)*{}="5",
(-6,0)*{}="6",
(-6,-3.5)*{}="6-u",
(-12,3.5)*{}="7",
(-12,0)*{}="8",
(-12,-3.6)*{}="8-u",
(-18,3.5)*{}="9",
(-18,0)*{}="10-a",
(-18,0)*{}="10",
(-18,-3.5)*{}="10-u",
(-24,-3.5)*{}="11-u",
(20,-3.5)*{}="0-u",
\ar@{-} "-1-u";"2-u"^{}
\ar@{->} "2-u";"0-u"^{}
\ar@{-} "10-u";"11-u"^{}
\ar@{-} "8";"10-a"^{}
\ar@{-} "10-u";"10-a"^{}
\ar@{-} "10-u";"2-u"^{}
\ar@{-} "8";"8-u"^{}
\ar@{-} "6";"6-u"^{}
\ar@{-} "2";"2-u"^{}
\ar@{-} "1";"1-u"^{}
\ar@{-} "1";"1-a"^{}
\ar@{-} "2-a";"1-a"^{}
\ar@{-} "1-a";"1-aa"^{}
\ar@{-} "2-aa";"1-aa"^{}
\ar@{-} "1";"2"^{}
\ar@{-} "2";"3"^{}
\ar@{-} "1";"6"^{}
\ar@{-} "6";"8"^{}
\end{xy}
\]
For $k\in I$ and a slot or block $S$ in a proper Young wall or truncated wall, we define
\begin{equation}\label{pka}
P_{k,S}^X(t):=
\begin{cases}
P_{\overline{T}^X_k}^X(t) & \text{if }S\text{ is in a proper Young wall or a truncated wall of supporting ground state}, \\
P_{\overline{\overline{T}}^X_k}^X(t) & \text{if }S\text{ is in a truncated wall of covering ground state} \\
\end{cases}
\end{equation}
and take $s\in\mathbb{Z}$.
Let $i\in\mathbb{Z}_{\geq0}$, $\ell\in\mathbb{Z}_{\geq1}$ and $S$ be the following slot or block in a proper Young wall or truncated wall: 
\[
S=
\begin{xy}
(-8,15) *{(-i-1,\ell+1)}="0000",
(17,15) *{(-i,\ell+1)}="000",
(15,-3) *{(-i,\ell)}="00",
(-7,-3) *{(-i-1,\ell)}="0",
(0,0) *{}="1",
(12,0)*{}="2",
(12,12)*{}="3",
(0,12)*{}="4",
\ar@{-} "1";"2"^{}
\ar@{-} "1";"4"^{}
\ar@{-} "2";"3"^{}
\ar@{-} "3";"4"^{}
\end{xy}
\]
When $S$ is colored by $t\in I$,
we set
\begin{equation}\label{LAdefad}
L^{A^{(1)}}_{s,k,{\rm ad}}(S):= x_{s+P_k^{A^{(1)}}(\ell-i)+{\rm min}\{i,\ell-k\},t},
\end{equation}
\begin{equation}\label{LAdefre}
L^{A^{(1)}}_{s,k,{\rm re}}(S):= x_{s+P_k^{A^{(1)}}(\ell-i)+{\rm min}\{i,\ell-k\}+1,t}
\end{equation}
and for $X(\neq A^{(1)}_{n-1})$, set
\begin{equation}\label{l1def1}
L^X_{s,k,{\rm ad}}(S):=
\begin{cases}
x_{s+P_{k,S}^X(\ell)+i,t} & \text{if }S\text{ is admissible,}
\\
0 & \text{otherwise,}
\end{cases}\quad
L^X_{s,k,{\rm re}}(S):=
\begin{cases}
x_{s+P_{k,S}^X(\ell)+i+1,t} & \text{if }S\text{ is removable,}
\\
0 & \text{otherwise.}
\end{cases}
\end{equation}
By $P_{k,S}^X(\ell)\geq0$, it is easy to see
\begin{equation}\label{non-neg1}
s+P_{k,S}^X(\ell)+i\geq s,\quad s+P_{k,S}^X(\ell)+i+1\geq s+1.
\end{equation}
Let $i\in\mathbb{Z}_{\geq0}$, $\ell\in\mathbb{Z}_{\geq1}$ and $S'$ be
a slot or block colored by $t=1$ or $t=n$ in a proper Young wall or truncated wall
such that the place is one of the following two: 
\[
S'=
\begin{xy}
(-8,10) *{(-i-1,\ell+\frac{1}{2})}="0000",
(17,10) *{(-i,\ell+\frac{1}{2})}="000",
(15,-3) *{(-i,\ell)}="00",
(-7,-3) *{(-i-1,\ell)}="0",
(0,0) *{}="1",
(12,0)*{}="2",
(12,6)*{}="3",
(0,6)*{}="4",
(32,3) *{{\rm or}}="or",
(52,10) *{(-i-1,\ell+1)}="a0000",
(77,10) *{(-i,\ell+1)}="a000",
(75,-3) *{(-i,\ell+\frac{1}{2})}="a00",
(54,-3) *{(-i-1,\ell+\frac{1}{2})}="a0",
(60,0) *{}="a1",
(72,0)*{}="a2",
(72,6)*{}="a3",
(60,6)*{}="a4",
\ar@{-} "1";"2"^{}
\ar@{-} "1";"4"^{}
\ar@{-} "2";"3"^{}
\ar@{-} "3";"4"^{}
\ar@{-} "a1";"a2"^{}
\ar@{-} "a1";"a4"^{}
\ar@{-} "a2";"a3"^{}
\ar@{-} "a3";"a4"^{}
\end{xy}
\]
Then we set
\begin{equation}\label{l1def2}
L^X_{s,k,{\rm ad}}(S'):=
\begin{cases}
x_{s+P_{k,S'}^X(\ell)+i,t} & \text{if }S'\text{ is admissible,}
\\
0 & \text{otherwise,}
\end{cases}\quad
L^X_{s,k,{\rm re}}(S'):=
\begin{cases}
x_{s+P_{k,S'}^X(\ell)+i+1,t} & \text{if }S'\text{ is removable,}
\\
0 & \text{otherwise.}
\end{cases}
\end{equation}
It is clear by $P_{k,S'}^X(\ell)\geq0$ that
\begin{equation}\label{non-neg2}
s+P_{k,S'}^X(\ell)+i\geq s,\quad s+P_{k,S'}^X(\ell)+i+1\geq s+1.
\end{equation}
Note that in the above settings, except for the case $X= A^{(1)}_{n-1}$, it holds
\[
t=\pi'_X(\ell).
\]
Let $i\in\mathbb{Z}_{\geq0}$, $\ell\in\mathbb{Z}_{\geq1}$ and $S''$ be a slot or block
\[
S''=
\begin{xy}
(17,15) *{(-i,\ell+1)}="000",
(15,-3) *{(-i,\ell)}="00",
(-7,-3) *{(-i-1,\ell)}="0",
(0,0) *{}="1",
(12,0)*{}="2",
(12,12)*{}="3",
(32,3) *{{\rm or}}="or",
(52,15) *{(-i-1,\ell+1)}="0000a",
(77,15) *{(-i,\ell+1)}="000a",
(53,-3) *{(-i-1,\ell)}="0a",
(60,0) *{}="1a",
(72,12)*{}="3a",
(60,12)*{}="4a",
\ar@{-} "1";"2"^{}
\ar@{-} "1";"3"^{}
\ar@{-} "2";"3"^{}
\ar@{-} "1a";"3a"^{}
\ar@{-} "1a";"4a"^{}
\ar@{-} "3a";"4a"^{}
\end{xy}
\]
in a proper Young wall or truncated wall.
If $S''$ is colored by $t=1$ or $t=n-1$ then it holds $t=\pi'_X(\ell)$ and we set
\begin{equation}\label{l1def3}
L^X_{s,k,{\rm ad}}(S''):=
\begin{cases}
x_{s+P_{k,S''}^X(\ell)+i,t} & \text{if }S''\text{ is admissible,}
\\
0 & \text{otherwise,}
\end{cases}
\quad
L^X_{s,k,{\rm re}}(S''):=
\begin{cases}
x_{s+P_{k,S''}^X(\ell)+i+1,t} & \text{if }S''\text{ is removable,}
\\
0 & \text{otherwise,}
\end{cases}
\end{equation}
and if colored by $t=2$ or $t=n$ then it holds $t=\pi'_X\left(\ell+\frac{1}{2}\right)$ and we set
\begin{equation}\label{l1def4}
L^X_{s,k,{\rm ad}}(S''):=
\begin{cases}
x_{s+P_{k,S''}^X\left(\ell+\frac{1}{2}\right)+i,t} & \text{if }S''\text{ is admissible,}
\\
0 & \text{otherwise,}
\end{cases}
\quad
L^X_{s,k,{\rm re}}(S''):=
\begin{cases}
x_{s+P_{k,S''}^X\left(\ell+\frac{1}{2}\right)+i+1,t} & \text{if }S''\text{ is removable,}
\\
0 & \text{otherwise.}
\end{cases}
\end{equation}
Writing $P_{k,S''}^X(\ell)=P_{j_1}^X(j_2)$, 
it holds $P_{j_1}^X(j_2)=-1$
only if $\{j_1,j_2\}=\{1,\frac{3}{2}\}$ or $\{j_1,j_2\}=\{n-2,n-\frac{3}{2}\}$.
Then if $S''$ is admissible or removable then we have $i>0$. Therefore, we get $P_{k,S''}^X(\ell)+i\geq0$.
Similarly, we get $P_{k,S''}^X\left(\ell+\frac{1}{2}\right)+i\geq0$. Hence,
\begin{equation}\label{non-neg3}
s+P_{k,S''}^X(\ell)+i\geq s,\quad s+P_{k,S''}^X(\ell)+i+1\geq s+1,\quad
s+P_{k,S''}^X\left(\ell+\frac{1}{2}\right)+i\geq s,\quad s+P_{k,S''}^X\left(\ell+\frac{1}{2}\right)+i+1\geq s+1.
\end{equation}
\begin{defn}\label{position1}
For above $S$, $S'$ or $S''$, we say the {\it position} of $S$, $S'$ or $S''$ is $(-i,\ell)$.
\end{defn}
Let $Y$ be a proper Young wall (resp. a proper truncated wall) of ground state
(resp. supporting or covering ground state) $\Lambda_k$ of type $X=A^{(1)}_{n-1}$, $B^{(1)}_{n-1}$, $C^{(1)}_{n-1}$, $D^{(1)}_{n-1}$, $A^{(2)\dagger}_{2n-2}$, $A^{(2)}_{2n-3}$ or
$D^{(2)}_{n}$. We define
\begin{eqnarray}
L^X_{s,k,\iota}(Y)&:=&
\sum_{t\in I} \left(
\sum_{P:\text{single }t\text{-admissible slot}} L^X_{s,k,{\rm ad}}(P)
-\sum_{P:\text{single removable }t\text{-block}} L^X_{s,k,{\rm re}}(P)\right)\nonumber\\
& & + 
\sum_{P:\text{double }t\text{-admissible slot}} 2L^X_{s,k,{\rm ad}}(P)
-\sum_{P:\text{double removable }t\text{-block}} 2L^X_{s,k,{\rm re}}(P) \label{L11-def}
\end{eqnarray}
for $s\in \mathbb{Z}$.

\begin{defn}\label{YWk-def}
\begin{enumerate}
\item
For $k\in I$ in class $1$,
we define ${\rm YW}^{X}_{k}$ as
the set of all proper Young walls of ground state $\Lambda_k$ of type $X$.
\item
For $k\in I$ in class $2$,
one defines $\overline{{\rm YW}^{X}_k}$ 
(resp. $\overline{\overline{{\rm YW}^{X}_k}}$)
as
the set of all proper truncated walls of supporting (resp. covering) ground state $\Lambda_k$ of type $X$.
\item
For $k\in I$ in class $2$,
\[
{\rm YW}^{X}_{k}:=\left\{(\overline{Y},\overline{\overline{Y}})\in \overline{{\rm YW}^{X}_k}\times \overline{\overline{{\rm YW}^{X}_k}} |
h_j=\frac{1}{2} \text{ if and only if }h_j'=\frac{1}{2}\ \text{for }j\in\mathbb{Z}_{\geq1}
\right\},
\]
where $h_j$ (resp. $h_j'$) is the height of $j$-th column of $\overline{Y}$ (resp. $\overline{\overline{Y}}$) from the right.
We set
\[
Y_{\Lambda_k}:=(
\begin{xy}
(-15.5,1.5) *{\ \cdots}="0000",
(-9.5,1.5) *{\ k}="000",
(-3.5,1.5) *{\ k}="00",
(1.5,1.5) *{\ \ k}="0",
(0,0) *{}="1",
(6,0)*{}="2",
(6,3.5)*{}="3",
(0,3.5)*{}="4",
(-6,3.5)*{}="5",
(-6,0)*{}="6",
(-12,3.5)*{}="7",
(-12,0)*{}="8",
\ar@{-} "1";"2"^{}
\ar@{-} "1";"4"^{}
\ar@{-} "2";"3"^{}
\ar@{-} "3";"4"^{}
\ar@{-} "5";"6"^{}
\ar@{-} "5";"4"^{}
\ar@{-} "1";"6"^{}
\ar@{-} "7";"8"^{}
\ar@{-} "7";"5"^{}
\ar@{-} "6";"8"^{}
\end{xy},
\begin{xy}
(-15.5,1.5) *{\ \cdots}="0000",
(-9.5,1.5) *{\ k}="000",
(-3.5,1.5) *{\ k}="00",
(1.5,1.5) *{\ \ k}="0",
(0,0) *{}="1",
(6,0)*{}="2",
(6,3.5)*{}="3",
(0,3.5)*{}="4",
(-6,3.5)*{}="5",
(-6,0)*{}="6",
(-12,3.5)*{}="7",
(-12,0)*{}="8",
\ar@{-} "1";"2"^{}
\ar@{-} "1";"4"^{}
\ar@{-} "2";"3"^{}
\ar@{-} "3";"4"^{}
\ar@{-} "5";"6"^{}
\ar@{-} "5";"4"^{}
\ar@{-} "1";"6"^{}
\ar@{-} "7";"8"^{}
\ar@{-} "7";"5"^{}
\ar@{-} "6";"8"^{}
\end{xy})\in YW^X_k.
\]
\end{enumerate}
\end{defn}
For example, for $\ovl{Y}_1$, $\ovl{\ovl{Y}}_1$ in (\ref{y1yy1}), it holds $(\ovl{Y}_1,\ovl{\ovl{Y}}_1)\in{\rm YW}^{B^{(1)}_3}_{3}$.
\begin{defn}
We take $k\in I$ in class $2$.
\begin{enumerate}
\item
For $Y=(\overline{Y},\overline{\overline{Y}})\in {\rm YW}^{X}_{k}$,
we suppose that
$Y'=(\overline{Y}',\overline{\overline{Y}}')\in {\rm YW}^{X}_{k}$ is obtained from $Y$
by removing one $k$-block with half-unit height
\[
\begin{xy}
(1.5,1.5) *{\ \ k}="0",
(0,0) *{}="1",
(6,0)*{}="2",
(6,3.5)*{}="3",
(0,3.5)*{}="4",
\ar@{-} "1";"2"^{}
\ar@{-} "1";"4"^{}
\ar@{-} "2";"3"^{}
\ar@{-} "3";"4"^{}
\end{xy}
\]
from the top of $j$-th column of $\overline{Y}$ and one from
the top of $j$-th column of $\overline{\overline{Y}}$ for some $j\in \mathbb{Z}_{\geq1}$.
Then the pair of these two $k$-blocks are said to be a {\it removable} $k$-{\it pair}.
\item
For $Y=(\overline{Y},\overline{\overline{Y}})\in {\rm YW}^{X}_{k}$,
we suppose that
$Y''=(\overline{Y}'',\overline{\overline{Y}}'')\in {\rm YW}^{X}_{k}$ is obtained from $Y$
by adding one $k$-block with half-unit height
\[
\begin{xy}
(1.5,1.5) *{\ \ k}="0",
(0,0) *{}="1",
(6,0)*{}="2",
(6,3.5)*{}="3",
(0,3.5)*{}="4",
\ar@{-} "1";"2"^{}
\ar@{-} "1";"4"^{}
\ar@{-} "2";"3"^{}
\ar@{-} "3";"4"^{}
\end{xy}
\]
to the top of $j$-th column of $\overline{Y}$ and one to 
the top of $j$-th column of
$\overline{\overline{Y}}$ for some $j\in \mathbb{Z}_{\geq1}$. 
Then the pair of these two $k$-slots are said to be a $k$-{\it admissible pair}.
\end{enumerate}
\end{defn}

Let $P$ be a pair in some $Y\in {\rm YW}^{X}_{k}$ in the following form for $k\in I$ in class $2$:
\[
P=
\begin{xy}
(-8,10) *{(-i-1,\overline{T}_k+1)}="0000",
(17,10) *{(-i,\overline{T}_k+1)}="000",
(15,-3) *{(-i,\overline{T}_k+\frac{1}{2})}="00",
(-7,-3) *{(-i-1,\overline{T}_k+\frac{1}{2})}="0",
(6,3) *{k}="0",
(0,0) *{}="1",
(12,0)*{}="2",
(12,6)*{}="3",
(0,6)*{}="4",
(32,3) *{}="or",
(52,10) *{(-i-1,\overline{\overline{T}}_k+1)}="a0000",
(80,10) *{(-i,\overline{\overline{T}}_k+1)}="a000",
(80,-3) *{(-i,\overline{\overline{T}}_k+\frac{1}{2})}="a00",
(54,-3) *{(-i-1,\overline{\overline{T}}_k+\frac{1}{2})}="a0",
(60,0) *{}="a1",
(72,0)*{}="a2",
(72,6)*{}="a3",
(60,6)*{}="a4",
(66,3) *{k}="0",
\ar@{-} "1";"2"^{}
\ar@{-} "1";"4"^{}
\ar@{-} "2";"3"^{}
\ar@{-} "3";"4"^{}
\ar@{-} "a1";"a2"^{}
\ar@{-} "a1";"a4"^{}
\ar@{-} "a2";"a3"^{}
\ar@{-} "a3";"a4"^{}
\end{xy}
\]
Then we set
\begin{equation}\label{l1def5}
L^X_{s,k,{\rm ad}}(P):=
\begin{cases}
x_{s+i,k} & \text{if }P\text{ is a }k\text{-admissible pair}, \\
0 & \text{otherwise,}
\end{cases}
\quad L^X_{s,k,{\rm re}}(P):=
\begin{cases}
x_{s+i+1,k} & \text{if }P\text{ is a removable }k\text{-pair}, \\
0 & \text{otherwise.}
\end{cases}
\end{equation}
\begin{defn}\label{position2}
For above $P$,
we say the {\it position} of $P$ is $(-i,\overline{T}_k)$.
\end{defn}
For $k\in I$ in class $2$ and
$Y=(\overline{Y},\overline{\overline{Y}})\in {\rm YW}^{X}_{k}$, 
combining (\ref{L11-def}) and (\ref{l1def5}),
we
define
\begin{eqnarray}
L^X_{s,k,\iota}(Y)&:=& L^X_{s,k,\iota}(\overline{Y})+L^X_{s,k,\iota}(\overline{\overline{Y}})\nonumber\\
& &+
\sum_{P: k\text{-admissible pair in }Y} L^X_{s,k,{\rm ad}}(P)
-\sum_{P: \text{removable }k\text{-pair in }Y} L^X_{s,k,{\rm re}}(P). \label{L-pair-def}
\end{eqnarray}
We set
${\rm COMB}^{X}_{\iota}[\infty]
:=
\{
L^X_{s,k,\iota}(Y) | s\in\mathbb{Z}_{\geq1},\ k\in I,\ Y\in {\rm YW}^X_k
\}$.
\begin{thm}\label{mainthm1}
Let $\mathfrak{g}$ be of type $X^L$ ($=A^{(1)}_{n-1}$, $B^{(1)}_{n-1}$, $C^{(1)}_{n-1}$, $D^{(1)}_{n-1}$, $A^{(2)}_{2n-2}$, $A^{(2)}_{2n-3}$ or
$D^{(2)}_{n}$). Here, $X^L$ is defined in (\ref{Langtype}).
We suppose $\iota$ is adapted. Then $\iota$ satisfies the $\Xi'$-positivity condition and it holds
\[
{\rm Im}(\Psi_{\iota})=
\{
\mathbf{a}\in\mathbb{Z}^{\infty} | \varphi(\mathbf{a})\geq0\ \text{for all }\varphi\in {\rm COMB}^{X}_{\iota}[\infty]
\}.
\]
\end{thm}

\begin{ex}\label{ex-1}
We consider the case $\mathfrak{g}$ is of type $D^{(2)}_{3}$, $\iota=(\cdots,1,2,3,1,2,3)$. 
We see that
\[
P_{\ovl{T}_1}^{C^{(1)}}(r)=P_{1}^{C^{(1)}}(r)\ (r\in\mathbb{Z}_{\geq1}),\quad P_{\ovl{\ovl{T}}_1}^{C^{(1)}}(r)=P_{5}^{C^{(1)}}(r)\ (r\in\mathbb{Z}_{\geq5})
\]
and
\[
P_{1}^{C^{(1)}}(1)=0,\ P_{1}^{C^{(1)}}(2)=1,\ P_{1}^{C^{(1)}}(3)=P_{1}^{C^{(1)}}(4)=P_{1}^{C^{(1)}}(5)=2, P_{1}^{C^{(1)}}(6)=3,\cdots,
\]
\[
P_{5}^{C^{(1)}}(5)=0,\ P_{5}^{C^{(1)}}(6)=1,\ P_{5}^{C^{(1)}}(7)=P_{5}^{C^{(1)}}(8)=P_{5}^{C^{(1)}}(9)=2, P_{5}^{C^{(1)}}(10)=3,\cdots.
\]
Taking
$Y_{\Lambda_1}$, $Y_1$, $Y_2$, $Y_3$, $Y_4\in {\rm YW}^{(D^{(2)}_3)^L}_1={\rm YW}^{C^{(1)}_2}_1$ as
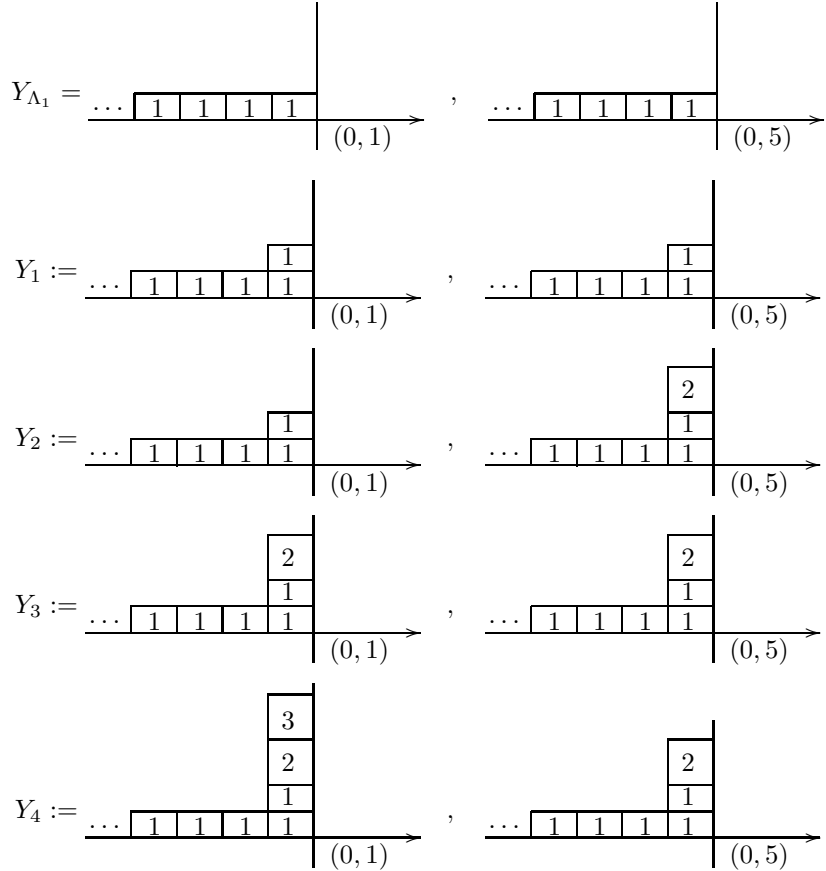
\[
Y_{\Lambda_1}=
\begin{xy}
(-15.5,-2) *{\ 1}="000-1",
(-9.5,-2) *{\ 1}="00-1",
(-3.5,-2) *{\ 1}="0-1",
(1.5,-2) *{\ \ 1}="0-1",
(-21.5,-2) *{\dots}="00000",
(12,-6) *{(0,1)}="origin",
(0,0) *{}="1",
(0,-3.5) *{}="1-u",
(6,0)*{}="2",
(6,-3.5)*{}="2-u",
(6,-7.5)*{}="-1-u",
(6,12)*{}="3",
(-6,3.5)*{}="5",
(-6,0)*{}="6",
(-6,-3.5)*{}="6-u",
(-12,3.5)*{}="7",
(-12,0)*{}="8",
(-12,-3.6)*{}="8-u",
(-18,3.5)*{}="9",
(-18,0)*{}="10-a",
(-18,0)*{}="10",
(-18,-3.5)*{}="10-u",
(-24,-3.5)*{}="11-u",
(20,-3.5)*{}="0-u",
\ar@{-} "-1-u";"2-u"^{}
\ar@{->} "2-u";"0-u"^{}
\ar@{-} "10-u";"11-u"^{}
\ar@{-} "8";"10-a"^{}
\ar@{-} "10-u";"10-a"^{}
\ar@{-} "10-u";"2-u"^{}
\ar@{-} "8";"8-u"^{}
\ar@{-} "6";"6-u"^{}
\ar@{-} "2";"2-u"^{}
\ar@{-} "1";"1-u"^{}
\ar@{-} "1";"2"^{}
\ar@{-} "2";"3"^{}
\ar@{-} "1";"6"^{}
\ar@{-} "6";"8"^{}
\end{xy}\quad,\quad
\begin{xy}
(-15.5,-2) *{\ 1}="000-1",
(-9.5,-2) *{\ 1}="00-1",
(-3.5,-2) *{\ 1}="0-1",
(1.5,-2) *{\ \ 1}="0-1",
(-21.5,-2) *{\dots}="00000",
(12,-6) *{(0,5)}="origin",
(0,0) *{}="1",
(0,-3.5) *{}="1-u",
(6,0)*{}="2",
(6,-3.5)*{}="2-u",
(6,-7.5)*{}="-1-u",
(6,12)*{}="3",
(-6,3.5)*{}="5",
(-6,0)*{}="6",
(-6,-3.5)*{}="6-u",
(-12,3.5)*{}="7",
(-12,0)*{}="8",
(-12,-3.6)*{}="8-u",
(-18,3.5)*{}="9",
(-18,0)*{}="10-a",
(-18,0)*{}="10",
(-18,-3.5)*{}="10-u",
(-24,-3.5)*{}="11-u",
(20,-3.5)*{}="0-u",
\ar@{-} "-1-u";"2-u"^{}
\ar@{->} "2-u";"0-u"^{}
\ar@{-} "10-u";"11-u"^{}
\ar@{-} "8";"10-a"^{}
\ar@{-} "10-u";"10-a"^{}
\ar@{-} "10-u";"2-u"^{}
\ar@{-} "8";"8-u"^{}
\ar@{-} "6";"6-u"^{}
\ar@{-} "2";"2-u"^{}
\ar@{-} "1";"1-u"^{}
\ar@{-} "1";"2"^{}
\ar@{-} "2";"3"^{}
\ar@{-} "1";"6"^{}
\ar@{-} "6";"8"^{}
\end{xy}
\]
\[
Y_1:=
\begin{xy}
(-15.5,-2) *{\ 1}="000-1",
(-9.5,-2) *{\ 1}="00-1",
(-3.5,-2) *{\ 1}="0-1",
(1.5,-2) *{\ \ 1}="0-1",
(1.5,2) *{\ \ 1}="0-2",
(-21.5,-2) *{\dots}="00000",
(12,-6) *{(0,1)}="origin",
(0,0) *{}="1",
(0,-3.5) *{}="1-u",
(0,3.5) *{}="1-a",
(6,0)*{}="2",
(6,-3.5)*{}="2-u",
(6,3.5)*{}="2-a",
(6,-7.5)*{}="-1-u",
(6,12)*{}="3",
(-6,3.5)*{}="5",
(-6,0)*{}="6",
(-6,-3.5)*{}="6-u",
(-6,3.5)*{}="6-a",
(-12,3.5)*{}="7",
(-12,0)*{}="8",
(-12,-3.6)*{}="8-u",
(-18,3.5)*{}="9",
(-18,0)*{}="10-a",
(-18,0)*{}="10",
(-18,-3.5)*{}="10-u",
(-24,-3.5)*{}="11-u",
(20,-3.5)*{}="0-u",
\ar@{-} "-1-u";"2-u"^{}
\ar@{->} "2-u";"0-u"^{}
\ar@{-} "10-u";"11-u"^{}
\ar@{-} "8";"10-a"^{}
\ar@{-} "10-u";"10-a"^{}
\ar@{-} "10-u";"2-u"^{}
\ar@{-} "8";"8-u"^{}
\ar@{-} "6";"6-u"^{}
\ar@{-} "2";"2-u"^{}
\ar@{-} "1";"1-u"^{}
\ar@{-} "1";"2"^{}
\ar@{-} "1";"1-a"^{}
\ar@{-} "2-a";"1-a"^{}
\ar@{-} "2";"3"^{}
\ar@{-} "1";"6"^{}
\ar@{-} "6";"8"^{}
\end{xy}\quad,\quad
\begin{xy}
(-15.5,-2) *{\ 1}="000-1",
(-9.5,-2) *{\ 1}="00-1",
(-3.5,-2) *{\ 1}="0-1",
(1.5,-2) *{\ \ 1}="0-1",
(1.5,2) *{\ \ 1}="0-2",
(-21.5,-2) *{\dots}="00000",
(12,-6) *{(0,5)}="origin",
(0,0) *{}="1",
(0,3.5) *{}="1-a",
(0,-3.5) *{}="1-u",
(6,0)*{}="2",
(6,-3.5)*{}="2-u",
(6,3.5)*{}="2-a",
(6,-7.5)*{}="-1-u",
(6,12)*{}="3",
(-6,3.5)*{}="5",
(-6,0)*{}="6",
(-6,-3.5)*{}="6-u",
(-12,3.5)*{}="7",
(-12,0)*{}="8",
(-12,-3.6)*{}="8-u",
(-18,3.5)*{}="9",
(-18,0)*{}="10-a",
(-18,0)*{}="10",
(-18,-3.5)*{}="10-u",
(-24,-3.5)*{}="11-u",
(20,-3.5)*{}="0-u",
\ar@{-} "-1-u";"2-u"^{}
\ar@{->} "2-u";"0-u"^{}
\ar@{-} "10-u";"11-u"^{}
\ar@{-} "8";"10-a"^{}
\ar@{-} "10-u";"10-a"^{}
\ar@{-} "10-u";"2-u"^{}
\ar@{-} "8";"8-u"^{}
\ar@{-} "6";"6-u"^{}
\ar@{-} "2";"2-u"^{}
\ar@{-} "1";"1-u"^{}
\ar@{-} "1";"1-a"^{}
\ar@{-} "2-a";"1-a"^{}
\ar@{-} "1";"2"^{}
\ar@{-} "2";"3"^{}
\ar@{-} "1";"6"^{}
\ar@{-} "6";"8"^{}
\end{xy}
\]
\[
Y_2:=
\begin{xy}
(-15.5,-2) *{\ 1}="000-1",
(-9.5,-2) *{\ 1}="00-1",
(-3.5,-2) *{\ 1}="0-1",
(1.5,-2) *{\ \ 1}="0-1",
(1.5,2) *{\ \ 1}="0-2",
(-21.5,-2) *{\dots}="00000",
(12,-6) *{(0,1)}="origin",
(0,0) *{}="1",
(0,-3.5) *{}="1-u",
(0,3.5) *{}="1-a",
(6,0)*{}="2",
(6,-3.5)*{}="2-u",
(6,3.5)*{}="2-a",
(6,-7.5)*{}="-1-u",
(6,12)*{}="3",
(-6,3.5)*{}="5",
(-6,0)*{}="6",
(-6,-3.5)*{}="6-u",
(-6,3.5)*{}="6-a",
(-12,3.5)*{}="7",
(-12,0)*{}="8",
(-12,-3.6)*{}="8-u",
(-18,3.5)*{}="9",
(-18,0)*{}="10-a",
(-18,0)*{}="10",
(-18,-3.5)*{}="10-u",
(-24,-3.5)*{}="11-u",
(20,-3.5)*{}="0-u",
\ar@{-} "-1-u";"2-u"^{}
\ar@{->} "2-u";"0-u"^{}
\ar@{-} "10-u";"11-u"^{}
\ar@{-} "8";"10-a"^{}
\ar@{-} "10-u";"10-a"^{}
\ar@{-} "10-u";"2-u"^{}
\ar@{-} "8";"8-u"^{}
\ar@{-} "6";"6-u"^{}
\ar@{-} "2";"2-u"^{}
\ar@{-} "1";"1-u"^{}
\ar@{-} "1";"2"^{}
\ar@{-} "1";"1-a"^{}
\ar@{-} "2-a";"1-a"^{}
\ar@{-} "2";"3"^{}
\ar@{-} "1";"6"^{}
\ar@{-} "6";"8"^{}
\end{xy}\quad,\quad
\begin{xy}
(-15.5,-2) *{\ 1}="000-1",
(-9.5,-2) *{\ 1}="00-1",
(-3.5,-2) *{\ 1}="0-1",
(1.5,-2) *{\ \ 1}="0-1",
(1.5,2) *{\ \ 1}="0-2",
(1.5,6.5) *{\ \ 2}="0-3",
(-21.5,-2) *{\dots}="00000",
(12,-6) *{(0,5)}="origin",
(0,0) *{}="1",
(0,3.5) *{}="1-a",
(0,9.5) *{}="1-aa",
(0,-3.5) *{}="1-u",
(6,0)*{}="2",
(6,-3.5)*{}="2-u",
(6,3.5)*{}="2-a",
(6,9.5)*{}="2-aa",
(6,-7.5)*{}="-1-u",
(6,12)*{}="3",
(-6,3.5)*{}="5",
(-6,0)*{}="6",
(-6,-3.5)*{}="6-u",
(-12,3.5)*{}="7",
(-12,0)*{}="8",
(-12,-3.6)*{}="8-u",
(-18,3.5)*{}="9",
(-18,0)*{}="10-a",
(-18,0)*{}="10",
(-18,-3.5)*{}="10-u",
(-24,-3.5)*{}="11-u",
(20,-3.5)*{}="0-u",
\ar@{-} "-1-u";"2-u"^{}
\ar@{->} "2-u";"0-u"^{}
\ar@{-} "10-u";"11-u"^{}
\ar@{-} "8";"10-a"^{}
\ar@{-} "10-u";"10-a"^{}
\ar@{-} "10-u";"2-u"^{}
\ar@{-} "8";"8-u"^{}
\ar@{-} "6";"6-u"^{}
\ar@{-} "2";"2-u"^{}
\ar@{-} "1";"1-u"^{}
\ar@{-} "1";"1-a"^{}
\ar@{-} "2-a";"1-a"^{}
\ar@{-} "1-a";"1-aa"^{}
\ar@{-} "2-aa";"1-aa"^{}
\ar@{-} "1";"2"^{}
\ar@{-} "2";"3"^{}
\ar@{-} "1";"6"^{}
\ar@{-} "6";"8"^{}
\end{xy}
\]
\[
Y_3:=
\begin{xy}
(-15.5,-2) *{\ 1}="000-1",
(-9.5,-2) *{\ 1}="00-1",
(-3.5,-2) *{\ 1}="0-1",
(1.5,-2) *{\ \ 1}="0-1",
(1.5,2) *{\ \ 1}="0-2",
(1.5,6.5) *{\ \ 2}="0-2",
(-21.5,-2) *{\dots}="00000",
(12,-6) *{(0,1)}="origin",
(0,0) *{}="1",
(0,-3.5) *{}="1-u",
(0,3.5) *{}="1-a",
(0,9.5) *{}="1-aa",
(0,-3.5) *{}="1-u",
(6,0)*{}="2",
(6,-3.5)*{}="2-u",
(6,3.5)*{}="2-a",
(6,9.5)*{}="2-aa",
(6,-7.5)*{}="-1-u",
(6,12)*{}="3",
(-6,3.5)*{}="5",
(-6,0)*{}="6",
(-6,-3.5)*{}="6-u",
(-6,3.5)*{}="6-a",
(-12,3.5)*{}="7",
(-12,0)*{}="8",
(-12,-3.6)*{}="8-u",
(-18,3.5)*{}="9",
(-18,0)*{}="10-a",
(-18,0)*{}="10",
(-18,-3.5)*{}="10-u",
(-24,-3.5)*{}="11-u",
(20,-3.5)*{}="0-u",
\ar@{-} "-1-u";"2-u"^{}
\ar@{->} "2-u";"0-u"^{}
\ar@{-} "10-u";"11-u"^{}
\ar@{-} "8";"10-a"^{}
\ar@{-} "10-u";"10-a"^{}
\ar@{-} "10-u";"2-u"^{}
\ar@{-} "8";"8-u"^{}
\ar@{-} "6";"6-u"^{}
\ar@{-} "2";"2-u"^{}
\ar@{-} "1";"1-u"^{}
\ar@{-} "1";"2"^{}
\ar@{-} "1";"1-a"^{}
\ar@{-} "2-a";"1-a"^{}
\ar@{-} "2";"3"^{}
\ar@{-} "1";"6"^{}
\ar@{-} "6";"8"^{}
\ar@{-} "1-a";"1-aa"^{}
\ar@{-} "2-aa";"1-aa"^{}
\end{xy}\quad,\quad
\begin{xy}
(-15.5,-2) *{\ 1}="000-1",
(-9.5,-2) *{\ 1}="00-1",
(-3.5,-2) *{\ 1}="0-1",
(1.5,-2) *{\ \ 1}="0-1",
(1.5,2) *{\ \ 1}="0-2",
(1.5,6.5) *{\ \ 2}="0-3",
(-21.5,-2) *{\dots}="00000",
(12,-6) *{(0,5)}="origin",
(0,0) *{}="1",
(0,3.5) *{}="1-a",
(0,9.5) *{}="1-aa",
(0,-3.5) *{}="1-u",
(6,0)*{}="2",
(6,-3.5)*{}="2-u",
(6,3.5)*{}="2-a",
(6,9.5)*{}="2-aa",
(6,-7.5)*{}="-1-u",
(6,12)*{}="3",
(-6,3.5)*{}="5",
(-6,0)*{}="6",
(-6,-3.5)*{}="6-u",
(-12,3.5)*{}="7",
(-12,0)*{}="8",
(-12,-3.6)*{}="8-u",
(-18,3.5)*{}="9",
(-18,0)*{}="10-a",
(-18,0)*{}="10",
(-18,-3.5)*{}="10-u",
(-24,-3.5)*{}="11-u",
(20,-3.5)*{}="0-u",
\ar@{-} "-1-u";"2-u"^{}
\ar@{->} "2-u";"0-u"^{}
\ar@{-} "10-u";"11-u"^{}
\ar@{-} "8";"10-a"^{}
\ar@{-} "10-u";"10-a"^{}
\ar@{-} "10-u";"2-u"^{}
\ar@{-} "8";"8-u"^{}
\ar@{-} "6";"6-u"^{}
\ar@{-} "2";"2-u"^{}
\ar@{-} "1";"1-u"^{}
\ar@{-} "1";"1-a"^{}
\ar@{-} "2-a";"1-a"^{}
\ar@{-} "1-a";"1-aa"^{}
\ar@{-} "2-aa";"1-aa"^{}
\ar@{-} "1";"2"^{}
\ar@{-} "2";"3"^{}
\ar@{-} "1";"6"^{}
\ar@{-} "6";"8"^{}
\end{xy}
\]
\[
Y_4:=
\begin{xy}
(-15.5,-2) *{\ 1}="000-1",
(-9.5,-2) *{\ 1}="00-1",
(-3.5,-2) *{\ 1}="0-1",
(1.5,-2) *{\ \ 1}="0-1",
(1.5,2) *{\ \ 1}="0-2",
(1.5,6.5) *{\ \ 2}="0-22",
(1.5,12) *{\ \ 3}="0-222",
(-21.5,-2) *{\dots}="00000",
(12,-6) *{(0,1)}="origin",
(0,0) *{}="1",
(0,-3.5) *{}="1-u",
(0,3.5) *{}="1-a",
(0,9.5) *{}="1-aa",
(0,15.5) *{}="1-aaa",
(0,-3.5) *{}="1-u",
(6,0)*{}="2",
(6,-3.5)*{}="2-u",
(6,3.5)*{}="2-a",
(6,9.5)*{}="2-aa",
(6,15.5)*{}="2-aaa",
(6,-7.5)*{}="-1-u",
(6,17)*{}="3",
(-6,3.5)*{}="5",
(-6,0)*{}="6",
(-6,-3.5)*{}="6-u",
(-6,3.5)*{}="6-a",
(-12,3.5)*{}="7",
(-12,0)*{}="8",
(-12,-3.6)*{}="8-u",
(-18,3.5)*{}="9",
(-18,0)*{}="10-a",
(-18,0)*{}="10",
(-18,-3.5)*{}="10-u",
(-24,-3.5)*{}="11-u",
(20,-3.5)*{}="0-u",
\ar@{-} "-1-u";"2-u"^{}
\ar@{->} "2-u";"0-u"^{}
\ar@{-} "10-u";"11-u"^{}
\ar@{-} "8";"10-a"^{}
\ar@{-} "10-u";"10-a"^{}
\ar@{-} "10-u";"2-u"^{}
\ar@{-} "8";"8-u"^{}
\ar@{-} "6";"6-u"^{}
\ar@{-} "2";"2-u"^{}
\ar@{-} "1";"1-u"^{}
\ar@{-} "1";"2"^{}
\ar@{-} "1";"1-a"^{}
\ar@{-} "2-a";"1-a"^{}
\ar@{-} "2";"3"^{}
\ar@{-} "1";"6"^{}
\ar@{-} "6";"8"^{}
\ar@{-} "1-a";"1-aa"^{}
\ar@{-} "2-aa";"1-aa"^{}
\ar@{-} "1-aa";"1-aaa"^{}
\ar@{-} "2-aaa";"1-aaa"^{}
\end{xy}\quad,\quad
\begin{xy}
(-15.5,-2) *{\ 1}="000-1",
(-9.5,-2) *{\ 1}="00-1",
(-3.5,-2) *{\ 1}="0-1",
(1.5,-2) *{\ \ 1}="0-1",
(1.5,2) *{\ \ 1}="0-2",
(1.5,6.5) *{\ \ 2}="0-3",
(-21.5,-2) *{\dots}="00000",
(12,-6) *{(0,5)}="origin",
(0,0) *{}="1",
(0,3.5) *{}="1-a",
(0,9.5) *{}="1-aa",
(0,-3.5) *{}="1-u",
(6,0)*{}="2",
(6,-3.5)*{}="2-u",
(6,3.5)*{}="2-a",
(6,9.5)*{}="2-aa",
(6,-7.5)*{}="-1-u",
(6,12)*{}="3",
(-6,3.5)*{}="5",
(-6,0)*{}="6",
(-6,-3.5)*{}="6-u",
(-12,3.5)*{}="7",
(-12,0)*{}="8",
(-12,-3.6)*{}="8-u",
(-18,3.5)*{}="9",
(-18,0)*{}="10-a",
(-18,0)*{}="10",
(-18,-3.5)*{}="10-u",
(-24,-3.5)*{}="11-u",
(20,-3.5)*{}="0-u",
\ar@{-} "-1-u";"2-u"^{}
\ar@{->} "2-u";"0-u"^{}
\ar@{-} "10-u";"11-u"^{}
\ar@{-} "8";"10-a"^{}
\ar@{-} "10-u";"10-a"^{}
\ar@{-} "10-u";"2-u"^{}
\ar@{-} "8";"8-u"^{}
\ar@{-} "6";"6-u"^{}
\ar@{-} "2";"2-u"^{}
\ar@{-} "1";"1-u"^{}
\ar@{-} "1";"1-a"^{}
\ar@{-} "2-a";"1-a"^{}
\ar@{-} "1-a";"1-aa"^{}
\ar@{-} "2-aa";"1-aa"^{}
\ar@{-} "1";"2"^{}
\ar@{-} "2";"3"^{}
\ar@{-} "1";"6"^{}
\ar@{-} "6";"8"^{}
\end{xy}
\]
Let us compute several inequalities
defining ${\rm Im}(\Psi_{\iota})$.
In $Y_{\Lambda_1}$, there is the unique $1$-admissible pair (denoted by $P_1$) at position $(0,1)$ and there are neither admissible slots nor
removable blocks so that
$L^{C^{(1)}_{2}}_{s,1,\iota}(Y_{\Lambda_1})=L^{C^{(1)}_{2}}_{s,1,{\rm ad}}(P_1)=x_{s,1}$ by (\ref{l1def5}), (\ref{L-pair-def}).
In $Y_1$, there are two $2$-admissible slots $S_1$ and $S_2$ at positions $(0,2)$ and $(0,6)$
and removable $1$-pair $P_2$ at position $(0,1)$. Thus,
$L^{C^{(1)}_{2}}_{s,1,\iota}(Y_1)=L^{C^{(1)}_{2}}_{s,1,{\rm ad}}(S_1)+L^{C^{(1)}_{2}}_{s,1,{\rm ad}}(S_2)-L^{C^{(1)}_{2}}_{s,1,{\rm re}}(P_2)
=x_{s+P_{\ovl{T}_1}^{C^{(1)}}(2),2}+x_{s+P_{\ovl{\ovl{T}}_1}^{C^{(1)}}(6),2}-x_{s+1,1}
=2x_{s+1,2}-x_{s+1,1}$.
Similarly, we get
$L^{C^{(1)}_{2}}_{s,1,\iota}(Y_2)=x_{s+1,2}+x_{s+2,3}-x_{s+2,2}$,
$L^{C^{(1)}_{2}}_{s,1,\iota}(Y_3)=x_{s+1,1}+2x_{s+2,3}-2x_{s+2,2}$ and
$L^{C^{(1)}_{2}}_{s,1,\iota}(Y_4)=x_{s+1,1}+x_{s+2,3}-x_{s+3,3}$.
Therefore, for all $s\in\mathbb{Z}_{\geq1}$,\ we see that $\mathbf{a}=(a_{m,j})_{m\in\mathbb{Z}_{\geq1},j\in I}\in{\rm Im}(\Psi_{\iota})$ satisfies
\begin{multline}\label{ex-1-eq1}
a_{s,1}\geq0,\ 2a_{s+1,2}-a_{s+1,1}\geq0,\ a_{s+1,2}+a_{s+2,3}-a_{s+2,2}\geq0,\\
 a_{s+1,1}+2a_{s+2,3}-2a_{s+2,2}\geq0,
a_{s+1,1}+a_{s+2,3}-a_{s+3,3}\geq0.
\end{multline}

\end{ex}

\begin{ex}\label{ex-2}
We consider the case $\mathfrak{g}$ is of type $A^{(2)}_{5}$, $\iota=(\cdots,1,3,4,2,1,3,4,2)$. It holds
\[
P_{\ovl{T}_3}^{B^{(1)}}(r)=P_{2}^{B^{(1)}}(r)\ (r\in D_{B^{(1)}},r\geq2),\quad P_{\ovl{\ovl{T}}_3}^{B^{(1)}}(r)=P_{4}^{B^{(1)}}(r)\ (r\in D_{B^{(1)}},r\geq4)
\]
and
\[
P_{2}^{B^{(1)}}(2)=0,\ P_{2}^{B^{(1)}}(3)=1,\ P_{2}^{B^{(1)}}(4)=1,\ P_{2}^{B^{(1)}}(5)=1,\ P_{2}^{B^{(1)}}\left(\frac{11}{2}\right)=2,\cdots,
\]
\[
P_{4}^{B^{(1)}}(4)=0,\ P_{4}^{B^{(1)}}(5)=0,\ P_{4}^{B^{(1)}}\left(\frac{11}{2}\right)=1,\ P_{4}^{B^{(1)}}(6)=1,\ P_{4}^{B^{(1)}}(7)=2,\cdots.
\]
Taking
$Y_{\Lambda_3}$, $Y_1$, $Y_2$, $Y_3$, $Y_4\in {\rm YW}^{(A^{(2)}_5)^L}_3={\rm YW}^{B^{(1)}_3}_3$ as
\[
Y_{\Lambda_3}=
\begin{xy}
(-15.5,-2) *{\ 3}="000-1",
(-9.5,-2) *{\ 3}="00-1",
(-3.5,-2) *{\ 3}="0-1",
(1.5,-2) *{\ \ 3}="0-1",
(-21.5,-2) *{\dots}="00000",
(12,-6) *{(0,2)}="origin",
(0,0) *{}="1",
(0,-3.5) *{}="1-u",
(6,0)*{}="2",
(6,-3.5)*{}="2-u",
(6,-7.5)*{}="-1-u",
(6,12)*{}="3",
(-6,3.5)*{}="5",
(-6,0)*{}="6",
(-6,-3.5)*{}="6-u",
(-12,3.5)*{}="7",
(-12,0)*{}="8",
(-12,-3.6)*{}="8-u",
(-18,3.5)*{}="9",
(-18,0)*{}="10-a",
(-18,0)*{}="10",
(-18,-3.5)*{}="10-u",
(-24,-3.5)*{}="11-u",
(20,-3.5)*{}="0-u",
\ar@{-} "-1-u";"2-u"^{}
\ar@{->} "2-u";"0-u"^{}
\ar@{-} "10-u";"11-u"^{}
\ar@{-} "8";"10-a"^{}
\ar@{-} "10-u";"10-a"^{}
\ar@{-} "10-u";"2-u"^{}
\ar@{-} "8";"8-u"^{}
\ar@{-} "6";"6-u"^{}
\ar@{-} "2";"2-u"^{}
\ar@{-} "1";"1-u"^{}
\ar@{-} "1";"2"^{}
\ar@{-} "2";"3"^{}
\ar@{-} "1";"6"^{}
\ar@{-} "6";"8"^{}
\end{xy}\quad,\quad
\begin{xy}
(-15.5,-2) *{\ 3}="000-1",
(-9.5,-2) *{\ 3}="00-1",
(-3.5,-2) *{\ 3}="0-1",
(1.5,-2) *{\ \ 3}="0-1",
(-21.5,-2) *{\dots}="00000",
(12,-6) *{(0,4)}="origin",
(0,0) *{}="1",
(0,-3.5) *{}="1-u",
(6,0)*{}="2",
(6,-3.5)*{}="2-u",
(6,-7.5)*{}="-1-u",
(6,12)*{}="3",
(-6,3.5)*{}="5",
(-6,0)*{}="6",
(-6,-3.5)*{}="6-u",
(-12,3.5)*{}="7",
(-12,0)*{}="8",
(-12,-3.6)*{}="8-u",
(-18,3.5)*{}="9",
(-18,0)*{}="10-a",
(-18,0)*{}="10",
(-18,-3.5)*{}="10-u",
(-24,-3.5)*{}="11-u",
(20,-3.5)*{}="0-u",
\ar@{-} "-1-u";"2-u"^{}
\ar@{->} "2-u";"0-u"^{}
\ar@{-} "10-u";"11-u"^{}
\ar@{-} "8";"10-a"^{}
\ar@{-} "10-u";"10-a"^{}
\ar@{-} "10-u";"2-u"^{}
\ar@{-} "8";"8-u"^{}
\ar@{-} "6";"6-u"^{}
\ar@{-} "2";"2-u"^{}
\ar@{-} "1";"1-u"^{}
\ar@{-} "1";"2"^{}
\ar@{-} "2";"3"^{}
\ar@{-} "1";"6"^{}
\ar@{-} "6";"8"^{}
\end{xy}
\]
\[
Y_1:=
\begin{xy}
(-15.5,-2) *{\ 3}="000-1",
(-9.5,-2) *{\ 3}="00-1",
(-3.5,-2) *{\ 3}="0-1",
(1.5,-2) *{\ \ 3}="0-1",
(1.5,2) *{\ \ 3}="0-2",
(-21.5,-2) *{\dots}="00000",
(12,-6) *{(0,2)}="origin",
(0,0) *{}="1",
(0,-3.5) *{}="1-u",
(0,3.5) *{}="1-a",
(6,0)*{}="2",
(6,-3.5)*{}="2-u",
(6,3.5)*{}="2-a",
(6,-7.5)*{}="-1-u",
(6,12)*{}="3",
(-6,3.5)*{}="5",
(-6,0)*{}="6",
(-6,-3.5)*{}="6-u",
(-6,3.5)*{}="6-a",
(-12,3.5)*{}="7",
(-12,0)*{}="8",
(-12,-3.6)*{}="8-u",
(-18,3.5)*{}="9",
(-18,0)*{}="10-a",
(-18,0)*{}="10",
(-18,-3.5)*{}="10-u",
(-24,-3.5)*{}="11-u",
(20,-3.5)*{}="0-u",
\ar@{-} "-1-u";"2-u"^{}
\ar@{->} "2-u";"0-u"^{}
\ar@{-} "10-u";"11-u"^{}
\ar@{-} "8";"10-a"^{}
\ar@{-} "10-u";"10-a"^{}
\ar@{-} "10-u";"2-u"^{}
\ar@{-} "8";"8-u"^{}
\ar@{-} "6";"6-u"^{}
\ar@{-} "2";"2-u"^{}
\ar@{-} "1";"1-u"^{}
\ar@{-} "1";"2"^{}
\ar@{-} "1";"1-a"^{}
\ar@{-} "2-a";"1-a"^{}
\ar@{-} "2";"3"^{}
\ar@{-} "1";"6"^{}
\ar@{-} "6";"8"^{}
\end{xy}\quad,\quad
\begin{xy}
(-15.5,-2) *{\ 3}="000-1",
(-9.5,-2) *{\ 3}="00-1",
(-3.5,-2) *{\ 3}="0-1",
(1.5,-2) *{\ \ 3}="0-1",
(1.5,2) *{\ \ 3}="0-2",
(-21.5,-2) *{\dots}="00000",
(12,-6) *{(0,4)}="origin",
(0,0) *{}="1",
(0,3.5) *{}="1-a",
(0,-3.5) *{}="1-u",
(6,0)*{}="2",
(6,-3.5)*{}="2-u",
(6,3.5)*{}="2-a",
(6,-7.5)*{}="-1-u",
(6,12)*{}="3",
(-6,3.5)*{}="5",
(-6,0)*{}="6",
(-6,-3.5)*{}="6-u",
(-12,3.5)*{}="7",
(-12,0)*{}="8",
(-12,-3.6)*{}="8-u",
(-18,3.5)*{}="9",
(-18,0)*{}="10-a",
(-18,0)*{}="10",
(-18,-3.5)*{}="10-u",
(-24,-3.5)*{}="11-u",
(20,-3.5)*{}="0-u",
\ar@{-} "-1-u";"2-u"^{}
\ar@{->} "2-u";"0-u"^{}
\ar@{-} "10-u";"11-u"^{}
\ar@{-} "8";"10-a"^{}
\ar@{-} "10-u";"10-a"^{}
\ar@{-} "10-u";"2-u"^{}
\ar@{-} "8";"8-u"^{}
\ar@{-} "6";"6-u"^{}
\ar@{-} "2";"2-u"^{}
\ar@{-} "1";"1-u"^{}
\ar@{-} "1";"1-a"^{}
\ar@{-} "2-a";"1-a"^{}
\ar@{-} "1";"2"^{}
\ar@{-} "2";"3"^{}
\ar@{-} "1";"6"^{}
\ar@{-} "6";"8"^{}
\end{xy}
\]
\[
Y_2:=
\begin{xy}
(-15.5,-2) *{\ 3}="000-1",
(-9.5,-2) *{\ 3}="00-1",
(-3.5,-2) *{\ 3}="0-1",
(1.5,-2) *{\ \ 3}="0-1",
(1.5,2) *{\ \ 3}="0-2",
(1.5,5.5) *{\ \ 4}="0-3",
(-21.5,-2) *{\dots}="00000",
(12,-6) *{(0,2)}="origin",
(0,0) *{}="1",
(0,-3.5) *{}="1-u",
(0,3.5) *{}="1-a",
(0,7) *{}="1-aa",
(6,0)*{}="2",
(6,-3.5)*{}="2-u",
(6,3.5)*{}="2-a",
(6,7)*{}="2-aa",
(6,-7.5)*{}="-1-u",
(6,12)*{}="3",
(-6,3.5)*{}="5",
(-6,0)*{}="6",
(-6,-3.5)*{}="6-u",
(-6,3.5)*{}="6-a",
(-12,3.5)*{}="7",
(-12,0)*{}="8",
(-12,-3.6)*{}="8-u",
(-18,3.5)*{}="9",
(-18,0)*{}="10-a",
(-18,0)*{}="10",
(-18,-3.5)*{}="10-u",
(-24,-3.5)*{}="11-u",
(20,-3.5)*{}="0-u",
\ar@{-} "-1-u";"2-u"^{}
\ar@{->} "2-u";"0-u"^{}
\ar@{-} "10-u";"11-u"^{}
\ar@{-} "8";"10-a"^{}
\ar@{-} "10-u";"10-a"^{}
\ar@{-} "10-u";"2-u"^{}
\ar@{-} "8";"8-u"^{}
\ar@{-} "6";"6-u"^{}
\ar@{-} "2";"2-u"^{}
\ar@{-} "1";"1-u"^{}
\ar@{-} "1";"2"^{}
\ar@{-} "1";"1-a"^{}
\ar@{-} "1-aa";"1-a"^{}
\ar@{-} "1-aa";"2-aa"^{}
\ar@{-} "2-a";"1-a"^{}
\ar@{-} "2";"3"^{}
\ar@{-} "1";"6"^{}
\ar@{-} "6";"8"^{}
\end{xy}\quad,\quad
\begin{xy}
(-15.5,-2) *{\ 3}="000-1",
(-9.5,-2) *{\ 3}="00-1",
(-3.5,-2) *{\ 3}="0-1",
(1.5,-2) *{\ \ 3}="0-1",
(1.5,2) *{\ \ 3}="0-2",
(-21.5,-2) *{\dots}="00000",
(12,-6) *{(0,4)}="origin",
(0,0) *{}="1",
(0,3.5) *{}="1-a",
(0,-3.5) *{}="1-u",
(6,0)*{}="2",
(6,-3.5)*{}="2-u",
(6,3.5)*{}="2-a",
(6,-7.5)*{}="-1-u",
(6,12)*{}="3",
(-6,3.5)*{}="5",
(-6,0)*{}="6",
(-6,-3.5)*{}="6-u",
(-12,3.5)*{}="7",
(-12,0)*{}="8",
(-12,-3.6)*{}="8-u",
(-18,3.5)*{}="9",
(-18,0)*{}="10-a",
(-18,0)*{}="10",
(-18,-3.5)*{}="10-u",
(-24,-3.5)*{}="11-u",
(20,-3.5)*{}="0-u",
\ar@{-} "-1-u";"2-u"^{}
\ar@{->} "2-u";"0-u"^{}
\ar@{-} "10-u";"11-u"^{}
\ar@{-} "8";"10-a"^{}
\ar@{-} "10-u";"10-a"^{}
\ar@{-} "10-u";"2-u"^{}
\ar@{-} "8";"8-u"^{}
\ar@{-} "6";"6-u"^{}
\ar@{-} "2";"2-u"^{}
\ar@{-} "1";"1-u"^{}
\ar@{-} "1";"1-a"^{}
\ar@{-} "2-a";"1-a"^{}
\ar@{-} "1";"2"^{}
\ar@{-} "2";"3"^{}
\ar@{-} "1";"6"^{}
\ar@{-} "6";"8"^{}
\end{xy}
\]
\[
Y_3:=
\begin{xy}
(-15.5,-2) *{\ 3}="000-1",
(-9.5,-2) *{\ 3}="00-1",
(-3.5,-2) *{\ 3}="0-1",
(1.5,-2) *{\ \ 3}="0-1",
(1.5,2) *{\ \ 3}="0-2",
(-21.5,-2) *{\dots}="00000",
(12,-6) *{(0,2)}="origin",
(0,0) *{}="1",
(0,-3.5) *{}="1-u",
(0,3.5) *{}="1-a",
(6,0)*{}="2",
(6,-3.5)*{}="2-u",
(6,3.5)*{}="2-a",
(6,-7.5)*{}="-1-u",
(6,12)*{}="3",
(-6,3.5)*{}="5",
(-6,0)*{}="6",
(-6,-3.5)*{}="6-u",
(-6,3.5)*{}="6-a",
(-12,3.5)*{}="7",
(-12,0)*{}="8",
(-12,-3.6)*{}="8-u",
(-18,3.5)*{}="9",
(-18,0)*{}="10-a",
(-18,0)*{}="10",
(-18,-3.5)*{}="10-u",
(-24,-3.5)*{}="11-u",
(20,-3.5)*{}="0-u",
\ar@{-} "-1-u";"2-u"^{}
\ar@{->} "2-u";"0-u"^{}
\ar@{-} "10-u";"11-u"^{}
\ar@{-} "8";"10-a"^{}
\ar@{-} "10-u";"10-a"^{}
\ar@{-} "10-u";"2-u"^{}
\ar@{-} "8";"8-u"^{}
\ar@{-} "6";"6-u"^{}
\ar@{-} "2";"2-u"^{}
\ar@{-} "1";"1-u"^{}
\ar@{-} "1";"2"^{}
\ar@{-} "1";"1-a"^{}
\ar@{-} "2-a";"1-a"^{}
\ar@{-} "2";"3"^{}
\ar@{-} "1";"6"^{}
\ar@{-} "6";"8"^{}
\end{xy}\quad,\quad
\begin{xy}
(-15.5,-2) *{\ 3}="000-1",
(-9.5,-2) *{\ 3}="00-1",
(-3.5,-2) *{\ 3}="0-1",
(1.5,-2) *{\ \ 3}="0-1",
(1.5,2) *{\ \ 3}="0-2",
(3,4.8) *{\ \ 1}="0-3",
(-21.5,-2) *{\dots}="00000",
(12,-6) *{(0,4)}="origin",
(0,0) *{}="1",
(0,3.5) *{}="1-a",
(0,-3.5) *{}="1-u",
(6,0)*{}="2",
(6,-3.5)*{}="2-u",
(6,3.5)*{}="2-a",
(6,9.5)*{}="2-aa",
(6,-7.5)*{}="-1-u",
(6,12)*{}="3",
(-6,3.5)*{}="5",
(-6,0)*{}="6",
(-6,-3.5)*{}="6-u",
(-12,3.5)*{}="7",
(-12,0)*{}="8",
(-12,-3.6)*{}="8-u",
(-18,3.5)*{}="9",
(-18,0)*{}="10-a",
(-18,0)*{}="10",
(-18,-3.5)*{}="10-u",
(-24,-3.5)*{}="11-u",
(20,-3.5)*{}="0-u",
\ar@{-} "-1-u";"2-u"^{}
\ar@{->} "2-u";"0-u"^{}
\ar@{-} "10-u";"11-u"^{}
\ar@{-} "8";"10-a"^{}
\ar@{-} "10-u";"10-a"^{}
\ar@{-} "10-u";"2-u"^{}
\ar@{-} "8";"8-u"^{}
\ar@{-} "6";"6-u"^{}
\ar@{-} "2";"2-u"^{}
\ar@{-} "1";"1-u"^{}
\ar@{-} "1";"1-a"^{}
\ar@{-} "2-a";"1-a"^{}
\ar@{-} "2-aa";"1-a"^{}
\ar@{-} "1";"2"^{}
\ar@{-} "2";"3"^{}
\ar@{-} "1";"6"^{}
\ar@{-} "6";"8"^{}
\end{xy}
\]
\[
Y_4:=
\begin{xy}
(-15.5,-2) *{\ 3}="000-1",
(-9.5,-2) *{\ 3}="00-1",
(-3.5,-2) *{\ 3}="0-1",
(1.5,-2) *{\ \ 3}="0-1",
(1.5,2) *{\ \ 3}="0-2",
(1.5,5.5) *{\ \ 4}="0-3",
(-21.5,-2) *{\dots}="00000",
(12,-6) *{(0,2)}="origin",
(0,0) *{}="1",
(0,-3.5) *{}="1-u",
(0,3.5) *{}="1-a",
(0,7) *{}="1-aa",
(6,0)*{}="2",
(6,-3.5)*{}="2-u",
(6,3.5)*{}="2-a",
(6,7)*{}="2-aa",
(6,-7.5)*{}="-1-u",
(6,12)*{}="3",
(-6,3.5)*{}="5",
(-6,0)*{}="6",
(-6,-3.5)*{}="6-u",
(-6,3.5)*{}="6-a",
(-12,3.5)*{}="7",
(-12,0)*{}="8",
(-12,-3.6)*{}="8-u",
(-18,3.5)*{}="9",
(-18,0)*{}="10-a",
(-18,0)*{}="10",
(-18,-3.5)*{}="10-u",
(-24,-3.5)*{}="11-u",
(20,-3.5)*{}="0-u",
\ar@{-} "-1-u";"2-u"^{}
\ar@{->} "2-u";"0-u"^{}
\ar@{-} "10-u";"11-u"^{}
\ar@{-} "8";"10-a"^{}
\ar@{-} "10-u";"10-a"^{}
\ar@{-} "10-u";"2-u"^{}
\ar@{-} "8";"8-u"^{}
\ar@{-} "6";"6-u"^{}
\ar@{-} "2";"2-u"^{}
\ar@{-} "1";"1-u"^{}
\ar@{-} "1";"2"^{}
\ar@{-} "1";"1-a"^{}
\ar@{-} "1-aa";"1-a"^{}
\ar@{-} "1-aa";"2-aa"^{}
\ar@{-} "2-a";"1-a"^{}
\ar@{-} "2";"3"^{}
\ar@{-} "1";"6"^{}
\ar@{-} "6";"8"^{}
\end{xy}\quad,\quad
\begin{xy}
(-15.5,-2) *{\ 3}="000-1",
(-9.5,-2) *{\ 3}="00-1",
(-3.5,-2) *{\ 3}="0-1",
(1.5,-2) *{\ \ 3}="0-1",
(1.5,2) *{\ \ 3}="0-2",
(3,4.8) *{\ \ 1}="0-3",
(-21.5,-2) *{\dots}="00000",
(12,-6) *{(0,4)}="origin",
(0,0) *{}="1",
(0,3.5) *{}="1-a",
(0,-3.5) *{}="1-u",
(6,0)*{}="2",
(6,-3.5)*{}="2-u",
(6,3.5)*{}="2-a",
(6,9.5)*{}="2-aa",
(6,-7.5)*{}="-1-u",
(6,12)*{}="3",
(-6,3.5)*{}="5",
(-6,0)*{}="6",
(-6,-3.5)*{}="6-u",
(-12,3.5)*{}="7",
(-12,0)*{}="8",
(-12,-3.6)*{}="8-u",
(-18,3.5)*{}="9",
(-18,0)*{}="10-a",
(-18,0)*{}="10",
(-18,-3.5)*{}="10-u",
(-24,-3.5)*{}="11-u",
(20,-3.5)*{}="0-u",
\ar@{-} "-1-u";"2-u"^{}
\ar@{->} "2-u";"0-u"^{}
\ar@{-} "10-u";"11-u"^{}
\ar@{-} "8";"10-a"^{}
\ar@{-} "10-u";"10-a"^{}
\ar@{-} "10-u";"2-u"^{}
\ar@{-} "8";"8-u"^{}
\ar@{-} "6";"6-u"^{}
\ar@{-} "2";"2-u"^{}
\ar@{-} "1";"1-u"^{}
\ar@{-} "1";"1-a"^{}
\ar@{-} "2-a";"1-a"^{}
\ar@{-} "2-aa";"1-a"^{}
\ar@{-} "1";"2"^{}
\ar@{-} "2";"3"^{}
\ar@{-} "1";"6"^{}
\ar@{-} "6";"8"^{}
\end{xy}
\]
Let us compute several inequalities
defining ${\rm Im}(\Psi_{\iota})$.
In $Y_{\Lambda_3}$, there is the unique $3$-admissible pair at position $(0,2)$ and there are neither admissible slots nor
removable blocks so that
$L^{B^{(1)}_{3}}_{s,3,\iota}(Y_{\Lambda_3})=x_{s,3}$. In $Y_1$, there are $1$-admissible slot and $2$-admissible slot at position $(0,5)$,
a double $4$-admissible slot at position $(0,3)$ and removable $3$-pair at position $(0,2)$. Thus,
$L^{B^{(1)}_{3}}_{s,3,\iota}(Y_1)=2x_{s+1,4}+x_{s,1}+x_{s+1,2}-x_{s+1,3}$.
Similarly, we get
$L^{B^{(1)}_{3}}_{s,3,\iota}(Y_2)=x_{s+1,4}+x_{s,1}+x_{s+1,2}-x_{s+2,4}$,
$L^{B^{(1)}_{3}}_{s,3,\iota}(Y_3)=2x_{s+1,4}+x_{s+1,2}-x_{s+1,1}$ and
$L^{B^{(1)}_{3}}_{s,3,\iota}(Y_4)=x_{s+1,4}+x_{s+1,3}+x_{s+1,2}-x_{s+1,1}-x_{s+2,4}$. Hence,
for all $s\in\mathbb{Z}_{\geq1}$, we see that $\mathbf{a}=(a_{m,j})_{m\in\mathbb{Z}_{\geq1},j\in I}\in{\rm Im}(\Psi_{\iota})$ satisfies
\begin{multline}\label{ex-2-eq2}
a_{s,3}\geq0,\ 2a_{s+1,4}+a_{s,1}+a_{s+1,2}-a_{s+1,3}\geq0,\ a_{s+1,4}+a_{s,1}+a_{s+1,2}-a_{s+2,4}\geq0,\\
2a_{s+1,4}+a_{s+1,2}-a_{s+1,1}\geq0,
a_{s+1,4}+a_{s+1,3}+a_{s+1,2}-a_{s+1,1}-a_{s+2,4}\geq0.
\end{multline}

\end{ex}

\subsection{Descriptions for $B(\lambda)$}\label{blamsec}

Let $\ell\in\mathbb{Z}_{\geq1}$.
We
define
\begin{equation}\label{r-tilde-A0}
\fbox{$r$}^{A^{(1)}}_{\ell}:=x_{P^{A^{(1)}}_{\ell}(r),\pi_{A^{(1)}}(r)}
-x_{1+P^{A^{(1)}}_{\ell}(r-1),\pi_{A^{(1)}}(r-1)} \qquad (r\in\mathbb{Z}_{\geq \ell+1}),
\end{equation}
\begin{equation}\label{r-tilde-A1}
\fbox{$\tilde{r}$}_{\ell}^{A^{(1)}}:=x_{P^{A^{(1)}}_{\ell}(r-1),\pi_{A^{(1)}}(r-1)}
-x_{1+P^{A^{(1)}}_{\ell}(r),\pi_{A^{(1)}}(r)} \qquad (r\in\mathbb{Z}_{\leq \ell}).
\end{equation}
For $X= C^{(1)}_{n-1}$, $A^{(2)\dagger}_{2n-2}$ or $D^{(2)}_{n}$ and $r\in\mathbb{Z}_{\geq \ell+1}$, we set
\begin{equation}\label{r-tilde-D0}
\fbox{$r$}^{X}_{\ell}:=
c_X(r)x_{P^X_{\ell}(r),\pi'_X(r)}
-c_X(r-1)x_{1+P^X_{\ell}(r-1),\pi'_X(r-1)},
\end{equation}
where
\[
c_X(r)=
\begin{cases}
2 & \text{ if }(X,\pi'_X(r))=(A^{(2)\dagger}_{2n-2},1),\ (D^{(2)}_{n},1) \text{ or }(D^{(2)}_{n},n),\\
1 & \text{otherwise.}
\end{cases}
\]
For $X=A^{(2)}_{2n-3}$ or $B^{(1)}_{n-1}$ and $r\in D_X$ such that $r\geq \ell+1$, we define
\begin{equation}\label{r-tilde-A2-odd0}
\fbox{$r$}^{X}_{\ell}:=
\begin{cases}
x_{P^X_{\ell}(r),1}
+x_{P^X_{\ell}\left(r+\frac{1}{2}\right),2}-x_{1+P^X_{\ell}(r-1),3}
 & \text{if }\pi'_X(r)=1,\\
x_{P^X_{\ell}(r),2}
-x_{1+P^X_{\ell}\left(r-\frac{1}{2}\right),1} & \text{if }\pi'_X(r)=2,\\
x_{P^X_{\ell}(r),3}
-x_{1+P^X_{\ell}\left(r-\frac{1}{2}\right),2}
-x_{1+P^X_{\ell}(r-1),1} & \text{if }\pi_X'(r)=3\text{ and }\pi_X'(r-1)=1,\\
c_X(r)x_{P^X_{\ell}(r),\pi'_X(r)}
-c_X(r-1)x_{1+P^X_{\ell}(r-1),\pi'_X(r-1)} & \text{otherwise},
\end{cases}
\end{equation}
where
\[
c_X(r)=
\begin{cases}
2 & \text{ if }(X,\pi'_X(r))=(B^{(1)}_{n-1},n),\\
1 & \text{otherwise},
\end{cases}
\]
and for $r$ such that $\pi'_X(r)=2$, define
\begin{equation}\label{r-tilde-A2-odd1}
\fbox{$\tilde{r}$}^{X}_{\ell}:=
x_{P^X_{\ell}\left(r-\frac{1}{2}\right),1}
-x_{1+P^X_{\ell}(r),2}.
\end{equation}
For $X$ and $r\in\mathbb{Z}_{\geq \ell+1}$ such that $(X,\pi'_X(r))=(A^{(2)\dagger}_{2n-2},1),\ (B^{(1)}_{n-1},n),\ (D^{(2)}_{n},1) \text{ or }(D^{(2)}_{n},n)$, we also set
\[
\fbox{$r+\frac{1}{2}$}^{X}_{\ell}:=
x_{P^X_{\ell}(r),\pi'_X(r)}
-x_{1+P^X_{\ell}(r),\pi'_X(r)}
\]
and for $\ell\in \mathbb{Z}_{\geq1}$,
\[
H_{X,\ell}^{\lambda}:=
\{
\fbox{$r+\frac{1}{2}$}^{X}_{\ell}
+\langle h_{\pi'_X(\ell)},\lambda \rangle| r\in \mathbb{Z}_{\geq \ell+1},\ 
(X,\pi'_X(r))=(A^{(2)\dagger}_{2n-2},1),\ (B^{(1)}_{n-1},n),\ (D^{(2)}_{n},1) \text{ or }(D^{(2)}_{n},n)
\}.
\]
For $X=D^{(1)}_{n-1}$ and $r\in D_X$ such that $r\geq \ell+1$,
\begin{equation}\label{boxr-D1}
\fbox{$r$}^{X}_{\ell}:=
\begin{cases}
x_{P^X_{\ell}(r),1}
+x_{P^X_{\ell}\left(r+\frac{1}{2}\right),2}-x_{1+P^X_{\ell}(r-1),3}
 & \text{if }\pi'_X(r)=1,\\
x_{P^X_{\ell}(r),2}
-x_{1+P^X_{\ell}\left(r-\frac{1}{2}\right),1} & \text{if }\pi'_X(r)=2,\\
x_{P^X_{\ell}(r),3}
-x_{1+P^X_{\ell}\left(r-\frac{1}{2}\right),2}
-x_{1+P^X_{\ell}(r-1),1} & \text{if }\pi_X'(r)=3\text{ and }\pi_X'(r-1)=1,\\
x_{P^X_{\ell}(r),n-2}
-x_{1+P^X_{\ell}\left(r-\frac{1}{2}\right),n}
-x_{1+P^X_{\ell}(r-1),n-1} & \text{if }\pi_X'(r)=n-2\text{ and }\pi_X'(r-1)=n-1,
\\
x_{P^X_{\ell}(r),n-1}
+x_{P^X_{\ell}\left(r+\frac{1}{2}\right),n}-x_{1+P^X_{\ell}(r-1),n-2}
 & \text{if }\pi'_X(r)=n-1,\\
x_{P^X_{\ell}(r),n}
-x_{1+P^X_{\ell}\left(r-\frac{1}{2}\right),n-1} & \text{if }\pi'_X(r)=n,\\
x_{P^X_{\ell}(r),\pi'_X(r)}
-x_{1+P^X_{\ell}(r-1),\pi'_X(r-1)} & \text{otherwise},
\end{cases}
\end{equation}
and for $r$ such that $\pi'_X(r)=2$ or $n$, define
\begin{equation}\label{r-tilde-D1}
\fbox{$\tilde{r}$}^{X}_{\ell}:=
\begin{cases}
x_{P^X_{\ell}\left(r-\frac{1}{2}\right),1}
-x_{1+P^X_{\ell}(r),2} & \text{ if }\pi'_X(r)=2,\\
x_{P^X_{\ell}\left(r-\frac{1}{2}\right),n-1}
-x_{1+P^X_{\ell}(r),n} & \text{ if }\pi'_X(r)=n.
\end{cases}
\end{equation}
For $X=A^{(2)}_{2n-3}$, $B^{(1)}_{n-1}$ or $D^{(1)}_{n-1}$ and $\ell\in \mathbb{Z}_{\geq1}$, we set
\[
\tilde{H}_{X,\ell}^{\lambda}:=
\{
\fbox{$\tilde{r}$}^{X}_{\ell}+\langle h_{\pi'_X(\ell)},\lambda \rangle | r\in D_{X},\ r\geq \ell+1,\ 
(X,\pi'_X(r))=(A^{(2)}_{2n-3},2),\ (B^{(1)}_{n-1},2),\ (D^{(1)}_{n-1},2) \text{ or }(D^{(1)}_{n-1},n)
\}.
\]

For $k\in I$,
we also set
\begin{equation}\label{comb-lam-1}
{\rm COMB}_{k,\iota}^{A^{(1)}}[\lambda]:=
\begin{cases}
\{-x_{1,k}+\langle h_k,\lambda \rangle\} & \text{if }(1,k)<(1,\pi_{A^{(1)}}(k+1)),\ (1,k)<(1,\pi_{A^{(1)}}(k-1)),\\
\{\fbox{$\tilde{r}$}_k^{A^{(1)}}+\langle h_k,\lambda \rangle\ | r\in\mathbb{Z}_{\leq k}\} & \text{if }(1,k)<(1,\pi_{A^{(1)}}(k+1)),\ (1,k)>(1,\pi_{A^{(1)}}(k-1)),\\
\{\fbox{$r$}^{A^{(1)}}_k+\langle h_k,\lambda \rangle\ | r\in\mathbb{Z}_{\geq k+1}\} & \text{if }(1,k)>(1,\pi_{A^{(1)}}(k+1)),\ (1,k)<(1,\pi_{A^{(1)}}(k-1)),\\
\{
L^{A^{(1)}}_{0,k,\iota}(T)+\langle h_k,\lambda \rangle | T\in{\rm YW}^{A^{(1)}}_{k}\setminus\{Y_{\Lambda_k}\}
\} & \text{if }(1,k)>(1,\pi_{A^{(1)}}(k+1)),\ (1,k)>(1,\pi_{A^{(1)}}(k-1))\\ 
\end{cases}
\end{equation}
and for $X=C^{(1)}_{n-1}, A^{(2)\dagger}_{2n-2}$ or $D^{(2)}_{n}$, set
\begin{equation}\label{comb-lam-2}
{\rm COMB}_{k,\iota}^{X}[\lambda]:=
\begin{cases}
\{-x_{1,k}+\langle h_k,\lambda \rangle\} & \text{if }(1,k)<(1,\pi_{X}'(k+1)),\ (1,k)<(1,\pi_{X}'(k-1)),\\
\{\fbox{$r$}^{X}_{\overline{\overline{T}}^X_k}+\langle h_k,\lambda \rangle\ | r\in \mathbb{Z}_{\geq \overline{\overline{T}}^X_k+1}\}
\cup
H_{X,\overline{\overline{T}}^X_k}^{\lambda}
 & \text{if }(1,k)<(1,\pi_{X}'(k+1)),\ (1,k)>(1,\pi_{X}'(k-1)),\\
\{\fbox{$r$}_{\overline{T}^X_k}^{X}+\langle h_k,\lambda \rangle\ | r\in \mathbb{Z}_{\geq \overline{T}^X_k+1}\}
\cup
H_{X,\overline{T}^X_k}^{\lambda}
 & \text{if }(1,k)>(1,\pi_{X}'(k+1)),\ (1,k)<(1,\pi_{X}'(k-1)),\\
\{
L^{X}_{0,k,\iota}(T)+\langle h_k,\lambda \rangle | T\in{\rm YW}^{X}_{k}\setminus\{Y_{\Lambda_k}\}
\} & \text{if }(1,k)>(1,\pi_{X}'(k+1)),\ (1,k)>(1,\pi_{X}'(k-1)),\\ 
\end{cases}
\end{equation}
where we understand $H_{C^{(1)}_{n-1},\overline{\overline{T}}^X_k}^{\lambda}=H_{C^{(1)}_{n-1},\overline{T}^X_k}^{\lambda}=\phi$.

Let $(a_{i,j})_{i,j\in I}$ be the generalized Cartan matrix of type $X$.
For $X=B^{(1)}_{n-1}$, $A^{(2)}_{2n-3}$ $(n\geq4)$ or $X=D^{(1)}_{n-1}$ $(n\geq6)$ and $k\in I$, we define ${\rm COMB}_{k,\iota}^{X}[\lambda]$ as follows (as for
$X=D^{(1)}_{4}$, we will define ${\rm COMB}_{k,\iota}^{X}[\lambda]$ in Appendix B):
\vspace{2mm}

\hspace{-7mm}
\underline{When $a_{k,j}=a_{j,k}=-1$ and $a_{k,m}=0$ ($m\in I\setminus\{k,j\}$) for some $j\in I$}

\begin{equation}\label{comb-lam-3}
{\rm COMB}_{k,\iota}^{X}[\lambda]:=
\begin{cases}
\{-x_{1,k}+\langle h_k,\lambda \rangle\} & \text{if }(1,k)<(1,j),\\
\{
L^{X}_{0,k,\iota}(T)+\langle h_k,\lambda \rangle | T\in{\rm YW}^{X}_{k}\setminus\{Y_{\Lambda_k}\}
\} & \text{if }(1,k)>(1,j).
\end{cases}
\end{equation}
Note that it holds
\begin{equation}\label{BAD}
k\in\{1,2\}\ \text{when }X=B^{(1)}_{n-1}\text{ or }A^{(2)}_{2n-3},\ k\in\{1,2,n-1,n\}\ \text{when }X=D^{(1)}_{n-1}.
\end{equation}

\hspace{-7mm}
\underline{When (\ref{BAD}) does not hold and $a_{k,\pi_{X}'(k-2)}<0$, $a_{k,\pi_{X}'(k)}<0$ and $a_{k,m}=0$ ($m\in I\setminus\{k,\pi_{X}'(k-2),\pi_{X}'(k)\}$)}

\vspace{2mm}
\nd
Note that in this case, it holds $\pi_{X}'(k-2)=k-1$ and $\pi_{X}'(k)=k+1$ when $k<n$ and $\pi_{X}'(k-2)=\pi_{X}'(k)=k-1=n-1$ when $k=n$.

\begin{equation}\label{comb-lam-3.5}
{\rm COMB}_{k,\iota}^{X}[\lambda]:=
\begin{cases}
\{-x_{1,k}+\langle h_k,\lambda \rangle\} & \text{if }(1,k)<(1,\pi_{X}'(k)),\ (1,k)<(1,\pi_{X}'(k-2)),\\
\{\fbox{$r$}^{X}_{\overline{\overline{T}}^X_k}+\langle h_k,\lambda \rangle\ | r\in D_{X\geq \overline{\overline{T}}^X_k+1}\}
\cup
H_{X,\overline{\overline{T}}^X_k}^{\lambda}
\cup
\tilde{H}_{X,\overline{\overline{T}}^X_k}^{\lambda}
 & \text{if }(1,k)<(1,\pi_{X}'(k)),\ (1,k)>(1,\pi_{X}'(k-2)),\\
\{\fbox{$r$}_{\overline{T}^X_k}^{X}+\langle h_k,\lambda \rangle\ | r\in D_{X\geq \overline{T}^X_k+1}\}
\cup
H_{X,\overline{T}^X_k}^{\lambda}
\cup
\tilde{H}_{X,\overline{T}^X_k}^{\lambda}
 & \text{if }(1,k)>(1,\pi_{X}'(k)),\ (1,k)<(1,\pi_{X}'(k-2)),\\
\{
L^{X}_{0,k,\iota}(T)+\langle h_k,\lambda \rangle | T\in{\rm YW}^{X}_{k}\setminus\{Y_{\Lambda_k}\}
\} & \text{if }(1,k)>(1,\pi_{X}'(k)),\ (1,k)>(1,\pi_{X}'(k-2)).\\ 
\end{cases}
\end{equation}

\vspace{2mm}

\hspace{-7mm}
\underline{When $a_{k,j_1}=a_{j_1,k}=a_{k,j_2}=a_{j_2,k}=-1$, $a_{k,j_3}<0$,
$j_1$ and $j_2$ are in class $1$ with distinct $j_1,j_2,j_3\in I$}

\vspace{2mm}
\nd
It holds $k=3$ or $k=n-2$.
For $t\in\{j_1,j_2,j_3\}$, we consider the following conditions:

\vspace{4mm}

\nd
$C_{t}$ : $(1,t)<(1,k)$ holds,\qquad $\ovl{C}_{t}$ : $(1,t)>(1,k)$ holds.
\begin{equation}\label{comb-lam-4}
{\rm COMB}_{k,\iota}^{X}[\lambda]:=
\begin{cases}
\{ -x_{1,k} + \lan h_{k},\lambda\ran \} & {\rm if}\ 
\ovl{C}_{j_1}, \ovl{C}_{j_2}, \ovl{C}_{j_3},
\\
\{
 -x_{1,k}+x_{1,j_1}+\lan h_{k},\lambda\ran,\ -x_{2,j_1}+\lan h_{k},\lambda\ran
\}
& {\rm if}\ 
C_{j_1}, \ovl{C}_{j_2},
\ovl{C}_{j_3},
\\
\{
 -x_{1,k}+x_{1,j_2}+\lan h_{k},\lambda\ran,\ -x_{2,j_2}+\lan h_{k},\lambda\ran
\}
& {\rm if}\ 
\ovl{C}_{j_1}, C_{j_2},
\ovl{C}_{j_3},
\\
\{\fbox{$r$}^{X}_{\overline{T}^X_{k}}+\langle h_k,\lambda \rangle\ | r\in D_{X\geq \overline{T}^X_{k}+1}\}
\cup H_{X,\overline{T}^X_{k}}^{\lambda}\cup \tilde{H}_{X,\overline{T}^X_{k}}^{\lambda}
 & {\rm if}\ 
\ovl{C}_{j_1}, \ovl{C}_{j_2}, C_{j_3}\text{ and }k=3,
\\
\{\fbox{$r$}^{X}_{\overline{T}^X_{k}}+\langle h_k,\lambda \rangle\ | r\in D_{X\geq \overline{T}^X_{k}+1}\}
\cup H_{X,\overline{T}^X_{k}}^{\lambda}\cup \tilde{H}_{X,\overline{T}^X_{k}}^{\lambda}
 & {\rm if}\ 
C_{j_1}, C_{j_2}, \ovl{C}_{j_3}\text{ and }k=n-2,
\\
\{\fbox{$r$}_{\overline{\overline{T}}^X_k}^{X}+\langle h_k,\lambda \rangle\ | r\in D_{X\geq \overline{\overline{T}}^X_k+1}\}
\cup H_{X,\overline{\overline{T}}^X_k}^{\lambda}\cup \tilde{H}_{X,\overline{\overline{T}}^X_k}^{\lambda}
& {\rm if}\ 
C_{j_1}, C_{j_2}, \ovl{C}_{j_3}\text{ and }k=3,
\\
\{\fbox{$r$}_{\overline{\overline{T}}^X_k}^{X}+\langle h_k,\lambda \rangle\ | r\in D_{X\geq \overline{\overline{T}}^X_k+1}\}
\cup H_{X,\overline{\overline{T}}^X_k}^{\lambda}\cup \tilde{H}_{X,\overline{\overline{T}}^X_k}^{\lambda}
& {\rm if}\ 
\ovl{C}_{j_1}, \ovl{C}_{j_2}, C_{j_3}\text{ and }k=n-2,
\\
\{
L^{X}_{-1,j_2,\iota}(T)+\langle h_{k},\lambda \rangle | T\in{\rm YW}^{X}_{j_2}\setminus\{Y_{\Lambda_{j_2}},Y_{j_2}^1\}
\}
& {\rm if}\ 
 C_{j_1}, \ovl{C}_{j_2}, C_{j_3}, 
\\
\{
L^{X}_{-1,j_1,\iota}(T)+\langle h_{k},\lambda \rangle | T\in{\rm YW}^{X}_{j_1}\setminus\{Y_{\Lambda_{j_1}},Y_{j_1}^1\}
\}
& {\rm if}\ 
\ovl{C}_{j_1}, C_{j_2}, C_{j_3}, 
\\
\{
L^{X}_{0,k,\iota}(T)+\langle h_{k},\lambda \rangle | T\in{\rm YW}^{X}_{k}\setminus\{Y_{\Lambda_k}\}
\}
& {\rm if}\ 
C_{j_1}, C_{j_2}, C_{j_3},
\end{cases}
\end{equation}
where $D_{X\geq j}$ is defined as $D_{X\geq j}:=\{r\in D_{X}| r\geq j\}$ for $j\in \mathbb{Z}$ and $D_{X}$
is defined in Definition \ref{periodic map}, and
$Y_{j_{r}}^1$ is obtained from $Y_{\Lambda_{j_r}}$ by adding a $j_r$-block to the $j_r$-admissible slot ($r=1,2$).
We understand $H_{X,\overline{\overline{T}}_k}^{\lambda}=H_{X,\overline{T}_k}^{\lambda}=\phi$ when $X=D^{(1)}_{n-1}$ or 
$X=A^{(2)}_{2n-3}$.

\begin{thm}\label{mainthm2}
Let $\mathfrak{g}$ be of type $X^L$ and $\lambda\in P^+$.
We suppose that $\iota$ is adapted. Then the pair $(\iota,\lambda)$ satisfies the $\Xi'$-ample condition and
\[
{\rm Im}(\Psi_{\iota}^{(\lambda)})
=\left\{ \mathbf{a}\in\mathbb{Z}^{\infty}
\left|
\varphi(\mathbf{a})\geq0\\
\text{ for any }
\varphi\in
{\rm COMB}^{X}_{\iota}[\infty]\cup
\bigcup_{k\in I}
{\rm COMB}^{X}_{k,\iota}[\lambda]
\right.
\right\}.
\]
\end{thm}

\begin{ex}\label{ex-3}
We consider the same setting as in Example \ref{ex-1}. It holds
${\rm COMB}_{3,\iota}^{C^{(1)}_2}[\lambda]=\{-x_{1,3}+\lan h_3,\lambda\ran\}$
and
${\rm COMB}_{2,\iota}^{C^{(1)}_2}[\lambda]=\{\fbox{$r$}^{C^{(1)}_2}_2+\lan h_2,\lambda\ran | r\in\mathbb{Z}_{\geq3}\}$.
We explicitly compute
$\fbox{$3$}^{C^{(1)}_2}_2=x_{1,3}-x_{1,2}$,
$\fbox{$4$}^{C^{(1)}_2}_2=x_{1,2}-x_{2,3}$,
$\fbox{$5$}^{C^{(1)}_2}_2=x_{1,1}-x_{2,2}$, $\cdots$.
Let us compute several inequalities in ${\rm COMB}_{1,\iota}^{C^{(1)}_2}[\lambda]
=\{
L^{C^{(1)}_2}_{0,1,\iota}(T)+\langle h_1,\lambda \rangle | T\in{\rm YW}^{C^{(1)}_2}_{1}\setminus\{Y_{\Lambda_1}\}
\}$.
For $Y_1$, $Y_2$, $Y_3$ and $Y_4$ in Example \ref{ex-1}, it holds
$L^{C^{(1)}_{2}}_{0,1,\iota}(Y_1)=2x_{1,2}-x_{1,1}$.
Similarly, we get
$L^{C^{(1)}_{2}}_{0,1,\iota}(Y_2)=x_{1,2}+x_{2,3}-x_{2,2}$,
$L^{C^{(1)}_{2}}_{0,1,\iota}(Y_3)=x_{1,1}+2x_{2,3}-2x_{2,2}$ and
$L^{C^{(1)}_{2}}_{0,1,\iota}(Y_4)=x_{1,1}+x_{2,3}-x_{3,3}$.
Therefore, $\mathbf{a}=(a_{m,j})_{m\in\mathbb{Z}_{\geq1},j\in I}\in{\rm Im}(\Psi_{\iota}^{(\lambda)})$ satisfies the inequalities in
(\ref{ex-1-eq1}) and
\begin{multline*}
-a_{1,3}+\lan h_3,\lambda\ran \geq0,\ a_{1,3}-a_{1,2}+\lan h_2,\lambda\ran\geq0,\ a_{1,2}-a_{2,3}+\lan h_2,\lambda\ran\geq0,\ a_{1,1}-a_{2,2}+\lan h_2,\lambda\ran\geq0,\cdots\\
2a_{1,2}-a_{1,1}+\lan h_1,\lambda\ran\geq0,\ a_{1,2}+a_{2,3}-a_{2,2}+\lan h_1,\lambda\ran\geq0,\ 
a_{1,1}+2a_{2,3}-2a_{2,2}+\lan h_1,\lambda\ran\geq0,\ a_{1,1}+a_{2,3}-a_{3,3}+\lan h_1,\lambda\ran\geq0,\cdots
\end{multline*}

\end{ex}

\begin{ex}\label{ex-4}
We consider the same setting as in Example \ref{ex-2}. It holds
${\rm COMB}_{2,\iota}^{B^{(1)}_3}[\lambda]=\{-x_{1,2}+\lan h_2,\lambda\ran\}$
and ${\rm COMB}_{4,\iota}^{B^{(1)}_3}[\lambda]=\{-x_{1,4}+\lan h_4,\lambda\ran\}$.
One obtains
\[
P_1^{B^{(1)}}(1)=0,\ P_1^{B^{(1)}}\left(\frac{3}{2}\right)=p_{3,1}-p_{3,2}=1,\ P_1^{B^{(1)}}(2)=1,\ P_1^{B^{(1)}}(3)=2.
\]
We put $X_1$, $X_2$, $X_3$, $X_4\in {\rm YW}^{B^{(1)}_3}_{1}\setminus\{Y_{\Lambda_1}\}$ as
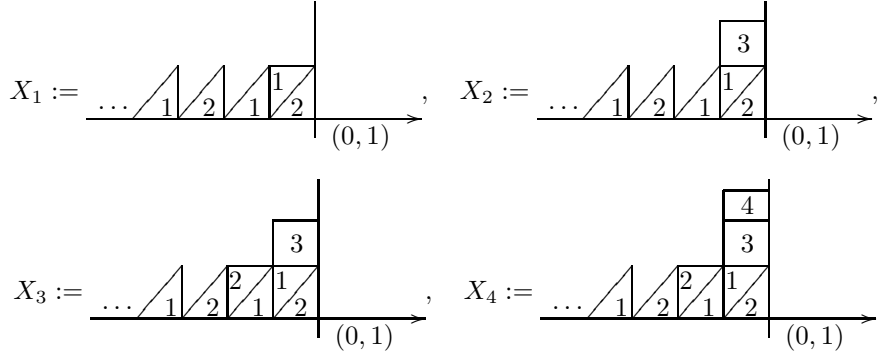
\[
X_1:=
\begin{xy}
(-14,-2) *{\ 1}="000-1",
(-8.5,-2) *{\ 2}="00-1",
(-2.5,-2) *{\ 1}="0-1",
(2.5,-2) *{\ \ 2}="0-11",
(0,1.5) *{\ \ 1}="0-111",
(-20.5,-2) *{\dots}="00000",
(12,-6) *{(0,1)}="origin",
(0,0) *{}="1",
(0,-3.5) *{}="1-u",
(6,0)*{}="2",
(6,-3.5)*{}="2-u",
(6,-6)*{}="2-uu",
(6,-7.5)*{}="-1-u",
(6,12)*{}="3",
(-6,3.5)*{}="5",
(-6,0)*{}="6",
(-6,-3.5)*{}="6-u",
(-12,3.5)*{}="7",
(-12,-3.5)*{}="7-u",
(-12,0)*{}="8",
(-12,-3.6)*{}="8-u",
(-18,3.5)*{}="9",
(-18,0)*{}="10-a",
(-18,0)*{}="10",
(-18,-3.5)*{}="10-u",
(-24,-3.5)*{}="11-u",
(20,-3.5)*{}="0-u",
(0,3.5) *{}="1-a",
(6,3.5) *{}="2-a",
\ar@{-} "2-a";"1-a"^{}
\ar@{-} "6-u";"1-a"^{}
\ar@{-} "1-u";"2-a"^{}
\ar@{-} "1-u";"1-a"^{}
\ar@{-} "6-u";"5"^{}
\ar@{-} "7-u";"5"^{}
\ar@{-} "7-u";"7"^{}
\ar@{-} "10-u";"7"^{}
\ar@{->} "2-u";"0-u"^{}
\ar@{-} "10-u";"11-u"^{}
\ar@{-} "10-u";"2-u"^{}
\ar@{-} "2";"2-u"^{}
\ar@{-} "2-uu";"2-u"^{}
\ar@{-} "2";"3"^{}
\end{xy},\quad
X_2:=
\begin{xy}
(-14,-2) *{\ 1}="000-1",
(-8.5,-2) *{\ 2}="00-1",
(-2.5,-2) *{\ 1}="0-1",
(2.5,-2) *{\ \ 2}="0-11",
(0,1.5) *{\ \ 1}="0-111",
(2,6.5) *{\ \ 3}="0-111-a",
(-20.5,-2) *{\dots}="00000",
(12,-6) *{(0,1)}="origin",
(0,0) *{}="1",
(0,-3.5) *{}="1-u",
(6,0)*{}="2",
(6,-3.5)*{}="2-u",
(6,-6)*{}="2-uu",
(6,-7.5)*{}="-1-u",
(6,12)*{}="3",
(-6,3.5)*{}="5",
(-6,0)*{}="6",
(-6,-3.5)*{}="6-u",
(-12,3.5)*{}="7",
(-12,-3.5)*{}="7-u",
(-12,0)*{}="8",
(-12,-3.6)*{}="8-u",
(-18,3.5)*{}="9",
(-18,0)*{}="10-a",
(-18,0)*{}="10",
(-18,-3.5)*{}="10-u",
(-24,-3.5)*{}="11-u",
(20,-3.5)*{}="0-u",
(0,3.5) *{}="1-a",
(6,3.5) *{}="2-a",
(0,9.5) *{}="1-aa",
(6,9.5) *{}="2-aa",
\ar@{-} "1-aa";"1-a"^{}
\ar@{-} "2-aa";"1-aa"^{}
\ar@{-} "2-a";"1-a"^{}
\ar@{-} "6-u";"1-a"^{}
\ar@{-} "1-u";"2-a"^{}
\ar@{-} "1-u";"1-a"^{}
\ar@{-} "6-u";"5"^{}
\ar@{-} "7-u";"5"^{}
\ar@{-} "7-u";"7"^{}
\ar@{-} "10-u";"7"^{}
\ar@{->} "2-u";"0-u"^{}
\ar@{-} "10-u";"11-u"^{}
\ar@{-} "10-u";"2-u"^{}
\ar@{-} "2";"2-u"^{}
\ar@{-} "2-uu";"2-u"^{}
\ar@{-} "2";"3"^{}
\end{xy},\]
\[
X_3:=
\begin{xy}
(-14,-2) *{\ 1}="000-1",
(-8.5,-2) *{\ 2}="00-1",
(-2.5,-2) *{\ 1}="0-1",
(2.5,-2) *{\ \ 2}="0-11",
(0,1.5) *{\ \ 1}="0-111",
(2,6.5) *{\ \ 3}="0-111-a",
(-5.5,1.5) *{\ 2}="0-12",
(-20.5,-2) *{\dots}="00000",
(12,-6) *{(0,1)}="origin",
(0,0) *{}="1",
(0,-3.5) *{}="1-u",
(6,0)*{}="2",
(6,-3.5)*{}="2-u",
(6,-6)*{}="2-uu",
(6,-7.5)*{}="-1-u",
(6,15)*{}="3",
(-6,3.5)*{}="5",
(-6,0)*{}="6",
(-6,-3.5)*{}="6-u",
(-12,3.5)*{}="7",
(-12,-3.5)*{}="7-u",
(-12,0)*{}="8",
(-12,-3.6)*{}="8-u",
(-18,3.5)*{}="9",
(-18,0)*{}="10-a",
(-18,0)*{}="10",
(-18,-3.5)*{}="10-u",
(-24,-3.5)*{}="11-u",
(20,-3.5)*{}="0-u",
(0,3.5) *{}="1-a",
(6,3.5) *{}="2-a",
(0,9.5) *{}="1-aa",
(6,9.5) *{}="2-aa",
\ar@{-} "1-aa";"1-a"^{}
\ar@{-} "2-aa";"1-aa"^{}
\ar@{-} "2-a";"1-a"^{}
\ar@{-} "6-u";"1-a"^{}
\ar@{-} "1-u";"2-a"^{}
\ar@{-} "1-u";"1-a"^{}
\ar@{-} "1-a";"5"^{}
\ar@{-} "6-u";"5"^{}
\ar@{-} "7-u";"5"^{}
\ar@{-} "7-u";"7"^{}
\ar@{-} "10-u";"7"^{}
\ar@{->} "2-u";"0-u"^{}
\ar@{-} "10-u";"11-u"^{}
\ar@{-} "10-u";"2-u"^{}
\ar@{-} "2";"2-u"^{}
\ar@{-} "2-uu";"2-u"^{}
\ar@{-} "2";"3"^{}
\end{xy},\quad
X_4:=
\begin{xy}
(-14,-2) *{\ 1}="000-1",
(-8.5,-2) *{\ 2}="00-1",
(-2.5,-2) *{\ 1}="0-1",
(2.5,-2) *{\ \ 2}="0-11",
(0,1.5) *{\ \ 1}="0-111",
(2,6.5) *{\ \ 3}="0-111-a",
(2,11.5) *{\ \ 4}="0-111-aa",
(-5.5,1.5) *{\ 2}="0-12",
(-20.5,-2) *{\dots}="00000",
(12,-6) *{(0,1)}="origin",
(0,0) *{}="1",
(0,-3.5) *{}="1-u",
(6,0)*{}="2",
(6,-3.5)*{}="2-u",
(6,-6)*{}="2-uu",
(6,-7.5)*{}="-1-u",
(6,15)*{}="3",
(-6,3.5)*{}="5",
(-6,0)*{}="6",
(-6,-3.5)*{}="6-u",
(-12,3.5)*{}="7",
(-12,-3.5)*{}="7-u",
(-12,0)*{}="8",
(-12,-3.6)*{}="8-u",
(-18,3.5)*{}="9",
(-18,0)*{}="10-a",
(-18,0)*{}="10",
(-18,-3.5)*{}="10-u",
(-24,-3.5)*{}="11-u",
(20,-3.5)*{}="0-u",
(0,3.5) *{}="1-a",
(6,3.5) *{}="2-a",
(0,9.5) *{}="1-aa",
(6,9.5) *{}="2-aa",
(0,13.5) *{}="1-aaa",
(6,13.5) *{}="2-aaa",
\ar@{-} "2-aaa";"1-aaa"^{}
\ar@{-} "1-aaa";"1-aa"^{}
\ar@{-} "2-aaa";"2-aa"^{}
\ar@{-} "1-aa";"1-a"^{}
\ar@{-} "2-aa";"1-aa"^{}
\ar@{-} "2-a";"1-a"^{}
\ar@{-} "6-u";"1-a"^{}
\ar@{-} "1-u";"2-a"^{}
\ar@{-} "1-u";"1-a"^{}
\ar@{-} "1-a";"5"^{}
\ar@{-} "6-u";"5"^{}
\ar@{-} "7-u";"5"^{}
\ar@{-} "7-u";"7"^{}
\ar@{-} "10-u";"7"^{}
\ar@{->} "2-u";"0-u"^{}
\ar@{-} "10-u";"11-u"^{}
\ar@{-} "10-u";"2-u"^{}
\ar@{-} "2";"2-u"^{}
\ar@{-} "2-uu";"2-u"^{}
\ar@{-} "2";"3"^{}
\end{xy}
\]
We see that
$
{\rm COMB}_{1,\iota}^{B^{(1)}_3}[\lambda]=\{
L^{B^{(1)}_3}_{0,1,\iota}(T)+\langle h_1,\lambda \rangle | T\in{\rm YW}^{B^{(1)}_3}_{1}\setminus\{Y_{\Lambda_1}\}
\}$,
$
{\rm COMB}_{3,\iota}^{B^{(1)}_3}[\lambda]=\{
L^{B^{(1)}_3}_{-1,1,\iota}(T)+\langle h_3,\lambda \rangle | T\in{\rm YW}^{B^{(1)}_3}_{1}\setminus\{Y_{\Lambda_1},X_1\}
\}$ and
$L^{B^{(1)}_3}_{0,1,\iota}(X_1)=x_{1,3}-x_{1,1}$, $L^{B^{(1)}_3}_{0,1,\iota}(X_2)=x_{2,2}+2x_{2,4}-x_{2,3}$, 
$L^{B^{(1)}_3}_{0,1,\iota}(X_3)=2x_{2,4}-x_{3,2}$, $L^{B^{(1)}_3}_{0,1,\iota}(X_4)=x_{2,4}+x_{2,3}-x_{3,2}-x_{3,4}$ and
$L^{B^{(1)}_3}_{-1,1,\iota}(X_2)=x_{1,2}+2x_{1,4}-x_{1,3}$, 
$L^{B^{(1)}_3}_{-1,1,\iota}(X_3)=2x_{1,4}-x_{2,2}$, $L^{B^{(1)}_3}_{-1,1,\iota}(X_4)=x_{1,4}+x_{1,3}-x_{2,2}-x_{2,4}$.
Therefore, $\mathbf{a}=(a_{m,j})_{m\in\mathbb{Z}_{\geq1},j\in I}\in{\rm Im}(\Psi_{\iota}^{(\lambda)})$ satisfies the inequalities in
(\ref{ex-2-eq2}) and
\begin{multline*}
-a_{1,2}+\lan h_2,\lambda\ran \geq0,\ 
-a_{1,4}+\lan h_4,\lambda\ran \geq0,\ 
a_{1,3}-a_{1,1}+\lan h_1,\lambda\ran \geq0,\ 
a_{2,2}+2a_{2,4}-a_{2,3}+\lan h_1,\lambda\ran \geq0,\\ 
2a_{2,4}-a_{3,2}+\lan h_1,\lambda\ran \geq0,\ 
a_{2,4}+a_{2,3}-a_{3,2}-a_{3,4}+\lan h_1,\lambda\ran \geq0,\cdots,\\
a_{1,2}+2a_{1,4}-a_{1,3}+\lan h_3,\lambda\ran\geq0,\ 
2a_{1,4}-a_{2,2}+\lan h_3,\lambda\ran\geq0,\ 
a_{1,4}+a_{1,3}-a_{2,2}-a_{2,4}+\lan h_3,\lambda\ran\geq0,\cdots.
\end{multline*}

\end{ex}

\subsection{Combinatorial descriptions of $\varepsilon_k^*$}

Let us give a combinatorial descriptions $\varepsilon_k^*$
in Example \ref{star-ex}.

\begin{thm}\label{mainthm3}
We suppose that $\iota$ is adapted and $\mathfrak{g}$ is of type $X^L$. 
Then
$\iota$ satisfies $\Xi'$-strict positivity condition and
for $k\in I$ and $\mathbf{a}\in{\rm Im}(\Psi_{\iota})$, we have
\[
\varepsilon_k^*(\mathbf{a})
={\rm max}\{-\varphi(\mathbf{a}) | \varphi\in{\rm COMB}_{k,\iota}^{X}[0]\}.
\]
\end{thm}

\begin{ex}\label{ex-5}
We consider the same setting as in Example \ref{ex-1}, \ref{ex-3}.
By a computation in Example \ref{ex-3} and above theorem, one obtains
\[
\varepsilon_3^*(\mathbf{a})
={\rm max}\{-\varphi(\mathbf{a}) | \varphi\in{\rm COMB}_{3,\iota}^{C^{(1)}_2}[0]=\{-x_{1,3}\}\}
=a_{1,3}
\]
for $\mathbf{a}=(a_{s,j})_{s\in\mathbb{Z}_{\geq1},j\in I}\in{\rm Im}(\Psi_{\iota})$. We also obtain
\[
\varepsilon_2^*(\mathbf{a})
={\rm max}\{-\varphi(\mathbf{a}) | \varphi\in{\rm COMB}_{2,\iota}^{C^{(1)}_2}[0]=\{\fbox{$r$}^{C^{(1)}_2}_2 | r\in\mathbb{Z}_{\geq3}\}\}.
\]
Note that there exists $(s,j)\in\mathbb{Z}_{\geq1}\times I$ such that $a_{s',j'}=0$ for $(s',j')>(s,j)$ so that
one can explicitly compute $\varepsilon_2^*(\mathbf{a})$. For example,
if $a_{s',j'}=0$ for all $(s',j')>(2,1)$ then $\fbox{$r$}^{C^{(1)}_2}_2(\mathbf{a})=0$ for $r>6$
so that
$\varepsilon_2^*(\mathbf{a})={\rm max}\{-\fbox{$3$}^{C^{(1)}_2}_2(\mathbf{a}),
-\fbox{$4$}^{C^{(1)}_2}_2(\mathbf{a}),
-\fbox{$5$}^{C^{(1)}_2}_2(\mathbf{a}),
-\fbox{$6$}^{C^{(1)}_2}_2(\mathbf{a})
\}={\rm max}\{a_{1,2}-a_{1,3},a_{2,3}-a_{1,2},a_{2,2}-a_{1,1},a_{2,1}-a_{2,2},0\}$. 
Similarly,
\[
\varepsilon_1^*(\mathbf{a})
=
{\rm max}\{-\varphi(\mathbf{a}) | \varphi\in
\{
L^{C^{(1)}_2}_{0,1,\iota}(T)+\langle h_1,\lambda \rangle | T\in{\rm YW}^{C^{(1)}_2}_{1}\setminus\{Y_{\Lambda_1}\}
\}\}.
\]

\end{ex}

\begin{rem}
In case of $X=A^{(1)}_{n-1}$, proper Young walls of ground state $\Lambda_k$ can be regarded as
extended Young diagrams via reflection over the line $y=x+k$.
By this identification, $t$-admissible slots and removable $t$-blocks with position $(-i,\ell)$ correspond to
concave corners $(\ell-k,k-i)$ and convex corners $(\ell+1-k,k-i-1)$ coloring by $t$ in the terminology of \cite{Ka2}, respectively.
The definitions (\ref{LAdefad}) and (\ref{LAdefre}) are obtained from (4.7) of \cite{Ka2}
by changing $i\mapsto \ell-k$, $j\mapsto k-i$ and $i\mapsto \ell+1-k$, $j\mapsto k-i-1$, respectively.
Thus, Theorem \ref{mainthm1}, \ref{mainthm2} and \ref{mainthm3}
follow from Theorem 4.3 in \cite{Ka} and Theorem 4.3, 4.9 in \cite{Ka2}.
Hence, in what follows, we assume $X\neq A^{(1)}_{n-1}$.
\end{rem}

\section{Action of $S$, $\what{S}$}

Let $X=B^{(1)}_{n-1}$, $C^{(1)}_{n-1}$, $D^{(1)}_{n-1}$, $A^{(2)}_{2n-2}$, $A^{(2)}_{2n-3}$ or
$D^{(2)}_{n}$.

\begin{prop}\label{prop-closedness1}
We assume $\mathfrak{g}$ is of type $X$ and $k\in I$.
Let $t\in I$ and
we assume that a proper Young wall (resp. proper truncated wall) $Y$ of type $X^L$ of ground state (resp. supporting or covering ground state)
$\Lambda_k$ has a $t$-admissible slot
\[
\begin{xy}
(-8,15) *{(-i-1,\ell+1)}="0000",
(17,15) *{(-i,\ell+1)}="000",
(15,-3) *{(-i,\ell)}="00",
(-7,-3) *{(-i-1,\ell)}="0",
(0,0) *{}="1",
(12,0)*{}="2",
(12,12)*{}="3",
(0,12)*{}="4",
\ar@{--} "1";"2"^{}
\ar@{--} "1";"4"^{}
\ar@{--} "2";"3"^{}
\ar@{--} "3";"4"^{}
\end{xy}
\]
and the below block 
\[
\begin{xy}
(-8,15) *{(-i-1,\ell)}="0000",
(17,15) *{(-i,\ell)}="000",
(15,-3) *{(-i,\ell-1)}="00",
(-7,-3) *{(-i-1,\ell-1)}="0",
(0,0) *{}="1",
(12,0)*{}="2",
(12,12)*{}="3",
(0,12)*{}="4",
\ar@{-} "1";"2"^{}
\ar@{-} "1";"4"^{}
\ar@{-} "2";"3"^{}
\ar@{-} "3";"4"^{}
\end{xy}
\]
is not
\[
\begin{xy}
(6,9) *{k}="00",
(6,3) *{k}="0",
(0,0) *{}="1",
(12,0)*{}="2",
(12,12)*{}="3",
(0,12)*{}="4",
(0,6)*{}="5",
(12,6)*{}="6",
\ar@{-} "1";"2"^{}
\ar@{-} "1";"4"^{}
\ar@{-} "2";"3"^{}
\ar@{-} "3";"4"^{}
\ar@{-} "5";"6"^{}
\end{xy}
\]
Let $Y'$ be a wall obtained from $Y$
by adding the $t$-block to the slot.
Then for $s\in\mathbb{Z}_{\geq-1}$
such that $s+P_{(k)}^{X^L}(\ell)+i\geq1$,
it follows
\[
L^{X^L}_{s,k,\iota}(Y')=L^{X^L}_{s,k,\iota}(Y)-\beta_{s+P_{(k)}^{X^L}(\ell)+i,t},
\]
where putting the $t$-admissible slot in our claim as $A$, we set $P_{(k)}^{X^L}(m):=P_{k,A}^{X^L}(m)$ (see (\ref{pka})).
\end{prop}

\nd{\it Proof.}

For a block or slot $P$ in $Y$, we define
\begin{equation}\label{rule1}
L_{s,k}^{X^L}(P):=
\begin{cases}
L_{s,k,{\rm ad}}^{X^L}(P) & \text{if }P\text{ is a single admissible slot,}\\
2L_{s,k,{\rm ad}}^{X^L}(P) & \text{if }P\text{ is a double admissible slot,}\\
-L_{s,k,{\rm re}}^{X^L}(P) & \text{if }P\text{ is a single removable block,}\\
-2L_{s,k,{\rm re}}^{X^L}(P) & \text{if }P\text{ is a double removable block,}\\
0 & \text{otherwise.}
\end{cases}
\end{equation}
If $P$ consists of $P_1$ and $P_2$ ($P_j$ is a block or slot with half-unit thickness or half-unit height for $j=1,2$) then we define 
\begin{equation}\label{rule2}
L_{s,k}^{X^L}(P):=L_{s,k}^{X^L}(P_1)+L_{s,k}^{X^L}(P_2).
\end{equation}
Let $A$ be the $t$-admissible slot in our claim. 
We name the blocks and slots around the slot $A$ in $Y$ as $B$, $C$, $D$ and $E$ as follows:
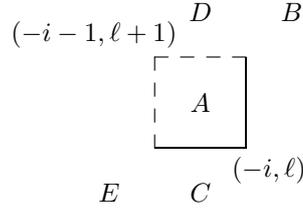
\[
\begin{xy}
(-8,15) *{(-i-1,\ell+1)}="0000",
(15,-3) *{(-i,\ell)}="00",
(6,6) *{A}="A",
(18,18) *{B}="B",
(6,-6) *{C}="C",
(6,18) *{D}="D",
(-6,-6) *{E}="E",
(0,0) *{}="1",
(12,0)*{}="2",
(12,12)*{}="3",
(0,12)*{}="4",
(12,12)*{}="5",
(0,0)*{}="6",
\ar@{-} "1";"2"^{}
\ar@{--} "1";"4"^{}
\ar@{-} "2";"3"^{}
\ar@{--} "3";"4"^{}
\ar@{-} "3";"5"^{}
\ar@{-} "1";"6"^{}
\end{xy}
\]
If $i=0$ then we consider
$B$ as a non-removable block. 
Let $A'$ be the removable $t$-block in $Y'$ and
we also name the blocks and slots around the slot $A'$ in $Y'$ as $B'$, $C'$, $D$ and $E'$ as follows:
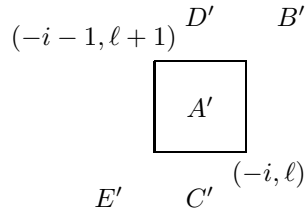
\[
\begin{xy}
(-8,15) *{(-i-1,\ell+1)}="0000",
(15,-3) *{(-i,\ell)}="00",
(6,6) *{A'}="A",
(18,18) *{B'}="B",
(6,-6) *{C'}="C",
(6,18) *{D'}="D",
(-6,-6) *{E'}="E",
(0,0) *{}="1",
(12,0)*{}="2",
(12,12)*{}="3",
(0,12)*{}="4",
(12,12)*{}="5",
(0,0)*{}="6",
\ar@{-} "1";"2"^{}
\ar@{-} "1";"4"^{}
\ar@{-} "2";"3"^{}
\ar@{-} "3";"4"^{}
\ar@{-} "3";"5"^{}
\ar@{-} "1";"6"^{}
\end{xy}
\]
The admissibility and removability of slots and blocks in $Y$ other than $A$, $B$, $C$, $D$ and $E$ are same as
those in $Y'$.
We consider the value of
\[
L^{X^L}_{s,k}(C')+L^{X^L}_{s,k}(E')- (L^{X^L}_{s,k}(C)+L^{X^L}_{s,k}(E)).
\]

If $C$ does not consist of
two blocks with unit width, unit height and half-unit thickness colored by $1$ and $2$ or $n-1$ and $n$
\begin{equation}\label{close-pr1}
\begin{xy}
(50,6) *{\{j_1,j_2\}=\{1,2\}\text{ or }\{j_1,j_2\}=\{n-1,n\}}="000",
(3,9) *{j_1}="00",
(9,3) *{j_2}="0",
(0,0) *{}="1",
(12,0)*{}="2",
(12,12)*{}="3",
(0,12)*{}="4",
(0,0)*{}="5",
(12,12)*{}="6",
\ar@{-} "1";"2"^{}
\ar@{-} "1";"4"^{}
\ar@{-} "2";"3"^{}
\ar@{-} "3";"4"^{}
\ar@{-} "5";"6"^{}
\end{xy}
\end{equation}
then
one can show
\begin{equation}\label{close-pr3}
L^{X^L}_{s,k}(C')+L^{X^L}_{s,k}(E')- (L^{X^L}_{s,k}(C)+L^{X^L}_{s,k}(E))
=c_{X^L}(\ell-1)x_{s+P^{X^L}_{(k)}(\ell)+p_{\pi'_{X^L}(\ell-1),t}+i,\pi'_{X^L}(\ell-1)}.
\end{equation}
 by the same way as in the proof of Proposition 5.7 (i) in \cite{Ka}. Here, just as in subsection \ref{blamsec}, we set
\[
c_X(r)=
\begin{cases}
2 & \text{ if }(X,\pi'_X(r))=(A^{(2)\dagger}_{2n-2},1),\ (B^{(1)}_{n-1},n),\ (D^{(2)}_{n},1) \text{ or }(D^{(2)}_{n},n),\\
1 & \text{otherwise.}
\end{cases}
\]
We consider the case $C$ is in the form (\ref{close-pr1})
and name the $j_1$-block and $j_2$-block as $C_1$ and $C_2$ respectively. 
By Definition \ref{def-YW1}, \ref{def-YW2} and \ref{def-YW3}, there are three cases for the form of $E$:
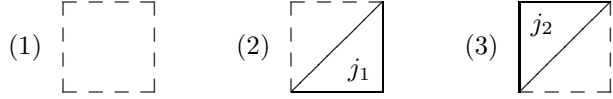
\[
\begin{xy}
(-5,6) *{(1)}="000",
(25,6) *{(2)}="0000",
(55,6) *{(3)}="00000",
(3,9) *{}="00",
(9,3) *{}="0",
(0,0) *{}="1",
(12,0)*{}="2",
(12,12)*{}="3",
(0,12)*{}="4",
(0,0)*{}="5",
(33,9) *{}="00a",
(39,3) *{j_1}="0a",
(30,0) *{}="1a",
(42,0)*{}="2a",
(42,12)*{}="3a",
(30,12)*{}="4a",
(30,0)*{}="5a",
(63,9) *{j_2}="00b",
(69,3) *{}="0b",
(60,0) *{}="1b",
(72,0)*{}="2b",
(72,12)*{}="3b",
(60,12)*{}="4b",
(60,0)*{}="5b",
\ar@{--} "1";"2"^{}
\ar@{--} "1";"4"^{}
\ar@{--} "2";"3"^{}
\ar@{--} "3";"4"^{}
\ar@{-} "1a";"2a"^{}
\ar@{--} "1a";"4a"^{}
\ar@{-} "2a";"3a"^{}
\ar@{--} "3a";"4a"^{}
\ar@{-} "3a";"5a"^{}
\ar@{--} "1b";"2b"^{}
\ar@{-} "1b";"4b"^{}
\ar@{--} "2b";"3b"^{}
\ar@{-} "3b";"4b"^{}
\ar@{-} "3b";"5b"^{}
\end{xy}
\]
We name the $j_1$-block or slot in $E$ as $E_1$ and $j_2$-block or slot in $E$ as $E_2$.
The admissibility and removability of these slots and blocks are as follows:

\begin{table}[H]
  \begin{tabular}{|c|l|l|} \hline
    slot or block & in $Y$ & in $Y'$ \\ \hline
        $C_1$ & (1), (2) $j_1$-removable & normal \\
	& (3) normal & \\ 
    $C_2$ & (1), (3) $j_2$-removable & normal \\
	& (2) normal & \\
    $E_1$ & (1) same as in $Y'$ & (1) same as in $Y$ \\
     & (2) same as in $Y'$ & (2) same as in $Y$  \\
     & (3) normal & (3) $j_1$-admissible \\
    $E_2$ & (1)  same as in $Y'$ & (1) same as in $Y$ \\
    & (2)  normal & (2) $j_2$-admissible \\
     & (3) same as in $Y'$ & (3) same as in $Y$  \\
    \hline
  \end{tabular}
\end{table}
We may assume that $j_2\in\{2,n\}$.
By a direct calculation, it holds
\begin{eqnarray}\label{close-pr4}
L^{X^L}_{s,k}(C')+L^{X^L}_{s,k}(E')- (L^{X^L}_{s,k}(C)+L^{X^L}_{s,k}(E))
&=&x_{s+P^{X^L}_{(k)}(\ell-1)+i+1,j_1}+x_{s+P^{X^L}_{(k)}\left(\ell-\frac{1}{2}\right)+i+1,j_2}\\
&=&x_{s+P^{X^L}_{(k)}(\ell)+i+p_{j_1,t},j_1}+x_{s+P^{X^L}_{(k)}(\ell)+i+p_{j_2,t},j_2}.\nonumber 
\end{eqnarray}
Here, we use $p_{j_1,t}+p_{t,j_1}=1$, $p_{j_2,t}+p_{t,j_2}=1$ and $P^{X^L}_{(k)}(\ell)=P^{X^L}_{(k)}(\ell-1)+p_{t,j_1}=P^{X^L}_{(k)}\left(\ell-\frac{1}{2}\right)+p_{t,j_2}$ in the second equality.
Next, we
consider the value of
\[
L^{X^L}_{s,k}(B')+L^{X^L}_{s,k}(D')- (L^{X^L}_{s,k}(B)+L^{X^L}_{s,k}(D)).
\]
If the slot $D$ does not consist of
two slots with unit width, unit height and half-unit thickness colored by $1$ and $2$ or $n-1$ and $n$
\begin{equation}\label{close-pr2}
\begin{xy}
(50,6) *{\{j_1,j_2\}=\{1,2\}\text{ or }\{j_1,j_2\}=\{n-1,n\}}="000",
(3,9) *{j_1}="00",
(9,3) *{j_2}="0",
(0,0) *{}="1",
(12,0)*{}="2",
(12,12)*{}="3",
(0,12)*{}="4",
(0,0)*{}="5",
(12,12)*{}="6",
\ar@{--} "1";"2"^{}
\ar@{--} "1";"4"^{}
\ar@{--} "2";"3"^{}
\ar@{--} "3";"4"^{}
\ar@{--} "5";"6"^{}
\end{xy}
\end{equation}
then
one can show
\begin{equation}\label{close-pr5}
L^{X^L}_{s,k}(B')+L^{X^L}_{s,k}(D')- (L^{X^L}_{s,k}(B)+L^{X^L}_{s,k}(D))
=c_{X^L}(\ell+1)x_{s+P^{X^L}_{(k)}(\ell)+p_{\pi'_X(\ell+1),t}+i,\pi'_X(\ell+1)}.
\end{equation}
by the same way as in the proof of Proposition 5.7 (i) in \cite{Ka}.
We assume $D$ is in the form (\ref{close-pr2}) and name $j_1$-slot and $j_2$-slot as $D_1$ and $D_2$ respectively.
By Definition \ref{def-YW1}, \ref{def-YW2} and \ref{def-YW3}, there are three cases for the form of $B$:
\[
\begin{xy}
(10,6) *{\text{or }i=0}="000000",
(-15,6) *{(1)}="000",
(25,6) *{(2)}="0000",
(55,6) *{(3)}="00000",
(-7,9) *{j_2}="00",
(-1,3) *{j_1}="0",
(-10,0) *{}="1",
(2,0)*{}="2",
(2,12)*{}="3",
(-10,12)*{}="4",
(-10,0)*{}="5",
(33,9) *{}="00a",
(39,3) *{j_1}="0a",
(30,0) *{}="1a",
(42,0)*{}="2a",
(42,12)*{}="3a",
(30,12)*{}="4a",
(30,0)*{}="5a",
(63,9) *{j_2}="00b",
(69,3) *{}="0b",
(60,0) *{}="1b",
(72,0)*{}="2b",
(72,12)*{}="3b",
(60,12)*{}="4b",
(60,0)*{}="5b",
\ar@{-} "1";"2"^{}
\ar@{-} "1";"4"^{}
\ar@{-} "2";"3"^{}
\ar@{-} "3";"4"^{}
\ar@{-} "1";"3"^{}
\ar@{-} "1a";"2a"^{}
\ar@{--} "1a";"4a"^{}
\ar@{-} "2a";"3a"^{}
\ar@{--} "3a";"4a"^{}
\ar@{-} "3a";"5a"^{}
\ar@{--} "1b";"2b"^{}
\ar@{-} "1b";"4b"^{}
\ar@{--} "2b";"3b"^{}
\ar@{-} "3b";"4b"^{}
\ar@{-} "3b";"5b"^{}
\end{xy}
\]
Let $B_1$ and $B_2$ be the $j_1$-block and $j_2$-block in $B$ respectively.
The admissibility and removability of these slots and blocks are as follows:
\begin{table}[H]
  \begin{tabular}{|c|l|l|} \hline
    slot or block & in $Y$ & in $Y'$ \\ \hline
        $D_1$ & normal & (1), (3) $j_1$-admissible \\
	&  & (2) normal \\ 
    $D_2$ & normal & (1), (2) $j_2$-admissible \\
	&  & (3) normal \\
    $B_1$ & (1) same as in $Y'$ & (1) same as in $Y$ \\
     & (2) removable $j_1$ & (2) normal  \\
     & (3) same as in $Y'$ & (3) same as in $Y$ \\
    $B_2$ & (1) same as in $Y'$ & (1) same as in $Y$ \\
    & (2) same as in $Y'$ & (2) same as in $Y$ \\
     & (3) removable $j_2$ & (3) normal  \\
    \hline
  \end{tabular}
\end{table}

By a direct calculation, it holds
\begin{eqnarray}\label{close-pr6}
L^{X^L}_{s,k}(B')+L^{X^L}_{s,k}(D')- (L^{X^L}_{s,k}(B)+L^{X^L}_{s,k}(D))
&=&x_{s+P^{X^L}_{(k)}(\ell+1)+i,j_1}+x_{s+P^{X^L}_{(k)}\left(\ell+\frac{3}{2}\right)+i,j_2}\\
&=&x_{s+P^{X^L}_{(k)}(\ell)+i+p_{j_1,t},j_1}+x_{s+P^{X^L}_{(k)}(\ell)+i+p_{j_2,t},j_2}. \nonumber
\end{eqnarray}
It is easy to check
\[
-L^{X^L}_{s,k,{\rm re}}(A')-L^{X^L}_{s,k,{\rm ad}}(A)=-x_{s+P^{X^L}_{(k)}(\ell)+i+1,t}-x_{s+P^{X^L}_{(k)}(\ell)+i,t}.
\]
Combining with (\ref{close-pr3}), (\ref{close-pr4}),  (\ref{close-pr5}) and (\ref{close-pr6}),  
we get our claim. \qed

\begin{prop}\label{prop-closedness2}
We assume $\mathfrak{g}$ is of type $X$ and $k\in I$.
For $t\in \{1,n\}$ with $t\neq k$, we assume that a proper Young wall (resp. proper truncated wall) $Y$ of type $X^L$ of ground state (resp. supporting or covering ground state)
$\Lambda_k$ has a $t$-admissible slot
\[
\begin{xy}
(-8,9) *{(-i-1,\ell+\frac{1}{2})}="0000",
(17,9) *{(-i,\ell+\frac{1}{2})}="000",
(15,-3) *{(-i,\ell)}="00",
(-7,-3) *{(-i-1,\ell)}="0",
(0,0) *{}="1",
(12,0)*{}="2",
(12,6)*{}="3",
(0,6)*{}="4",
(38,3)*{{\rm or}}="or",
(60,9) *{(-i-1,\ell+1)}="a0000",
(85,9) *{(-i,\ell+1)}="a000",
(83,-3) *{(-i,\ell+\frac{1}{2})}="a00",
(61,-3) *{(-i-1,\ell+\frac{1}{2})}="a0",
(68,0) *{}="a1",
(80,0)*{}="a2",
(80,6)*{}="a3",
(68,6)*{}="a4",
\ar@{--} "1";"2"^{}
\ar@{--} "1";"4"^{}
\ar@{--} "2";"3"^{}
\ar@{--} "3";"4"^{}
\ar@{--} "a1";"a2"^{}
\ar@{--} "a1";"a4"^{}
\ar@{--} "a2";"a3"^{}
\ar@{--} "a3";"a4"^{}
\end{xy}
\]
and the below block 
\[
\begin{xy}
(-8,15) *{(-i-1,\ell)}="0000",
(17,15) *{(-i,\ell)}="000",
(15,-3) *{(-i,\ell-1)}="00",
(-7,-3) *{(-i-1,\ell-1)}="0",
(0,0) *{}="1",
(12,0)*{}="2",
(12,12)*{}="3",
(0,12)*{}="4",
\ar@{-} "1";"2"^{}
\ar@{-} "1";"4"^{}
\ar@{-} "2";"3"^{}
\ar@{-} "3";"4"^{}
\end{xy}
\]
is not
\[
\begin{xy}
(6,9) *{k}="00",
(6,3) *{k}="0",
(0,0) *{}="1",
(12,0)*{}="2",
(12,12)*{}="3",
(0,12)*{}="4",
(0,6)*{}="5",
(12,6)*{}="6",
\ar@{-} "1";"2"^{}
\ar@{-} "1";"4"^{}
\ar@{-} "2";"3"^{}
\ar@{-} "3";"4"^{}
\ar@{-} "5";"6"^{}
\end{xy}
\]
Let $Y'$ be a wall obtained from $Y$
by adding the $t$-block to the slot.
For $s\in\mathbb{Z}_{\geq-1}$ such that $s+P_{(k)}^{X^L}(\ell)+i\geq1$,
it follows
\[
L^{X^L}_{s,k,\iota}(Y')=L^{X^L}_{s,k,\iota}(Y)-\beta_{s+P_{(k)}^{X^L}(\ell)+i,t},
\]
where we define $P_{(k)}^{X^L}(m)$ just as in the previous proposition.
\end{prop}

\nd
{\it Proof.}

We can show the claim by the same way as in the proof of Proposition 5.7 (ii) in \cite{Ka}. \qed

\begin{prop}\label{prop-closedness3}
We assume $\mathfrak{g}$ is of type $X$ and $k\in I$.
For $t\in \{1,2,n-1,n\}$, we assume that a proper Young wall (resp. proper truncated wall) $Y$ of type $X^L$ of ground state (resp. supporting or covering ground state)
$\Lambda_k$ has a $t$-admissible slot
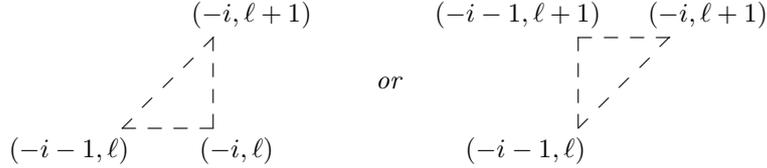
\[
\begin{xy}
(35,6) *{\text{or}}="0",
(17,15) *{(-i,\ell+1)}="000",
(15,-3) *{(-i,\ell)}="00",
(-7,-3) *{(-i-1,\ell)}="0",
(0,0) *{}="1",
(12,0)*{}="2",
(12,12)*{}="3",
(52,15) *{(-i-1,\ell+1)}="0000a",
(77,15) *{(-i,\ell+1)}="000a",
(53,-3) *{(-i-1,\ell)}="0a",
(60,0) *{}="1a",
(72,0)*{}="2a",
(72,12)*{}="3a",
(60,12)*{}="4a",
\ar@{--} "1";"2"^{}
\ar@{--} "1";"3"^{}
\ar@{--} "2";"3"^{}
\ar@{--} "1a";"4a"^{}
\ar@{--} "1a";"3a"^{}
\ar@{--} "3a";"4a"^{}
\end{xy}
\]
and the below block 
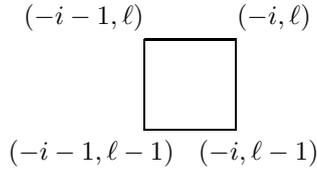
\[
\begin{xy}
(-8,15) *{(-i-1,\ell)}="0000",
(17,15) *{(-i,\ell)}="000",
(15,-3) *{(-i,\ell-1)}="00",
(-7,-3) *{(-i-1,\ell-1)}="0",
(0,0) *{}="1",
(12,0)*{}="2",
(12,12)*{}="3",
(0,12)*{}="4",
\ar@{-} "1";"2"^{}
\ar@{-} "1";"4"^{}
\ar@{-} "2";"3"^{}
\ar@{-} "3";"4"^{}
\end{xy}
\]
is not
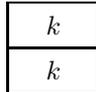
\[
\begin{xy}
(6,9) *{k}="00",
(6,3) *{k}="0",
(0,0) *{}="1",
(12,0)*{}="2",
(12,12)*{}="3",
(0,12)*{}="4",
(0,6)*{}="5",
(12,6)*{}="6",
\ar@{-} "1";"2"^{}
\ar@{-} "1";"4"^{}
\ar@{-} "2";"3"^{}
\ar@{-} "3";"4"^{}
\ar@{-} "5";"6"^{}
\end{xy}
\]
Let $Y'$ be a wall obtained from $Y$
by adding the $t$-block to the slot.
Then for $s\in\mathbb{Z}_{\geq-1}$
such that $s+P_k^{X^L}(\ell)+i\geq1$ (resp. $s+P_{(k)}^{X^L}\left(\ell+\frac{1}{2}\right)+i\geq1$) when $t\in\{1,n-1\}$ (resp. $t\in\{2,n\}$),
it follows
\[
L^{X^L}_{s,k,\iota}(Y')=L^{X^L}_{s,k,\iota}(Y)-\beta_{s+P_{(k)}^{X^L}(\ell)+i,t}\quad (\text{if }t\in\{1,n-1\}),
\]
\[
L^{X^L}_{s,k,\iota}(Y')=L^{X^L}_{s,k,\iota}(Y)-\beta_{s+P_{(k)}^{X^L}\left(\ell+\frac{1}{2}\right)+i,t}\quad (\text{if }t\in\{2,n\}),
\]
where we define $P_{(k)}^{X^L}(m)$ just as in the previous proposition.
\end{prop}

\nd{\it Proof.}

We only consider the first $t$-admissible slot in our claim
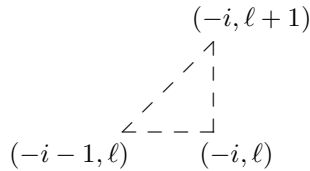
\[
\begin{xy}
(17,15) *{(-i,\ell+1)}="000",
(15,-3) *{(-i,\ell)}="00",
(-7,-3) *{(-i-1,\ell)}="0",
(0,0) *{}="1",
(12,0)*{}="2",
(12,12)*{}="3",
\ar@{--} "1";"2"^{}
\ar@{--} "1";"3"^{}
\ar@{--} "2";"3"^{}
\end{xy}
\]
since one can similarly treat the second admissible slot.
Let $A_1$ be this $t$-admissible slot. One can take $u\in I$ such that $\{u,t\}=\{1,2\}$ or $\{u,t\}=\{n-1,n\}$.
We name the blocks and slots around the slot $A_1$ in $Y$ as follows:
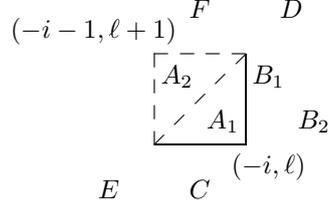
\[
\begin{xy}
(-8,15) *{(-i-1,\ell+1)}="0000",
(15,-3) *{(-i,\ell)}="00",
(9,3) *{A_1}="A",
(3,9) *{A_2}="AA",
(21,3) *{B_2}="B",
(15,9) *{B_1}="BB",
(6,18) *{F}="F",
(18,18) *{D}="D",
(6,-6) *{C}="C",
(-6,-6) *{E}="E",
(0,0) *{}="1",
(12,0)*{}="2",
(12,12)*{}="3",
(0,12)*{}="4",
(12,12)*{}="5",
(0,0)*{}="6",
\ar@{-} "1";"2"^{}
\ar@{--} "1";"4"^{}
\ar@{-} "2";"3"^{}
\ar@{--} "3";"4"^{}
\ar@{--} "3";"1"^{}
\ar@{-} "3";"5"^{}
\ar@{-} "1";"6"^{}
\end{xy}
\]
If $i=0$ then we consider
$B_1$, $B_2$, $D$ as non-removable blocks. 
By Definition \ref{def-YW1}, \ref{def-YW2} and \ref{def-YW3},
$B_2$ is an $u$-block.
We can consider the following five cases:

Case 1 : $B_1$ is a block, $D$ is a block, $A_2$ is a block, $E$ is a slot,

Case 2 : $B_1$ is a block, $D$ is a block, $A_2$ is a block, $E$ is a block,

Case 3 : $B_1$ is a block, $D$ is a block, $A_2$ is a slot,

Case 4 : $B_1$ is a block, $D$ is a slot, $A_2$ is a slot,

Case 5 : $B_1$ is a slot, $D$ is a slot, $A_2$ is a slot.

\vspace{3mm}

\underline{Case 1 : $B_1$ is a block, $D$ is a block, $A_2$ is a block, $E$ is a slot}

\vspace{2mm}

In this case, the admissibility and removability of $A_2$, $B_1$, $B_2$, $C$ and $E$ in $Y$ are same
as those in $Y'$. Those of $A_1$, $D$ and $F$ are as follows:

\begin{table}[H]
  \begin{tabular}{|c|l|l|} \hline
    slot or block & in $Y$ & in $Y'$ \\ \hline
        $A_1$ & $t$-admissible & removable $t$ \\
        $D$ & normal (resp. removable $\pi'_{X^L}(\ell+1)$) & normal \\
    $F$ & normal & $\pi'_{X^L}(\ell+1)$-admissible (resp. normal) \\
    \hline
  \end{tabular}
\end{table}
\nd
The admissibility and removability of other slots and blocks in $Y$ are same
as those in $Y'$. 
By a direct calculation, we obtain
\[
L^{X^L}_{s,k,\iota}(Y')=L^{X^L}_{s,k,\iota}(Y)-\beta_{s+P_{(k)}^{X^L}(\ell)+i,t}\quad (\text{if }t\in\{1,n-1\}),
\]
\[
L^{X^L}_{s,k,\iota}(Y')=L^{X^L}_{s,k,\iota}(Y)-\beta_{s+P_{(k)}^{X^L}\left(\ell+\frac{1}{2}\right)+i,t}\quad (\text{if }t\in\{2,n\}).
\]
We can similarly prove above formulas in other cases so that
write only the table of admissibility and removability for each case.
The admissibility and removability of $A_1$ are same as in Case 1.
\vspace{1mm}

\underline{Case 2 : $B_1$ is a block, $D$ is a block, $A_2$ is a block, $E$ is a block}

\vspace{2mm}

We name two blocks as $G_1$, $G_2$ as follows:
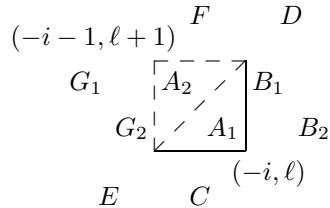
\[
\begin{xy}
(-8,15) *{(-i-1,\ell+1)}="0000",
(15,-3) *{(-i,\ell)}="00",
(9,3) *{A_1}="A",
(3,9) *{A_2}="AA",
(-3,3) *{G_2}="G",
(-9,9) *{G_1}="GG",
(21,3) *{B_2}="B",
(15,9) *{B_1}="BB",
(6,18) *{F}="F",
(18,18) *{D}="D",
(6,-6) *{C}="C",
(-6,-6) *{E}="E",
(0,0) *{}="1",
(12,0)*{}="2",
(12,12)*{}="3",
(0,12)*{}="4",
(12,12)*{}="5",
(0,0)*{}="6",
\ar@{-} "1";"2"^{}
\ar@{--} "1";"4"^{}
\ar@{-} "2";"3"^{}
\ar@{--} "3";"4"^{}
\ar@{--} "3";"1"^{}
\ar@{-} "3";"5"^{}
\ar@{-} "1";"6"^{}
\end{xy}
\]
By Definition \ref{def-YW1} (iv) and \ref{def-YW3} (iii), we see that $G_2$ is a slot.
In this case, the admissibility and removability of $B_1$, $B_2$, $C$ and $E$ in $Y$ are same
as those in $Y'$. Those of $A_1$, $D$ and $F$ are same as Case 1. 
There are two patterns for $G_1$ : (1) $G_1$ is a block, (2) $G_1$ is a slot.
Those of $A_2$, $G_1$ and $G_2$ are as follows:

\begin{table}[H]
  \begin{tabular}{|c|l|l|} \hline
    slot or block & in $Y$ & in $Y'$ \\ \hline
        $A_2$ & normal & (1) normal \\
         & & (2) removable $u$ \\
        $G_1$ & same as in $Y'$ & same as in $Y$ \\
        $G_2$ &normal & (1) normal \\
         & & (2) $u$-admissible \\
    \hline
  \end{tabular}
\end{table}

\underline{Case 3 : $B_1$ is a block, $D$ is a block, $A_2$ is a slot}

\vspace{2mm}

In this case, the admissibility and removability of $A_2$, $B_1$, $B_2$, $D$ and $F$ in $Y$ are same
as those in $Y'$. Those of $C$ and $E$ are as follows:

\begin{table}[H]
  \begin{tabular}{|c|l|l|} \hline
    slot or block & in $Y$ & in $Y'$ \\ \hline
        $C$ & normal (resp. removable $\pi'_{X^L}(\ell-1)$) & normal \\
    $E$ & normal & $\pi'_{X^L}(\ell-1)$-admissible (resp. normal) \\
    \hline
  \end{tabular}
\end{table}

\vspace{1mm}

\underline{Case 4 : $B_1$ is a block, $D$ is a slot, $A_2$ is a slot}

\vspace{2mm}

In this case, the admissibility and removability of $B_1$, $D$ and $F$ in $Y$ are same
as those in $Y'$. Those of $A_2$, $B_2$, $C$ and $E$ are as follows:

\begin{table}[H]
  \begin{tabular}{|c|l|l|} \hline
    slot or block & in $Y$ & in $Y'$ \\ \hline
        $A_2$ & $u$-admissible & normal \\
        $B_2$ & removable $u$ & normal \\
        $C$ & normal (resp. removable $\pi'_{X^L}(\ell-1)$) & normal \\
    $E$ & normal & $\pi'_{X^L}(\ell-1)$-admissible (resp. normal) \\
    \hline
  \end{tabular}
\end{table}

\vspace{1mm}

\underline{Case 5 : $B_1$ is a slot, $D$ is a slot, $A_2$ is a slot}

\vspace{2mm}

In this case, the admissibility and removability of $A_2$, $B_1$, $B_2$, $D$ and $F$ in $Y$ are same
as those in $Y'$. Those of $C$ and $E$ are same as in Case 3. \qed

\begin{prop}\label{prop-closedness4}
We assume $\mathfrak{g}$ is of type $X$. Let $k\in I$ in class $2$. 
Let $t\in I$ and $Y=(\overline{Y},\overline{\overline{Y}})\in {\rm YW}_k^{X^L}$. 
We assume that $\overline{Y}$ or $\overline{\overline{Y}}$
has a $t$-admissible slot
\[
\begin{xy}
(-8,15) *{(-i-1,\ell+1)}="0000",
(17,15) *{(-i,\ell+1)}="000",
(15,-3) *{(-i,\ell)}="00",
(-7,-3) *{(-i-1,\ell)}="0",
(0,0) *{}="1",
(12,0)*{}="2",
(12,12)*{}="3",
(0,12)*{}="4",
\ar@{--} "1";"2"^{}
\ar@{--} "1";"4"^{}
\ar@{--} "2";"3"^{}
\ar@{--} "3";"4"^{}
\end{xy}
\]
and the below block is
\[
\begin{xy}
(-8,15) *{(-i-1,\ell)}="0000a",
(17,15) *{(-i,\ell)}="000a",
(15,-3) *{(-i,\ell-1)}="00a",
(-7,-3) *{(-i-1,\ell-1)}="0a",
(6,9) *{k}="00",
(6,3) *{k}="0",
(0,0) *{}="1",
(12,0)*{}="2",
(12,12)*{}="3",
(0,12)*{}="4",
(0,6)*{}="5",
(12,6)*{}="6",
\ar@{-} "1";"2"^{}
\ar@{-} "1";"4"^{}
\ar@{-} "2";"3"^{}
\ar@{-} "3";"4"^{}
\ar@{-} "5";"6"^{}
\end{xy}
\]
Let $Y'\in {\rm YW}^{X^L}_k$ be the pair of proper truncated walls obtained from $Y$
by adding the $t$-block to the slot.
Then for $s\in\mathbb{Z}_{\geq-1}$
such that $s+P_{(k)}^{X^L}(\ell)+i\geq1$,
it follows
\[
L^{X^L}_{s,k,\iota}(Y')=L^{X^L}_{s,k,\iota}(Y)-\beta_{s+P_{(k)}^{X^L}(\ell)+i,t},
\]
where $P_{(k)}^{X^L}(m)$ is defined as in the previous proposition.
\end{prop}

\nd{\it Proof.}

We only consider the case the $t$-admissible slot in our claim is in $\ovl{Y}$ since one can similarly show the case
the $t$-admissible slot is in $\ovl{Y}$.
Let $\ovl{A}$ be the $t$-admissible slot. 
Since the slot $\ovl{A}$ is $t$-admissible, by Definition \ref{def-YW1}, \ref{def-YW2} and \ref{def-YW3}, blocks and slots around $\ovl{A}$ are as follows:

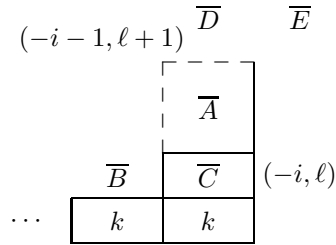
\[
\begin{xy}
(-8,15) *{(-i-1,\ell+1)}="0000",
(18,-3) *{(-i,\ell)}="00",
(6,6) *{\ovl{A}}="A",
(18,18) *{\ovl{E}}="E",
(6,-3) *{\ovl{C}}="C",
(-6,-3) *{\ovl{B}}="B",
(6,-9) *{k}="CC",
(6,18) *{\ovl{D}}="D",
(-6,-9) *{k}="BB",
(-18,-9) *{\cdots}="BBB",
(0,0) *{}="1",
(12,0)*{}="2",
(12,12)*{}="3",
(0,12)*{}="4",
(12,12)*{}="5",
(0,0)*{}="6",
(0,-12) *{}="11",
(12,-12) *{}="22",
(0,-6) *{}="111",
(12,-6) *{}="222",
(-12,-6) *{}="1111",
(-12,-12) *{}="11111",
\ar@{-} "11";"22"^{}
\ar@{-} "1";"11"^{}
\ar@{-} "22";"2"^{}
\ar@{-} "1";"2"^{}
\ar@{-} "111";"222"^{}
\ar@{-} "1111";"11111"^{}
\ar@{-} "1111";"111"^{}
\ar@{-} "11";"11111"^{}
\ar@{--} "1";"4"^{}
\ar@{-} "2";"3"^{}
\ar@{--} "3";"4"^{}
\ar@{-} "3";"5"^{}
\ar@{-} "1";"6"^{}
\end{xy}
\]
Here, we named blocks and slots as above and
$\ovl{C}$ is the $k$-block
with unit width, half-unit height and unit thickness.
Let $\ovl{\ovl{B}}$ (resp. $\ovl{\ovl{C}}$) be a slot (resp. block) with half-unit thickness such that the position of the pair $B:=(\ovl{B},\ovl{\ovl{B}})$
(resp. $C:=(\ovl{C},\ovl{\ovl{C}})$) is $(-i-1,\ovl{T}^{X^L}_k)$ (resp. $(-i,\ovl{T}^{X^L}_k)$).
The admissibility and removability of $B$ and $C$ are as follows:

\begin{table}[H]
  \begin{tabular}{|c|l|l|} \hline
    slot or block & in $Y$ & in $Y'$ \\ \hline
        $B$ & normal & $k$-admissible pair (resp. normal) \\
    $C$ & normal (resp. removable $k$-pair) & normal \\
    \hline
  \end{tabular}
\end{table}
\nd
Hence, we have
\[
-L^{X^L}_{s,k,{\rm re}}(C')+L^{X^L}_{s,k,{\rm ad}}(B')- (-L^{X^L}_{s,k,{\rm re}}(C)+L^{X^L}_{s,k,{\rm ad}}(B))
=x_{s+i+1,k}=x_{s+i+P_{(k)}^{X^L}(\ell)+p_{k,t},k},
\]
where $B'$ and $C'$ are corresponding pairs in $Y'$.
Let $\ovl{A}'$, $\ovl{D}'$ and $\ovl{E}'$ be the corresponding blocks and slots in $Y'$ to $\ovl{A}$, $\ovl{D}$ and $\ovl{E}$.
One can compute
\[
L^{X^L}_{s,k}(\ovl{E}')+L^{X^L}_{s,k}(\ovl{D})- (L^{X^L}_{s,k}(\ovl{E})+L^{X^L}_{s,k}(\ovl{D}))
\]
by the same way as in (\ref{close-pr5}) and (\ref{close-pr6}). Here, we used the same notation as in (\ref{rule1}).
It is easy to check
\[
-L^{X^L}_{s,k,{\rm re}}(\ovl{A}')-L^{X^L}_{s,k,{\rm ad}}(\ovl{A})=-x_{s+P^{X^L}_{(k)}(\ell)+i+1,t}-x_{s+P^{X^L}_{(k)}(\ell)+i,t}.
\]
Combining the above calculations, we get our claim. \qed

\begin{prop}\label{prop-closedness5}
We assume $\mathfrak{g}$ is of type $X$. Let $k\in I$ in class $2$. 
Let $t\in I$ and $Y=(\overline{Y},\overline{\overline{Y}})\in {\rm YW}_k^{X^L}$. 
We assume that $\overline{Y}$ or $\overline{\overline{Y}}$
has a $t$-admissible slot
\[
\begin{xy}
(-8,9) *{(-i-1,\ell+\frac{1}{2})}="0000",
(17,9) *{(-i,\ell+\frac{1}{2})}="000",
(15,-3) *{(-i,\ell)}="00",
(-7,-3) *{(-i-1,\ell)}="0",
(0,0) *{}="1",
(12,0)*{}="2",
(12,6)*{}="3",
(0,6)*{}="4",
(38,3)*{{\rm or}}="or",
(60,9) *{(-i-1,\ell+1)}="a0000",
(85,9) *{(-i,\ell+1)}="a000",
(83,-3) *{(-i,\ell+\frac{1}{2})}="a00",
(61,-3) *{(-i-1,\ell+\frac{1}{2})}="a0",
(68,0) *{}="a1",
(80,0)*{}="a2",
(80,6)*{}="a3",
(68,6)*{}="a4",
\ar@{--} "1";"2"^{}
\ar@{--} "1";"4"^{}
\ar@{--} "2";"3"^{}
\ar@{--} "3";"4"^{}
\ar@{--} "a1";"a2"^{}
\ar@{--} "a1";"a4"^{}
\ar@{--} "a2";"a3"^{}
\ar@{--} "a3";"a4"^{}
\end{xy}
\]
and the below block is
\[
\begin{xy}
(-8,15) *{(-i-1,\ell)}="0000a",
(17,15) *{(-i,\ell)}="000a",
(15,-3) *{(-i,\ell-1)}="00a",
(-7,-3) *{(-i-1,\ell-1)}="0a",
(6,9) *{k}="00",
(6,3) *{k}="0",
(0,0) *{}="1",
(12,0)*{}="2",
(12,12)*{}="3",
(0,12)*{}="4",
(0,6)*{}="5",
(12,6)*{}="6",
\ar@{-} "1";"2"^{}
\ar@{-} "1";"4"^{}
\ar@{-} "2";"3"^{}
\ar@{-} "3";"4"^{}
\ar@{-} "5";"6"^{}
\end{xy}
\]
Let $Y'\in {\rm YW}^{X^L}_k$ be the pair of proper truncated walls obtained from $Y$
by adding the $t$-block to the slot.
Then for $s\in\mathbb{Z}_{\geq-1}$
such that $s+P_{(k)}^{X^L}(\ell)+i\geq1$,
it follows
\[
L^{X^L}_{s,k,\iota}(Y')=L^{X^L}_{s,k,\iota}(Y)-\beta_{s+P_{(k)}^{X^L}(\ell)+i,t},
\]
where $P_{(k)}^{X^L}(\ell)$ is same as in the previous proposition.
\end{prop}

\nd{\it Proof.}

We only consider the case the slot in our claim is in $\ovl{Y}$.
First, we consider the case the $t$-admissible slot is
\[
\begin{xy}
(-8,9) *{(-i-1,\ell+\frac{1}{2})}="0000",
(17,9) *{(-i,\ell+\frac{1}{2})}="000",
(15,-3) *{(-i,\ell)}="00",
(-7,-3) *{(-i-1,\ell)}="0",
(0,0) *{}="1",
(12,0)*{}="2",
(12,6)*{}="3",
(0,6)*{}="4",
\ar@{--} "1";"2"^{}
\ar@{--} "1";"4"^{}
\ar@{--} "2";"3"^{}
\ar@{--} "3";"4"^{}
\end{xy}
\]
and name this slot as $\ovl{A}_1$. Blocks and slots around $\ovl{A}_1$ are as follows:
\[
\begin{xy}
(-10,9) *{(-i-1,\ell+\frac{1}{2})}="0000",
(18,-3) *{(-i,\ell)}="00",
(6,9) *{\ovl{A}_2}="AA",
(6,3) *{\ovl{A}_1}="A",
(18,18) *{\ovl{D}_3}="DDD",
(18,9) *{\ovl{D}_2}="DD",
(18,3) *{\ovl{D}_1}="D",
(6,-3) *{\ovl{C}}="C",
(-6,-3) *{\ovl{B}}="B",
(6,-9) *{k}="CC",
(-6,-9) *{k}="BB",
(-18,-9) *{\cdots}="BBB",
(0,0) *{}="1",
(12,0)*{}="2",
(12,6)*{}="3",
(0,6)*{}="4",
(12,6)*{}="5",
(12,6)*{}="55",
(24,6)*{}="555",
(0,0)*{}="6",
(0,-12) *{}="11",
(12,-12) *{}="22",
(0,-6) *{}="111",
(12,-6) *{}="222",
(-12,-6) *{}="1111",
(-12,-12) *{}="11111",
\ar@{-} "11";"22"^{}
\ar@{-} "1";"11"^{}
\ar@{-} "22";"2"^{}
\ar@{-} "1";"2"^{}
\ar@{-} "111";"222"^{}
\ar@{-} "1111";"11111"^{}
\ar@{-} "1111";"111"^{}
\ar@{-} "11";"11111"^{}
\ar@{--} "1";"4"^{}
\ar@{-} "2";"3"^{}
\ar@{--} "3";"4"^{}
\ar@{-} "3";"5"^{}
\ar@{-} "55";"5"^{}
\ar@{-} "55";"555"^{}
\ar@{-} "1";"6"^{}
\end{xy}
\]
Here, we named blocks and slots as above.
Let $\ovl{\ovl{B}}$ (resp. $\ovl{\ovl{C}}$) be a slot (resp. block) with half-unit thickness such that the position of the pair $B:=(\ovl{B},\ovl{\ovl{B}})$
(resp. $C:=(\ovl{C},\ovl{\ovl{C}})$) is $(-i-1,\ovl{T}^{X^L}_k)$ (resp. $(-i,\ovl{T}^{X^L}_k)$).
There are three cases (1) $\ovl{D}_2$ and $\ovl{D}_3$ are slots, (2) $\ovl{D}_2$ is a block and $\ovl{D}_3$ is a slot, (3) $\ovl{D}_2$ and $\ovl{D}_3$ are blocks.
The admissibility and removability of $\ovl{A}_1$, $\ovl{A}_2$, $\ovl{D}_1$, $\ovl{D}_2$ and $\ovl{D}_3$ for each case are as follows:
\begin{table}[H]
  \begin{tabular}{|c|l|l|} \hline
    slot or block & in $Y$ & in $Y'$ \\ \hline
        $\ovl{A}_1$ & (1), (2) single $t$-admissible & single removable $t$ \\
	& (3) double $t$-admissible & \\ 
    $\ovl{A}_2$ & normal & (1), (2) normal \\
	&  & (3) single $t$-admissible \\
    $\ovl{D}_1$ & normal & normal \\
    $\ovl{D}_2$ & (1) same as in $Y'$ & (1) same as in $Y$ \\
    & (2) single removable $t$ & (2) single removable $t$ \\
     & (3) normal & (3) normal  \\
    $\ovl{D}_3$ & (1) normal & (1) normal  \\
    & (2), (3) same as in $Y'$ & (2), (3) same as in $Y$ \\
    \hline
  \end{tabular}
\end{table}
The admissibility and removability of $B$ and $C$ are same as in the proof of Proposition 5.4:
\begin{table}[H]
  \begin{tabular}{|c|l|l|} \hline
    slot or block & in $Y$ & in $Y'$ \\ \hline
        $B$ & normal & $k$-admissible pair (resp. normal) \\
    $C$ & normal (resp. removable $k$-pair) & normal \\
    \hline
  \end{tabular}
\end{table}
\nd
Thus, our claim is a consequence of a direct calculation.

Next, we consider the case our $t$-admissible slot is
\[
\begin{xy}
(60,9) *{(-i-1,\ell+1)}="a0000",
(85,9) *{(-i,\ell+1)}="a000",
(83,-3) *{(-i,\ell+\frac{1}{2})}="a00",
(61,-3) *{(-i-1,\ell+\frac{1}{2})}="a0",
(68,0) *{}="a1",
(80,0)*{}="a2",
(80,6)*{}="a3",
(68,6)*{}="a4",
\ar@{--} "a1";"a2"^{}
\ar@{--} "a1";"a4"^{}
\ar@{--} "a2";"a3"^{}
\ar@{--} "a3";"a4"^{}
\end{xy}
\]
and name this slot as $\ovl{A}_2$. Blocks and slots around $\ovl{A}_2$ are as follows:
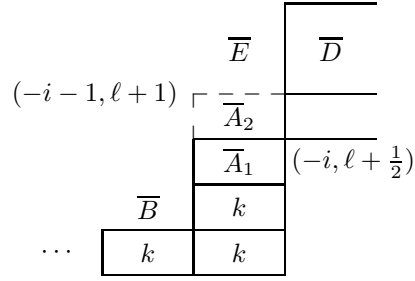
\[
\begin{xy}
(-13,12) *{(-i-1,\ell+1)}="0000",
(21,3) *{(-i,\ell+\frac{1}{2})}="00",
(6,9) *{\ovl{A}_2}="AA",
(6,3) *{\ovl{A}_1}="A",
(18,18) *{\ovl{D}}="D",
(6,18) *{\ovl{E}}="E",
(6,-3) *{k}="C",
(-6,-3) *{\ovl{B}}="B",
(6,-9) *{k}="CC",
(-6,-9) *{k}="BB",
(-18,-9) *{\cdots}="BBB",
(0,0) *{}="1",
(12,0)*{}="2",
(12,6)*{}="3",
(12,12)*{}="3a",
(24,12)*{}="3aa",
(12,24)*{}="3aaa",
(24,24)*{}="3aaaa",
(0,6)*{}="4",
(0,12)*{}="4a",
(12,12)*{}="4aa",
(12,6)*{}="5",
(12,6)*{}="55",
(24,6)*{}="555",
(0,0)*{}="6",
(0,-12) *{}="11",
(12,-12) *{}="22",
(0,-6) *{}="111",
(12,-6) *{}="222",
(-12,-6) *{}="1111",
(-12,-12) *{}="11111",
\ar@{-} "3";"3a"^{}
\ar@{-} "3a";"3aa"^{}
\ar@{-} "3a";"3aaa"^{}
\ar@{-} "3aaa";"3aaaa"^{}
\ar@{--} "4a";"4"^{}
\ar@{--} "4aa";"4a"^{}
\ar@{-} "11";"22"^{}
\ar@{-} "1";"11"^{}
\ar@{-} "22";"2"^{}
\ar@{-} "1";"2"^{}
\ar@{-} "111";"222"^{}
\ar@{-} "1111";"11111"^{}
\ar@{-} "1111";"111"^{}
\ar@{-} "11";"11111"^{}
\ar@{-} "1";"4"^{}
\ar@{-} "2";"3"^{}
\ar@{-} "3";"4"^{}
\ar@{-} "3";"5"^{}
\ar@{-} "55";"5"^{}
\ar@{-} "55";"555"^{}
\ar@{-} "1";"6"^{}
\end{xy}
\]
Let $\ovl{\ovl{B}}$ be a slot with half-unit thickness such that the position of the pair $B:=(\ovl{B},\ovl{\ovl{B}})$ is $(-i-1,\ovl{T}^{X^L}_k)$.
There are two cases : (1) $\ovl{B}$ is a slot, (2) $\ovl{B}$ is a block.
The admissibility and removability of $\ovl{A}_1$, $\ovl{A}_2$ and $B$ are as follows: 
\begin{table}[H]
  \begin{tabular}{|c|l|l|} \hline
    slot or block & in $Y$ & in $Y'$ \\ \hline
        $\ovl{A}_1$ & (1) single removable $t$ & normal \\
	& (2) normal & \\ 
    $\ovl{A}_2$ & single $t$-admissible & (1) double removable $t$ \\
	&  & (2) single removable $t$ \\
    $B$ & same as in $Y'$ & same as in $Y$ \\
    \hline
  \end{tabular}
\end{table}
The admissibility and removability of $\ovl{D}$ and $\ovl{E}$ are as follows:
\begin{table}[H]
  \begin{tabular}{|c|l|l|} \hline
    slot or block & in $Y$ & in $Y'$ \\ \hline
        $\ovl{D}$ & removable $\pi'_{X^L}(\ell+1)$ (resp. normal) & normal \\
    $\ovl{E}$ & normal & normal (resp. $\pi'_{X^L}(\ell+1)$-admissible) \\
    \hline
  \end{tabular}
\end{table}
Thus, we get
\[
-L^{X^L}_{s,k,{\rm re}}(D')+L^{X^L}_{s,k,{\rm ad}}(E')- (-L^{X^L}_{s,k,{\rm re}}(D)+L^{X^L}_{s,k,{\rm ad}}(E))
=x_{s+P^{X^L}_{(k)}(\ell)+p_{\pi'_{X^L}(\ell+1),\pi'_{X^L}(\ell)}+i,\pi'_{X^L}(\ell+1)},
\]
where $D'$ and $E'$ are
the block and slot in $Y'$ corresponding to $D$ and $E$, respectively.
A direct calculation yields our claim. \qed

\begin{prop}\label{prop-closedness6}
We assume $\mathfrak{g}$ is of type $X$. Let $k\in I$ in class $2$. 
Let $t\in I$ and $Y=(\overline{Y},\overline{\overline{Y}})\in {\rm YW}_k^{X^L}$. 
We assume that $\overline{Y}$ or $\overline{\overline{Y}}$
has a $t$-admissible slot
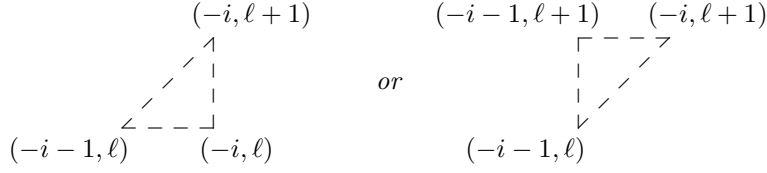
\[
\begin{xy}
(35,6) *{\text{or}}="0",
(17,15) *{(-i,\ell+1)}="000",
(15,-3) *{(-i,\ell)}="00",
(-7,-3) *{(-i-1,\ell)}="0",
(0,0) *{}="1",
(12,0)*{}="2",
(12,12)*{}="3",
(52,15) *{(-i-1,\ell+1)}="0000a",
(77,15) *{(-i,\ell+1)}="000a",
(53,-3) *{(-i-1,\ell)}="0a",
(60,0) *{}="1a",
(72,0)*{}="2a",
(72,12)*{}="3a",
(60,12)*{}="4a",
\ar@{--} "1";"2"^{}
\ar@{--} "1";"3"^{}
\ar@{--} "2";"3"^{}
\ar@{--} "1a";"4a"^{}
\ar@{--} "1a";"3a"^{}
\ar@{--} "3a";"4a"^{}
\end{xy}
\]
and the below block is
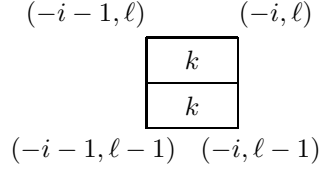
\[
\begin{xy}
(-8,15) *{(-i-1,\ell)}="0000a",
(17,15) *{(-i,\ell)}="000a",
(15,-3) *{(-i,\ell-1)}="00a",
(-7,-3) *{(-i-1,\ell-1)}="0a",
(6,9) *{k}="00",
(6,3) *{k}="0",
(0,0) *{}="1",
(12,0)*{}="2",
(12,12)*{}="3",
(0,12)*{}="4",
(0,6)*{}="5",
(12,6)*{}="6",
\ar@{-} "1";"2"^{}
\ar@{-} "1";"4"^{}
\ar@{-} "2";"3"^{}
\ar@{-} "3";"4"^{}
\ar@{-} "5";"6"^{}
\end{xy}
\]
Let $Y'\in {\rm YW}^{X^L}_k$ be the pair of proper truncated walls obtained from $Y$
by adding the $t$-block to the slot.
Then for $s\in\mathbb{Z}_{\geq-1}$
such that $s+P_{(k)}^{X^L}(\ell)+i\geq1$ (resp. $s+P_{(k)}^{X^L}\left(\ell+\frac{1}{2}\right)+i\geq1$) when $t\in\{1,n-1\}$ (resp. $t\in\{2,n\}$)
it follows
\[
L^{X^L}_{s,k,\iota}(Y')=L^{X^L}_{s,k,\iota}(Y)-\beta_{s+P_{(k)}^{X^L}(\ell)+i,t} \quad \text{if }t\in\{1,n-1\},
\]
\[
L^{X^L}_{s,k,\iota}(Y')=L^{X^L}_{s,k,\iota}(Y)-\beta_{s+P_{(k)}^{X^L}\left(\ell+\frac{1}{2}\right)+i,t} \quad \text{if }t\in\{2,n\},
\]
where $P_{(k)}^{X^L}$ is the same notation as in the previous proposition.
\end{prop}

\nd{\it Proof.}

We prove the case the $t$-admissible slot in our claim is in $\ovl{Y}$.
We only consider the first $t$-admissible slot
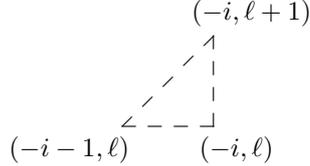
\[
\begin{xy}
(17,15) *{(-i,\ell+1)}="000",
(15,-3) *{(-i,\ell)}="00",
(-7,-3) *{(-i-1,\ell)}="0",
(0,0) *{}="1",
(12,0)*{}="2",
(12,12)*{}="3",
\ar@{--} "1";"2"^{}
\ar@{--} "1";"3"^{}
\ar@{--} "2";"3"^{}
\end{xy}
\]
since one can similarly treat the second slot. We name this slot as $\ovl{A}_1$.
Blocks and slots around $\ovl{A}_1$ are as follows:
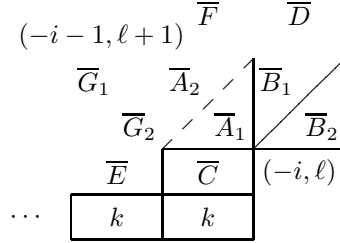
\[
\begin{xy}
(-8,15) *{(-i-1,\ell+1)}="0000",
(18,-3) *{(-i,\ell)}="00",
(9,3) *{\ovl{A}_1}="A",
(3,9) *{\ovl{A}_2}="AA",
(21,3) *{\ovl{B}_2}="B_2",
(15,9) *{\ovl{B}_1}="B_1",
(-3,3) *{\ovl{G}_2}="G_2",
(-9,9) *{\ovl{G}_1}="G_1",
(18,18) *{\ovl{D}}="D",
(6,-3) *{\ovl{C}}="C",
(-6,-3) *{\ovl{E}}="E",
(6,-9) *{k}="CC",
(6,18) *{\ovl{F}}="F",
(-6,-9) *{k}="BB",
(-18,-9) *{\cdots}="BBB",
(0,0) *{}="1",
(12,0)*{}="2",
(12,12)*{}="3",
(0,12)*{}="4",
(12,12)*{}="5",
(0,0)*{}="6",
(24,0)*{}="2-a",
(24,12)*{}="2-aa",
(0,-12) *{}="11",
(12,-12) *{}="22",
(0,-6) *{}="111",
(12,-6) *{}="222",
(-12,-6) *{}="1111",
(-12,-12) *{}="11111",
\ar@{-} "2";"2-a"^{}
\ar@{-} "2";"2-aa"^{}
\ar@{-} "2-aa";"2-a"^{}
\ar@{-} "11";"22"^{}
\ar@{-} "1";"11"^{}
\ar@{-} "22";"2"^{}
\ar@{-} "1";"2"^{}
\ar@{-} "111";"222"^{}
\ar@{-} "1111";"11111"^{}
\ar@{-} "1111";"111"^{}
\ar@{-} "11";"11111"^{}
\ar@{--} "1";"3"^{}
\ar@{-} "2";"3"^{}
\ar@{-} "3";"5"^{}
\ar@{-} "1";"6"^{}
\end{xy}
\]
Here, we named blocks and slots as above.
Let $\ovl{\ovl{C}}$ (resp. $\ovl{\ovl{E}}$) be a block (resp. slot) in $\ovl{\ovl{Y}}$ with half-unit thickness such that the position of the pair $C:=(\ovl{C},\ovl{\ovl{C}})$
(resp. $E:=(\ovl{E},\ovl{\ovl{E}})$) is $(-i,\ovl{T}^{X^L}_k)$ (resp. $(-i-1,\ovl{T}^{X^L}_k)$).
We take $u\in I$ such that $\{t,u\}=\{1,2\}$ or 
$\{t,u\}=\{n-1,n\}$.

First, we consider the case $\ovl{A}_2$ is a $u$-block. Note that $\ovl{B}_1$ and $\ovl{D}$ are blocks. 
The admissibility and removability of $\ovl{B}_1$, $\ovl{B}_2$ $C$ and $E$ in $Y$ are same as in $Y'$.
Those of $\ovl{A}_1$, $\ovl{D}$ and $\ovl{F}$ are same as in Case 1 of the proof of Proposition 5.3.
If $\ovl{E}$ is a slot then the admissibility and removability of $\ovl{A}_2$, $\ovl{G}_1$ and $\ovl{G}_2$ in $Y$ are same as those in $Y'$.
If $\ovl{E}$ is a block then those of $\ovl{A}_2$, $\ovl{G}_1$ and $\ovl{G}_2$ are same as in Case 2 of the proof of Proposition 5.3.

Next, we consider the case $\ovl{A}_2$ is a $u$-slot.
Then $\ovl{E}$ is a slot and $\ovl{G}_1$, $\ovl{G}_2$ and $\ovl{F}$ are normal slots in $Y$ and $Y'$.
As for $\ovl{A}_2$, $\ovl{B}_1$, $\ovl{B}_2$, $\ovl{D}$ and $\ovl{F}$, the admissibility and removability 
are same as in Case 3 - Case 5 in the proof of Proposition 5.3.
Note that the pair
$E$ is $k$-admissible in $Y'$ if and only if $C$ is not a removable $k$-pair in $Y$. 
Thus, we get our claim. \qed

\begin{prop}\label{prop-closedness7}
We assume $\mathfrak{g}$ is of type $X$. Let $k\in I$ in class $2$. 
Let $t\in I$ and $Y=(\overline{Y},\overline{\overline{Y}})\in {\rm YW}_k^{X^L}$. 
We assume that $Y'=(\overline{Y}',\overline{\overline{Y}}')\in {\rm YW}_k^{X^L}$
is obtained from $Y$ by adding the $k$-block
with unit width, half-unit height and unit thickness
to each top of $i+1$-th column of $\overline{Y}'$ and $\overline{Y}'$ from the right ($i\in\mathbb{Z}_{\geq0}$).
Then for $s\in\mathbb{Z}_{\geq0}$
such that $s+i\geq1$,
it follows
\[
L^{X^L}_{s,k,\iota}(Y')=L^{X^L}_{s,k,\iota}(Y)-\beta_{s+i,k}.
\]

\end{prop}

\nd
{\it Proof.}

We name the $k$-admissible slot in $\overline{Y}$ (resp. $\overline{\overline{Y}}$) as $\overline{A}$ (resp. $\overline{\overline{A}}$).
We name blocks and slots around $\overline{A}$ as follows:
\[
\begin{xy}
(6,6) *{\overline{B}}="B",
(6,-3) *{k}="C",
(-6,-3) *{\overline{A}}="AA",
(6,-9) *{k}="CC",
(-6,6) *{\overline{C}}="C",
(-6,-9) *{k}="BB",
(-18,-9) *{\cdots}="BBB",
(0,0) *{}="1",
(12,0)*{}="2",
(12,12)*{}="3",
(0,12)*{}="4",
(12,12)*{}="5",
(0,0)*{}="6",
(0,-12) *{}="11",
(12,-12) *{}="22",
(0,-6) *{}="111",
(-12,0) *{}="111a",
(12,-6) *{}="222",
(-12,-6) *{}="1111",
(-12,-12) *{}="11111",
\ar@{--} "111a";"1111"^{}
\ar@{--} "111a";"1"^{}
\ar@{-} "11";"22"^{}
\ar@{-} "1";"11"^{}
\ar@{-} "22";"2"^{}
\ar@{-} "1";"2"^{}
\ar@{-} "111";"222"^{}
\ar@{-} "1111";"11111"^{}
\ar@{-} "1111";"111"^{}
\ar@{-} "11";"11111"^{}
\ar@{-} "3";"5"^{}
\ar@{-} "1";"6"^{}
\end{xy}
\]
We name blocks and slots around the slot $\overline{\overline{A}}$ as follows:
\[
\begin{xy}
(6,6) *{\overline{\overline{B}}}="B",
(6,-3) *{k}="C",
(-6,-3) *{\overline{\overline{A}}}="AA",
(6,-9) *{k}="CC",
(-6,6) *{\overline{\overline{C}}}="C",
(-6,-9) *{k}="BB",
(-18,-9) *{\cdots}="BBB",
(0,0) *{}="1",
(12,0)*{}="2",
(12,12)*{}="3",
(0,12)*{}="4",
(12,12)*{}="5",
(0,0)*{}="6",
(0,-12) *{}="11",
(12,-12) *{}="22",
(0,-6) *{}="111",
(-12,0) *{}="111a",
(12,-6) *{}="222",
(-12,-6) *{}="1111",
(-12,-12) *{}="11111",
\ar@{--} "111a";"1111"^{}
\ar@{--} "111a";"1"^{}
\ar@{-} "11";"22"^{}
\ar@{-} "1";"11"^{}
\ar@{-} "22";"2"^{}
\ar@{-} "1";"2"^{}
\ar@{-} "111";"222"^{}
\ar@{-} "1111";"11111"^{}
\ar@{-} "1111";"111"^{}
\ar@{-} "11";"11111"^{}
\ar@{-} "3";"5"^{}
\ar@{-} "1";"6"^{}
\end{xy}
\]
\nd
The position of $(\ovl{A},\ovl{\ovl{A}})$ is $(-i,\ovl{T}^{X^L}_k)$.
Let $(T,\tilde{Y},B,C)$ be $(\ovl{T}^{X^L}_k,\ovl{Y},\ovl{B},\ovl{C})$ or $(\ovl{\ovl{T}}^{X^L}_k,\ovl{\ovl{Y}},\ovl{\ovl{B}},\ovl{\ovl{C}})$ and
$B'$, $C'$ be the corresponding blocks or slots to $B$, $C$ in $Y'$.
Let us
compute $L^{X^L}_{s,k}(B')+L^{X^L}_{s,k}(C')-(L^{X^L}_{s,k}(B)+L^{X^L}_{s,k}(C))$, where we use the same notation as in (\ref{rule1}), (\ref{rule2}).
There are four cases:

Case 1 : $B$ is a unit cube, 

Case 2 : $B$ consists of two blocks with half-unit thickness,

Case 3 : $B$ consists of a block with half-unit thickness,

Case 4 : $B$ consists of one or two blocks with half-unit height.

\vspace{2mm}

\underline{Case 1 : $B$ is a unit cube}

\vspace{2mm}

Note that $B$ is a removable $\pi'_{X^L}(T+1)$-block in $\tilde{Y}$
if and only if
$C$
 is not a $\pi'_{X^L}(T+1)$-admissible slot in $\tilde{Y}'$.
Thus, we obtain
 $L^{X^L}_{s,k}(B')+L^{X^L}_{s,k}(C')-(L^{X^L}_{s,k}(B)+L^{X^L}_{s,k}(C))=x_{s+P^{X^L}_{T}(T+1)+i,\pi'_{X^L}(T+1)}=x_{s+p_{\pi'_{X^L}(T+1),k}+i,\pi'_{X^L}(T+1)}$.

\vspace{2mm}

\underline{Case 2 : $B$ consists of two blocks with half-unit thickness}

\vspace{2mm}

The removability of B in $\tilde{Y}$ are same as those in 
$\tilde{Y}'$. Writing $C$ as
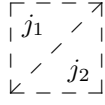
\[
\begin{xy}
(3,9) *{j_1}="00",
(9,3) *{j_2}="0",
(0,0) *{}="1",
(12,0)*{}="2",
(12,12)*{}="3",
(0,12)*{}="4",
(0,0)*{}="5",
(12,12)*{}="6",
\ar@{--} "1";"2"^{}
\ar@{--} "1";"4"^{}
\ar@{--} "2";"3"^{}
\ar@{--} "3";"4"^{}
\ar@{--} "5";"6"^{}
\end{xy}
\]
we see that $C$ is $j_1$-admissible and 
$j_2$-admissible in 
$\tilde{Y}'$.
Thus, we see that
$L^{X^L}_{s,k}(B')+L^{X^L}_{s,k}(C')-(L^{X^L}_{s,k}(B)+L^{X^L}_{s,k}(C))=x_{s+p_{j_1,k}+i,j_1}+x_{s+p_{j_2,k}+i,j_2}$.

\vspace{2mm}

\underline{Case 3 : $B$ consists of a block with half-unit thickness}

\vspace{2mm}

$B$ is in the form
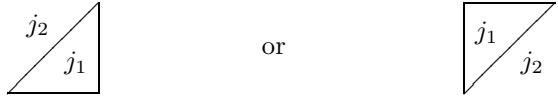
\[
\begin{xy}
(35,6) *{\text{or}}="0",
(9,4) *{j_1}="00",
(4,9) *{j_2}="0000",
(63,8) *{j_1}="000",
(69,4) *{j_2}="00000",
(0,0) *{}="1",
(12,0)*{}="2",
(12,12)*{}="3",
(60,0) *{}="1a",
(72,0)*{}="2a",
(72,12)*{}="3a",
(60,12)*{}="4a",
\ar@{-} "1";"2"^{}
\ar@{-} "1";"3"^{}
\ar@{-} "2";"3"^{}
\ar@{-} "1a";"4a"^{}
\ar@{-} "1a";"3a"^{}
\ar@{-} "3a";"4a"^{}
\end{xy}
\]
Let $B_1$ be the $j_1$-block and $B_2$ be the $j_2$-slot. Note that $C$ consists of $j_1$-slot $C_1$ and $j_2$-slot $C_2$.
We see that $B_1$ is a removable $j_1$-block (resp. normal)
in $\tilde{Y}$ (resp. $\tilde{Y}'$). The admissibilities of
$B_2$ and $C_1$ in $\tilde{Y}$ are same as those in $\tilde{Y}'$.
The slot $C_2$ is normal (resp. $j_2$-admissible) in 
$\tilde{Y}$ (resp. $\tilde{Y}'$).
Hence, we obtain $L^{X^L}_{s,k}(B')+L^{X^L}_{s,k}(C')-(L^{X^L}_{s,k}(B)+L^{X^L}_{s,k}(C))=x_{s+p_{j_1,k}+i,j_1}+x_{s+p_{j_2,k}+i,j_2}$.

\vspace{2mm}

\underline{Case 4 : $B$ consists of one or two blocks with half-unit height}

\vspace{2mm}

There are three cases:

(1) : $B$ consists of 
one block $B_1$ with half-unit height and one slot with half-unit height $B_2$,

(2) : $B$ consists of 
two blocks $B_1$, $B_2$ with half-unit height and the above of
$B_2$ is a slot,

(3) : $B$ consists of 
two blocks $B_1$, $B_2$ with half-unit height and the above of
$B_2$ is a block.

\begin{table}[H]
  \begin{tabular}{|c|l|l|} \hline
    slot or block & in $Y$ & in $Y'$ \\ \hline
    $B$ & (1) $B_1$ : single $\pi'_{X^L}(T+1)$-removable & (1) $B_1$ : normal \\
     & (1) $B_2$ : same as in $Y'$ & (1) $B_2$ : same as in $Y$ \\
     & (2) $B_1$ : normal & (2) $B_1$ : normal \\
     & (2) $B_2$ : double removable $\pi'_{X^L}(T+1)$ & (2) $B_2$ : single removable $\pi'_{X^L}(T+1)$  \\
     & (3) $B_1$, $B_2$ : normal & (3) $B_1$, $B_2$ : normal \\
    $C$ & normal & (1) single $\pi'_{X^L}(T+1)$-admissible \\
    & & (2) single $\pi'_{X^L}(T+1)$-admissible \\
     & & (3) double $\pi'_{X^L}(T+1)$-admissible \\
    \hline
  \end{tabular}
\end{table}
Hence,
$L^{X^L}_{s,k}(B')+L^{X^L}_{s,k}(C')-(L^{X^L}_{s,k}(B)+L^{X^L}_{s,k}(C))=2x_{s+p_{\pi'_{X^L}(T+1),k}+i,\pi'_{X^L}(T+1)}$.
In all cases, the pair $(\overline{A},\overline{\overline{A}})$ (resp. $(\overline{A}',\overline{\overline{A}}')$) is a $k$-admissible pair (resp. removable $k$-pair)
in $Y$ (resp. $Y'$) so that
\[
-L^{X^L}_{s,k,{\rm re}}(\overline{A}',\overline{\overline{A}}')-L^{X^L}_{s,k,{\rm ad}}(\overline{A},\overline{\overline{A}})=-x_{s+i+1,k}-x_{s+i,k},
\]
where $(\overline{A}',\overline{\overline{A}}')$ is the pair in $Y'$ at the same position as $(\overline{A},\overline{\overline{A}})$.
In conjunction with above calculations, one obtains
$L^{X^L}_{s,k,\iota}(Y')=L^{X^L}_{s,k,\iota}(Y)-\beta_{s+i,k}$.\qed

\section{Proof}

We assume that $\mathfrak{g}$ is of type $X=B^{(1)}_{n-1}$, $C^{(1)}_{n-1}$, $D^{(1)}_{n-1}$, $A^{(2)}_{2n-2}$, $A^{(2)}_{2n-3}$ or
$D^{(2)}_{n}$.

\subsection{Proof of Theorem \ref{mainthm1}}

For $s\in\mathbb{Z}_{\geq1}$ and $k\in I$, we set
\[
\Xi_{s,k,\io}' :=  \{S_{j_{\ell}}'\cd S_{j_2}'S_{j_1}'x_{s,k}\,|\,
\ell\in \mathbb{Z}_{\geq0},j_1,\cd,j_{\ell}\in\mathbb{Z}_{\geq1}\},
\]
and simply write
\begin{equation}\label{omit}
L_{s,k}:=L^{X^L}_{s,k,\iota},\ L_{s,k,{\rm ad}}:=L^{X^L}_{s,k,{\rm ad}},\ L_{s,k,{\rm re}}:=L^{X^L}_{s,k,{\rm re}}
\end{equation}
Let us prove $\Xi'_{s,k,\io}=L_{s,k}({\rm YW}^{X^L}_k)$. First, we show $\Xi'_{s,k,\io}\subset L_{s,k}({\rm YW}^{X^L}_k)$.
We take $Y_{\Lambda_k}\in {\rm YW}^{X^L}_k$
defined in the subsection \ref{pyws} and Definition \ref{YWk-def}.
One can verify that $L_{s,k}(Y_{\Lambda_k})=x_{s,k}$. Thus, we have $x_{s,k}\in L_{s,k}({\rm YW}^{X^L}_k)$
so that need to prove $L_{s,k}({\rm YW}^{X^L}_k)$ is closed under the action of 
$S'_{m,t}$ for any $(m,t)\in\mathbb{Z}_{\geq1}\times I$.
For any $Y\in {\rm YW}^{X^L}_k$ and
$(m,t)\in\mathbb{Z}_{\geq1}\times I$, if $x_{m,t}$ has a positive coefficient in $L_{s,k}(Y)$
then there is a $t$-admissible slot or $t$-admissible pair $S$
such as in Proposition \ref{prop-closedness1}-\ref{prop-closedness7}
in $Y$ and it holds $L_{s,k,{\rm ad}}(S)=x_{m,t}$. Here, writing the position of $S$ as $(-i,\ell)$, we have $(m,t)=(s+P_{(k)}^{X^L}(\ell)+i,\pi'_{X^L}(\ell))$
or $(m,t)=(s+P_{(k)}^{X^L}\left(\ell+\frac{1}{2}\right)+i,\pi'_{X^L}\left(\ell+\frac{1}{2}\right))$ and
$P_{(k)}^{X^L}(\ell)=P_{k,S}^{X^L}(\ell)$ is defined in Proposition \ref{prop-closedness1} (we understand $S$ belongs to a truncated wall of supporting ground state when $S$ is a pair
).
Let $Y'\in {\rm YW}^{X^L}_k$ be the element obtained from $Y$ by putting a $t$-block or pair of $t$-blocks to $S$.
By the propositions, we see that
\[
L_{s,k}(Y')=L_{s,k}(Y)-\beta_{m,t}=S'_{m,t}L_{s,k}(Y)
\]
and $S'_{m,t}L_{s,k}(Y)\in L_{s,k}({\rm YW}^{X^L}_k)$. Similarly, in case that $x_{m,t}$ has a negative coefficient in $L_{s,k}(Y)$,
there is a removable $t$-block or removable $t$-pair $B$ with position $(-i,\ell)$
such as in Proposition \ref{prop-closedness1}-\ref{prop-closedness7}
in $Y$ and it holds $L_{s,k,{\rm re}}(B)=x_{m,t}$. We have $(m,t)=(s+P_{(k)}^{X^L}(\ell)+i+1,\pi'_{X^L}(\ell))$
or $(m,t)=(s+P_{(k)}^{X^L}\left(\ell+\frac{1}{2}\right)+i+1,\pi'_{X^L}\left(\ell+\frac{1}{2}\right))$.
It is clear $m\geq2$ from (\ref{non-neg1}), (\ref{non-neg2}) and (\ref{non-neg3}).
Let $Y''\in {\rm YW}^{X^L}_k$ be the element obtained from $Y$ by removing $B$.
By the propositions, we see that
\[
L_{s,k}(Y)=L_{s,k}(Y'')-\beta_{m-1,t}
\]
so that $L_{s,k}({\rm YW}^{X^L}_k)\ni L_{s,k}(Y'')=L_{s,k}(Y)+\beta_{m-1,t}=S'_{m,t}L_{s,k}(Y)$.
In this way, we can show $\Xi'_{s,k,\io}\subset L_{s,k}({\rm YW}^{X^L}_k)$.

Next, we prove $\Xi'_{s,k,\io}\supset L_{s,k}({\rm YW}^{X^L}_k)$.
When we get a wall $Y$ by adding $m$ blocks to $Y_{\Lambda_k}$,
we say the number of blocks in $Y$ is $m$.
For any $Y\in {\rm YW}^{X^L}_k$, we prove $L_{s,k}(Y)\in \Xi'_{s,k,\iota}$
using the induction on the number of blocks in $Y$.
If the number of blocks is $0$ then it holds $Y=Y_{\Lambda_k}$ and
$L_{s,k}(Y)=x_{s,k}\in \Xi'_{s,k,\iota}$. Thus, we may assume $Y\neq Y_{\Lambda_k}$.
Using (\ref{non-neg1}), (\ref{non-neg2}), (\ref{non-neg3}) and
Proposition \ref{prop-closedness1}-\ref{prop-closedness7},
we see that $L_{s,k}(Y)$ can be expressed as
\[
L_{s,k}(Y)=x_{s,k}-\sum_{(r,j)\in\mathbb{Z}_{\geq s}\times I} c_{r,j}\beta_{r,j}
\]
with non-negative integers $\{c_{r,j}\}$. Except for finitely many $(r,j)$, it follows $c_{r,j}=0$.
Since we assumed
$Y\neq Y_{\Lambda_k}$, one can take $(m',t')\in\mathbb{Z}_{\geq s}\times I$ as
$(m',t')
={\rm max}\{(r,j)\in\mathbb{Z}_{\geq s}\times I |c_{r,j}>0 \}$. Here, we consider the order on $\mathbb{Z}_{\geq 1}\times I$
defined in the subsection \ref{s-d}.
Considering the definition of $\beta_{m',t'}$ and the adaptedness of $\iota$, it is easy to see 
the coefficient of $x_{m'+1,t'}$ is negative in $L_{s,k}(Y)$.
Hence,
we can also take $(m'',t'')$ as
\[
(m'',t'')={\rm min}\{(r,j)\in\mathbb{Z}_{\geq 1}\times I | \text{the coefficient of }x_{r,j}\text{ in }L_{s,k}(Y)\text{ is negative}\}.
\]
There exists a removable $t''$-block or removable $t''$-pair $B$ in $Y$ such that $L_{s,k,{\rm re}}(B)=x_{m'',t''}$. 
Let $Y''\in{\rm YW}^{X^L}_k$ be the element obtained from $Y$ by removing $B$. 
Then we can prove 
\begin{equation}\label{A2YW-pr2}
L_{s,k}(Y)=L_{s,k}(Y'')-\beta_{m''-1,t''}=S'_{m''-1,t''}L_{s,k}(Y'')
\end{equation}
 by (\ref{non-neg1}), (\ref{non-neg2}), (\ref{non-neg3}) and
Proposition \ref{prop-closedness1}-\ref{prop-closedness7}.
Since the number of boxes in $Y''$ is smaller than $Y$, by induction assumption,
we obtain $L_{s,k}(Y'')\in \Xi'_{s,k,\iota}$, which yields $L_{s,k}(Y)\in \Xi'_{s,k,\iota}$ by (\ref{A2YW-pr2}).
Therefore, we get the inclusion $\Xi'_{s,k,\io}\supset L_{s,k}({\rm YW}^{X^L}_k)$. Therefore, we obtain
$\Xi'_{\iota}={\rm COMB}^{X^L}_{\iota}[\infty]$. It follows from (\ref{non-neg1}), (\ref{non-neg2}), (\ref{non-neg3}) and (\ref{l1def5})
that $\iota$ satisfies the $\Xi'$-positivity condition. Our claim is an easy consequence
of Theorem \ref{inf-thm}.
\qed

\subsection{Proof of Theorem \ref{mainthm2}}

For simplicity, we assume $X\neq D^{(1)}_4$. The proof for $X=D^{(1)}_4$
is essentially same as other cases.
Let us prove $\Xi'_{k,\iota}[\lambda]:=\{\what{S}'_{j_{\ell}}\cdots \what{S}'_{j_1} \lambda^{(k)} | \ell\in\mathbb{Z}_{\geq0},\ j_1,\cdots,j_{\ell}\in \mathbb{Z}_{\geq1} \}={\rm COMB}_{k,\iota}^{X^L}[\lambda]\cup\{0\}$.
We consider only $k\in I$ such that either

(1) for any $j\in I$ such that $a_{k,j}<0$, it holds $(1,k)>(1,j)$, or

(2) there exists distinct $j_1,j_2,j_3\in I$ such that $j_1$ and $j_2$ are in class $1$,
$a_{k,j_1}=a_{j_1,k}=a_{k,j_2}=a_{j_2,k}=-1$ and $a_{k,j_3}<0$ and
$(1,k)>(1,j_1)$, $(1,k)>(1,j_3)$ and $(1,k)<(1,j_2)$,

\hspace{-7mm}
since other cases are trivial via direct calculations or Appendix A. 

First, we consider the case (1).
Putting 
$Z:=\{L_{0,k}(Y)+\langle h_k,\lambda \rangle | Y\in{\rm YW}^{X^L}_{k}\setminus\{Y_{\Lambda_k}\}\}\cup\{0\}$,
we show $\{\what{S}'_{j_{\ell}}\cdots \what{S}'_{j_1} \lambda^{(k)} | \ell\in\mathbb{Z}_{\geq0},\ j_1,\cdots,j_{\ell}\in \mathbb{Z}_{\geq1} \}\subset Z$.
Here, we use the same notation as in (\ref{omit}).
Let $Y_k^1\in {\rm YW}^{X^L}_k$ be the element obtained from $Y_{\Lambda_k}$ by adding $k$-block to the unique $k$-admissible slot (or $k$-admissible
pair). One can confirm $L_{0,k}(Y_k^1)+\lan h_k,\lambda \ran=\lambda^{(k)}$ so that $\lambda^{(k)}\in Z$.
Let us prove $Z$ is closed under the action of $\what{S}'_{m,t}$ for all $(m,t)\in\mathbb{Z}_{\geq1}\times I$. We take $Y\in{\rm YW}^{X^L}_{k}\setminus\{Y_{\Lambda_k}\}$.
If the coefficient of $x_{m,t}$ in $L_{0,k}(Y)$ is positive then there is a $t$-admissible slot or $t$-admissible pair $S$ in $Y$ such that
$L_{0,k,{\rm ad}}(S)=x_{m,t}$. Let $(-i,\ell)$ be its position. We have $(m,t)=(P_{(k)}^{X^L}(\ell)+i,\pi'_{X^L}(\ell))$ or $(m,t)=(P_{(k)}^{X^L}\left(\ell+\frac{1}{2}\right)+i,\pi'_{X^L}\left(\ell+\frac{1}{2}\right))$. 
Let $Y'\in {\rm YW}^{X^L}_k$ be the element obtained from $Y$ by putting a $t$-block to the slot or pair.
By Proposition \ref{prop-closedness1}-\ref{prop-closedness7}, we see that
\[
L_{0,k}(Y')=L_{0,k}(Y)-\beta_{m,t}=\what{S}'_{m,t}L_{0,k}(Y)
\]
and $\what{S}'_{m,t}L_{0,k}(Y)\in L_{0,k}({\rm YW}^{X^L}_k\setminus\{Y_{\Lambda_k}\})$.
If the coefficient of $x_{m,t}$ in $L_{0,k}(Y)$ is negative then there is a removable $t$-block or removable $t$-pair $S$ in $Y$
such that $L_{0,k,{\rm re}}(S)=x_{m,t}$. Let $(-i,\ell)$ be its position. 
Putting
\[
T:=
\begin{cases}
\ovl{T}_k & \text{if }S\text{ is in a proper Young wall},\\
\ovl{T}_k & \text{if }S\text{ is a pair or is in a truncated wall of supporting ground state},\\
\ovl{\ovl{T}}_k & \text{if }S\text{ is in a truncated wall of covering ground state},
\end{cases}
\]
we have $(m,t)=(P^{X^L}_{T}(\ell)+i+1,\pi'_{X^L}(\ell))$ or
$(m,t)=(P^{X^L}_{T}\left(\ell+\frac{1}{2}\right)+i+1,\pi'_{X^L}\left(\ell+\frac{1}{2}\right))$.
If $Y=Y_k^1$ then 
$S$ is the removable $k$-pair and $(-i,\ell)=(0,\overline{T}_k)$ so that $(m,t)=(1,k)$ and
$\what{S}'_{1,k}(L_{0,k}(Y_k^1)+\lan h_k,\lambda \ran)=\what{S}'(\lambda^{(k)})=0$. If $Y\neq Y_k^1$ then it holds
$i>0$ or 
$\ell\geq T+1$, which implies 
\begin{equation}\label{D-pr-7}
P^{X^L}_{T}(\ell)+i+1\geq2\quad (\text{or }P^{X^L}_{T}\left(\ell+\frac{1}{2}\right)+i+1\geq2).
\end{equation}
Let $Y''\in {\rm YW}^{X^L}_k\setminus\{Y_{\Lambda_k}\}$ be the element obtained from $Y$
by removing $S$. 
Using Proposition \ref{prop-closedness1}-\ref{prop-closedness7}, we get
\[
L_{0,k}(Y'')=L_{0,k}(Y)+\beta_{m-1,t}=\what{S}'_{m,t}L_{0,k}(Y)
\]
and $\what{S}'_{m,t}L_{0,k}(Y)=L_{0,k}(Y'')\in L_{0,k}({\rm YW}^{X^L}_k\setminus\{Y_{\Lambda_k}\})$. Therefore, the set $Z$ is closed under the action of $\what{S}'_{m,t}$ so that
the inclusion $\{\what{S}'_{j_{\ell}}\cdots \what{S}'_{j_1} \lambda^{(k)} | \ell\in\mathbb{Z}_{\geq0},\ j_1,\cdots,j_{\ell}\in \mathbb{Z}_{\geq1} \}\subset Z$ follows.

Next, to prove another inclusion, let us take any $Y\in {\rm YW}^{X^L}_{k}\setminus\{Y_{\Lambda_k}\}$. We show $L_{0,k}(Y)+\lan h_k,\lambda\ran\in
\{\what{S}'_{j_{\ell}}\cdots \what{S}'_{j_1} \lambda^{(k)} | \ell\in\mathbb{Z}_{\geq0},\ j_1,\cdots,j_{\ell}\in \mathbb{Z}_{\geq1} \}$ by induction on the number of boxes in $Y$.
The minimum value of number of boxes is $2$, which occurs only if $Y=Y_k^1$. By
$L_{0,k}(Y_k^1)+\lan h_k,\lambda \ran=\lambda^{(k)}\in \{\what{S}'_{j_{\ell}}\cdots \what{S}'_{j_1} \lambda^{(k)} | \ell\in\mathbb{Z}_{\geq0},\ j_1,\cdots,j_{\ell}\in \mathbb{Z}_{\geq1} \}$,
we may assume $Y\neq Y_k^1$.
By
Proposition \ref{prop-closedness1}-\ref{prop-closedness7},
$L_{0,k}(Y)$ can be expressed as
\[
L_{0,k}(Y)=(\lambda^{(k)}-\lan h_k,\lambda \ran)-\sum_{(r,j)\in\mathbb{Z}_{\geq 1}\times I} c_{r,j}\beta_{r,j}
\]
with non-negative integers $\{c_{r,j}\}$. Except for finitely many $(r,j)$, it follows $c_{r,j}=0$. 
It follows from the assumption $Y\neq Y_k^1$ that one can take
$(m',t')\in\mathbb{Z}_{\geq 1}\times I$ as
$(m',t')
={\rm max}\{(r,j)\in\mathbb{Z}_{\geq 1}\times I |c_{r,j}>0 \}$.
The coefficient of $x_{m'+1,t'}$ is negative in $L_{0,k}(Y)$ by the definition of $\beta_{m',t'}$. We set
\[
(m'',t''):={\rm min}\{(r,j)\in\mathbb{Z}_{\geq 1}\times I | \text{the coefficient of }x_{r,j}\text{ in }L_{0,k}(Y)\text{ is negative}\}
\]
and take a removable block $B$ such that $L_{0,k,{\rm re}}(B)=x_{m'',t''}$. By (\ref{D-pr-7}), it holds $m''\geq2$.
Let $Y''\in{\rm YW}^{X^L}_k\setminus\{Y_{\Lambda_k}\}$ be the element obtained from $Y$ by removing $B$. 
Then we get 
\begin{equation}\label{D-pr-8}
L_{0,k}(Y)=L_{0,k}(Y'')-\beta_{m''-1,t''}=\what{S}'_{m''-1,t''}L_{0,k}(Y'')
\end{equation}
by 
Proposition \ref{prop-closedness1}-\ref{prop-closedness7}.
Since the number of boxes in $Y''$ is smaller than $Y$, by induction assumption,
we obtain $L_{0,k}(Y'')+\lan h_k,\lambda \ran$ belongs
to
$\{\what{S}'_{j_{\ell}}\cdots \what{S}'_{j_1} \lambda^{(k)} | \ell\in\mathbb{Z}_{\geq0},\ j_1,\cdots,j_{\ell}\in \mathbb{Z}_{\geq1} \}$
and so is $L_{0,k}(Y)+\lan h_k,\lambda \ran$ by (\ref{D-pr-8}).
Therefore, we get the inclusion $\{\what{S}'_{j_{\ell}}\cdots \what{S}'_{j_1} \lambda^{(k)} | \ell\in\mathbb{Z}_{\geq0},\ j_1,\cdots,j_{\ell}\in \mathbb{Z}_{\geq1} \}\supset Z$.

Next, we turn to the case (2). 
One obtains $\lambda^{(k)}=-x_{1,k}+x_{1,j_1}+cx_{1,j_3}+\lan h_k,\lambda \ran$, where
$c=2$ when $X^L=B^{(1)}_{n-1}$ and $j_3=n$ and $c=1$ otherwise. Let $Y_{k}^2\in {\rm YW}_{j_2}^{X^L}$ be the wall obtained
from $Y_{\Lambda_k}$ by adding a $j_2$-block and $k$-block, that is,
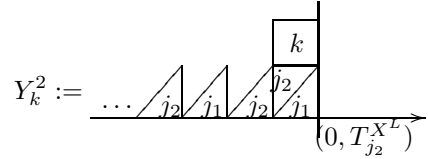
\[
Y_{k}^2:=
\begin{xy}
(-14,-2) *{\ j_2}="000-1",
(-8.5,-2) *{\ j_1}="00-1",
(-2.5,-2) *{\ j_2}="0-1",
(2.5,-2) *{\ \ j_1}="0-11",
(0,1.5) *{\ \ j_2}="0-111",
(2,6.5) *{\ \ k}="0-111-a",
(-20.5,-2) *{\dots}="00000",
(12,-6) *{(0,T^{X^L}_{j_2})}="origin",
(0,0) *{}="1",
(0,-3.5) *{}="1-u",
(6,0)*{}="2",
(6,-3.5)*{}="2-u",
(6,-6)*{}="2-uu",
(6,-7.5)*{}="-1-u",
(6,12)*{}="3",
(-6,3.5)*{}="5",
(-6,0)*{}="6",
(-6,-3.5)*{}="6-u",
(-12,3.5)*{}="7",
(-12,-3.5)*{}="7-u",
(-12,0)*{}="8",
(-12,-3.6)*{}="8-u",
(-18,3.5)*{}="9",
(-18,0)*{}="10-a",
(-18,0)*{}="10",
(-18,-3.5)*{}="10-u",
(-24,-3.5)*{}="11-u",
(20,-3.5)*{}="0-u",
(0,3.5) *{}="1-a",
(6,3.5) *{}="2-a",
(0,9.5) *{}="1-aa",
(6,9.5) *{}="2-aa",
\ar@{-} "1-aa";"1-a"^{}
\ar@{-} "2-aa";"1-aa"^{}
\ar@{-} "2-a";"1-a"^{}
\ar@{-} "6-u";"1-a"^{}
\ar@{-} "1-u";"2-a"^{}
\ar@{-} "1-u";"1-a"^{}
\ar@{-} "6-u";"5"^{}
\ar@{-} "7-u";"5"^{}
\ar@{-} "7-u";"7"^{}
\ar@{-} "10-u";"7"^{}
\ar@{->} "2-u";"0-u"^{}
\ar@{-} "10-u";"11-u"^{}
\ar@{-} "10-u";"2-u"^{}
\ar@{-} "2";"2-u"^{}
\ar@{-} "2-uu";"2-u"^{}
\ar@{-} "2";"3"^{}
\end{xy}
\]
Then one can verify
$L_{-1,j_2}(Y_{k}^2)+\lan h_k,\lambda \ran=\lambda^{(k)}$. Note that any $Y\in {\rm YW}_{j_2}^{X^L}\setminus\{Y_{\Lambda_{j_2}},Y_{j_2}^1\}$
is obtained from $Y_{k}^2$ by adding several boxes to admissible slots. 
When $Y\neq Y_{k}^2$,
for an admissible slot or removable block with position $(-i,\ell)$ in $Y$,
we have either $i>0$ or $\ell\geq T_{j_2}+2$. From this, one can deduce
\begin{equation}\label{non-neg4}
-1+P^{X^L}_{\ovl{T}_{j_2}}(\ell)+i \geq1,\quad -1+P^{X^L}_{\ovl{T}_{j_2}}(\ell)+i+1 \geq2.
\end{equation}
In the same manner as in (1),
we can show
$\{\what{S}'_{j_{\ell}}\cdots \what{S}'_{j_1} \lambda^{(k)} | \ell\in\mathbb{Z}_{\geq0},\ j_1,\cdots,j_{\ell}\in \mathbb{Z}_{\geq1} \}
=\{
L^{X^L}_{-1,j_2,\iota}(Y)+\langle h_{k},\lambda \rangle | Y\in{\rm YW}^{X^L}_{j_2}\setminus\{Y_{\Lambda_{j_2}},Y_{j_2}^1\}
\}\cup\{0\}$. 

By (\ref{twoS}) and the proof in the previous subsection,
we see that $\{\what{S}'_{j_{\ell}}\cdots \what{S}'_{j_1} x_{j_0} | \ell\in\mathbb{Z}_{\geq0},\ j_0,j_1,\cdots,j_{\ell}\in \mathbb{Z}_{\geq1} \}
=\{S'_{j_{\ell}}\cdots S'_{j_1} x_{j_0} | \ell\in\mathbb{Z}_{\geq0},\ j_0,j_1,\cdots,j_{\ell}\in \mathbb{Z}_{\geq1} \}={\rm COMB}_{\iota}^{X^L}[\infty]$.
Combining with $\{\what{S}'_{j_{\ell}}\cdots \what{S}'_{j_1} \lambda^{(k)} | \ell\in\mathbb{Z}_{\geq0},\ j_1,\cdots,j_{\ell}\in \mathbb{Z}_{\geq1} \}
={\rm COMB}_{k,\iota}^{X^L}[\lambda]\cup\{0\}$, it is easy to see $(\iota,\lambda)$ satisfies $\Xi'$-ample condition by the definition
of ${\rm COMB}_{k,\iota}^{X^L}[\lambda]$. Therefore, our claim holds from Theorem \ref{Nthm1}.
\qed

\subsection{Proof of Theorem \ref{mainthm3}}

By Theorem \ref{mainthm1}, the sequence $\iota$ satisfies the $\Xi'$-positivity condition.
Therefore, it holds
\[
\text{if }\ell^{(-)}=0\text{ then }\varphi_{\ell}\geq0\text{ for any }
\vp=\sum_{\ell\geq1} \varphi_{\ell}x_{\ell}\in \Xi'_{\iota}.
\]
Let us show
for any $k\in I$ and $\vp=\sum_{\ell\geq1} \varphi_{\ell}x_{\ell}\in \Xi'^{(k)}_{\iota}\setminus\{\xi^{(k)}\}$,
if $\ell^{(-)}=0$ then $\varphi_{\ell}\geq0$.
In the previous subsection, we proved $\{\what{S}'_{j_{\ell}}\cdots \what{S}'_{j_1} \lambda^{(k)} | \ell\in\mathbb{Z}_{\geq0},\ j_1,\cdots,j_{\ell}\in \mathbb{Z}_{\geq1} \}={\rm COMB}_{k,\iota}^{X^L}[\lambda]\cup\{0\}$.
Thus, taking $\lambda$ as $0$, we obtain
$\Xi'^{(k)}_{\iota}={\rm COMB}_{k,\iota}^{X^L}[0]\cup\{0\}$ by (2.24) of \cite{Ka2}.
Using (\ref{D-pr-7}), (\ref{non-neg4}) and Appendix A, 
we see that $\ell^{(-)}=0$ then $\varphi_{\ell}\geq0$. Hence,
$\iota$ satisfies the $\Xi'$-strict positivity condition. 
As a consequence of Theorem \ref{thm2} and (2.24) of \cite{Ka2}, we obtain
\begin{eqnarray*}
\varepsilon_k^*(x)
&=&{\rm max}\{-\varphi(x)|\varphi\in \Xi'^{(k)}_{\iota}\}\\
&=&
{\rm max}\{-\varphi(x)|\varphi\in 
\{S'_{j_{\ell}}\cdots S'_{j_1} 0^{(k)} | \ell\in\mathbb{Z}_{\geq0},\ j_1,\cdots,j_{\ell}\in \mathbb{Z}_{\geq1} \}\}\\
&=&
{\rm max}\{-\varphi(x)|\varphi\in 
\{\what{S}'_{j_{\ell}}\cdots \what{S}'_{j_1} 0^{(k)} | \ell\in\mathbb{Z}_{\geq0},\ j_1,\cdots,j_{\ell}\in \mathbb{Z}_{\geq1} \}\setminus\{0\}\}\\
&=&
{\rm max}\{-\varphi(x)|\varphi\in {\rm Comb}_{k,\iota}^{X^L}[0]\}.
\end{eqnarray*}
\qed

\appendix

\section{Appendix A}


We explicitly give the set
\[
\{\what{S}'_{j_{\ell}}\cdots \what{S}'_{j_1} \lambda^{(k)} | \ell\in\mathbb{Z}_{\geq0},\ j_1,\cdots,j_{\ell}\in \mathbb{Z}_{\geq1} \}
\]
and action of $\what{S}'_{j}$ on elements when the set is expressed by boxes. We assume that
$1<k<n$ and $k\in I$ is in class $2$.

\vspace{2mm}

\hspace{-7mm}
\underline{When $\mathfrak{g}$ is of type $A^{(1)}_{n-1}$}

If $(1,k)>(1,k-1)$ and $(1,k)<(1,k+1)$ then
we get
$\lambda^{(k)}=-x_{1,k}+x_{1,k-1}+\langle h_k,\lambda \rangle=
\fbox{$\tilde{k}$}^{A^{(1)}}_k
+\langle h_k,\lambda \rangle$.
Following the definition,
we get the following diagram of actions of $\what{S}'$ on boxes:
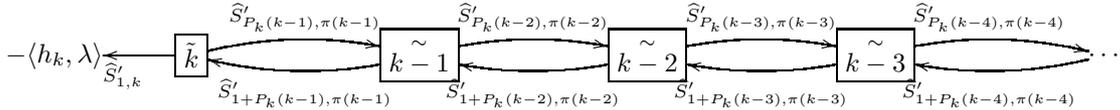
\[
\begin{xy}
(-8,0) *{-\langle h_k,\lambda \rangle}="00",
(0,0) *{\ }="0",
(10,0) *{\fbox{$\tilde{k}$}}="r-1r-1",
(40,0)*{\fbox{$\overset{\sim}{k-1}$}}="r-2r-2",
(70,0)*{\fbox{$\overset{\sim}{k-2}$}}="r-3r-3",
(100,0)*{\fbox{$\overset{\sim}{k-3}$}}="r-4r-4",
(130,0)*{\cdots}="dot",
\ar@{->} "r-1r-1";"00"^{\what{S}'_{1,k}}
\ar@/^/ @{->} "r-1r-1";"r-2r-2"^{\what{S}'_{P_k(k-1),\pi(k-1)}}
\ar@/^/ @{->} "r-2r-2";"r-3r-3"^{\what{S}'_{P_k(k-2),\pi(k-2)}}
\ar@/^/ @{->} "r-3r-3";"r-4r-4"^{\what{S}'_{P_k(k-3),\pi(k-3)}}
\ar@/^/ @{->} "r-4r-4";"dot"^{\what{S}'_{P_k(k-4),\pi(k-4)}}
\ar@/_/ @{<-} "r-1r-1";"r-2r-2"_{\what{S}'_{1+P_k(k-1),\pi(k-1)}}
\ar@/_/ @{<-} "r-2r-2";"r-3r-3"_{\what{S}'_{1+P_k(k-2),\pi(k-2)}}
\ar@/_/ @{<-} "r-3r-3";"r-4r-4"_{\what{S}'_{1+P_k(k-3),\pi(k-3)}}
\ar@/_/ @{<-} "r-4r-4";"dot"_{\what{S}'_{1+P_k(k-4),\pi(k-4)}}
\end{xy}
\]
Here, we simply write $\fbox{$\tilde{a}$}^{A^{(1)}}_k$, $P^{A^{(1)}}_k$ and $\pi_{A^{(1)}}$ as $\fbox{$\tilde{a}$}$, $P_k$ and $\pi$, respectively.

If $(1,k)<(1,k-1)$ and $(1,k)>(1,k+1)$ then
we get
$\lambda^{(k)}=-x_{1,k}+x_{1,k+1}+\langle h_k,\lambda \rangle
=\fbox{$k+1$}^{A^{(1)}}_k+\langle h_k,\lambda \rangle$
and
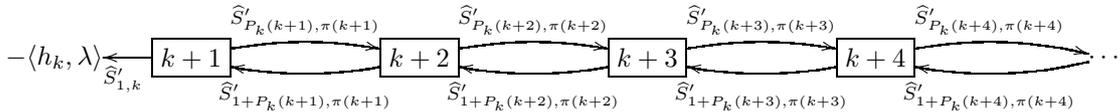
\[
\begin{xy}
(-8,0) *{-\langle h_k,\lambda \rangle}="00",
(0,0) *{\ }="0",
(10,0) *{\fbox{$k+1$}}="r-1r-1",
(40,0)*{\fbox{$k+2$}}="r-2r-2",
(70,0)*{\fbox{$k+3$}}="r-3r-3",
(100,0)*{\fbox{$k+4$}}="r-4r-4",
(130,0)*{\cdots}="dot",
\ar@{->} "r-1r-1";"00"^{\what{S}'_{1,k}}
\ar@/^/ @{->} "r-1r-1";"r-2r-2"^{\what{S}'_{P_k(k+1),\pi(k+1)}}
\ar@/^/ @{->} "r-2r-2";"r-3r-3"^{\what{S}'_{P_k(k+2),\pi(k+2)}}
\ar@/^/ @{->} "r-3r-3";"r-4r-4"^{\what{S}'_{P_k(k+3),\pi(k+3)}}
\ar@/^/ @{->} "r-4r-4";"dot"^{\what{S}'_{P_k(k+4),\pi(k+4)}}
\ar@/_/ @{<-} "r-1r-1";"r-2r-2"_{\what{S}'_{1+P_k(k+1),\pi(k+1)}}
\ar@/_/ @{<-} "r-2r-2";"r-3r-3"_{\what{S}'_{1+P_k(k+2),\pi(k+2)}}
\ar@/_/ @{<-} "r-3r-3";"r-4r-4"_{\what{S}'_{1+P_k(k+3),\pi(k+3)}}
\ar@/_/ @{<-} "r-4r-4";"dot"_{\what{S}'_{1+P_k(k+4),\pi(k+4)}}
\end{xy}
\]
\hspace{-7mm}
\underline{When $\mathfrak{g}$ is of type $B^{(1)}_{n-1}$}

We put $K:=\overline{T}^{A^{(2)}_{2n-3}}_k$ (resp. $K:=\overline{\overline{T}}^{A^{(2)}_{2n-3}}_k)$ and suppose that either
$k\neq3$, $(1,k)<(1,k-1)$ and $(1,k)>(1,k+1)$ or $k=3$, $(1,3)<(1,2), (1,1)$ and $(1,3)>(1,4)$
(resp. either $k\neq3$, $(1,k)>(1,k-1)$ and $(1,k)<(1,k+1)$ or $k=3$, $(1,3)>(1,2),(1,1)$ and $(1,3)<(1,4)$)
then
we get
$\lambda^{(k)}=\fbox{$K+1$}^{A^{(2)}_{2n-3}}_{K}+\langle h_k,\lambda \rangle$
and
\[
\begin{xy}
(-8,0) *{-\langle h_k,\lambda \rangle}="00",
(0,0) *{\ }="0",
(10,0) *{\fbox{$K+1$}}="r-1r-1",
(40,0)*{\fbox{$K+2$}}="r-2r-2",
(70,0)*{\fbox{$K+3$}}="r-3r-3",
(100,0)*{\fbox{$K+4$}}="r-4r-4",
(130,0)*{\cdots}="dot",
\ar@{->} "r-1r-1";"00"^{\what{S}'_{1,k}}
\ar@/^/ @{->} "r-1r-1";"r-2r-2"^{\what{S}'_{P_K(K+1),\pi(K+1)}}
\ar@/^/ @{->} "r-2r-2";"r-3r-3"^{\what{S}'_{P_K(K+2),\pi(K+2)}}
\ar@/^/ @{->} "r-3r-3";"r-4r-4"^{\what{S}'_{P_K(K+3),\pi(K+3)}}
\ar@/^/ @{->} "r-4r-4";"dot"^{\what{S}'_{P_K(K+4),\pi(K+4)}}
\ar@/_/ @{<-} "r-1r-1";"r-2r-2"_{\what{S}'_{1+P_K(K+1),\pi(K+1)}}
\ar@/_/ @{<-} "r-2r-2";"r-3r-3"_{\what{S}'_{1+P_K(K+2),\pi(K+2)}}
\ar@/_/ @{<-} "r-3r-3";"r-4r-4"_{\what{S}'_{1+P_K(K+3),\pi(K+3)}}
\ar@/_/ @{<-} "r-4r-4";"dot"_{\what{S}'_{1+P_K(K+4),\pi(K+4)}}
\end{xy}
\]
\[
\begin{xy}
(-20,0) *{\cdots}="rr",
(10,0) *{\fbox{$\overline{\overline{T}}_1$}}="r-1r-1",
(40,15)*{\fbox{$\overline{\overline{T}}_1+\frac{1}{2}$}}="r-2r-2",
(40,-15)*{\fbox{$\widetilde{\overline{\overline{T}}_1+\frac{1}{2}}$}}="r-2r-2t",
(70,0)*{\fbox{$\overline{\overline{T}}_1+1$}}="r-3r-3",
(100,0)*{\fbox{$\overline{\overline{T}}_1+2$}}="r-4r-4",
(130,0)*{\cdots}="dot",
\ar@/^/ @{->} "rr";"r-1r-1"^{\what{S}'_{P_K(\overline{\overline{T}}_1-1),3}}
\ar@/^/ @{->} "r-1r-1";"r-2r-2"^{\what{S}'_{P_K(\overline{\overline{T}}_1),1}}
\ar@/^/ @{->} "r-1r-1";"r-2r-2t"^{\what{S}'_{P_K(\overline{\overline{T}}_1+\frac{1}{2}),2}}
\ar@/^/ @{->} "r-2r-2";"r-3r-3"^{\what{S}'_{P_K(\overline{\overline{T}}_1+\frac{1}{2}),2}}
\ar@/^/ @{->} "r-2r-2t";"r-3r-3"^{\what{S}'_{P_K(\overline{\overline{T}}_1),1}}
\ar@/^/ @{->} "r-3r-3";"r-4r-4"^{\what{S}'_{P_K(\overline{\overline{T}}_1+1),3}}
\ar@/^/ @{->} "r-4r-4";"dot"^{\what{S}'_{P_K(\overline{\overline{T}}_1+2),4}}
\ar@/_/ @{<-} "rr";"r-1r-1"_{\what{S}'_{1+P_K(\overline{\overline{T}}_1-1),3}}
\ar@/_/ @{<-} "r-1r-1";"r-2r-2"_{\what{S}'_{1+P_K(\overline{\overline{T}}_1),1}}
\ar@/_/ @{<-} "r-1r-1";"r-2r-2t"_{\what{S}'_{1+P_K(\overline{\overline{T}}_1+\frac{1}{2}),2}}
\ar@/_/ @{<-} "r-2r-2";"r-3r-3"_{\what{S}'_{1+P_K(\overline{\overline{T}}_1+\frac{1}{2}),2}}
\ar@/_/ @{<-} "r-2r-2t";"r-3r-3"_{\what{S}'_{1+P_K(\overline{\overline{T}}_1),1}}
\ar@/_/ @{<-} "r-3r-3";"r-4r-4"_{\what{S}'_{1+P_K(\overline{\overline{T}}_1+1),3}}
\ar@/_/ @{<-} "r-4r-4";"dot"_{\what{S}'_{1+P_K(\overline{\overline{T}}_1+2),4}}
\end{xy}
\]
Here, we simply write $\fbox{$a$}^{A^{(2)}_{2n-3}}_K$, $P^{A^{(2)}_{2n-3}}_K$ and $\pi'_{A^{(2)}_{2n-3}}$ as $\fbox{$a$}$, $P_K$ and $\pi$, respectively.

\hspace{-7mm}
\underline{When $\mathfrak{g}$ is of type $C^{(1)}_{n-1}$}

We put $K:=\overline{T}^{D^{(2)}}_k$, (resp. $K:=\overline{\overline{T}}^{D^{(2)}}_k)$
and $T:=\overline{T}^{D^{(2)}}_n$, $T':=\overline{\overline{T}}^{D^{(2)}}_1$ (resp. $T:=\overline{\overline{T}}^{D^{(2)}}_1$, $T':=\overline{\overline{T}}^{D^{(2)}}_n$).
If $(1,k)<(1,k-1)$ and $(1,k)>(1,k+1)$
(resp. $(1,k)>(1,k-1)$ and $(1,k)<(1,k+1)$)
then
we get
$\lambda^{(k)}=\fbox{$K+1$}^{D^{(2)}}_{K}+\langle h_k,\lambda \rangle$
and
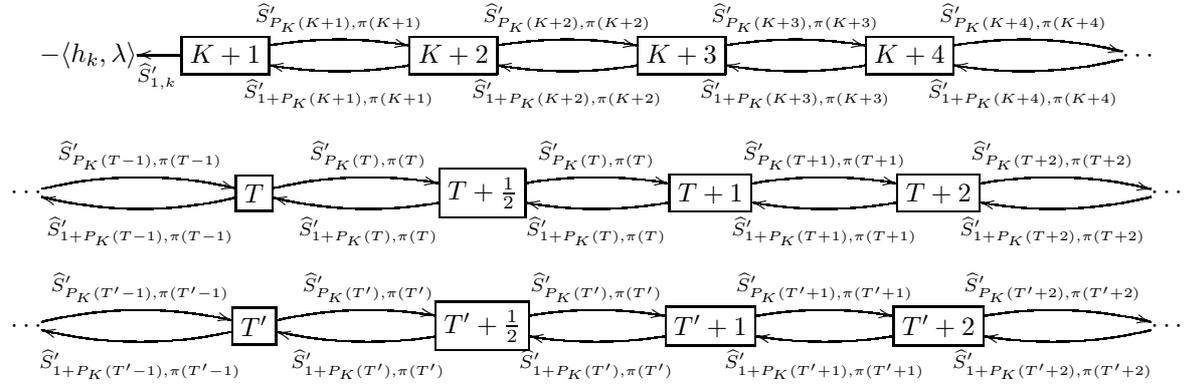
\[
\begin{xy}
(-8,0) *{-\langle h_k,\lambda \rangle}="00",
(0,0) *{\ }="0",
(10,0) *{\fbox{$K+1$}}="r-1r-1",
(40,0)*{\fbox{$K+2$}}="r-2r-2",
(70,0)*{\fbox{$K+3$}}="r-3r-3",
(100,0)*{\fbox{$K+4$}}="r-4r-4",
(130,0)*{\cdots}="dot",
\ar@{->} "r-1r-1";"00"^{\what{S}'_{1,k}}
\ar@/^/ @{->} "r-1r-1";"r-2r-2"^{\what{S}'_{P_K(K+1),\pi(K+1)}}
\ar@/^/ @{->} "r-2r-2";"r-3r-3"^{\what{S}'_{P_K(K+2),\pi(K+2)}}
\ar@/^/ @{->} "r-3r-3";"r-4r-4"^{\what{S}'_{P_K(K+3),\pi(K+3)}}
\ar@/^/ @{->} "r-4r-4";"dot"^{\what{S}'_{P_K(K+4),\pi(K+4)}}
\ar@/_/ @{<-} "r-1r-1";"r-2r-2"_{\what{S}'_{1+P_K(K+1),\pi(K+1)}}
\ar@/_/ @{<-} "r-2r-2";"r-3r-3"_{\what{S}'_{1+P_K(K+2),\pi(K+2)}}
\ar@/_/ @{<-} "r-3r-3";"r-4r-4"_{\what{S}'_{1+P_K(K+3),\pi(K+3)}}
\ar@/_/ @{<-} "r-4r-4";"dot"_{\what{S}'_{1+P_K(K+4),\pi(K+4)}}
\end{xy}
\]
\[
\begin{xy}
(-20,0) *{\cdots}="rr",
(10,0) *{\fbox{$T$}}="r-1r-1",
(40,0)*{\fbox{$T+\frac{1}{2}$}}="r-2r-2",
(70,0)*{\fbox{$T+1$}}="r-3r-3",
(100,0)*{\fbox{$T+2$}}="r-4r-4",
(130,0)*{\cdots}="dot",
\ar@/^/ @{->} "rr";"r-1r-1"^{\what{S}'_{P_K(T-1),\pi(T-1)}}
\ar@/^/ @{->} "r-1r-1";"r-2r-2"^{\what{S}'_{P_K(T),\pi(T)}}
\ar@/^/ @{->} "r-2r-2";"r-3r-3"^{\what{S}'_{P_K(T),\pi(T)}}
\ar@/^/ @{->} "r-3r-3";"r-4r-4"^{\what{S}'_{P_K(T+1),\pi(T+1)}}
\ar@/^/ @{->} "r-4r-4";"dot"^{\what{S}'_{P_K(T+2),\pi(T+2)}}
\ar@/_/ @{<-} "rr";"r-1r-1"_{\what{S}'_{1+P_K(T-1),\pi(T-1)}}
\ar@/_/ @{<-} "r-1r-1";"r-2r-2"_{\what{S}'_{1+P_K(T),\pi(T)}}
\ar@/_/ @{<-} "r-2r-2";"r-3r-3"_{\what{S}'_{1+P_K(T),\pi(T)}}
\ar@/_/ @{<-} "r-3r-3";"r-4r-4"_{\what{S}'_{1+P_K(T+1),\pi(T+1)}}
\ar@/_/ @{<-} "r-4r-4";"dot"_{\what{S}'_{1+P_K(T+2),\pi(T+2)}}
\end{xy}
\]
\[
\begin{xy}
(-20,0) *{\cdots}="rr",
(10,0) *{\fbox{$T'$}}="r-1r-1",
(40,0)*{\fbox{$T'+\frac{1}{2}$}}="r-2r-2",
(70,0)*{\fbox{$T'+1$}}="r-3r-3",
(100,0)*{\fbox{$T'+2$}}="r-4r-4",
(130,0)*{\cdots}="dot",
\ar@/^/ @{->} "rr";"r-1r-1"^{\what{S}'_{P_K(T'-1),\pi(T'-1)}}
\ar@/^/ @{->} "r-1r-1";"r-2r-2"^{\what{S}'_{P_K(T'),\pi(T')}}
\ar@/^/ @{->} "r-2r-2";"r-3r-3"^{\what{S}'_{P_K(T'),\pi(T')}}
\ar@/^/ @{->} "r-3r-3";"r-4r-4"^{\what{S}'_{P_K(T'+1),\pi(T'+1)}}
\ar@/^/ @{->} "r-4r-4";"dot"^{\what{S}'_{P_K(T'+2),\pi(T'+2)}}
\ar@/_/ @{<-} "rr";"r-1r-1"_{\what{S}'_{1+P_K(T'-1),\pi(T'-1)}}
\ar@/_/ @{<-} "r-1r-1";"r-2r-2"_{\what{S}'_{1+P_K(T'),\pi(T')}}
\ar@/_/ @{<-} "r-2r-2";"r-3r-3"_{\what{S}'_{1+P_K(T'),\pi(T')}}
\ar@/_/ @{<-} "r-3r-3";"r-4r-4"_{\what{S}'_{1+P_K(T'+1),\pi(T'+1)}}
\ar@/_/ @{<-} "r-4r-4";"dot"_{\what{S}'_{1+P_K(T'+2),\pi(T'+2)}}
\end{xy}
\]
Here, we simply write $\fbox{$a$}^{D^{(2)}}_K$, $P^{D^{(2)}}_K$ and $\pi'_{D^{(2)}}$ as $\fbox{$a$}$, $P_K$ and $\pi$, respectively.

\hspace{-7mm}
\underline{When $\mathfrak{g}$ is of type $D^{(1)}_{n-1}$}

We put $K:=\overline{T}^{D^{(1)}}_k$, (resp. $K:=\overline{\overline{T}}^{D^{(1)}}_k)$
and $T:=\overline{T}^{D^{(1)}}_{n-1}$, $T':=\overline{\overline{T}}^{D^{(1)}}_1$ (resp. $T:=\overline{\overline{T}}^{D^{(1)}}_1$, $T':=\overline{\overline{T}}^{D^{(1)}}_{n-1}$).
We assume that for all $t\in I$ such that $t<k$ and $a_{k,t}=-1$, it holds
$(1,k)<(1,t)$ (resp. $(1,k)>(1,t)$) and for all $t'\in I$ such that $t'>k$ and $a_{k,t'}=-1$, it holds
$(1,k)>(1,t')$ (resp. $(1,k)<(1,t')$).
Then one obtains
$\lambda^{(k)}=\fbox{$K+1$}^{D^{(1)}}_{K}+\langle h_k,\lambda \rangle$
and
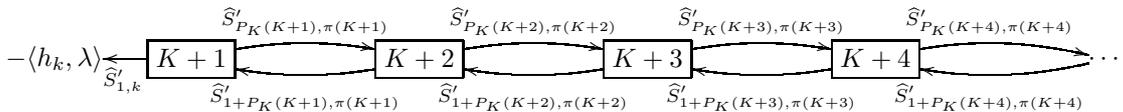
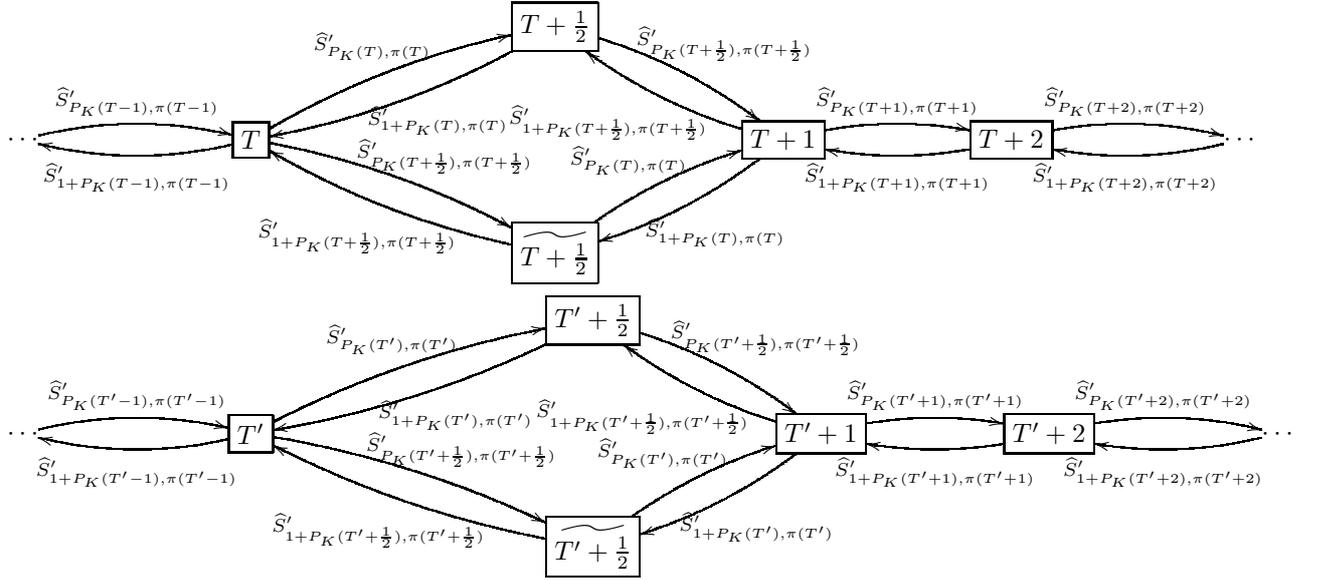
\[
\begin{xy}
(-8,0) *{-\langle h_k,\lambda \rangle}="00",
(0,0) *{\ }="0",
(10,0) *{\fbox{$K+1$}}="r-1r-1",
(40,0)*{\fbox{$K+2$}}="r-2r-2",
(70,0)*{\fbox{$K+3$}}="r-3r-3",
(100,0)*{\fbox{$K+4$}}="r-4r-4",
(130,0)*{\cdots}="dot",
\ar@{->} "r-1r-1";"00"^{\what{S}'_{1,k}}
\ar@/^/ @{->} "r-1r-1";"r-2r-2"^{\what{S}'_{P_K(K+1),\pi(K+1)}}
\ar@/^/ @{->} "r-2r-2";"r-3r-3"^{\what{S}'_{P_K(K+2),\pi(K+2)}}
\ar@/^/ @{->} "r-3r-3";"r-4r-4"^{\what{S}'_{P_K(K+3),\pi(K+3)}}
\ar@/^/ @{->} "r-4r-4";"dot"^{\what{S}'_{P_K(K+4),\pi(K+4)}}
\ar@/_/ @{<-} "r-1r-1";"r-2r-2"_{\what{S}'_{1+P_K(K+1),\pi(K+1)}}
\ar@/_/ @{<-} "r-2r-2";"r-3r-3"_{\what{S}'_{1+P_K(K+2),\pi(K+2)}}
\ar@/_/ @{<-} "r-3r-3";"r-4r-4"_{\what{S}'_{1+P_K(K+3),\pi(K+3)}}
\ar@/_/ @{<-} "r-4r-4";"dot"_{\what{S}'_{1+P_K(K+4),\pi(K+4)}}
\end{xy}
\]
\[
\begin{xy}
(-20,0) *{\cdots}="rr",
(10,0) *{\fbox{$T$}}="r-1r-1",
(50,15)*{\fbox{$T+\frac{1}{2}$}}="r-2r-2",
(50,-15)*{\fbox{$\widetilde{T+\frac{1}{2}}$}}="r-2r-2t",
(80,0)*{\fbox{$T+1$}}="r-3r-3",
(110,0)*{\fbox{$T+2$}}="r-4r-4",
(140,0)*{\cdots}="dot",
\ar@/^/ @{->} "rr";"r-1r-1"^{\what{S}'_{P_K(T-1),\pi(T-1)}}
\ar@/^/ @{->} "r-1r-1";"r-2r-2"^{\what{S}'_{P_K(T),\pi(T)}}
\ar@/^/ @{->} "r-1r-1";"r-2r-2t"^{\what{S}'_{P_K(T+\frac{1}{2}),\pi(T+\frac{1}{2})}}
\ar@/^/ @{->} "r-2r-2";"r-3r-3"^{\what{S}'_{P_K(T+\frac{1}{2}),\pi(T+\frac{1}{2})}}
\ar@/^/ @{->} "r-2r-2t";"r-3r-3"^{\what{S}'_{P_K(T),\pi(T)}}
\ar@/^/ @{->} "r-3r-3";"r-4r-4"^{\what{S}'_{P_K(T+1),\pi(T+1)}}
\ar@/^/ @{->} "r-4r-4";"dot"^{\what{S}'_{P_K(T+2),\pi(T+2)}}
\ar@/_/ @{<-} "rr";"r-1r-1"_{\what{S}'_{1+P_K(T-1),\pi(T-1)}}
\ar@/_/ @{<-} "r-1r-1";"r-2r-2"_{\what{S}'_{1+P_K(T),\pi(T)}}
\ar@/_/ @{<-} "r-1r-1";"r-2r-2t"_{\what{S}'_{1+P_K(T+\frac{1}{2}),\pi(T+\frac{1}{2})}}
\ar@/_/ @{<-} "r-2r-2";"r-3r-3"_{\what{S}'_{1+P_K(T+\frac{1}{2}),\pi(T+\frac{1}{2})}}
\ar@/_/ @{<-} "r-2r-2t";"r-3r-3"_{\what{S}'_{1+P_K(T),\pi(T)}}
\ar@/_/ @{<-} "r-3r-3";"r-4r-4"_{\what{S}'_{1+P_K(T+1),\pi(T+1)}}
\ar@/_/ @{<-} "r-4r-4";"dot"_{\what{S}'_{1+P_K(T+2),\pi(T+2)}}
\end{xy}
\]
\[
\begin{xy}
(-20,0) *{\cdots}="rr",
(10,0) *{\fbox{$T'$}}="r-1r-1",
(55,15)*{\fbox{$T'+\frac{1}{2}$}}="r-2r-2",
(55,-15)*{\fbox{$\widetilde{T'+\frac{1}{2}}$}}="r-2r-2t",
(85,0)*{\fbox{$T'+1$}}="r-3r-3",
(115,0)*{\fbox{$T'+2$}}="r-4r-4",
(145,0)*{\cdots}="dot",
\ar@/^/ @{->} "rr";"r-1r-1"^{\what{S}'_{P_K(T'-1),\pi(T'-1)}}
\ar@/^/ @{->} "r-1r-1";"r-2r-2"^{\what{S}'_{P_K(T'),\pi(T')}}
\ar@/^/ @{->} "r-1r-1";"r-2r-2t"^{\what{S}'_{P_K(T'+\frac{1}{2}),\pi(T'+\frac{1}{2})}}
\ar@/^/ @{->} "r-2r-2";"r-3r-3"^{\what{S}'_{P_K(T'+\frac{1}{2}),\pi(T'+\frac{1}{2})}}
\ar@/^/ @{->} "r-2r-2t";"r-3r-3"^{\what{S}'_{P_K(T'),\pi(T')}}
\ar@/^/ @{->} "r-3r-3";"r-4r-4"^{\what{S}'_{P_K(T'+1),\pi(T'+1)}}
\ar@/^/ @{->} "r-4r-4";"dot"^{\what{S}'_{P_K(T'+2),\pi(T'+2)}}
\ar@/_/ @{<-} "rr";"r-1r-1"_{\what{S}'_{1+P_K(T'-1),\pi(T'-1)}}
\ar@/_/ @{<-} "r-1r-1";"r-2r-2"_{\what{S}'_{1+P_K(T'),\pi(T')}}
\ar@/_/ @{<-} "r-1r-1";"r-2r-2t"_{\what{S}'_{1+P_K(T'+\frac{1}{2}),\pi(T'+\frac{1}{2})}}
\ar@/_/ @{<-} "r-2r-2";"r-3r-3"_{\what{S}'_{1+P_K(T'+\frac{1}{2}),\pi(T'+\frac{1}{2})}}
\ar@/_/ @{<-} "r-2r-2t";"r-3r-3"_{\what{S}'_{1+P_K(T'),\pi(T')}}
\ar@/_/ @{<-} "r-3r-3";"r-4r-4"_{\what{S}'_{1+P_K(T'+1),\pi(T'+1)}}
\ar@/_/ @{<-} "r-4r-4";"dot"_{\what{S}'_{1+P_K(T'+2),\pi(T'+2)}}
\end{xy}
\]
Here, we simply write $\fbox{$a$}^{D^{(1)}}_K$, $P^{D^{(1)}}_K$ and $\pi'_{D^{(1)}}$ as $\fbox{$a$}$, $P_K$ and $\pi$, respectively.

\hspace{-7mm}
\underline{When $\mathfrak{g}$ is of type $A^{(2)}_{2n-2}$}

We put $K:=\overline{T}^{A^{(2)\dagger}_{2n-2}}_k$ (resp. $K:=\overline{\overline{T}}^{A^{(2)\dagger}_{2n-2}}_k)$.
If
$(1,k)<(1,k-1)$ and $(1,k)>(1,k+1)$ 
(resp. $(1,k)>(1,k-1)$ and $(1,k)<(1,k+1)$)
then
we get
$\lambda^{(k)}=\fbox{$K+1$}^{A^{(2)\dagger}_{2n-2}}_{K}+\langle h_k,\lambda \rangle$
and
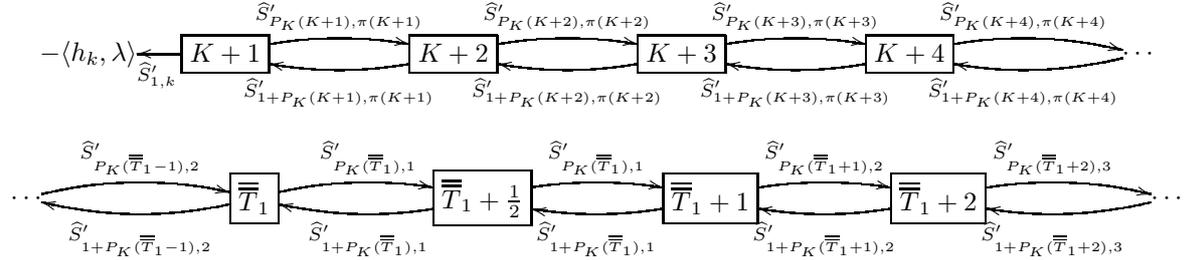
\[
\begin{xy}
(-8,0) *{-\langle h_k,\lambda \rangle}="00",
(0,0) *{\ }="0",
(10,0) *{\fbox{$K+1$}}="r-1r-1",
(40,0)*{\fbox{$K+2$}}="r-2r-2",
(70,0)*{\fbox{$K+3$}}="r-3r-3",
(100,0)*{\fbox{$K+4$}}="r-4r-4",
(130,0)*{\cdots}="dot",
\ar@{->} "r-1r-1";"00"^{\what{S}'_{1,k}}
\ar@/^/ @{->} "r-1r-1";"r-2r-2"^{\what{S}'_{P_K(K+1),\pi(K+1)}}
\ar@/^/ @{->} "r-2r-2";"r-3r-3"^{\what{S}'_{P_K(K+2),\pi(K+2)}}
\ar@/^/ @{->} "r-3r-3";"r-4r-4"^{\what{S}'_{P_K(K+3),\pi(K+3)}}
\ar@/^/ @{->} "r-4r-4";"dot"^{\what{S}'_{P_K(K+4),\pi(K+4)}}
\ar@/_/ @{<-} "r-1r-1";"r-2r-2"_{\what{S}'_{1+P_K(K+1),\pi(K+1)}}
\ar@/_/ @{<-} "r-2r-2";"r-3r-3"_{\what{S}'_{1+P_K(K+2),\pi(K+2)}}
\ar@/_/ @{<-} "r-3r-3";"r-4r-4"_{\what{S}'_{1+P_K(K+3),\pi(K+3)}}
\ar@/_/ @{<-} "r-4r-4";"dot"_{\what{S}'_{1+P_K(K+4),\pi(K+4)}}
\end{xy}
\]
\[
\begin{xy}
(-20,0) *{\cdots}="rr",
(10,0) *{\fbox{$\overline{\overline{T}}_1$}}="r-1r-1",
(40,0)*{\fbox{$\overline{\overline{T}}_1+\frac{1}{2}$}}="r-2r-2",
(70,0)*{\fbox{$\overline{\overline{T}}_1+1$}}="r-3r-3",
(100,0)*{\fbox{$\overline{\overline{T}}_1+2$}}="r-4r-4",
(130,0)*{\cdots}="dot",
\ar@/^/ @{->} "rr";"r-1r-1"^{\what{S}'_{P_K(\overline{\overline{T}}_1-1),2}}
\ar@/^/ @{->} "r-1r-1";"r-2r-2"^{\what{S}'_{P_K(\overline{\overline{T}}_1),1}}
\ar@/^/ @{->} "r-2r-2";"r-3r-3"^{\what{S}'_{P_K(\overline{\overline{T}}_1),1}}
\ar@/^/ @{->} "r-3r-3";"r-4r-4"^{\what{S}'_{P_K(\overline{\overline{T}}_1+1),2}}
\ar@/^/ @{->} "r-4r-4";"dot"^{\what{S}'_{P_K(\overline{\overline{T}}_1+2),3}}
\ar@/_/ @{<-} "rr";"r-1r-1"_{\what{S}'_{1+P_K(\overline{\overline{T}}_1-1),2}}
\ar@/_/ @{<-} "r-1r-1";"r-2r-2"_{\what{S}'_{1+P_K(\overline{\overline{T}}_1),1}}
\ar@/_/ @{<-} "r-2r-2";"r-3r-3"_{\what{S}'_{1+P_K(\overline{\overline{T}}_1),1}}
\ar@/_/ @{<-} "r-3r-3";"r-4r-4"_{\what{S}'_{1+P_K(\overline{\overline{T}}_1+1),2}}
\ar@/_/ @{<-} "r-4r-4";"dot"_{\what{S}'_{1+P_K(\overline{\overline{T}}_1+2),3}}
\end{xy}
\]
Here, we simply write $\fbox{$a$}^{A^{(2)\dagger}_{2n-2}}_K$, $P^{A^{(2)\dagger}_{2n-2}}_K$ and $\pi'_{A^{(2)\dagger}_{2n-2}}$ as $\fbox{$a$}$, $P_K$ and $\pi$, respectively.

\hspace{-7mm}
\underline{When $\mathfrak{g}$ is of type $A^{(2)}_{2n-3}$}

We suppose that either
$k\neq3$, $(1,k)<(1,k-1)$ and $(1,k)>(1,k+1)$ or $k=3$, $(1,3)<(1,2), (1,1)$ and $(1,3)>(1,4)$.
We write
$\overline{T}^{B^{(1)}}_m=\overline{T}_m$, $\overline{\overline{T}}^{B^{(1)}}_m=\overline{\overline{T}}_m$.
Then
$\lambda^{(k)}=\fbox{$\overline{T}_k+1$}^{B^{(1)}}_{\overline{T}_k}+\langle h_k,\lambda \rangle$
and
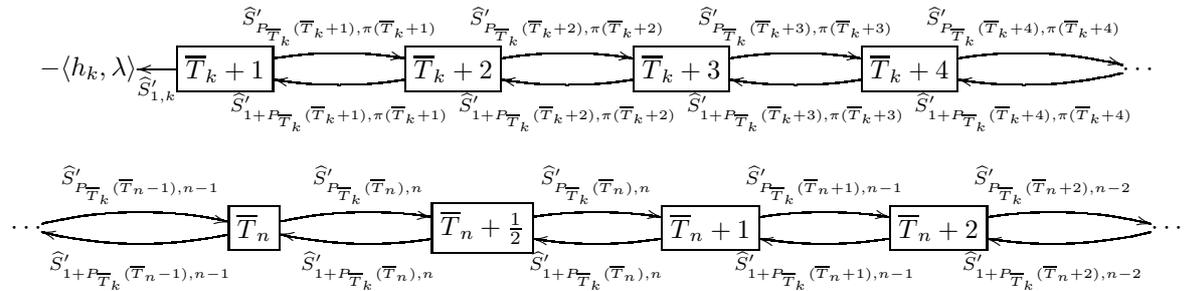
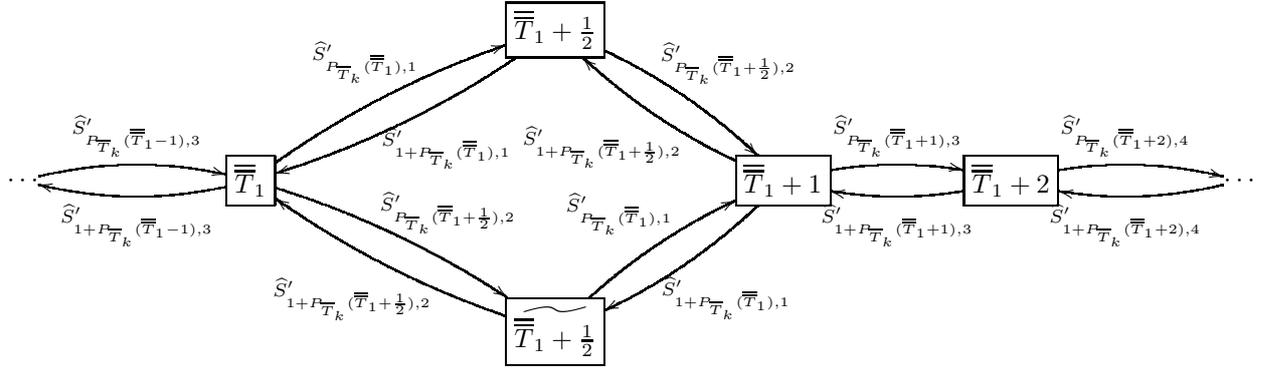
\[
\begin{xy}
(-8,0) *{-\langle h_k,\lambda \rangle}="00",
(0,0) *{\ }="0",
(10,0) *{\fbox{$\overline{T}_k+1$}}="r-1r-1",
(40,0)*{\fbox{$\overline{T}_k+2$}}="r-2r-2",
(70,0)*{\fbox{$\overline{T}_k+3$}}="r-3r-3",
(100,0)*{\fbox{$\overline{T}_k+4$}}="r-4r-4",
(130,0)*{\cdots}="dot",
\ar@{->} "r-1r-1";"00"^{\what{S}'_{1,k}}
\ar@/^/ @{->} "r-1r-1";"r-2r-2"^{\what{S}'_{P_{\overline{T}_k}(\overline{T}_k+1),\pi(\overline{T}_k+1)}}
\ar@/^/ @{->} "r-2r-2";"r-3r-3"^{\what{S}'_{P_{\overline{T}_k}(\overline{T}_k+2),\pi(\overline{T}_k+2)}}
\ar@/^/ @{->} "r-3r-3";"r-4r-4"^{\what{S}'_{P_{\overline{T}_k}(\overline{T}_k+3),\pi(\overline{T}_k+3)}}
\ar@/^/ @{->} "r-4r-4";"dot"^{\what{S}'_{P_{\overline{T}_k}(\overline{T}_k+4),\pi(\overline{T}_k+4)}}
\ar@/_/ @{<-} "r-1r-1";"r-2r-2"_{\what{S}'_{1+P_{\overline{T}_k}(\overline{T}_k+1),\pi(\overline{T}_k+1)}}
\ar@/_/ @{<-} "r-2r-2";"r-3r-3"_{\what{S}'_{1+P_{\overline{T}_k}(\overline{T}_k+2),\pi(\overline{T}_k+2)}}
\ar@/_/ @{<-} "r-3r-3";"r-4r-4"_{\what{S}'_{1+P_{\overline{T}_k}(\overline{T}_k+3),\pi(\overline{T}_k+3)}}
\ar@/_/ @{<-} "r-4r-4";"dot"_{\what{S}'_{1+P_{\overline{T}_k}(\overline{T}_k+4),\pi(\overline{T}_k+4)}}
\end{xy}
\]
\[
\begin{xy}
(-20,0) *{\cdots}="rr",
(10,0) *{\fbox{$\overline{T}_n$}}="r-1r-1",
(40,0)*{\fbox{$\overline{T}_n+\frac{1}{2}$}}="r-2r-2",
(70,0)*{\fbox{$\overline{T}_n+1$}}="r-3r-3",
(100,0)*{\fbox{$\overline{T}_n+2$}}="r-4r-4",
(130,0)*{\cdots}="dot",
\ar@/^/ @{->} "rr";"r-1r-1"^{\what{S}'_{P_{\overline{T}_k}(\overline{T}_n-1),n-1}}
\ar@/^/ @{->} "r-1r-1";"r-2r-2"^{\what{S}'_{P_{\overline{T}_k}(\overline{T}_n),n}}
\ar@/^/ @{->} "r-2r-2";"r-3r-3"^{\what{S}'_{P_{\overline{T}_k}(\overline{T}_n),n}}
\ar@/^/ @{->} "r-3r-3";"r-4r-4"^{\what{S}'_{P_{\overline{T}_k}(\overline{T}_n+1),n-1}}
\ar@/^/ @{->} "r-4r-4";"dot"^{\what{S}'_{P_{\overline{T}_k}(\overline{T}_n+2),n-2}}
\ar@/_/ @{<-} "rr";"r-1r-1"_{\what{S}'_{1+P_{\overline{T}_k}(\overline{T}_n-1),n-1}}
\ar@/_/ @{<-} "r-1r-1";"r-2r-2"_{\what{S}'_{1+P_{\overline{T}_k}(\overline{T}_n),n}}
\ar@/_/ @{<-} "r-2r-2";"r-3r-3"_{\what{S}'_{1+P_{\overline{T}_k}(\overline{T}_n),n}}
\ar@/_/ @{<-} "r-3r-3";"r-4r-4"_{\what{S}'_{1+P_{\overline{T}_k}(\overline{T}_n+1),n-1}}
\ar@/_/ @{<-} "r-4r-4";"dot"_{\what{S}'_{1+P_{\overline{T}_k}(\overline{T}_n+2),n-2}}
\end{xy}
\]
\[
\begin{xy}
(-20,0) *{\cdots}="rr",
(10,0) *{\fbox{$\overline{\overline{T}}_1$}}="r-1r-1",
(50,20)*{\fbox{$\overline{\overline{T}}_1+\frac{1}{2}$}}="r-2r-2",
(50,-20)*{\fbox{$\widetilde{\overline{\overline{T}}_1+\frac{1}{2}}$}}="r-2r-2t",
(80,0)*{\fbox{$\overline{\overline{T}}_1+1$}}="r-3r-3",
(110,0)*{\fbox{$\overline{\overline{T}}_1+2$}}="r-4r-4",
(140,0)*{\cdots}="dot",
\ar@/^/ @{->} "rr";"r-1r-1"^{\what{S}'_{P_{\overline{T}_k}(\overline{\overline{T}}_1-1),3}}
\ar@/^/ @{->} "r-1r-1";"r-2r-2"^{\what{S}'_{P_{\overline{T}_k}(\overline{\overline{T}}_1),1}}
\ar@/^/ @{->} "r-1r-1";"r-2r-2t"^{\what{S}'_{P_{\overline{T}_k}(\overline{\overline{T}}_1+\frac{1}{2}),2}}
\ar@/^/ @{->} "r-2r-2";"r-3r-3"^{\what{S}'_{P_{\overline{T}_k}(\overline{\overline{T}}_1+\frac{1}{2}),2}}
\ar@/^/ @{->} "r-2r-2t";"r-3r-3"^{\what{S}'_{P_{\overline{T}_k}(\overline{\overline{T}}_1),1}}
\ar@/^/ @{->} "r-3r-3";"r-4r-4"^{\what{S}'_{P_{\overline{T}_k}(\overline{\overline{T}}_1+1),3}}
\ar@/^/ @{->} "r-4r-4";"dot"^{\what{S}'_{P_{\overline{T}_k}(\overline{\overline{T}}_1+2),4}}
\ar@/_/ @{<-} "rr";"r-1r-1"_{\what{S}'_{1+P_{\overline{T}_k}(\overline{\overline{T}}_1-1),3}}
\ar@/_/ @{<-} "r-1r-1";"r-2r-2"_{\what{S}'_{1+P_{\overline{T}_k}(\overline{\overline{T}}_1),1}}
\ar@/_/ @{<-} "r-1r-1";"r-2r-2t"_{\what{S}'_{1+P_{\overline{T}_k}(\overline{\overline{T}}_1+\frac{1}{2}),2}}
\ar@/_/ @{<-} "r-2r-2";"r-3r-3"_{\what{S}'_{1+P_{\overline{T}_k}(\overline{\overline{T}}_1+\frac{1}{2}),2}}
\ar@/_/ @{<-} "r-2r-2t";"r-3r-3"_{\what{S}'_{1+P_{\overline{T}_k}(\overline{\overline{T}}_1),1}}
\ar@/_/ @{<-} "r-3r-3";"r-4r-4"_{\what{S}'_{1+P_{\overline{T}_k}(\overline{\overline{T}}_1+1),3}}
\ar@/_/ @{<-} "r-4r-4";"dot"_{\what{S}'_{1+P_{\overline{T}_k}(\overline{\overline{T}}_1+2),4}}
\end{xy}
\]
Here, we simply write
and $\fbox{$a$}^{B^{(1)}}_{\overline{T}_k}$, $P^{B^{(1)}}_{\overline{T}_k}$ and $\pi'_{B^{(1)}}$ as $\fbox{$a$}$, $P_{\overline{T}_k}$ and $\pi$, respectively.

We suppose that either
$k\neq3$, $(1,k)>(1,k-1)$ and $(1,k)<(1,k+1)$ or $k=3$, $(1,3)>(1,2), (1,1)$ and $(1,3)<(1,4)$.
Then
$\lambda^{(k)}=\fbox{$\overline{\overline{T}}_k+1$}^{B^{(1)}}_{\overline{\overline{T}}_k}+\langle h_k,\lambda \rangle$
and
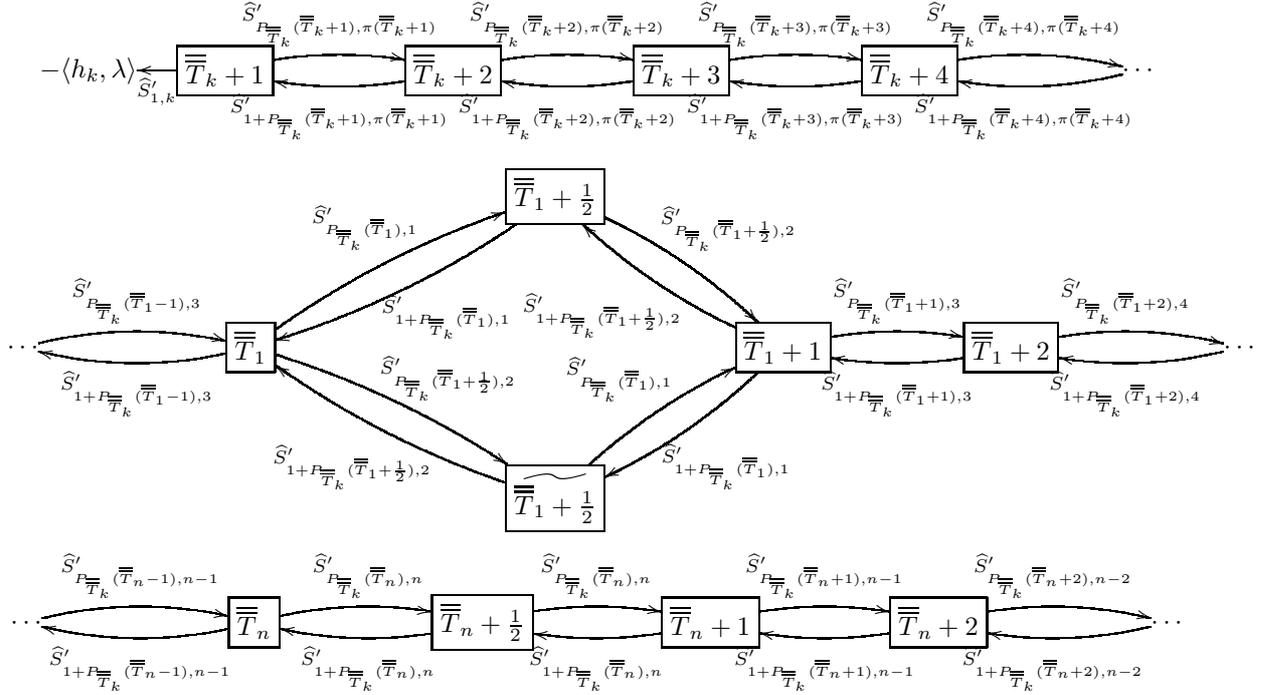
\[
\begin{xy}
(-8,0) *{-\langle h_k,\lambda \rangle}="00",
(0,0) *{\ }="0",
(10,0) *{\fbox{$\overline{\overline{T}}_k+1$}}="r-1r-1",
(40,0)*{\fbox{$\overline{\overline{T}}_k+2$}}="r-2r-2",
(70,0)*{\fbox{$\overline{\overline{T}}_k+3$}}="r-3r-3",
(100,0)*{\fbox{$\overline{\overline{T}}_k+4$}}="r-4r-4",
(130,0)*{\cdots}="dot",
\ar@{->} "r-1r-1";"00"^{\what{S}'_{1,k}}
\ar@/^/ @{->} "r-1r-1";"r-2r-2"^{\what{S}'_{P_{\overline{\overline{T}}_k}(\overline{\overline{T}}_k+1),\pi(\overline{\overline{T}}_k+1)}}
\ar@/^/ @{->} "r-2r-2";"r-3r-3"^{\what{S}'_{P_{\overline{\overline{T}}_k}(\overline{\overline{T}}_k+2),\pi(\overline{\overline{T}}_k+2)}}
\ar@/^/ @{->} "r-3r-3";"r-4r-4"^{\what{S}'_{P_{\overline{\overline{T}}_k}(\overline{\overline{T}}_k+3),\pi(\overline{\overline{T}}_k+3)}}
\ar@/^/ @{->} "r-4r-4";"dot"^{\what{S}'_{P_{\overline{\overline{T}}_k}(\overline{\overline{T}}_k+4),\pi(\overline{\overline{T}}_k+4)}}
\ar@/_/ @{<-} "r-1r-1";"r-2r-2"_{\what{S}'_{1+P_{\overline{\overline{T}}_k}(\overline{\overline{T}}_k+1),\pi(\overline{\overline{T}}_k+1)}}
\ar@/_/ @{<-} "r-2r-2";"r-3r-3"_{\what{S}'_{1+P_{\overline{\overline{T}}_k}(\overline{\overline{T}}_k+2),\pi(\overline{\overline{T}}_k+2)}}
\ar@/_/ @{<-} "r-3r-3";"r-4r-4"_{\what{S}'_{1+P_{\overline{\overline{T}}_k}(\overline{\overline{T}}_k+3),\pi(\overline{\overline{T}}_k+3)}}
\ar@/_/ @{<-} "r-4r-4";"dot"_{\what{S}'_{1+P_{\overline{\overline{T}}_k}(\overline{\overline{T}}_k+4),\pi(\overline{\overline{T}}_k+4)}}
\end{xy}
\]
\[
\begin{xy}
(-20,0) *{\cdots}="rr",
(10,0) *{\fbox{$\overline{\overline{T}}_1$}}="r-1r-1",
(50,20)*{\fbox{$\overline{\overline{T}}_1+\frac{1}{2}$}}="r-2r-2",
(50,-20)*{\fbox{$\widetilde{\overline{\overline{T}}_1+\frac{1}{2}}$}}="r-2r-2t",
(80,0)*{\fbox{$\overline{\overline{T}}_1+1$}}="r-3r-3",
(110,0)*{\fbox{$\overline{\overline{T}}_1+2$}}="r-4r-4",
(140,0)*{\cdots}="dot",
\ar@/^/ @{->} "rr";"r-1r-1"^{\what{S}'_{P_{\overline{\overline{T}}_k}(\overline{\overline{T}}_1-1),3}}
\ar@/^/ @{->} "r-1r-1";"r-2r-2"^{\what{S}'_{P_{\overline{\overline{T}}_k}(\overline{\overline{T}}_1),1}}
\ar@/^/ @{->} "r-1r-1";"r-2r-2t"^{\what{S}'_{P_{\overline{\overline{T}}_k}(\overline{\overline{T}}_1+\frac{1}{2}),2}}
\ar@/^/ @{->} "r-2r-2";"r-3r-3"^{\what{S}'_{P_{\overline{\overline{T}}_k}(\overline{\overline{T}}_1+\frac{1}{2}),2}}
\ar@/^/ @{->} "r-2r-2t";"r-3r-3"^{\what{S}'_{P_{\overline{\overline{T}}_k}(\overline{\overline{T}}_1),1}}
\ar@/^/ @{->} "r-3r-3";"r-4r-4"^{\what{S}'_{P_{\overline{\overline{T}}_k}(\overline{\overline{T}}_1+1),3}}
\ar@/^/ @{->} "r-4r-4";"dot"^{\what{S}'_{P_{\overline{\overline{T}}_k}(\overline{\overline{T}}_1+2),4}}
\ar@/_/ @{<-} "rr";"r-1r-1"_{\what{S}'_{1+P_{\overline{\overline{T}}_k}(\overline{\overline{T}}_1-1),3}}
\ar@/_/ @{<-} "r-1r-1";"r-2r-2"_{\what{S}'_{1+P_{\overline{\overline{T}}_k}(\overline{\overline{T}}_1),1}}
\ar@/_/ @{<-} "r-1r-1";"r-2r-2t"_{\what{S}'_{1+P_{\overline{\overline{T}}_k}(\overline{\overline{T}}_1+\frac{1}{2}),2}}
\ar@/_/ @{<-} "r-2r-2";"r-3r-3"_{\what{S}'_{1+P_{\overline{\overline{T}}_k}(\overline{\overline{T}}_1+\frac{1}{2}),2}}
\ar@/_/ @{<-} "r-2r-2t";"r-3r-3"_{\what{S}'_{1+P_{\overline{\overline{T}}_k}(\overline{\overline{T}}_1),1}}
\ar@/_/ @{<-} "r-3r-3";"r-4r-4"_{\what{S}'_{1+P_{\overline{\overline{T}}_k}(\overline{\overline{T}}_1+1),3}}
\ar@/_/ @{<-} "r-4r-4";"dot"_{\what{S}'_{1+P_{\overline{\overline{T}}_k}(\overline{\overline{T}}_1+2),4}}
\end{xy}
\]
\[
\begin{xy}
(-20,0) *{\cdots}="rr",
(10,0) *{\fbox{$\overline{\overline{T}}_n$}}="r-1r-1",
(40,0)*{\fbox{$\overline{\overline{T}}_n+\frac{1}{2}$}}="r-2r-2",
(70,0)*{\fbox{$\overline{\overline{T}}_n+1$}}="r-3r-3",
(100,0)*{\fbox{$\overline{\overline{T}}_n+2$}}="r-4r-4",
(130,0)*{\cdots}="dot",
\ar@/^/ @{->} "rr";"r-1r-1"^{\what{S}'_{P_{\overline{\overline{T}}_k}(\overline{\overline{T}}_n-1),n-1}}
\ar@/^/ @{->} "r-1r-1";"r-2r-2"^{\what{S}'_{P_{\overline{\overline{T}}_k}(\overline{\overline{T}}_n),n}}
\ar@/^/ @{->} "r-2r-2";"r-3r-3"^{\what{S}'_{P_{\overline{\overline{T}}_k}(\overline{\overline{T}}_n),n}}
\ar@/^/ @{->} "r-3r-3";"r-4r-4"^{\what{S}'_{P_{\overline{\overline{T}}_k}(\overline{\overline{T}}_n+1),n-1}}
\ar@/^/ @{->} "r-4r-4";"dot"^{\what{S}'_{P_{\overline{\overline{T}}_k}(\overline{\overline{T}}_n+2),n-2}}
\ar@/_/ @{<-} "rr";"r-1r-1"_{\what{S}'_{1+P_{\overline{\overline{T}}_k}(\overline{\overline{T}}_n-1),n-1}}
\ar@/_/ @{<-} "r-1r-1";"r-2r-2"_{\what{S}'_{1+P_{\overline{\overline{T}}_k}(\overline{\overline{T}}_n),n}}
\ar@/_/ @{<-} "r-2r-2";"r-3r-3"_{\what{S}'_{1+P_{\overline{\overline{T}}_k}(\overline{\overline{T}}_n),n}}
\ar@/_/ @{<-} "r-3r-3";"r-4r-4"_{\what{S}'_{1+P_{\overline{\overline{T}}_k}(\overline{\overline{T}}_n+1),n-1}}
\ar@/_/ @{<-} "r-4r-4";"dot"_{\what{S}'_{1+P_{\overline{\overline{T}}_k}(\overline{\overline{T}}_n+2),n-2}}
\end{xy}
\]
Here, we simply write $\fbox{$a$}^{B^{(1)}}_{\overline{\overline{T}}_k}$, $P^{B^{(1)}}_{\overline{\overline{T}}_k}$ and $\pi'_{B^{(1)}}$ as $\fbox{$a$}$, $P_{\overline{\overline{T}}_k}$ and $\pi$, respectively.

\hspace{-7mm}
\underline{When $\mathfrak{g}$ is of type $D^{(2)}_{n}$}

We put $K:=\overline{T}^{C^{(1)}}_k$ (resp. $K:=\overline{\overline{T}}^{C^{(1)}}_k)$.
If $(1,k)<(1,k-1)$ and $(1,k)>(1,k+1)$
(resp. $(1,k)>(1,k-1)$ and $(1,k)<(1,k+1)$)
then
we get
$\lambda^{(k)}=\fbox{$K+1$}^{C^{(1)}}_{K}+\langle h_k,\lambda \rangle$
and
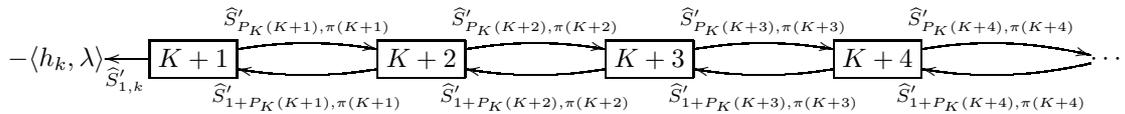
\[
\begin{xy}
(-8,0) *{-\langle h_k,\lambda \rangle}="00",
(0,0) *{\ }="0",
(10,0) *{\fbox{$K+1$}}="r-1r-1",
(40,0)*{\fbox{$K+2$}}="r-2r-2",
(70,0)*{\fbox{$K+3$}}="r-3r-3",
(100,0)*{\fbox{$K+4$}}="r-4r-4",
(130,0)*{\cdots}="dot",
\ar@{->} "r-1r-1";"00"^{\what{S}'_{1,k}}
\ar@/^/ @{->} "r-1r-1";"r-2r-2"^{\what{S}'_{P_K(K+1),\pi(K+1)}}
\ar@/^/ @{->} "r-2r-2";"r-3r-3"^{\what{S}'_{P_K(K+2),\pi(K+2)}}
\ar@/^/ @{->} "r-3r-3";"r-4r-4"^{\what{S}'_{P_K(K+3),\pi(K+3)}}
\ar@/^/ @{->} "r-4r-4";"dot"^{\what{S}'_{P_K(K+4),\pi(K+4)}}
\ar@/_/ @{<-} "r-1r-1";"r-2r-2"_{\what{S}'_{1+P_K(K+1),\pi(K+1)}}
\ar@/_/ @{<-} "r-2r-2";"r-3r-3"_{\what{S}'_{1+P_K(K+2),\pi(K+2)}}
\ar@/_/ @{<-} "r-3r-3";"r-4r-4"_{\what{S}'_{1+P_K(K+3),\pi(K+3)}}
\ar@/_/ @{<-} "r-4r-4";"dot"_{\what{S}'_{1+P_K(K+4),\pi(K+4)}}
\end{xy}
\]
Here, we simply write $\fbox{$a$}^{C^{(1)}}_K$, $P^{C^{(1)}}_K$ and $\pi'_{C^{(1)}}$ as $\fbox{$a$}$, $P_K$ and $\pi$, respectively.

\section{Appendix B}

When $\mathfrak{g}$ is of type $X=D_{4}^{(1)}$ and $k=3$, we give the definition of ${\rm COMB}^X_{k,\iota}[\lambda]$.
Let us take $a,b,c,d\in I$ such that $\{a,b,c,d\}=\{1,2,4,5\}$.

If $(1,3)<(1,a),(1,b),(1,c),(1,d)$ then ${\rm COMB}^{D^{(1)}}_{3,\iota}[\lambda]:=\{-x_{1,3}+\lan h_3,\lambda\ran\}$.

If $(1,3)<(1,b),(1,c),(1,d)$ and $(1,3)>(1,a)$ then
${\rm COMB}^{D^{(1)}}_{3,\iota}[\lambda]:=\{-x_{1,3}+x_{1,a}+\lan h_3,\lambda\ran, -x_{2,a}+\lan h_3,\lambda\ran\}$.

If $(1,3)<(1,c),(1,d)$ and $(1,3)>(1,a),(1,b)$ then we put
\[
\fbox{$r_1,r_2$}_s:=x_{s,r_1}-x_{s+1,r_2},\quad
\fbox{$r_1,r_2,3$}_s:=x_{s,r_1}+x_{s,r_2}-x_{s+p_{3,r_1},3},\quad
\fbox{$3,r_1,r_2$}_s:=x_{s+p_{3,r_1},3}-x_{s+1,r_1}-x_{s+1,r_2}
\]
for $s\in\mathbb{Z}_{\geq1}$ and $r_1,r_2\in\{1,2,4,5\}$ and set
${\rm COMB}^{D^{(1)}}_{3,\iota}[\lambda]:=
\{
\fbox{$r_1,r_2,3$}_s+\langle h_{3},\lambda \rangle ,
\fbox{$r_1,r_2$}_s+\langle h_{3},\lambda \rangle ,
\fbox{$r_2,r_1$}_s+\langle h_{3},\lambda \rangle ,
\fbox{$3,r_1,r_2$}_s+\langle h_{3},\lambda \rangle  | (r_1,r_2)\in\{(a,b),(c,d)\},s\in\mathbb{Z}_{\geq1},s:\text{odd}
\}$.

If $(1,3)<(1,d)$ and $(1,3)>(1,a),(1,b),(1,c)$ then we put
\[
{\rm COMB}^{D^{(1)}}_{3,\iota}[\lambda]:=
\{
L^{D^{(1)}}_{-1,d,\iota}(T)+\langle h_{3},\lambda \rangle | T\in{\rm YW}^{D^{(1)}}_{d}\setminus\{Y_{\Lambda_{d}},Y_{d}^1\} 
\},
\]
where $Y_{d}^1$ is obtained from $Y_{\Lambda_{d}}$ by adding a $d$-block to the $d$-admissible slot.

If $(1,3)>(1,a),(1,b),(1,c),(1,d)$ then
\[
{\rm COMB}^{D^{(1)}}_{3,\iota}[\lambda]:=
\{
L^{D^{(1)}}_{0,3,\iota}(T)+\langle h_{3},\lambda \rangle | T\in{\rm YW}^{D^{(1)}}_{3}\setminus\{Y_{\Lambda_3}\}
\}.
\]

\end{document}